\documentclass[12pt,a4paper]{article}
\usepackage{amsmath,amssymb,epic,eepic,float}
\usepackage[dvipdf]{graphicx}
\graphicspath{{Fig-EPSc/}}
\makeatletter
\def\section{\@startsection {section}{1}{\z@}{-3.5ex plus -1ex minus 
 -.2ex}{2.3ex plus .2ex}{\normalsize\bf}}
\def\@maketitle{\newpage
 \null
 \vskip 2em \begin{center}
 {\large \@title \par} \vskip 1.5em {\normalsize \lineskip .5em
\begin{tabular}[t]{c}\@author 
 \end{tabular}\par} 
 \end{center}
 \par
 \vskip 1.5em} 
\makeatother

\begin{document}

\title{\bf %\family{cmr}\series{bx}\shape{it} 
On the fusion algebras of bimodules \\
arising from Goodman-de la Harpe-Jones subfactors}

\author{
By Satoshi {\sc Goto} %~~(Sophia University)
\thanks{This manuscript is typeset by LATEX2e 
\newline
\hspace{5mm}
2010 {\it Mathematics Subject Classification.} 
~46L37.}
}

%%%%%%%%%%%%%%%%%%%%%%%%%%%%%%%%%%%%%%%%%%%%%%%%%%%%%%%%%%%%%%%%%%%%%%
\maketitle
%%%%%%%%%%%%%%%%%%%%%%%%%%%%%%%%%%%%%%%%%%%%%%%%%%%%%%%%%%%%%%%%%%%%%%

\def\coef{{{\rm coef}}}
\def\Edge{{{\rm Edge\;}}}
\def\End{{{\rm End}}}
\def\EssPath{{{\rm EssPath}}}
\def\gap{{{\rm gap}}}
\def\0-gap{{{\rm 0-gap}}}
\def\GHJ{{{\rm GHJ}}}
\def\Hom{{{\rm Hom}}}
\def\HPath{{{\rm HPath}}}
\def\mod{{{\rm mod}}}
\def\Path{{{\rm Path}}}
\def\Proj{{{\rm Proj}}}
\def\String{{{\rm String}}}
\def\Vert{{{\rm Vert}}}

\def\N{{\bf N}}
\def\Z{{\bf Z}}
\def\Q{{\bf Q}}
\def\R{{\bf R}}
\def\C{{\bf C}}

\def\A{{\cal A}}
\def\B{{\cal B}}
\def\E{{\cal E}}
\def\G{{\cal G}}
\def\H{{\cal H}}
\def\K{{\cal K}}
\def\T{{\cal T}}
\def\l{{\cal L}}
\def\M{{\cal M}}
\def\n{{\cal N}}
\def\p{{\cal P}}
\def\q{{\cal Q}}
\def\r{{\cal R}}
\def\U{{\cal U}}
\def\V{{\cal V}}
\def\W{{\cal W}}
\def\z{{\cal Z}}

\def\a{\alpha}
\def\be{\beta}
\def\de{\delta}
\def\De{\Delta}
\def\e{\varepsilon}
\def\epsilon{\varepsilon}
\def\ga{\gamma}
\def\Ga{\Gamma}
\def\la{\lambda}
\def\La{\Lambda}
\def\si{\sigma}
\def\th{\theta}
\def\t{\tau}
\def\om{\omega}
\def\Om{\Omega}

\newcommand{\Frac}[2]{\displaystyle{\frac{#1}{#2}}}

\def\scr{\scriptstyle}

\newtheorem{theorem}{Theorem}[section]
\newtheorem{lemma}[theorem]{Lemma}
\newtheorem{corollary}[theorem]{Corollary}
\newtheorem{proposition}[theorem]{Proposition}
\newtheorem{remark}[theorem]{Remark}
\newtheorem{definition}[theorem]{Definition}
\newtheorem{example}[theorem]{Example}
\newtheorem{axiom}{Axiom}

\begin{abstract}
By using Ocneanu's result 
on the classification of all irreducible connections on the Dynkin diagrams,
we show that the dual principal graphs as well as the fusion rules 
of bimodules arising from any Goodman-de la Harpe-Jones subfactors 
are obtained by a purely combinatorial method.  
In particular we obtain the dual principal graph and the fusion rule 
of bimodules arising from the Goodman-de la Harpe-Jones subfactor 
corresponding to the Dynkin diagram $E_8$. 
As an application, we also show some subequivalence 
among $A$-$D$-$E$ paragroups.
\end{abstract}

%%%%%%%%%%%%%%%%%%%%%%%%%%%%%%%%%%%%%%%%%%%%%%%%%%%%%%%%%%%%%%%%%%%%%%
\section{Introduction}
%%%%%%%%%%%%%%%%%%%%%%%%%%%%%%%%%%%%%%%%%%%%%%%%%%%%%%%%%%%%%%%%%%%%%%

Since V. F. R. Jones initiated the index theory for subfactors
in \cite{J1}, intensive studies on the classification of subfactors
have been made by many people.  
The classification of subfactors of the AFD type II$_1$ factor 
with index less than 4 has been completed by many people's contribution
(\cite{BiN, I1, I2, J1, Ka1, SV}, see also \cite{EK-book}) 
after A. Ocneanu's announcement \cite{O1}.

 Goodman-de la Harpe-Jones subfactors (abbreviated as GHJ subfactors) \cite{GHJ} 
are known as a series of interesting non-trivial examples 
of irreducible subfactors with indices greater than 4, 
though some of them have indices less than 4.
The indices of all GHJ subfactors are given in \cite{GHJ}.
They are constructed from the commuting squares arising from the embeddings 
of type $A$ string algebras into other string algebras of type $ADE$.
(See \cite[Chapter 11]{EK-book} for the construction of GHJ subfactors 
from a viewpoint of string algebra embedding.)
The principal graphs of these subfactors are easily obtained by
a simple method but the dual principal graphs as well as their 
fusion rules are much more difficult to compute. 
(Okamoto first computed their principal graphs in \cite{Ok}.)

One of the most important examples of GHJ subfactor has index $3+\sqrt 3$
and it is constructed from the embedding of the string algebra of $A_{11}$
into that of $E_6$. In this particular case 
it happens that it is not very difficult to compute 
the dual principal graph (see \cite{Ka2}, \cite[Section 11.6]{EK-book}).
But it is more difficult to determine its fusion rules. 
Actually D. Bisch has tried to compute the fusion rule just from the graph 
but there were five possibilities and it turned out that the fusion
rule cannot be determined from the graph only \cite{Bs}.
Some more information is needed and Y. Kawahigashi obtained the fusion 
rule as an application of paragroup actions in \cite{Ka2}.

In his lectures at The Fields Institute 
A. Ocneanu introduced a new algebra called 
{\sl double triangle algebra} by using the notion of essential paths
and extension of Temperley-Lieb recoupling theory of Kauffman-Lins \cite{Oc}.
He also announced a solution to the problem of determining 
the dual principal graphs and their fusion rules of the GHJ subfactors
as one of some applications of his theory.
But the details have not been published.

After A. Ocneanu's works, F. Xu and J. B\"ockenhauer-D. E. Evans 
have revealed a relation between the GHJ subfactors and conformal inclusions 
(\cite{Xu}, \cite{BE1}, \cite{BE2}, \cite{BE3}) and  
J. B\"ockenhauer, D. E. Evans and Y. Kawahigashi (\cite{BEK2}) 
obtained essentially 
the same fusion algebras of GHJ subfactors of type $D_{2n}, E_6, E_8$ 
by using conformal field theory and the Cappelli-Itzykson-Zuber's 
classification of modular invariant \cite{CIZ}.

In this paper we give detailed computations 
of the dual principal graphs and the fusion rules 
for any GHJ subfactors by a purely combinatorial method.
For this purpose we will make the most use of Ocneanu's result 
on the classification of all irreducible connections on the Dynkin diagrams
(See \cite{Oc}. Our method here is based on the observation in \cite{G}). 
Especially we will make use of Figures \ref{EP-A4}`\ref{qsym-e8}, 
which were first found by A. Ocneanu \cite{Oc}. 
Our result does not rely on either conformal field theory or
the classification of modular invariant.

\vspace{-3mm}
%%%%%%%%%%%%%%%%%%%%%%%%%%%%%%%%%%%%%%%%%%%%%%%%%%%%%%%%%%%%%%%%%%%%%%
\section{Correspondence between system of connections and
system of bimodules}
%%%%%%%%%%%%%%%%%%%%%%%%%%%%%%%%%%%%%%%%%%%%%%%%%%%%%%%%%%%%%%%%%%%%%%

Let $K$ and $L$ be two connected finite bipartite graphs.
A bi-unitary connection on four graphs is called a {\sl $K$-$L$ bi-unitary
connection} if it has the graph $K$ as an upper horizontal graph 
and the graph $L$ as a lower horizontal graph as in Figure \ref{comm-1}.

If we have a $K$-$L$ connection, 
we can construct a subfactor $N\subset M$
by choosing a distinguished vertex $*_K$ of the upper graph $K$ and 
applying string algebra construction to the connection.
(See \cite[Section 11]{EK-book}.)
This construction seems to depend on the choice of the vertex $*_K$.
But it is well-known that the subfactors constructed from this connection
does not depend on the choice of the vertex $*_K$, that is,  
they become all isomorphic because of the relative McDuff property
\cite{Bs2}.

%%%%%%%%%%%%%%%%%%%%%
\begin{figure}[H]
\unitlength 0.7mm
\thicklines
\begin{center}
\begin{picture}(80,30)(0,5)
\drawline(15,30)(30,30)(30,10)(10,10)(10,25)
\put(10,30){\makebox(0,0){$*_K$}}
\put(5,20){\makebox(0,0){}}
\put(20,35){\makebox(0,0){$K$}}
\put(35,20){\makebox(0,0){}}
\put(20,5){\makebox(0,0){$L$}}
\put(20,20){\makebox(0,0){$w$}}
\put(50,30){\makebox(0,0){$\cdots \quad \longrightarrow$}}
\put(50,10){\makebox(0,0){$\cdots \quad \longrightarrow$}}
\put(70,30){\makebox(0,0){$N$}}
\put(70,20){\makebox(0,0){$\cap$}}
\put(70,10){\makebox(0,0){$M$}}
\end{picture}
\caption{ }
\label{comm-1}
\end{center}
\end{figure}

On the one hand as a paragroup of the subfactor $N\subset M$ obtained 
from the connection $w$ as above,
we obtain the system of 4-kinds of bimodules, i.e. 
$N$-$N$, $N$-$M$, $M$-$N$, $M$-$M$ bimodules, by taking irreducible 
decomposition of alternating relative tensor products of ${}_NM_M$ and 
its conjugate bimodule ${}_MM_N$ as usual. (See \cite{EK-book} for details.)

On the other hand we also get the system of 4-kinds of connections,
i.e. $K$-$K$, $K$-$L$, $L$-$K$, $L$-$L$ bi-unitary connections, 
by taking irreducible decomposition of alternating compositions of 
the connection $w$ and its conjugate $L$-$K$ connection $\bar w$.

Now the problem is the relation between the system of bimodules and 
the system of connections obtained as above.
We can easily see that those two systems become the same paragroup 
for $N\subset M$ if the subfactor $N\subset M$ has finite depth.

To see this it is enough to see the relation among a usual paragroup
based on bimodules,
a system of generalized open string bimodules and a system of bi-unitary
connections. The details of these relations are found in \cite{AH}.
%%%
Note that when we consider a system of bi-unitary connections
forms a paragroup, we need the notion of intertwiners between two
connections. For this purpose, we need to fix distinguished vertices 
$*_K$ and $*_L$ of both even and odd part of the graphs $K$ and $L$, 
then we identify all the bi-unitary connections of the system as 
the generalized open string bimodules constructed from those connections.
Then we define the intertwiners between two connections by those 
between the corresponding two generalized open string bimodules. 
Now from the argument in \cite{AH}, 
the intertwiners between two connections can naturally be identified with 
the intertwiners between the correponding 4 kinds of bimodules, i.e.
$N$-$N$, $N$-$M$, $M$-$N$, $M$-$M$ bimodules arising 
from the usual paragroup.
See Theorem 4 in \cite{AH} for more details.

Hence we obtain the following theorem.

\begin{theorem}{\rm
If the subfactor $N\subset M$ constructed from a $K$-$L$ connection ${}_Kw_L$ 
has finite depth, the system of 4-kinds of connections obtained from ${}_Kw_L$ 
and the system of 4-kinds of bimodules obtained from the subfactor $N\subset M$
have the same fusion rules. 
Moreover these two systems defines the same paragroup for $N\subset M$
via the correspondence between connections and 
generalized open string bimodulesD
}\end{theorem}

\begin{remark}{\rm As we mentioned above, the subfactor constructed 
from a connection ${}_Kw_L$ does not depend on the choice of 
the distinguished vertex $*_K$. 
In the same way we need to fix two vertices $*_K$ and $*_L$ 
in order to construct a generalized open string bimodule from 
a connection. But the above theorem holds true for arbitrary 
choice of two distinguished vertices $*_K$ and $*_L$ of the graphs
$K$ and $L$ respectively.
}\end{remark}

The above theorem provide us a purely combinatorial method
to compute fusion rules for the subfactor obtained from 
a connection ${}_Kw_L$. 
Actually we can compute the fusion rules of a system of connections
by looking at the composition and decomposition of their vertical graphs.

%%%%%%%%%%%%%%%%%%%%%%%%%%%%%%%%%%%%%%%%%%%%%%%%%%%%%%%%%%%%%%%%%%%%%%
\section{The (dual) principal graphs and their fusion rules of the
Goodman-de la Harpe-Jones subfactor}
%%%%%%%%%%%%%%%%%%%%%%%%%%%%%%%%%%%%%%%%%%%%%%%%%%%%%%%%%%%%%%%%%%%%%%

\begin{figure}[H]
\centering\includegraphics[width=110mm,clip]{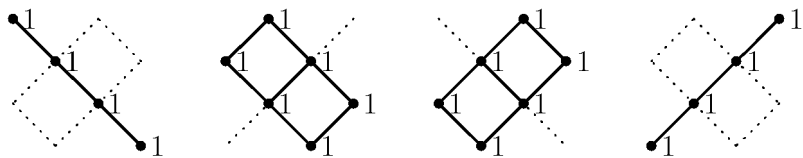}
\caption{Essential paths on the Coxeter graph $A_4$}
\label{EP-A4}
\end{figure}

\begin{figure}[H]
\centering\includegraphics[width=140mm,clip]{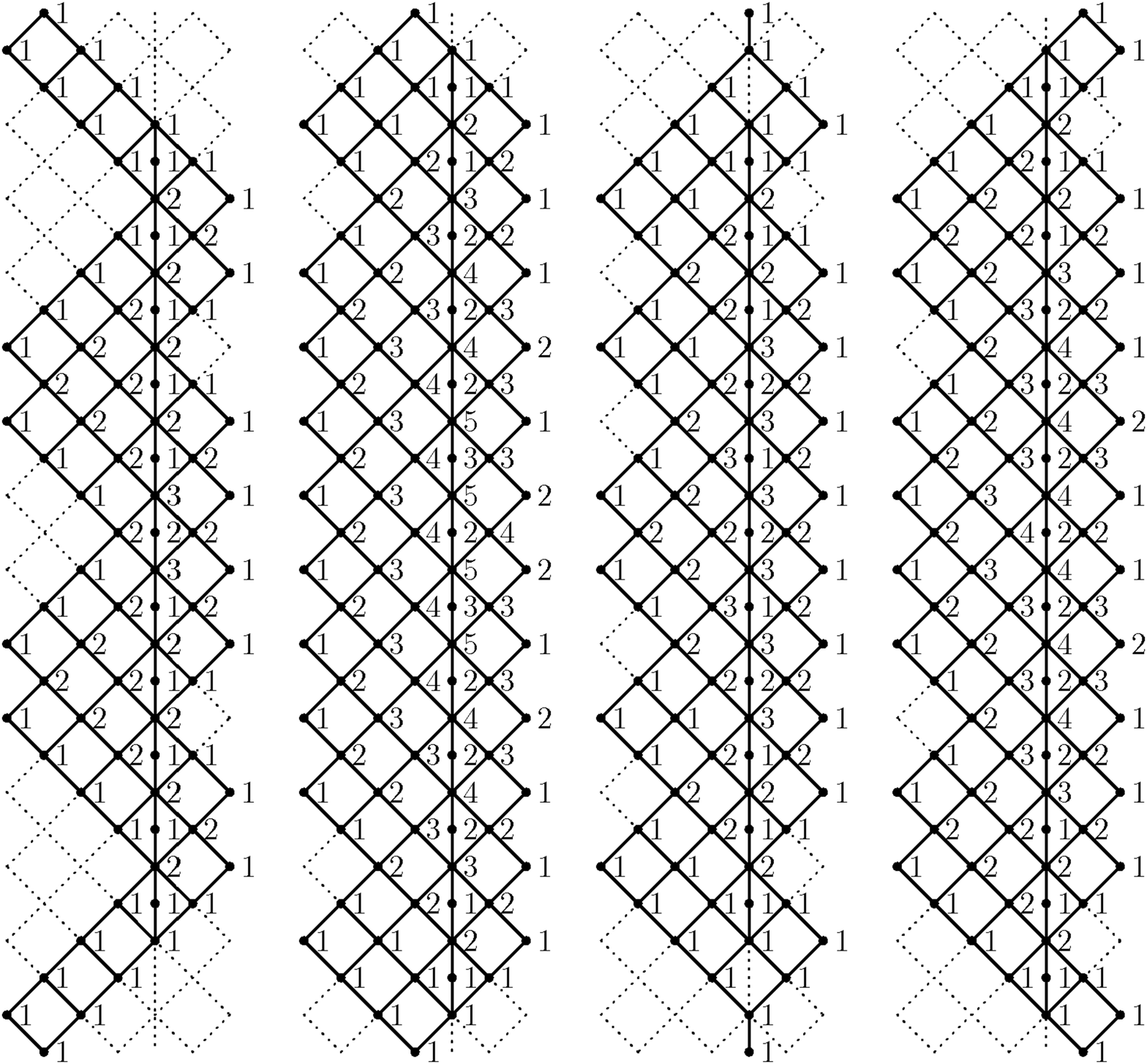}
\caption{Essential paths on the Coxeter graph $E_8$ (2)}
\label{EP-E8-2}
\end{figure}

Let $A$ be the Dynkin diagram $A_n$ and 
$K$ one of the $ADE$ Dynkin diagrams with the same Coxeter number.
The subfactors constructed from the commuting square 
as in Figure \ref{comm-2} are called the Goodman-de la Harpe-Jones subfactors
(abbreviated as GHJ subfactors).
Here the construction depends only on the graph $K$ and the vertex $*_K=x$.
(See \cite{GHJ} for details.)
We denote this subfactor $\GHJ(K,*_K=x)$.
We remark that the vertical graphs $G$ and $G'$ as in Figure \ref{comm-2} 
are easily obtained from the dimension of essential paths 
on the graph $K$ (Figures \ref{EP-A4}`\ref{EP-E8-2}). 
Here we note that the graphs $G$ and $G'$ may be disconnected. 

%%%%%%%%%%%%%%%%%%%%%
\begin{figure}[H]
\unitlength 0.7mm
\thicklines
\begin{center}
\begin{picture}(80,30)(0,5)
\drawline(15,30)(30,30)(30,10)(15,10)
\drawline(10,15)(10,25)
\put(10,30){\makebox(0,0){$*_A$}}
\put(10,10){\makebox(0,0){$*_K$}}
\put(5,20){\makebox(0,0){$G$}}
\put(20,35){\makebox(0,0){$A$}}
\put(35,20){\makebox(0,0){$G'$}}
\put(20,5){\makebox(0,0){$K$}}
\put(20,20){\makebox(0,0){$w$}}
\put(50,30){\makebox(0,0){$\cdots \quad \longrightarrow$}}
\put(50,10){\makebox(0,0){$\cdots \quad \longrightarrow$}}
\put(70,30){\makebox(0,0){$N$}}
\put(70,20){\makebox(0,0){$\cap$}}
\put(70,10){\makebox(0,0){$M$}}
\end{picture}
\caption{ }
\label{comm-2}
\end{center}
\end{figure}

We use the next two propositions to compute the fusion rule of 
the Goodman-de la Harpe-Jones subfactors.

\begin{proposition}{\rm (\cite[Proposition 5.6]{G})}
\label{AKconnection}
Let $A, K, G$ and $G'$ be the four graphs connected as in Figure \ref{4graphs-1}. 
Suppose there is a bi-unitary connection on the four graphs. 
Then the connecting vertical graphs $G$ and $G'$ are uniquely determined 
by the {\sl initial condition},
i.e., the condition of edges connected to the distinguished vertex 
of the graph $A$ (see Figure \ref{comm-3}). 
Moreover such a connection is unique up to vertical gauge choice.
\end{proposition}

%%%%%%%%%% 4graphs-1 %%%%%%%%%%%%%%%%%%%
\unitlength 0.7mm
\begin{figure}[H]
\begin{center}
\begin{picture}(45,16)
\path(10,0)(30,0)(30,20)(10,20)(10,0)
\put(10,0){\circle*{1}}
\put(30,0){\circle*{1}}
\put(30,20){\circle*{1}}
\put(10,20){\circle*{1}}
\put(5,10){\makebox(0,0){$G$}}
\put(20,25){\makebox(0,0){$A$}}
\put(35,10){\makebox(0,0){$G'$}}
\put(20,-5){\makebox(0,0){$K$}}
\end{picture}
\end{center}
\caption{}
\label{4graphs-1}
\end{figure} 
%%%%%%%%%%%%%%%%%%%%%%%%%%%%%%%%%%%%%%

%%%%%%%%%%%%%%%%%%%%%
\begin{figure}[H]
\unitlength 0.7mm
\thicklines
\begin{center}
\begin{picture}(40,30)(2,5)
\drawline(15,30)(30,30)(30,10)(15,10)
\drawline(10,15)(10,25)
\put(10,30){\makebox(0,0){$*_A$}}
\put(10,10){\makebox(0,0){$x$}}
\put(5,20){\makebox(0,0){}}
\put(20,35){\makebox(0,0){$A$}}
\put(35,20){\makebox(0,0){}}
\put(20,5){\makebox(0,0){$K$}}
\put(20,20){\makebox(0,0){$w$}}
%\put(50,30){\makebox(0,0){$\cdots \quad \longrightarrow$}}
%\put(50,10){\makebox(0,0){$\cdots \quad \longrightarrow$}}
%\put(70,30){\makebox(0,0){$N$}}
%\put(70,20){\makebox(0,0){$\cap$}}
%\put(70,10){\makebox(0,0){$M$}}
\end{picture}
\caption{ }
\label{comm-3}
\end{center}
\end{figure}
%%%%%%%%%%%%%%%%%%%%%

\begin{proposition}{\rm 
{\bf (Frobenius reciprocity)} (\cite[Proposition 3.21]{G})}
Let $K$, $L$ and $M$ be three connected finite bipartite graphs
with the same Perron-Frobenius eigenvalue. Let ${}_K\a_L$,
${}_L\be_M$ and ${}_K\ga_M$ be three irreducible bi-unitary 
connections which are $K$-$L$, $L$-$M$ and $K$-$M$ respectively.
If $\ga$ appears $n$ times in the composite connection
$\a \be$, then $\a$ appears $n$ times in $\ga \bar \be$
and $\be$ appears $n$ times in $\bar \a \ga$.
\end{proposition}

\vspace{3mm}
\noindent
%%%%%%%%%%%%%%%%%%%%%%%%%%%%%%%%%%%%%%%%%%%%%%%%%%%%%%%%%%%%%
%%%%% section 3.1 %%%%%%%%%%%%%%%%%%%%%%%%%%%%%%%%%%%%%%%%%%%
%%%%%%%%%%%%%%%%%%%%%%%%%%%%%%%%%%%%%%%%%%%%%%%%%%%%%%%%%%%%%
{\sl 3.1. The fusion rules of four kinds of connections arising from GHJ subfactors}
\medskip

\vspace{1mm}
%\noindent
The system of connections arising from a GHJ subfactor consists of four kinds
of connections, i.e. $A$-$A$, $A$-$K$, $K$-$A$ and $K$-$K$ connections.
So the fusion rules consist of the following 8 kinds of multiplication table.

\vspace{2mm}
\noindent
$(1)$~$A$-$A$ $\times$ $A$-$A$ $\longrightarrow$ $A$-$A$ \\
$(2)$~$A$-$A$ $\times$ $A$-$K$ $\longrightarrow$ $A$-$K$ \\
%\hspace{6.4mm}
$(2)'$~$K$-$A$ $\times$ $A$-$A$ $\longrightarrow$ $K$-$A$ 
\hspace{6.3mm}
$(2)''$~$A$-$K$ $\times$ $K$-$A$ $\longrightarrow$ $A$-$A$ \\
$(3)$~$A$-$K$ $\times$ $K$-$K$ $\longrightarrow$ $A$-$K$ \\
%\hspace{6mm}
$(3)'$~$K$-$K$ $\times$ $K$-$A$ $\longrightarrow$ $K$-$A$ 
\hspace{6mm}
$(3)''$~$K$-$A$ $\times$ $A$-$K$ $\longrightarrow$ $K$-$K$ \\
$(4)$~$K$-$K$ $\times$ $K$-$K$ $\longrightarrow$ $K$-$K$ 

\vspace{2mm}
Among these multiplication tables, $(2)'$ and $(3)'$ are 
obtained by taking conjugation of $(2)$ and $(3)$ respectively.
The tables $(2)''$ and $(3)''$ are also obtained 
from $(2)$ and $(3)$ respectively by Frobenius reciprocity.
So it is enough to determine four multiplication table 
$(1)$, $(2)$, $(3)$ and $(4)$.

\vspace{3mm}
\noindent
%%%%%%%%%%%%%%%%%%%%%%%%%%%%%%%%%%%%%%%%%%%%%%%%%%%%%%%%%%%%%
%%%%% section 3.1.1 %%%%%%%%%%%%%%%%%%%%%%%%%%%%%%%%%%%%%%%%%
%%%%%%%%%%%%%%%%%%%%%%%%%%%%%%%%%%%%%%%%%%%%%%%%%%%%%%%%%%%%%
{\sl 3.1.1. The fusion rules of 
~$(1)$~$A$-$A$ $\times$ $A$-$A$ $\longrightarrow$ $A$-$A$ ~and~
$(2)$~$A$-$A$ $\times$ $A$-$K$ $\longrightarrow$ $A$-$K$ and the principal graphs}
\medskip

We put the labels $0,1,2,\cdots,m-1$ of vertices of the Dynkin diagram $A_m$ 
as in Figure \ref{label-Am}.
We denote the unique irreducible $A$-$A$ connection with the ``initial edge''
connected to the vertex $n$ in the lower graph $A_m$ by ${}_An_A$
(Figure \ref{comm-4}).
We also denote the unique irreducible $A$-$K$ connection with the ``initial edge''
connected to the vertex $x$ in the lower graph $K$ by ${}_Ax_K$
(Figure \ref{comm-4}).

%%%%%%%%%%%%%%%%%%%%%%%%%%%%%%%%%%%%%%%%%%%%%%%%%%%%
\begin{figure}[H]
\unitlength 0.85mm
\thinlines
\begin{center}
\begin{picture}(80,17)(-20,-5)
\multiput(10,0)(20,0){2}{\circle*{2}}
\multiput(20,10)(20,0){4}{\circle*{2}}
\put(70,0){\circle*{2}}
\path(2,8)(10,0)(20,10)(30,0)(40,10)
\path(60,10)(70,0)(80,10) 
\put(0,10){\makebox(0,0){$*_A$}}
\put(0,14){\makebox(0,0){$0$}}
\put(10,-4){\makebox(0,0){$1$}}
\put(20,14){\makebox(0,0){$2$}}
\put(30,-4){\makebox(0,0){$3$}}
\put(40,14){\makebox(0,0){$4$}}
\put(50,5){\makebox(0,0){$\cdots$}}
\put(60,14){\makebox(0,0){$m-3$}}
\put(70,-4){\makebox(0,0){$m-2$}}
\put(80,14){\makebox(0,0){$m-1$}}
\put(-30,5){\makebox(0,0){Dynkin diagram $A_m$}}
%\put(80,14){\makebox(0,0){even}}
%\put(80,-4){\makebox(0,0){odd}}
\end{picture}
\caption{The label of vertices of the Dynkin diagram $A_m$.}
\label{label-Am}
\end{center}
\end{figure}
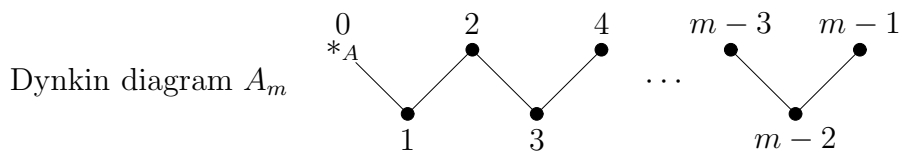
%%%%%%%%%%%%%%%%%%%%%%%%%%%%%%%%%%%%%%%%%%%%%%%%%%%%

%%%%%%%%%%%%%%%%%%%%%
\begin{figure}[H]
\unitlength 0.8mm
\thicklines
\begin{center}
\begin{picture}(80,17)(5,5)
\drawline(15,30)(30,30)(30,10)(15,10)
\drawline(10,15)(10,25)
\put(10,30){\makebox(0,0){$*_A$}}
\put(10,10){\makebox(0,0){$n$}}
\put(5,20){\makebox(0,0){}}
\put(20,35){\makebox(0,0){$A$}}
\put(35,20){\makebox(0,0){}}
\put(20,5){\makebox(0,0){$A$}}
\put(20,20){\makebox(0,0){${}_An_A$}}
\drawline(65,30)(80,30)(80,10)(65,10)
\drawline(60,15)(60,25)
\put(60,30){\makebox(0,0){$*_A$}}
\put(60,10){\makebox(0,0){$x$}}
\put(55,20){\makebox(0,0){}}
\put(70,35){\makebox(0,0){$A$}}
\put(85,20){\makebox(0,0){}}
\put(70,5){\makebox(0,0){$K$}}
\put(70,20){\makebox(0,0){${}_Ax_K$}}
\end{picture}
\caption{ }
\label{comm-4}
\end{center}
\end{figure}
%%%%%%%%%%%%%%%%%%%%%

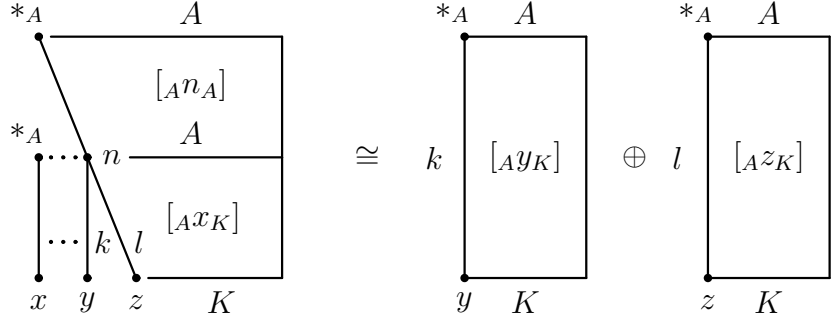
\begin{figure}[H]
\unitlength 0.8mm
\thicklines
\begin{center}
\begin{picture}(120,50)(5,5)
\drawline(0,50)(8,30)(8,10)
\put(0,50){\circle*{1}}
\put(0,30){\circle*{1}}
\put(8,30){\circle*{1}}
\multiput(0,10)(8,0){3}{\circle*{1}}
\drawline(8,30)(16,10)
\drawline(0,30)(0,10)
\multiput(2,30)(2,0){3}{\circle*{0.3}}
\multiput(2,16)(2,0){3}{\circle*{0.3}}
\put(0,6){\makebox(0,0){$x$}}
\put(8,6){\makebox(0,0){$y$}}
\put(16,6){\makebox(0,0){$z$}}
\put(12,30){\makebox(0,0){$n$}}
\put(10.5,16){\makebox(0,0){$k$}}
\put(16.5,16){\makebox(0,0){$l$}}
\put(-2,54){\makebox(0,0){$*_A$}}
\put(-2,34){\makebox(0,0){$*_A$}}
\drawline(2,50)(40,50)(40,30)(15,30)
\drawline(40,30)(40,10)(18,10)
\put(25,54){\makebox(0,0){$A$}}
\put(25,34){\makebox(0,0){$A$}}
\put(30,6){\makebox(0,0){$K$}}
\put(25,42){\makebox(0,0){$[{}_An_A]$}}
\put(27,20){\makebox(0,0){$[{}_Ax_K]$}}
\put(54,30){\makebox(0,0){$\cong$}}
\put(65,30){\makebox(0,0){$k$}}
\drawline(70,10)(70,50)(90,50)(90,10)(70,10)
\multiput(70,10)(0,40){2}{\circle*{1}}
\put(68,54){\makebox(0,0){$*_A$}}
\put(70,6){\makebox(0,0){$y$}}
\put(80,54){\makebox(0,0){$A$}}
\put(80,6){\makebox(0,0){$K$}}
\put(80,30){\makebox(0,0){$[{}_Ay_K]$}}
\put(98,30){\makebox(0,0){$\oplus$}}
\put(105,30){\makebox(0,0){$l$}}
\drawline(110,10)(110,50)(130,50)(130,10)(110,10)
\multiput(110,10)(0,40){2}{\circle*{1}}
\put(108,54){\makebox(0,0){$*_A$}}
\put(110,6){\makebox(0,0){$z$}}
\put(120,54){\makebox(0,0){$A$}}
\put(120,6){\makebox(0,0){$K$}}
\put(120,30){\makebox(0,0){$[{}_Az_K]$}}
\end{picture}
\caption{The fusion rule of $A$-$A$ $\times$ $A$-$K$ $\longrightarrow$ $A$-$K$ }
\label{fusion-AAAK-AK}
\end{center}
\end{figure}

Then the fusion rules of 
$(1)$~$A$-$A$ $\times$ $A$-$A$ $\longrightarrow$ $A$-$A$ and 
$(2)$~$A$-$A$ $\times$ $A$-$K$ $\longrightarrow$ $A$-$K$ can be 
obtained by composition and decomposition of 
the (left) vertical edges of the two connections ${}_An_A$ and ${}_Ax_K$
as in Figure \ref{fusion-AAAK-AK}.
So we have only to count the vertical edges of the connection
${}_Ax_K$ in order to determine the fusion tables of (1) and (2).

%%%%%%%%%%%%%%%%%%%%%
%%%%%%%%%%%%%%%%%%%%%

Because we need the notion of essential paths on graphs 
in order to describe these fusion rules, we review the definition
here for readers convenience. 
Please see \cite[section 32.2, page 254--256]{Oc} for more
details and the proof of the moderated Pascal rule.

\begin{definition}
A space of essential paths of a graph $G$ with length $n$ is defined
by $\EssPath^{(n)} G = p_n \cdot \HPath^{(n)} G$.
Here $p_n = 1- e_1 \vee e_2 \vee \cdots \vee e_{n-1}$ is 
the Wenzl projector and $e_k$ is the $k$-th Jones projection.  
We denote the space of essential paths of a graph $G$ with length $n$, 
with starting point $x$ and end point $y$ 
by $\EssPath_{x,y}^{(n)} G$.
\end{definition}

The dimensions of spaces of essential paths of length $n$
is easily obtained by using the following {\sl moderated Pascal rule}.
\[
\dim \EssPath_{a,x}^{(n+1)} G = \sum_{\xi \in \Edge G, r(\xi) =x}
       \dim \EssPath_{a,s(\xi)}^{(n)} G - \dim \EssPath_{a,x}^{(n-1)} G
\]

Now we continue the description of the fusion rules (1) and (2).
Because the connection ${}_Ax_K$ comes from the inclusion of 
the string algebras $\String_*A \subset \String_xK$, the number of 
vertical edges of this inclusion coincides with the dimension of 
essential paths from the vertex $x$ to $y$ of $K$ with length $n$.
(See Figures \ref{EP-A4}`\ref{EP-E8-2} for the dimension of 
essential paths.)
Hence we get the fusion tables of (1) and (2) as follows.  

\vspace{3mm}
%%%%%%%%%%
%\vspace{-2mm}
\hspace{38mm}
$\displaystyle{{}_An_A\cdot{}_Ax_K \cong \bigoplus_{y\in\Vert K}
(\dim\EssPath_{x,y}^{(n)}K)~{}_Ay_K}$ \\

%\vspace{-2mm}
\hspace{38mm}
$\displaystyle{{}_K\bar x_A\cdot{}_An_A \cong \bigoplus_{y\in\Vert K}
(\dim\EssPath_{x,y}^{(n)}K)~{}_K\bar y_A}$ \\

%\vspace{-2mm}
\hspace{38mm}
$\displaystyle{{}_Ay_K\cdot{}_K\bar x_A \cong \bigoplus_{n\in\Vert A}
(\dim\EssPath_{x,y}^{(n)}K)~{}_An_A}$
%%%%%%%%%%
\vspace{3mm}

Since the principal graph is obtained from the fusion rule of 
$A$-$A$ $\times$ $A$-$K$ $\longrightarrow$ $A$-$K$,
we can easily see that the principal graph of $\GHJ(K,*_K=x)$
coincides with the connected component of the vertical edges 
of the connection ${}_Ax_K$ including the distinguished vertex $*_A$.
This principal graph can be obtained easily by counting the dimension
of essential path.
It follows from this fact that the even vertices of the the principal graph 
of $\GHJ(K,*_K=x)$ coincides with (possibly a subset of) the even vertices of 
the Dynkin diagram $A_m$.

%\newpage
\vspace{3mm}
\noindent
%%%%%%%%%%%%%%%%%%%%%%%%%%%%%%%%%%%%%%%%%%%%%%%%%%%%%%%%%%%%%
%%%%% section 3.1.2 %%%%%%%%%%%%%%%%%%%%%%%%%%%%%%%%%%%%%%%%%
%%%%%%%%%%%%%%%%%%%%%%%%%%%%%%%%%%%%%%%%%%%%%%%%%%%%%%%%%%%%%
{\sl 3.1.2. The fusion rules of 
~$(3)$~$A$-$K$ $\times$ $K$-$K$ $\longrightarrow$ $A$-$K$ ~and 
the dual principal graphs}
\medskip

\begin{figure}[H]
\unitlength 0.8mm
\thicklines
\begin{center}
\begin{picture}(120,40)(5,5)
\drawline(0,50)(8,30)(0,10)
\put(0,50){\circle*{1}}
\put(0,10){\circle*{1}}
\put(16,10){\circle*{1}}
\put(8,30){\circle*{1}}
\drawline(8,30)(16,10)
\multiput(6,17)(2,0){3}{\circle*{0.3}}
\multiput(5,10)(3,0){3}{\circle*{0.3}}
\put(4,30){\makebox(0,0){$x$}}
\put(0,6){\makebox(0,0){$y$}}
\put(16,6){\makebox(0,0){$z$}}
\put(-1,17){\makebox(0,0){$k$}}
\put(16,17){\makebox(0,0){$l$}}
\put(-2,54){\makebox(0,0){$*_A$}}
%\put(-2,34){\makebox(0,0){$*_A$}}
\drawline(2,50)(40,50)(40,30)(10,30)
\drawline(40,30)(40,10)(18,10)
\put(25,54){\makebox(0,0){$A$}}
\put(25,34){\makebox(0,0){$K$}}
\put(30,6){\makebox(0,0){$K$}}
\put(25,42){\makebox(0,0){$[{}_Ax_K]$}}
\put(27,20){\makebox(0,0){$[{}_Kw_i{}_K]$}}
\put(54,30){\makebox(0,0){$\cong$}}
\put(65,30){\makebox(0,0){$k$}}
\drawline(70,10)(70,50)(90,50)(90,10)(70,10)
\multiput(70,10)(0,40){2}{\circle*{1}}
\put(68,54){\makebox(0,0){$*_A$}}
\put(70,6){\makebox(0,0){$y$}}
\put(80,54){\makebox(0,0){$A$}}
\put(80,6){\makebox(0,0){$K$}}
\put(80,30){\makebox(0,0){$[{}_Ay_K]$}}
\put(98,30){\makebox(0,0){$\oplus$}}
\put(105,30){\makebox(0,0){$l$}}
\drawline(110,10)(110,50)(130,50)(130,10)(110,10)
\multiput(110,10)(0,40){2}{\circle*{1}}
\put(108,54){\makebox(0,0){$*_A$}}
\put(110,6){\makebox(0,0){$z$}}
\put(120,54){\makebox(0,0){$A$}}
\put(120,6){\makebox(0,0){$K$}}
\put(120,30){\makebox(0,0){$[{}_Az_K]$}}
\end{picture}
\caption{The fusion rule of $A$-$K$ $\times$ $K$-$K$ $\longrightarrow$ $A$-$K$ }
\label{fusion-AKKK-AK}
\end{center}
\end{figure}
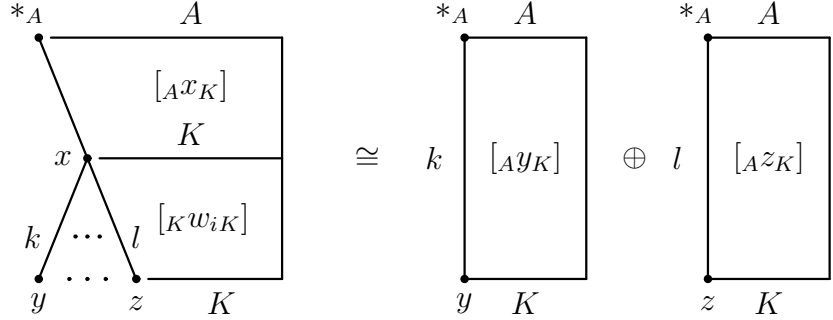

\begin{figure}[H]
%\centering\includegraphics[width=140mm,clip]{qsym-a.eps}
\centering\includegraphics[height=40mm,clip]{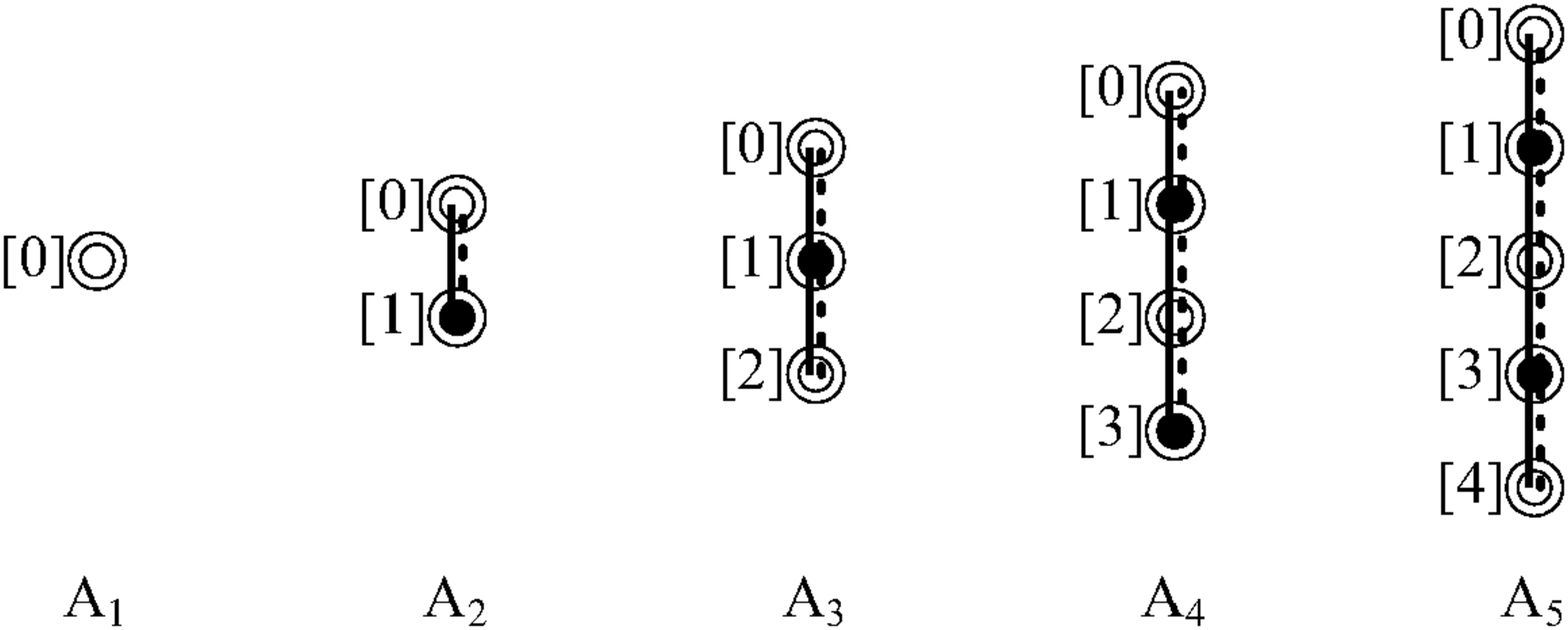}
\caption{Chiral symmetry for the Coxeter graph $A_n$}
\label{qsym-a}
\end{figure}

\begin{figure}[H]
%\centering\includegraphics[width=140mm,clip]{qsym-e8.eps}
%\centering\includegraphics[height=150mm,clip]{qsym-e8.eps}
\centering\includegraphics[height=78mm,clip]{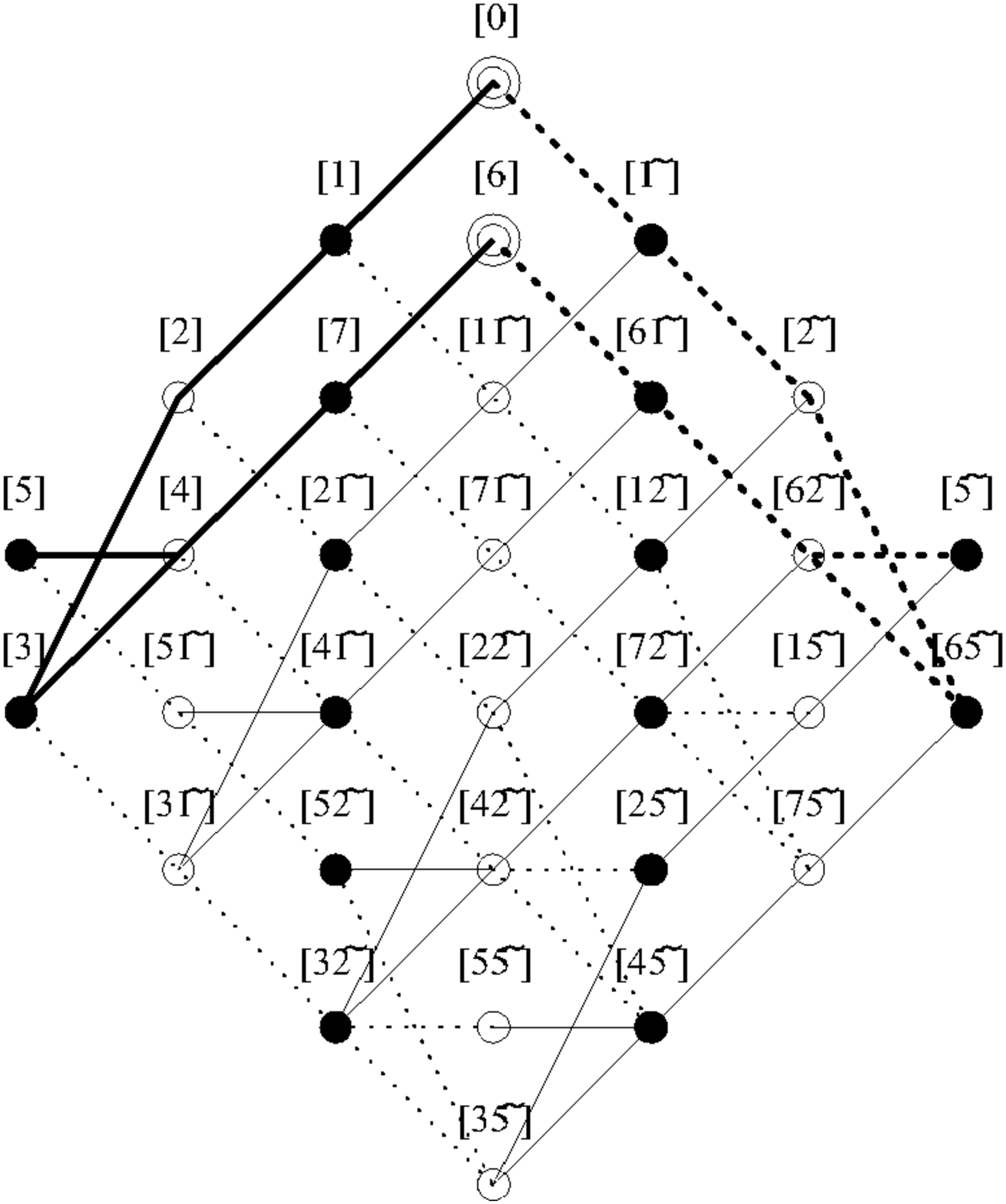}
\caption{Chiral symmetry for the Coxeter graph $E_8$}
\label{qsym-e8}
\end{figure}

\vspace{3mm}
We denote the unique irreducible $A$-$K$ connection with the ``initial edge''
connected to the vertex $x$ in the lower graph $K$ by ${}_Ax_K$ as before
and an irreducible $K$-$K$ connection by ${}_Kw_i{}_K$
(Figure \ref{comm-5}).
Here ${}_Kw_i{}_K$ is one of the connections of all $K$-$K$ connection system
(Figures \ref{qsym-a}`\ref{qsym-e8}).
(See \cite[section 5.3, pages 244--252]{G} for details.)
In this case the fusion rule of 
$(3)$~$A$-$K$ $\times$ $K$-$K$ $\longrightarrow$ $A$-$K$ is also obtained 
by composition and decomposition of the (left) vertical edges of 
the two connections ${}_Ax_K$ and ${}_Kw_i{}_K$
as in Figure \ref{fusion-AKKK-AK}.
We can get the fusion table of (3) by counting the vertical edges 
of the connection ${}_Kw_i{}_K$ in the same way as subsection 3.1.1.

%%%%%%%%%%%%%%%%%%%%%
\begin{figure}[H]
\unitlength 0.8mm
\thicklines
\begin{center}
\begin{picture}(80,30)(5,5)
\drawline(15,30)(30,30)(30,10)(15,10)
\drawline(10,15)(10,25)
\put(10,30){\makebox(0,0){$*_A$}}
\put(10,10){\makebox(0,0){$x$}}
\put(5,20){\makebox(0,0){}}
\put(20,35){\makebox(0,0){$A$}}
\put(35,20){\makebox(0,0){}}
\put(20,5){\makebox(0,0){$K$}}
\put(20,20){\makebox(0,0){${}_Ax_A$}}
\drawline(60,30)(80,30)(80,10)(60,10)
\drawline(60,10)(60,30)
\put(70,35){\makebox(0,0){$K$}}
\put(70,5){\makebox(0,0){$K$}}
\put(70,20){\makebox(0,0){${}_Kw_i{}_K$}}
\end{picture}
\caption{ }
\label{comm-5}
\end{center}
\end{figure}
%%%%%%%%%%%%%%%%%%%%%

%%%%%%%%%%%%%%%%%%%%%
%%%%%%%%%%%%%%%%%%%%%

This time the method of counting dimensions of essential paths 
does not work in order to get the vertical edges 
of the connection ${}_Kw_i{}_K$.
But we can compute them by using Ocneanu's classification 
of all irreducible $K$-$K$ connections and their fusion rules 
(\cite[section 5.3, pages 244--252]{G}).

\begin{figure}[H]
\centering\includegraphics[width=100mm,clip]{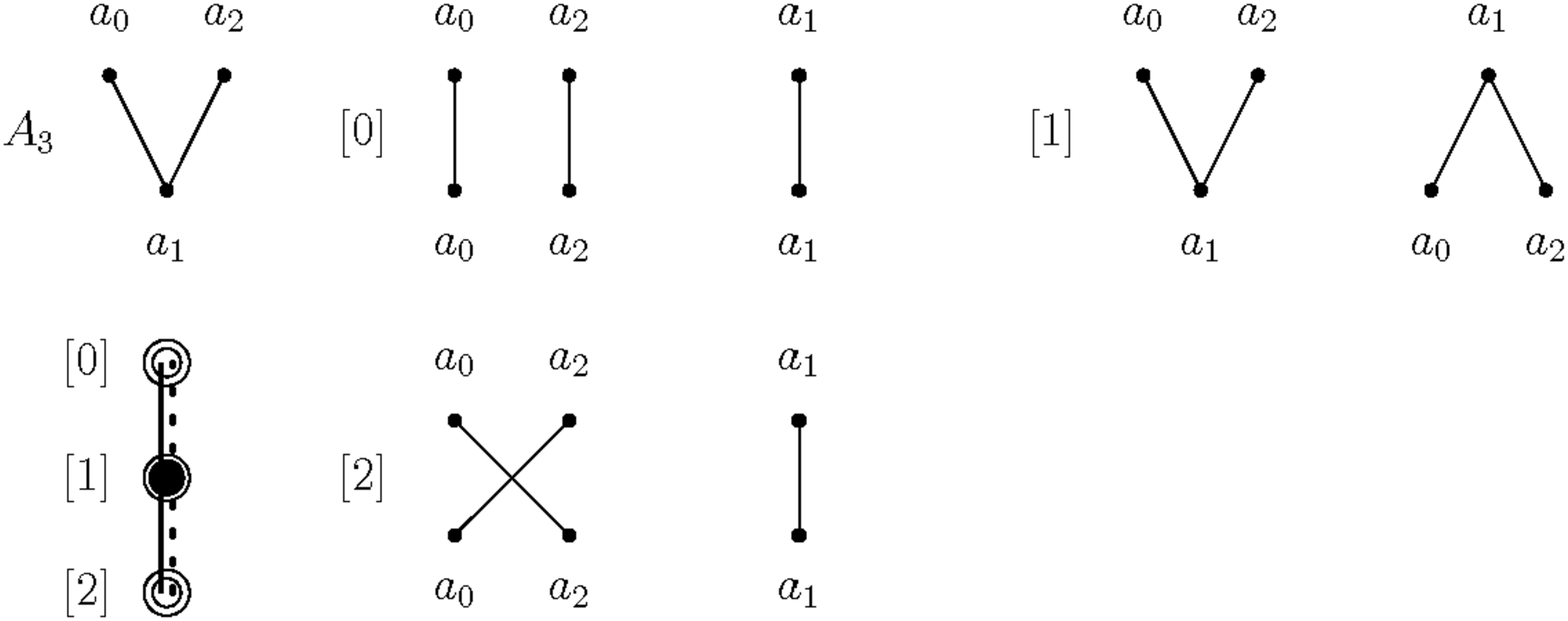}
\caption{Vertical graphs for connections on 
the Coxeter graph $A_3$}
\label{v-edge-A3}
\end{figure}

For example, the vertical edges of all $K$-$K$ connections are given 
in Figures \ref{v-edge-A3}`\ref{v-edge-E8odd} in the case of 
$K=A_3,A_4,A_5,A_6,D_4,D_5,D_6,E_6,E_7,E_8$.
Here in the case of $E_6,E_7,E_8$, we give list of incidence matrices 
of vertical graphs instead of graphs themselves because it is complicated
to draw them all.

Now we get the fusion rule of $A$-$K$ $\times$ $K$-$K$ $\longrightarrow$ $A$-$K$
as follows.

\vspace{4mm}
%%%%%%%%%%
%\vspace{-2mm}
\hspace{45mm}
$\displaystyle{{}_Ax_K\cdot{}_Kw_i{}_K \cong \bigoplus_{y\in\Vert K}
n(w_i)_{x,y}~{}_Ay_K}$ \\

%\vspace{-2mm}
\hspace{45mm}
$\displaystyle{{}_Kw_i{}_K\cdot{}_K\bar x_A \cong \bigoplus_{y\in\Vert K}
n(w_i)_{x,y}~{}_K\bar y_A}$ \\

%\vspace{-2mm}
\hspace{45mm}
$\displaystyle{{}_K\bar x_A\cdot{}_Ay_K \hspace{2.5mm} \cong \bigoplus_{w_i\in {}_KZ_K}
n(w_i)_{x,y}~{}_Kw_i{}_K}$
%%%%%%%%%%

\vspace{4mm}
Here ${}_KZ_K$ represents the system of all $K$-$K$ connections 
which is isomorphic to the fusion algebras of the center of 
$K$-$K$ double triangle algebra (\cite[Theorem 4.1, Corollary 4.5]{G}).
And $n(w_i)_{x,y}$ means the number of vertical edges of the $K$-$K$ connection
${}_Kw_i{}_K$ connecting the vertices $x$ and $y$. 

Now we can get the dual principal graph from the fusion rule of 
$(3)$~$A$-$K$ $\times$ $K$-$K$ $\longrightarrow$ $A$-$K$.
It is the connected component of the fusion graph of (3) which 
contains the connection ${}_Ax_K$.

\vspace{3mm}
\noindent
%%%%%%%%%%%%%%%%%%%%%%%%%%%%%%%%%%%%%%%%%%%%%%%%%%%%%%%%%%%%%
%%%%% section 3.1.3 %%%%%%%%%%%%%%%%%%%%%%%%%%%%%%%%%%%%%%%%%
%%%%%%%%%%%%%%%%%%%%%%%%%%%%%%%%%%%%%%%%%%%%%%%%%%%%%%%%%%%%%
{\sl 3.1.3. The fusion rules of 
~$(4)$~$K$-$K$ $\times$ $K$-$K$ $\longrightarrow$ $K$-$K$}
\medskip

\vspace{2mm}
This is the fusion rule of the system of all $K$-$K$ connections 
obtained by Ocneanu 
(Figures \ref{qsym-a}`\ref{qsym-e8}, 
\cite[section 5.3, pages 244--252]{G}). 
It is isomorphic to the fusion algebras of the center of 
$K$-$K$ double triangle algebra $({}_KZ_K,\cdot)$
with dot product (vertical product) ``$\cdot$''.
We know that this fusion algebra $({}_KZ_K,\cdot)$ 
is generated by chiral left part and chiral right part 
which are isomorphic to the fusion algebra of connections 
arising from corresponding $ADE$ subfactor 
and that the chiral left and right part are relatively commutative
\cite[Theorem 5.16]{G}.
So we can compute the fusion rule of $({}_KZ_K,\cdot)$ 
from the above facts. 

We remark that the commutativity of the chiral left and right part 
is proved at the same time when we draw the diagrams of all $K$-$K$ 
connections (Figures \ref{qsym-a}`\ref{qsym-e8}). 
The proof is based on coset decomposition, fusion rules of chiral
left (right) part and indices of irreducible connections.
We refer readers to \cite[section 5.3, pages 244--252]{G} for details.

\begin{figure}[H]
\centering
\includegraphics[width=112mm,clip]
{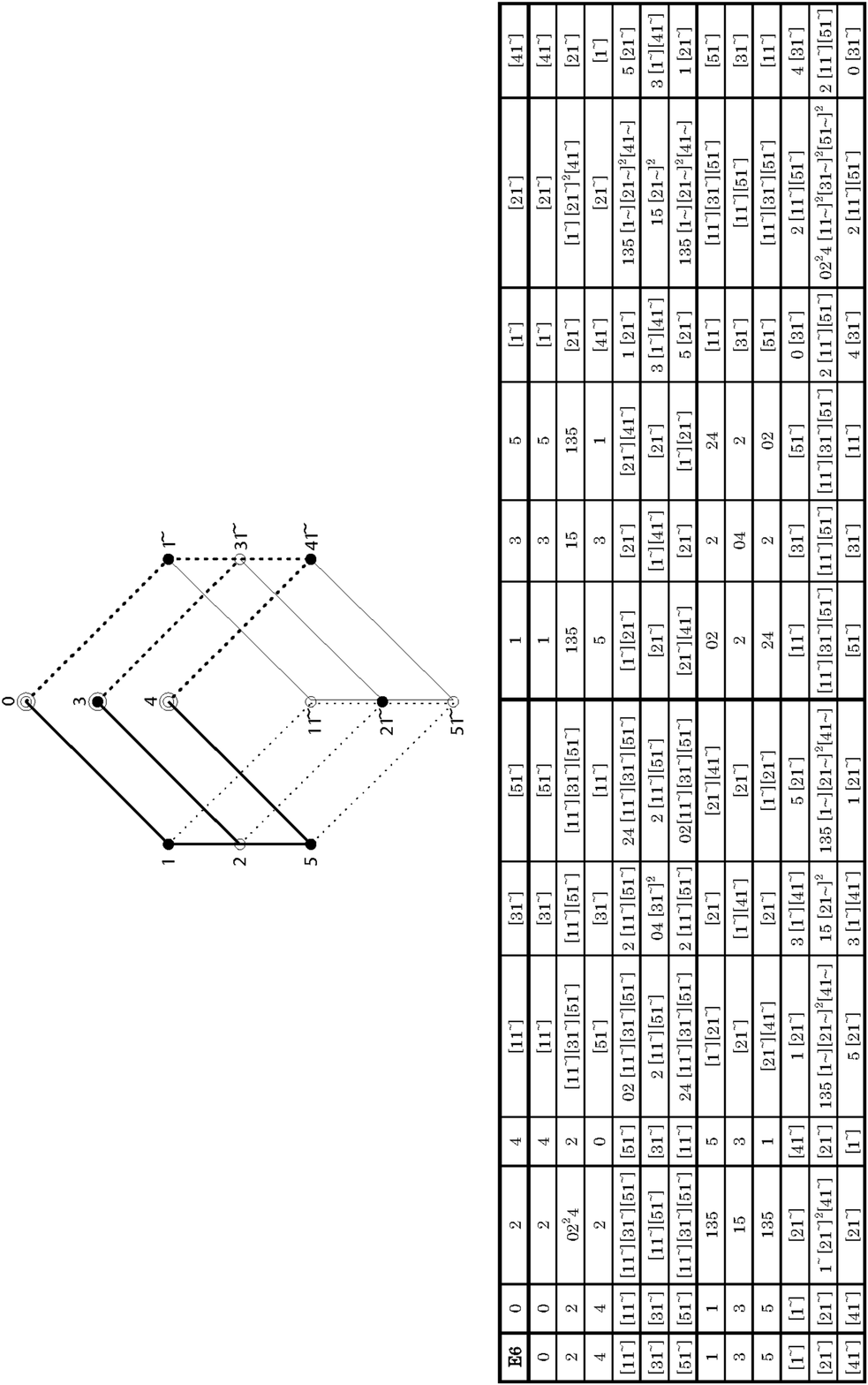}
\caption{The fusion table of $E_6$-$E_6$ connections}
\label{Fusion-Table-E6}
\end{figure}

\begin{figure}[H]
\centering
\includegraphics[width=107mm,clip]
{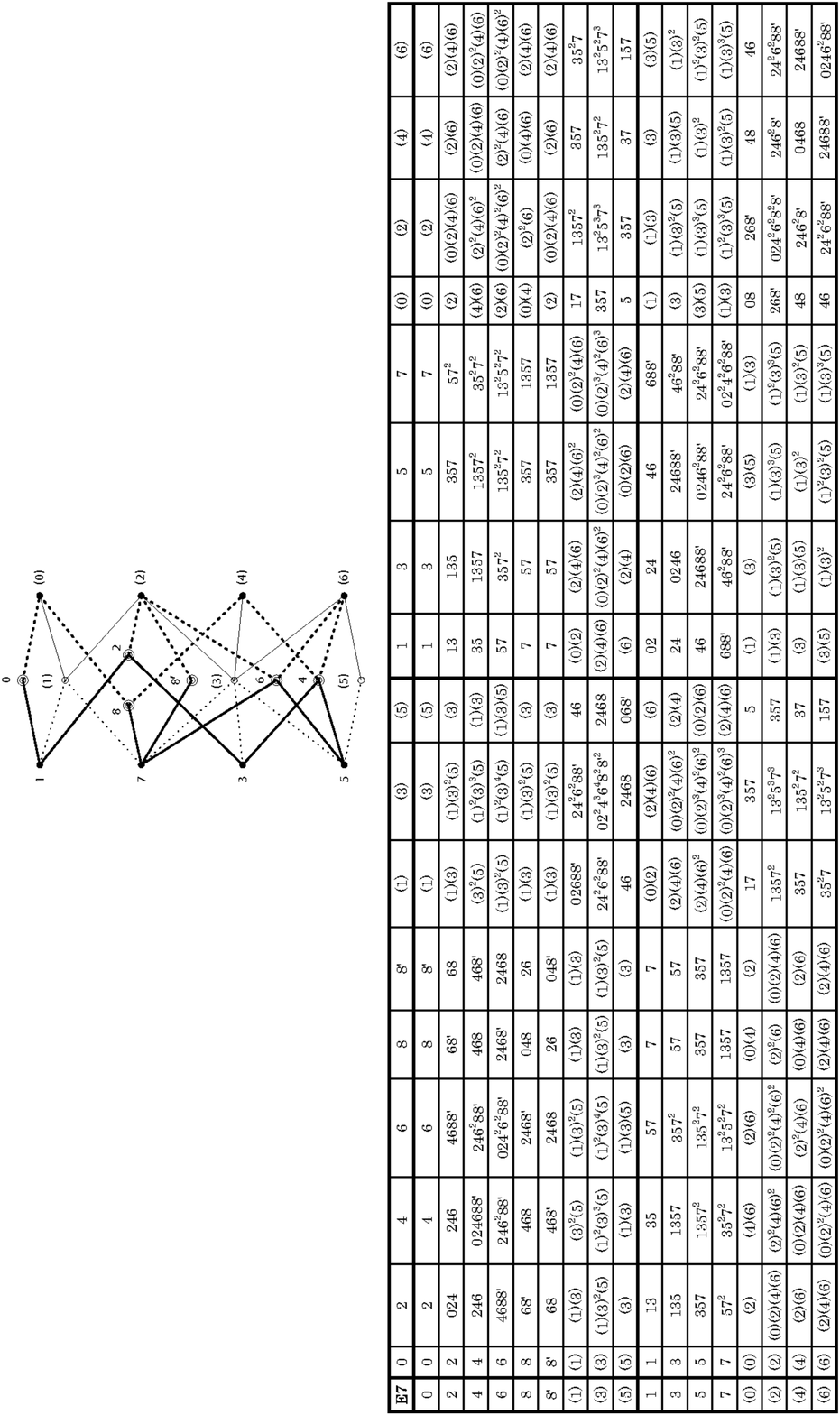}
\caption{The fusion table of $E_7$-$E_7$ connections}
\label{Fusion-Table-E7}
\end{figure}

\begin{figure}[H]
\centering
\includegraphics[width=110mm,clip]
{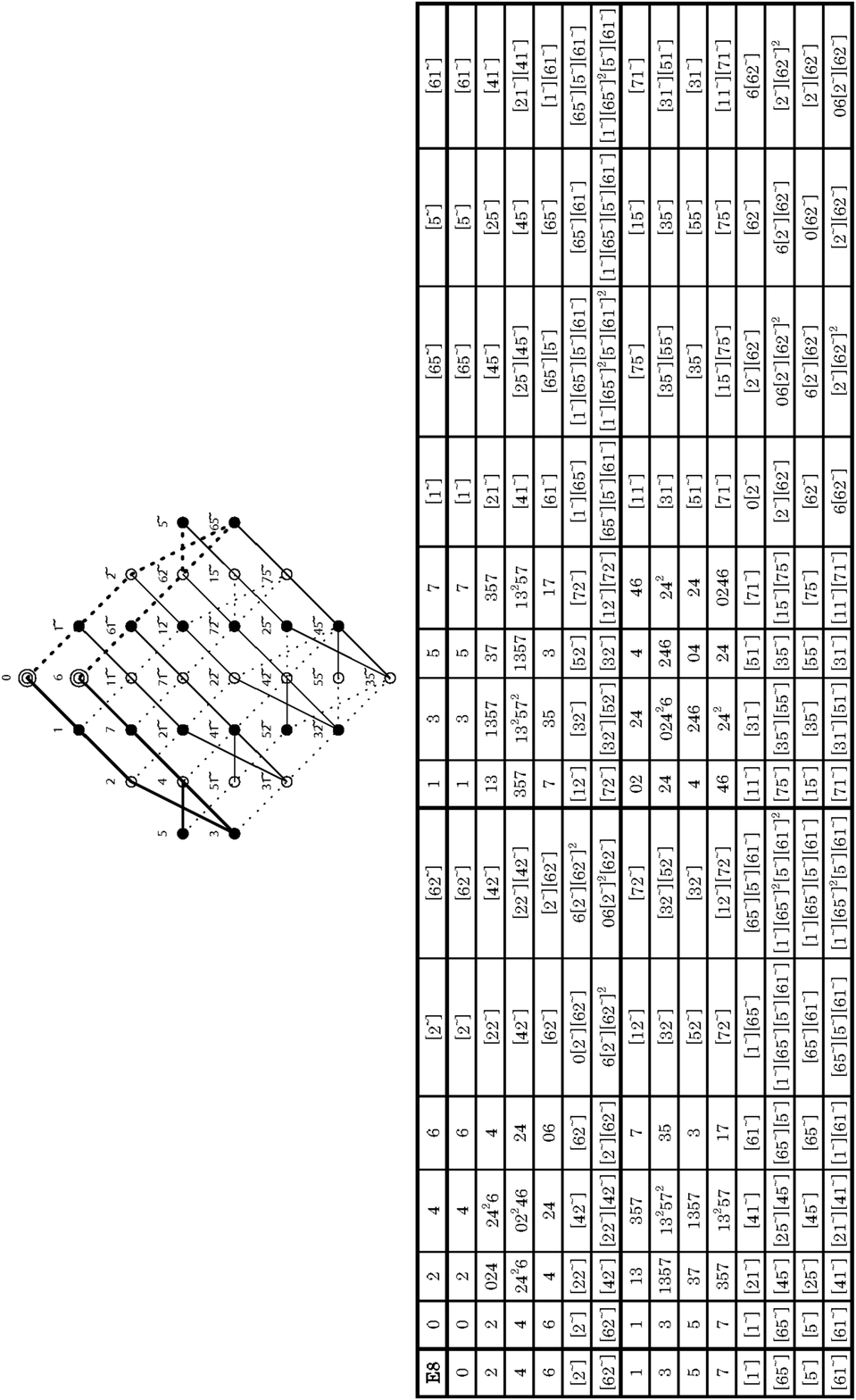}
\caption{A part of the fusion table of $E_8$-$E_8$ connections}
\label{Fusion-Table-E8}
\end{figure}

The fusion tables of $({}_KZ_K,\cdot)$, i.e. the system of all $K$-$K$ 
connections for $K$ = $E_6$, $E_7$ and (a part of) $E_8$ is given in
Figures \ref{Fusion-Table-E6}`\ref{Fusion-Table-E8}.
We note that these fusion tables is expressed in product form.
For example in the table \ref{Fusion-Table-E7}, 
we can read $(3)\cdot 4$=$(1)^2(3)^3(5)$, which means the fusion rule
$w_{(3)}\cdot w_4$=$2w_{(1)}+3w_{(3)}+w_{(5)}$ holds.

\vspace{3mm}
\noindent
%%%%%%%%%%%%%%%%%%%%%%%%%%%%%%%%%%%%%%%%%%%%%%%%%%%%%%%%%%%%%
%%%%% section 3.2 %%%%%%%%%%%%%%%%%%%%%%%%%%%%%%%%%%%%%%%%%%%
%%%%%%%%%%%%%%%%%%%%%%%%%%%%%%%%%%%%%%%%%%%%%%%%%%%%%%%%%%%%%
{\sl 3.2. The fusion rules of even vertices of 
the (dual) principal graphs of $\GHJ(K,*_K=x)$ }
\medskip

\vspace{2mm}
Let $N\subset M$ be the Goodman-de la Harpe-Jones subfactor $\GHJ(K,*_K=x)$.
Here we will compute the fusion rules of even vertices of 
the (dual) principal graphs of $\GHJ(K,*_K=x)$, that is, 
the fusion rules of $N$-$N$ bimodules and $M$-$M$ bimodules of 
the subfactor $N\subset M$.

The system of $N$-$N$ bimodules are isomorphic to the system of 
$A$-$A$ connections generated by ${}_Ax_K$ and this is the same
as ${}_AZ^{\rm even}_A$, i.e. the fusion algebra of even part of ${}_AZ_A$.
So the fusion algebra of $N$-$N$ bimodules are isomorphic to 
the fusion algebra $A^{\rm even}$, i.e. the fusion algebra of even vertices 
of the Jones' type $A$ subfactor.
Hence it turns out that the fusion algebra of $N$-$N$ bimodules are always 
commutative for any GHJ subfactors.

The system of $M$-$M$ bimodules are similarly isomorphic to the system of 
$K$-$K$ connections generated by ${}_Ax_K$ and this is the same
as (a part of) ${}_KZ^{\rm even}_K$, 
i.e. the fusion algebra of even part of ${}_KZ_K$. 
So we have only to compute the fusion rule of ${}_KZ^{\rm even}_K$.

%\vspace{2mm}
Here the fusion rule of ${}_KZ^{\rm even}_K$ and the vertical edges
of irreducible $K$-$K$ connections can be summarized as 
in the Table \ref{fusion-KZK}.
As we mentioned above, we can compute the fusion rule of ${}_KZ^{\rm even}_K$
in detail from the fusion graph of all $K$-$K$ connections 
as in Figures \ref{qsym-a}`\ref{qsym-e8} and 
\ref{Fusion-Table-E6}`\ref{Fusion-Table-E8}.

\begin{figure}[H]
\centering\includegraphics[width=130mm,clip]{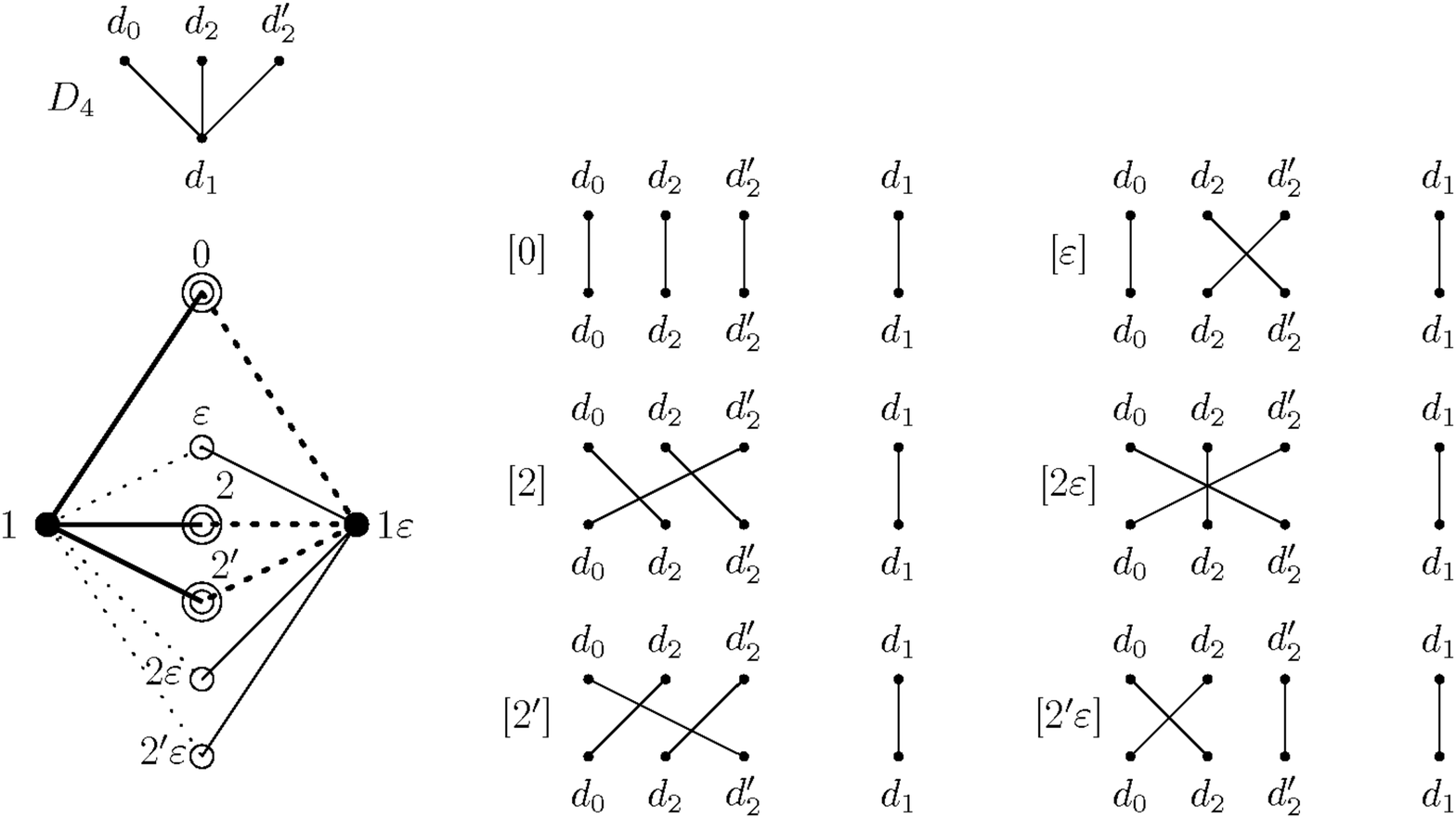}
\caption{Vertical graphs for connections on 
the Coxeter graph $D_4$}
\label{v-edge-D4}
\end{figure}

In the following table, $\e$ represents the index 1 $D_{2n}$-$D_{2n}$
connection which corresponds to the flip of two tails of $D_{2n}$.
Because ${}_{D_{2n}}Z_{D_{2n}}$ has coset decomposition 
$D_{2n} \cup D_{2n}\cdot \e$ and $\e^2=id$ 
as shown in Figures \ref{qsym-de}, \ref{v-edge-D4} and \ref{v-edge-D6}, 
we can easily compute the fusion rule for ${}_{D_{2n}}Z_{D_{2n}}$. 

\vspace{-1mm}
\begin{table}[H]
\begin{center}
\begin{tabular}{@{\vrule width 1pt\ }c|c|c@{\ \vrule width 1pt}}
\noalign{\hrule height 1pt}
{Graph $K$} & { fusion rule of ${}_KZ^{\rm even}_K$} & { vertical edges 
of $K$-$K$ connections  }\\ 
\noalign{\hrule height 1pt}
$A_n$ & commutative & $\EssPath A_n$~
(Figures \ref{v-edge-A3}`\ref{v-edge-A6}) \\ \hline
$D_{2n}$ & non-commutative & $\EssPath D_{2n} + \e$~
(Figures \ref{v-edge-D4} and \ref{v-edge-D6}) \\ \hline
$D_{2n+1}$ & commutative & $\EssPath D_{2n+1}$~
(Figure \ref{v-edge-D5}) \\ \hline
$E_6$ & commutative  & Figure \ref{v-edge-E6}  \\ \hline
$E_7$ & commutative & Figure \ref{v-edge-E7}  \\ \hline
$E_8$ & commutative & Figures \ref{v-edge-E8even},~\ref{v-edge-E8odd} \\ \hline
\noalign{\hrule height 1pt}
\end{tabular}
\end{center}

\vspace{-4mm}
\caption{The fusion rule of ${}_KZ^{\rm even}_K$ and
vertical edges of $K$-$K$ connections }
\label{fusion-KZK}
\end{table}

\vspace{-7mm}
%%%%%%%%%%%%%%%%%%%%%%%%%%%%%%%%%%%%%%%%%%%%%%%%%%%%%%%%%%%%%%%%%%%%%%
\newpage
%%%%%%%%%%%%%%%%%%%%%%%%%%%%%%%%%%%%%%%%%%%%%%%%%%%%%%%%%%%%%%%%%%%%%%
\section{The structure of Goodman-de la Harpe-Jones subfactors}
%%%%%%%%%%%%%%%%%%%%%%%%%%%%%%%%%%%%%%%%%%%%%%%%%%%%%%%%%%%%%%%%%%%%%%

\noindent
%%%%%%%%%%%%%%%%%%%%%%%%%%%%%%%%%%%%%%%%%%%%%%%%%%%%%%%%%%%%%
%%%%% section 4.1 %%%%%%%%%%%%%%%%%%%%%%%%%%%%%%%%%%%%%%%%%%%
%%%%%%%%%%%%%%%%%%%%%%%%%%%%%%%%%%%%%%%%%%%%%%%%%%%%%%%%%%%%%
{\sl 4.1. Goodman-de la Harpe-Jones subfactors of type $A_n$}
\medskip

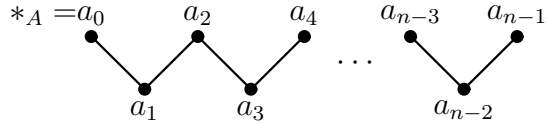
\begin{figure}[H]
%\unitlength 0.85mm
\unitlength 0.7mm
\thicklines
%\thinlines
\begin{center}
\begin{picture}(80,18)(0,-5)
\multiput(10,0)(20,0){2}{\circle*{2}}
\multiput(0,10)(20,0){5}{\circle*{2}}
\put(70,0){\circle*{2}}
\path(0,10)(10,0)(20,10)(30,0)(40,10)
\path(60,10)(70,0)(80,10) 
\put(-9,14){\makebox(0,0){$*_A=$}}
\put(0,14){\makebox(0,0){$a_0$}}
\put(10,-4){\makebox(0,0){$a_1$}}
\put(20,14){\makebox(0,0){$a_2$}}
\put(30,-4){\makebox(0,0){$a_3$}}
\put(40,14){\makebox(0,0){$a_4$}}
\put(50,5){\makebox(0,0){$\cdots$}}
\put(60,14){\makebox(0,0){$a_{n-3}$}}
\put(70,-4){\makebox(0,0){$a_{n-2}$}}
\put(80,14){\makebox(0,0){$a_{n-1}$}}
%\put(-40,5){\makebox(0,0){Dynkin diagram $A_n$}}
%\put(80,14){\makebox(0,0){even}}
%\put(80,-4){\makebox(0,0){odd}}
\end{picture}
\caption{The label of vertices of the Dynkin diagram $A_n$.}
\label{label-An}
\end{center}
\end{figure}

Let $N\subset M$ be the Jones' subfactor of type $A_n$ and 
$N\subset M\subset M_1\subset M_2\subset \cdots\subset M_k\subset$
be the Jones tower. 
We label the vertices of the Dynkin diagram $A_n$ 
by $a_0,a_1,\cdots, a_{n-1}$ as in Figure \ref{label-An}.
Then the Goodman-de la Harpe-Jones subfactor $\GHJ(A_n,*=a_m)$ 
is isomorphic to $pN\subset pM_{m-1}p$, where $p$ is 
a minimal projection in $\Proj(N'\cap M_{m-1})$ corresponding to
the vertex $a_m$.
Hence in this case the principal graph and the dual principal graph
coincide and fusion rule of even vertices of both graphs becomes
$A^{\rm even}_n$.

%%%%%%%%%%%%%%%%%%%%%%%%%%%%%%%%%%%%%%%%%%%%%%%%%%%%
%%%%%%%%%%%%%%%%%%%%%%%%%%%%%%%%%%%%%%%%%%%%%%%%%%%%

\vspace{-7mm}
\noindent
%%%%%%%%%%%%%%%%%%%%%%%%%%%%%%%%%%%%%%%%%%%%%%%%%%%%%%%%%%%%%
%%%%% section 4.2 %%%%%%%%%%%%%%%%%%%%%%%%%%%%%%%%%%%%%%%%%%%
%%%%%%%%%%%%%%%%%%%%%%%%%%%%%%%%%%%%%%%%%%%%%%%%%%%%%%%%%%%%%
{\sl 4.2. Goodman-de la Harpe-Jones subfactors of type $D_{2n+1}$}
\medskip

We label the vertices of the Dynkin diagram $D_{2n+1}$ 
by $d_0,d_1,d_2,\cdots,d_{2n-2}$, \\
$d_{2n-1},d_{2n-1}'$
as in Figure \ref{label-D2n+1}.

The Goodman-de la Harpe-Jones subfactor $\GHJ(D_{2n+1},*_K=d_0)$ 
is isomorphic to the unique index 2 subfactor 
$N\subset N\rtimes\Z_2$. 

\begin{figure}[H]
%\centering\includegraphics[width=140mm,clip]{qsym-do.eps}
\centering\includegraphics[height=83mm,clip]{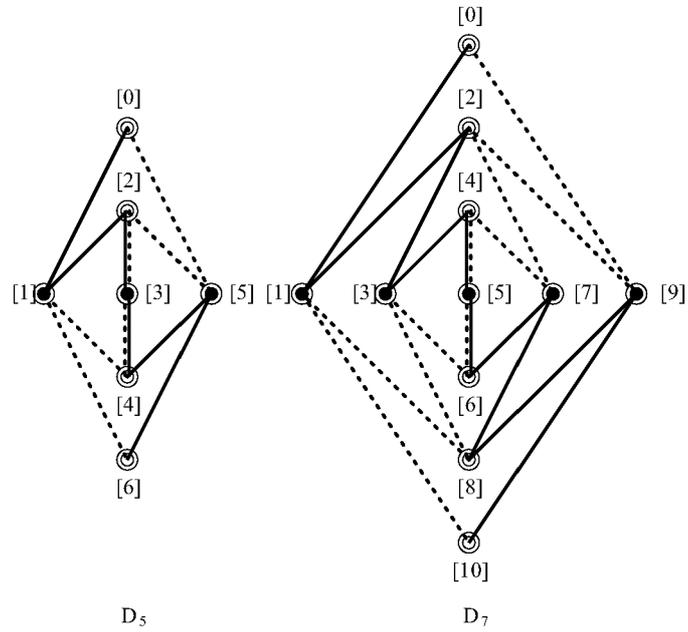}
\caption{Chiral symmetry for the Coxeter graph $D_{odd}$}
\label{qsym-do}
\end{figure}

\begin{figure}[H]
\centering\includegraphics[width=125mm,clip]{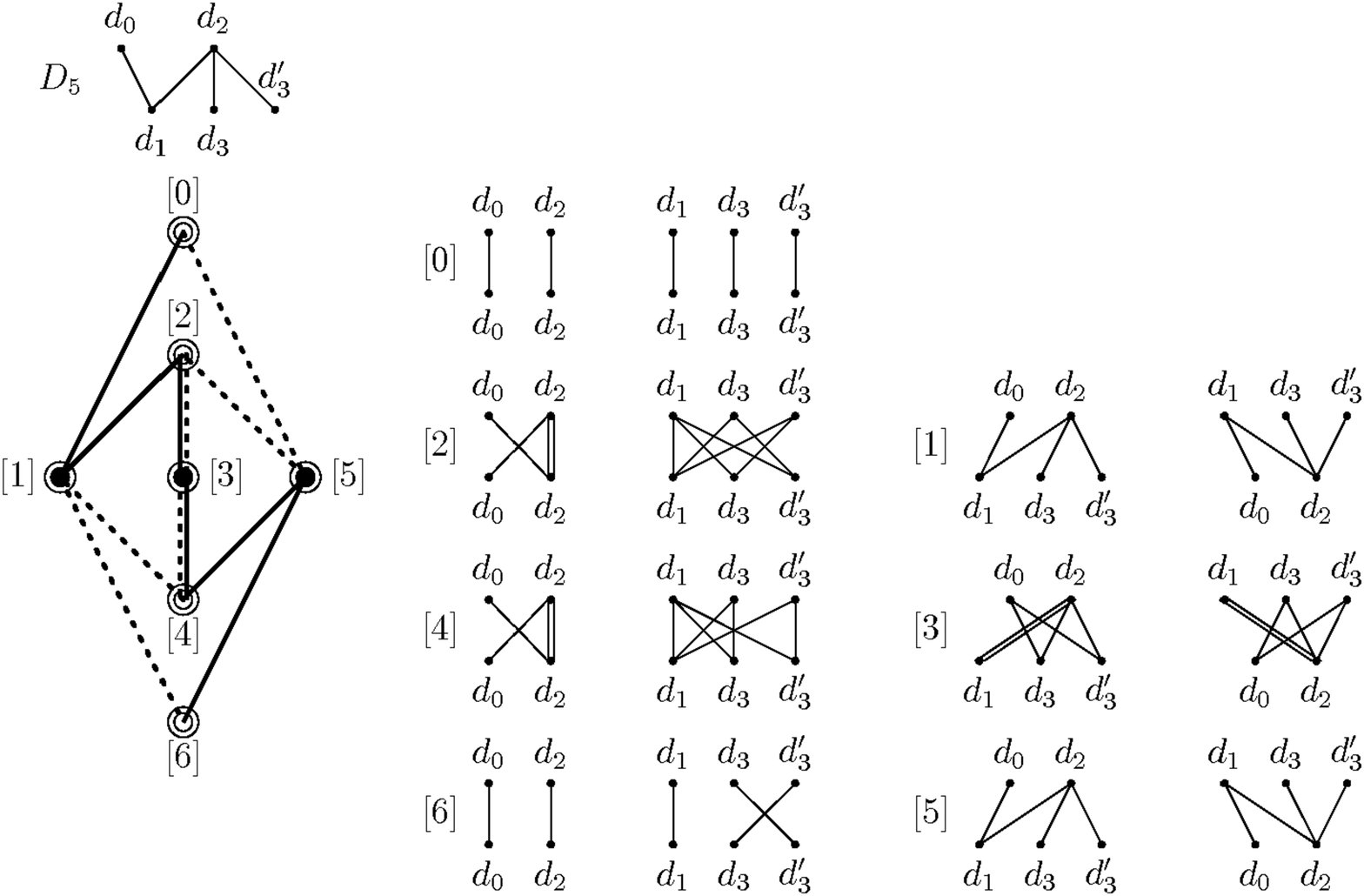}
\caption{Vertical graphs for connections on 
the Coxeter graph $D_5$}
\label{v-edge-D5}
\end{figure}

If the vertex $*_K\neq d_0,d_{2n-1},d_{2n-1}'$, 
$\GHJ(D_{2n+1},*_K)$ has nontrivial intermediate subfactor 
as in Figure \ref{GHJ-D2n+1} 
because we have the decomposition of connections 
${}_Ad_k{}_D={}_Ad_0{}_D\cdot {}_D[k]{}_D$ for $k=1,2,\dots, 2n-2$. 
Here ${}_D[k]{}_D$ is the $D_{2n+1}$-$D_{2n+1}$ connection corresponding
to the vertex $[k]$ as in Figures \ref{qsym-do} and \ref{v-edge-D5}.

\begin{figure}[H]
\centering
\includegraphics[width=140mm,clip]{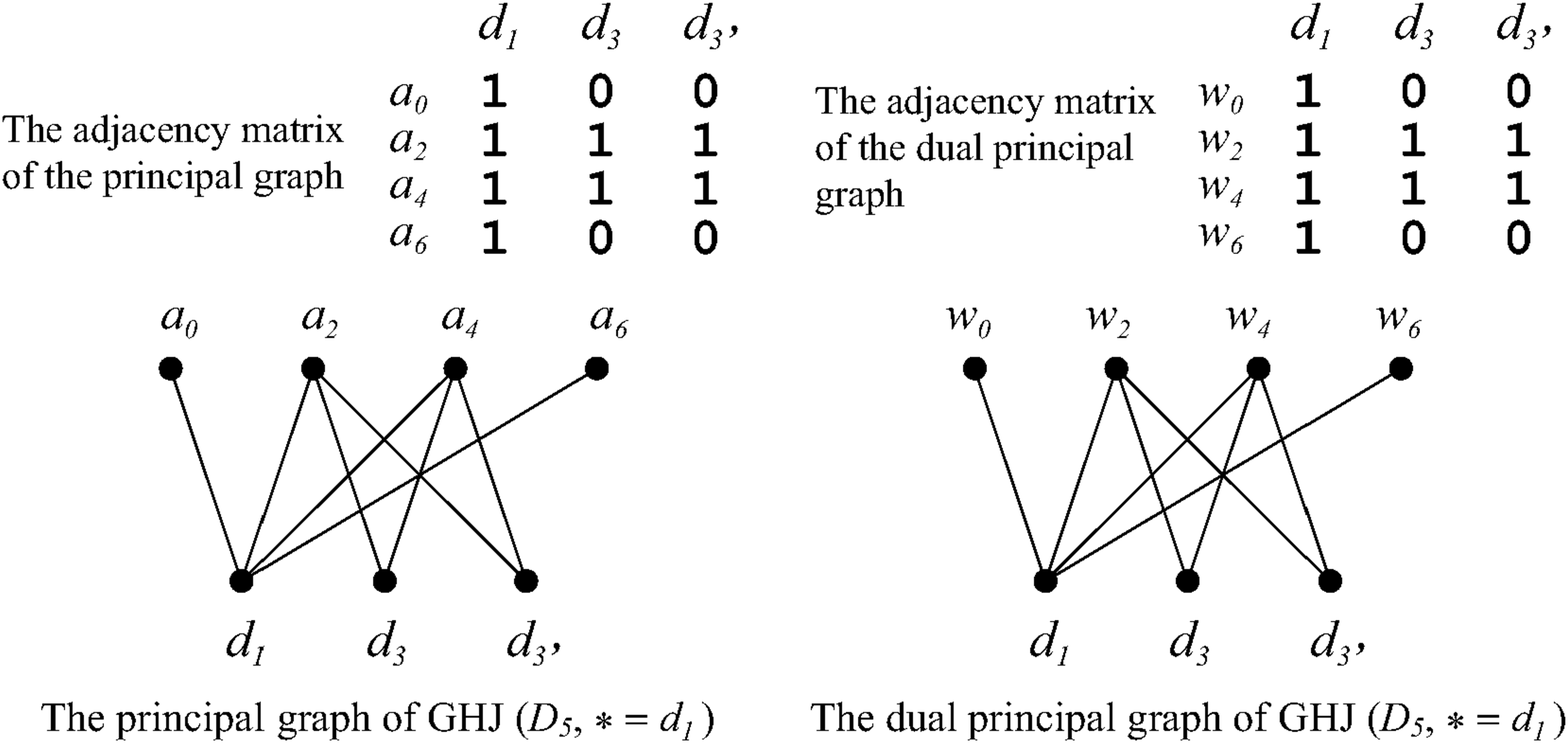}
\caption{The (dual) principal graph of GHJ$(D_5, *=d_1)$.}
\label{GHJ(D5-d1)}
\end{figure}

\begin{figure}[H]
\centering
\includegraphics[width=130mm,clip]{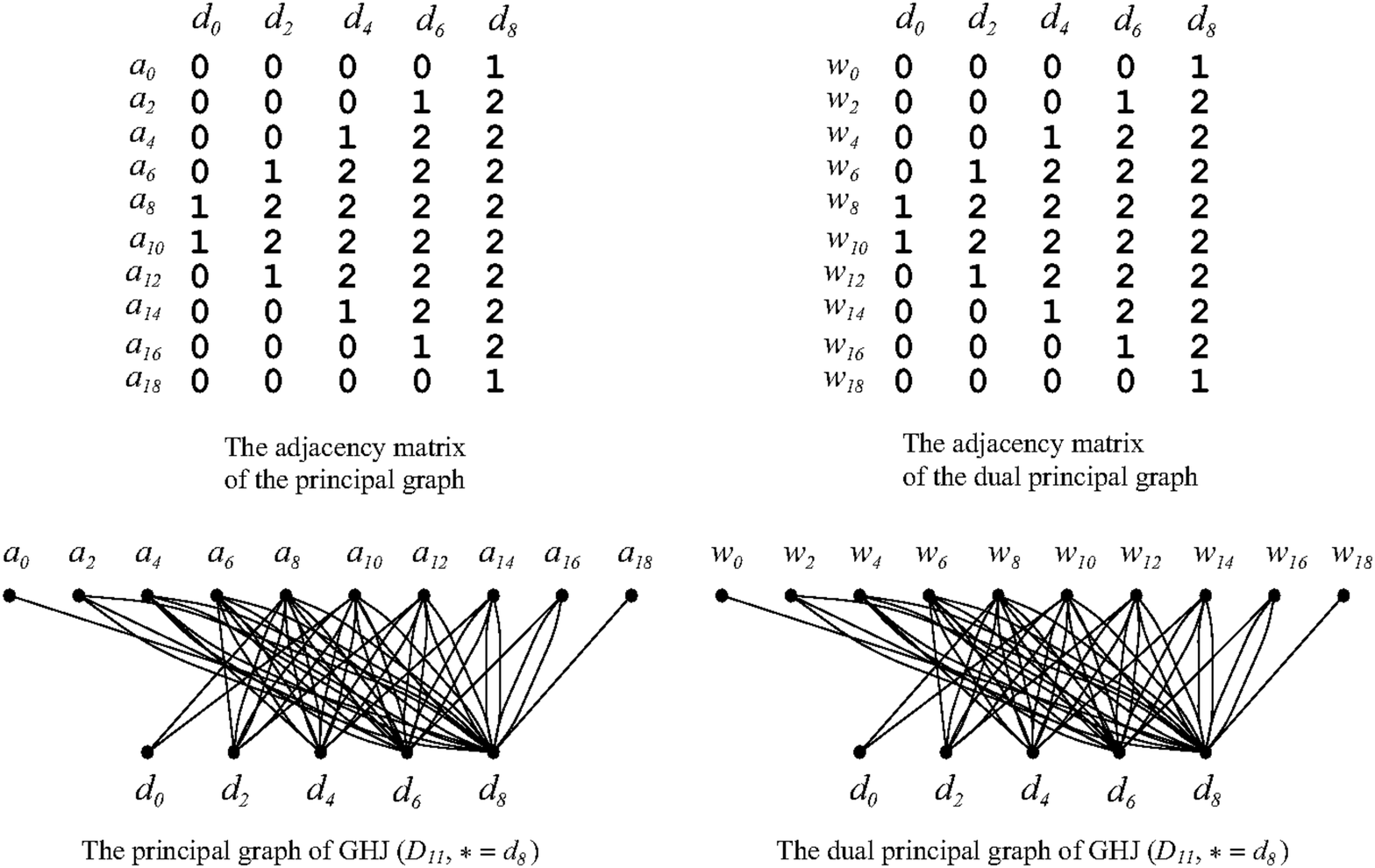}
\caption{The (dual) principal graph of GHJ$(D_{11}, *=d_8)$.}
\label{GHJ(D11-d8)}
\end{figure}

The (dual) principal graphs of $\GHJ(D_{2n+1},*_K)$ are given 
in Figures \ref{GHJ(D5-d1)}`\ref{GHJ(D11-d8)} for $n=2,3,4,5$.

\begin{figure}[H]
\centering
\includegraphics[width=150mm,clip]{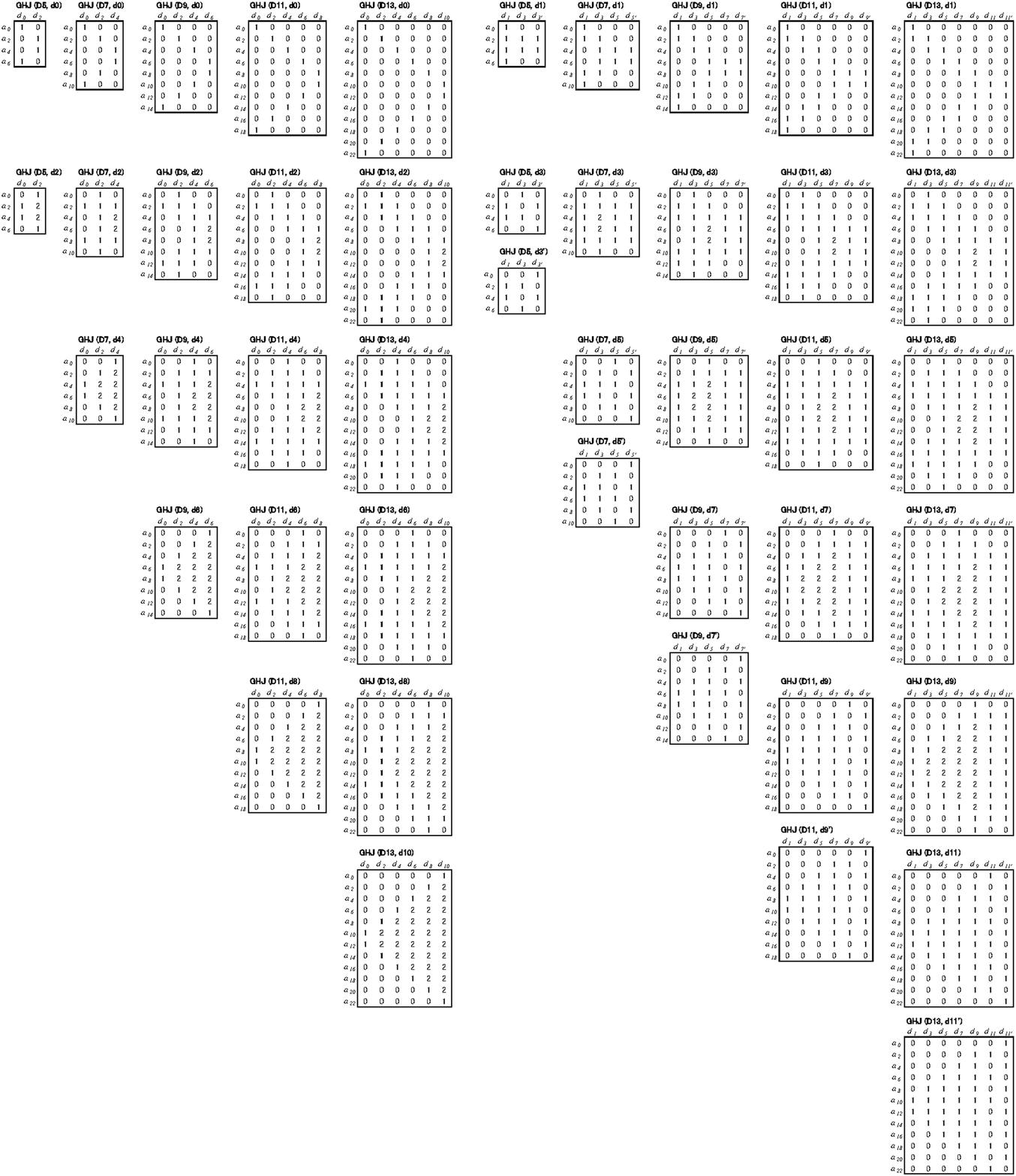}
\caption{The incidence matrices of 
the (dual) principal graphs of GHJ$(D_{odd})$.}
\label{GHJ(D_odd)-PG=DPG}
\end{figure}

The incidence matrices of the (dual) principal graphs of 
GHJ$(D_{odd},*_K)$ are also given 
in Figure \ref{GHJ(D_odd)-PG=DPG}.

%If the vertex $*_D$ is not among $d_0$, $d_{2n-1}$, and $d_{2n-1}'$, 
%If the vertex $*_D$ is except $d_0$, $d_{2n-1}$, and $d_{2n-1}'$,

\vspace{-1mm}
%%%%%%%%%%%%%%%%%%%%
%%    D_{2n+1}    %%
%%%%%%%%%%%%%%%%%%%%
\begin{figure}[H]
\unitlength 0.75mm
\thicklines
\begin{center}
\begin{picture}(90,20)(0,-5)
\multiput(10,0)(20,0){2}{\circle*{2}}
\multiput(0,10)(20,0){5}{\circle*{2}}
\multiput(70,0)(10,0){3}{\circle*{2}}
\path(0,10)(10,0)(20,10)(30,0)(40,10)
\path(60,10)(70,0)(80,10)(80,0) 
\path(80,10)(90,0)
%\put(-9,14){\makebox(0,0){$*_A=$}}
\put(0,14){\makebox(0,0){$d_0$}}
\put(10,-4){\makebox(0,0){$d_1$}}
\put(20,14){\makebox(0,0){$d_2$}}
\put(30,-4){\makebox(0,0){$d_3$}}
\put(40,14){\makebox(0,0){$d_4$}}
\put(50,5){\makebox(0,0){$\cdots$}}
\put(60,14){\makebox(0,0){$d_{2n-4}$}}
\put(69,-5){\makebox(0,0){$d_{2n-3}$}}
\put(80,14){\makebox(0,0){$d_{2n-2}$}}
\put(82,-5){\makebox(0,0){$d_{2n-1}$}}
\put(95,-5){\makebox(0,0){$d_{2n-1}'$}}
\end{picture}
\caption{The label of vertices of the Dynkin diagram $D_{2n+1}$.}
\label{label-D2n+1}
\end{center}
%%%%%%%%%%%%%%%%%%%%
\begin{center}
\unitlength 0.75mm
\begin{picture}(80,80)(30,0)
\path(-20,-3)(-20,80)(150,80)(150,-3)(-20,-3)
\drawline(3,70)(20,70)(20,30)(0,30)(0,67)
\drawline(0,50)(20,50)
\multiput(0,30)(0,20){2}{\circle*{1}}
\put(0,70){\makebox(0,0){$*$}}
\put(10,74){\makebox(0,0){$A_{4n-1}$}}
\put(10,54){\makebox(0,0){$D_{2n+1}$}}
\put(10,26){\makebox(0,0){$D_{2n+1}$}}
\put(-5,50){\makebox(0,0){$d_0$}}
\put(-5,30){\makebox(0,0){$d_k$}}
\put(40,70){\makebox(0,0){$\cdots \quad \longrightarrow$}}
\put(40,50){\makebox(0,0){$\cdots \quad \longrightarrow$}}
\put(40,30){\makebox(0,0){$\cdots \quad \longrightarrow$}}
\put(60,70){\makebox(0,0){$N$}}
\put(60,60){\makebox(0,0){$\cap$}}
\put(59,50){\makebox(0,0){${}^{\exists} P$}}
\put(60,40){\makebox(0,0){$\cap$}}
\put(60,30){\makebox(0,0){$M$}}
\put(110,63){\makebox(0,0){$(N\subset P) \cong (N\subset N\rtimes\Z_2) $}}
\put(110,57){\makebox(0,0){index~=~2}}
\put(110,40){\makebox(0,0){$P\neq M$}}
\put(50,15){\makebox(0,0){The principal graph}}
\put(80,15){\makebox(0,0){$\cong$}}
\put(115,15){\makebox(0,0){The dual principal graph}}
\put(5,8){\makebox(0,0){The fusion rule}}
\put(6,3){\makebox(0,0){of even vertices}}
\put(60,5){\makebox(0,0){$A_{4n-1}^{\rm even}$}}
\put(80,5){\makebox(0,0){=}}
\put(102,5){\makebox(0,0){$A_{4n-1}^{\rm even}$}}
\end{picture}
\caption{ }
\label{GHJ-D2n+1}
\end{center}
\end{figure}
%%%%%%%%%%%%%%%%%%%%  

%\vspace{3mm}
\noindent
%%%%%%%%%%%%%%%%%%%%%%%%%%%%%%%%%%%%%%%%%%%%%%%%%%%%%%%%%%%%%
%%%%% section 4.3 %%%%%%%%%%%%%%%%%%%%%%%%%%%%%%%%%%%%%%%%%%%
%%%%%%%%%%%%%%%%%%%%%%%%%%%%%%%%%%%%%%%%%%%%%%%%%%%%%%%%%%%%%
{\sl 4.3. Goodman-de la Harpe-Jones subfactors of type $D_{2n}$}
\medskip

We label the vertices of the Dynkin diagram $D_{2n}$ 
by $d_0,d_1,d_2,\cdots,d_{2n-3}$, \\
$d_{2n-2},d_{2n-2}'$
as in Figure \ref{label-D2n}.

The Goodman-de la Harpe-Jones subfactor $\GHJ(D_{2n},*_K=d_0)$,  
$\GHJ(D_4,*_K=d_2)$ and $\GHJ(D_4,*_K=d_2')$ are isomorphic to 
the unique index 2 subfactor $N\subset N\rtimes\Z_2$. 

\begin{figure}[H]
%\centering\includegraphics[width=140mm,clip]{qsym-de.eps}
\centering\includegraphics[height=115mm,clip]{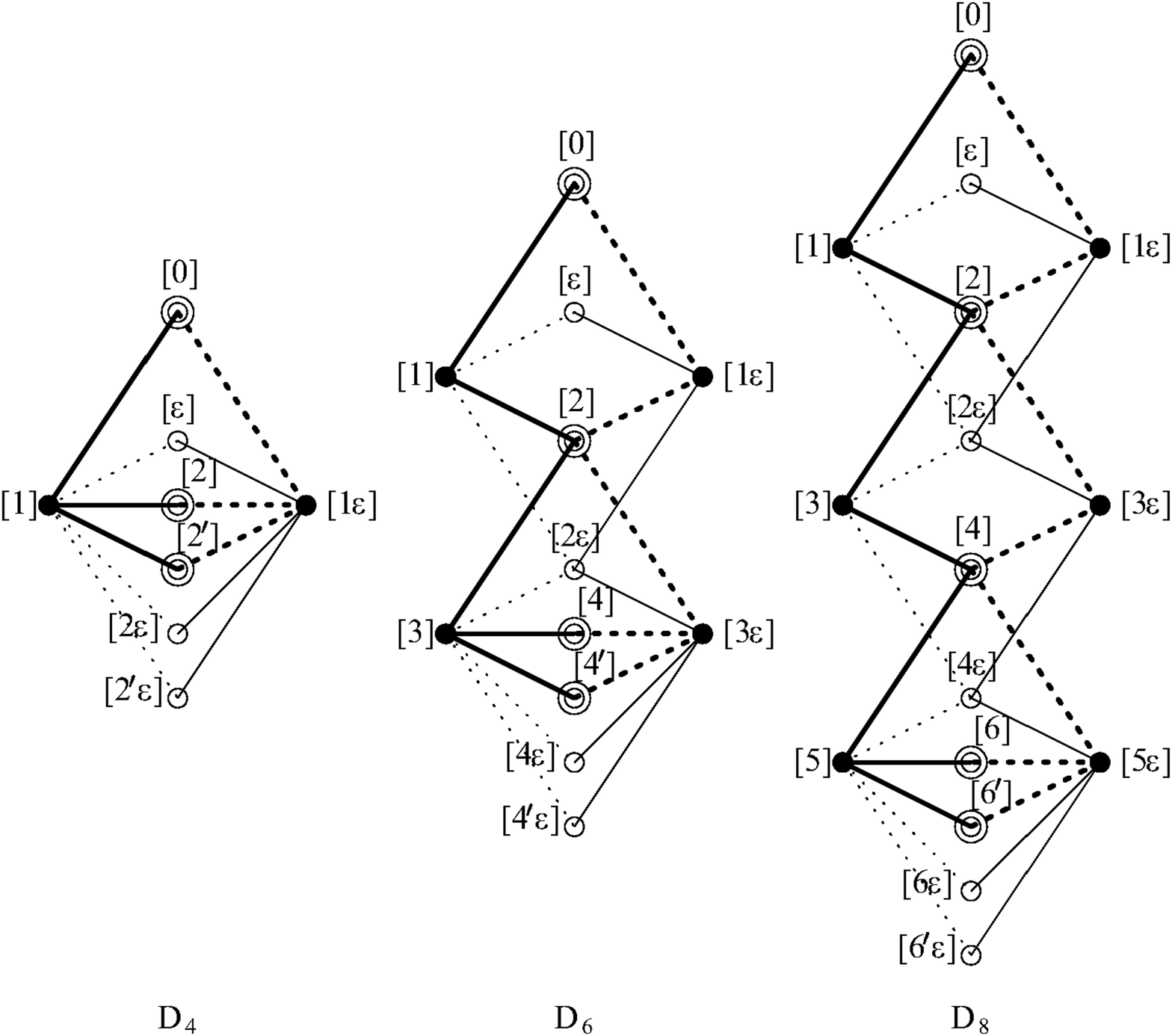}
\caption{Chiral symmetry for the Coxeter graph $D_{even}$}
\label{qsym-de}
\end{figure}

\begin{figure}[H]
\centering\includegraphics[width=125mm,clip]{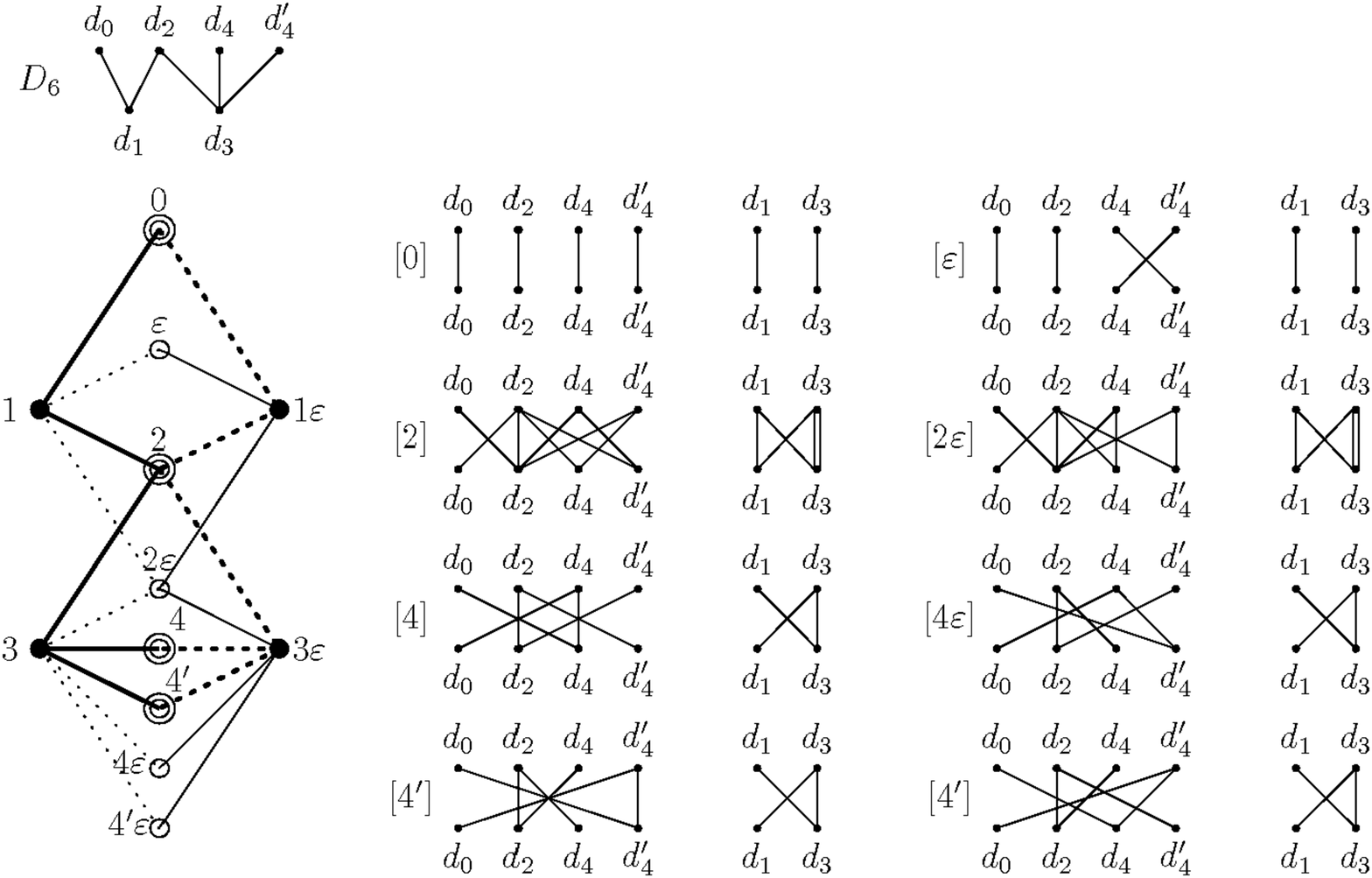}
\caption{Vertical graphs for connections on 
the Coxeter graph $D_6$}
\label{v-edge-D6}
\end{figure}

If $n>2$ and the vertex $*_K\neq d_0$, 
$\GHJ(D_{2n},*_K)$ has nontrivial intermediate subfactor 
as in  Figure \ref{GHJ-D2n+1}
because we have the decomposition of connections 
${}_Ad_k{}_D={}_Ad_0{}_D\cdot {}_D[k]{}_D$ for $k\neq 0$. 
Here ${}_D[k]{}_D$ is the $D_{2n}$-$D_{2n}$ connection corresponding
to the vertex $[k]$ as in Figures \ref{qsym-de} and \ref{v-edge-D6}.

\begin{figure}[H]
\centering
\includegraphics[width=130mm,clip]{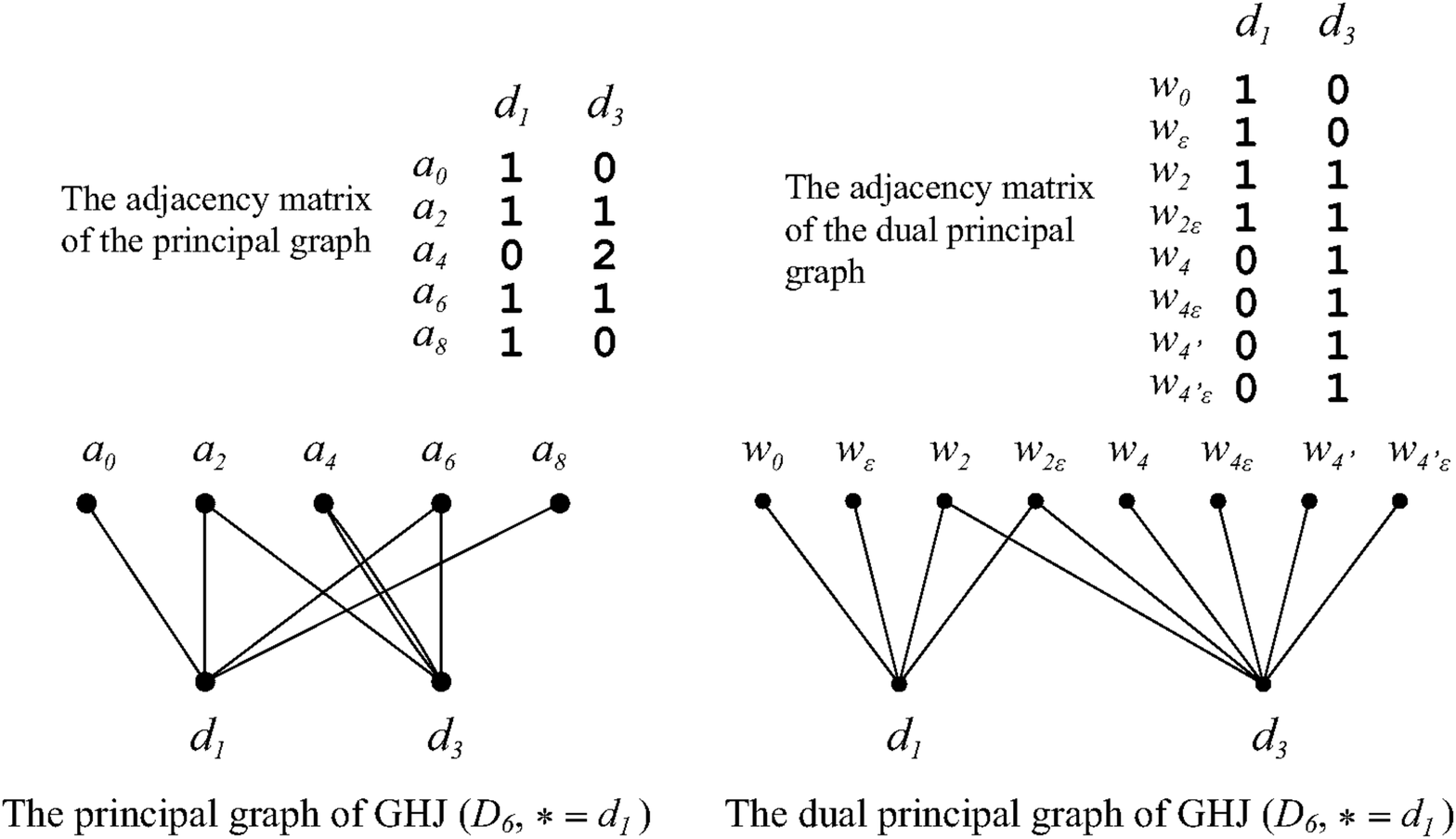}
\caption{The (dual) principal graph of GHJ$(D_6, *=d_1)$.}
\label{GHJ(D6-d1)}
\end{figure}

\begin{figure}[H]
\centering
\includegraphics[width=130mm,clip]{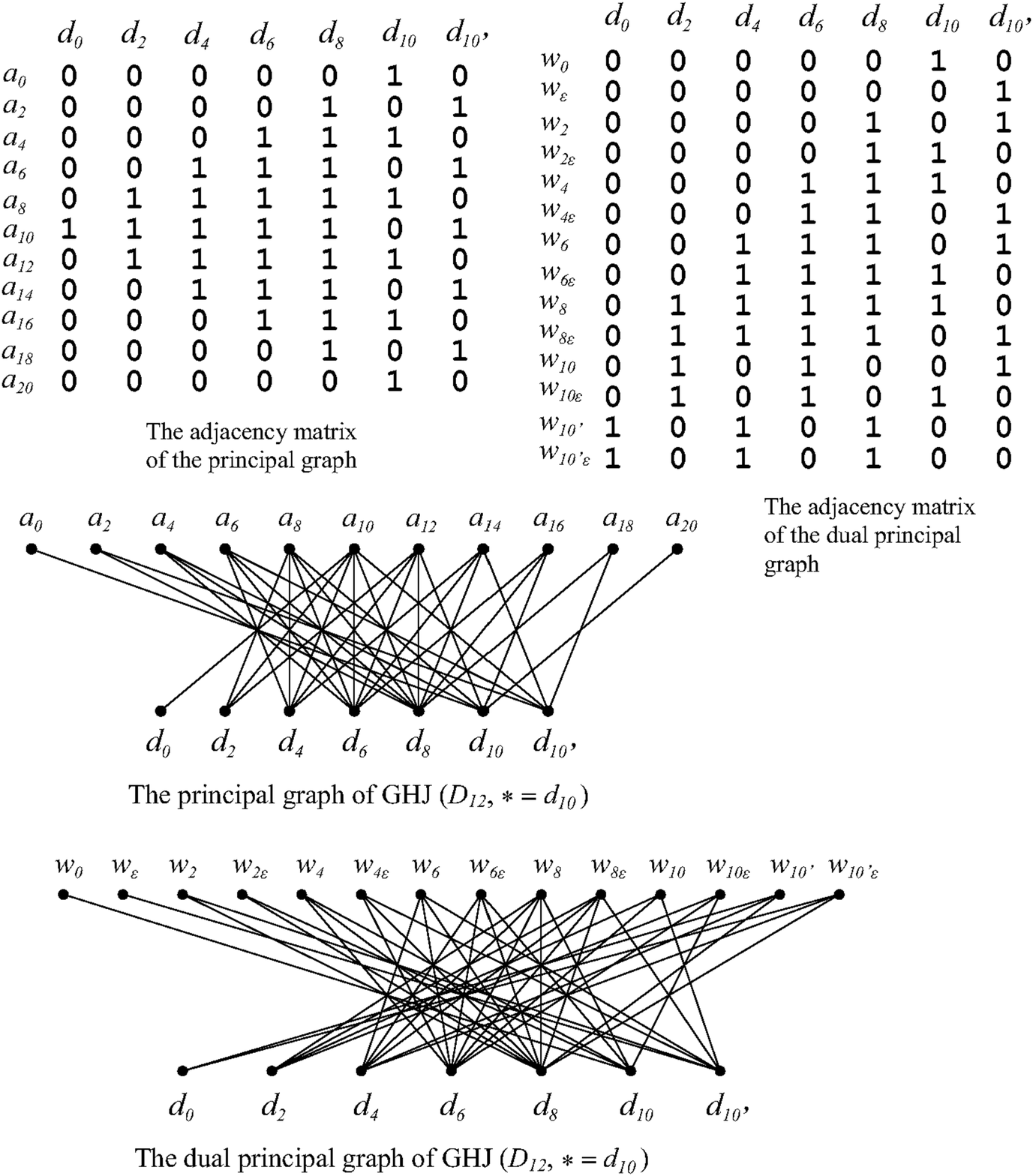}
\caption{The (dual) principal graph of GHJ$(D_{12}, *=d_{10})$.}
\label{GHJ(D12-d10)}
\end{figure}

The (dual) principal graphs of $\GHJ(D_{2n},*_K)$ are given 
in Figures \ref{GHJ(D6-d1)}`\ref{GHJ(D12-d10)} for $n=3,4,5,6$.

\begin{figure}[H]
\centering
\includegraphics[width=130mm,clip]{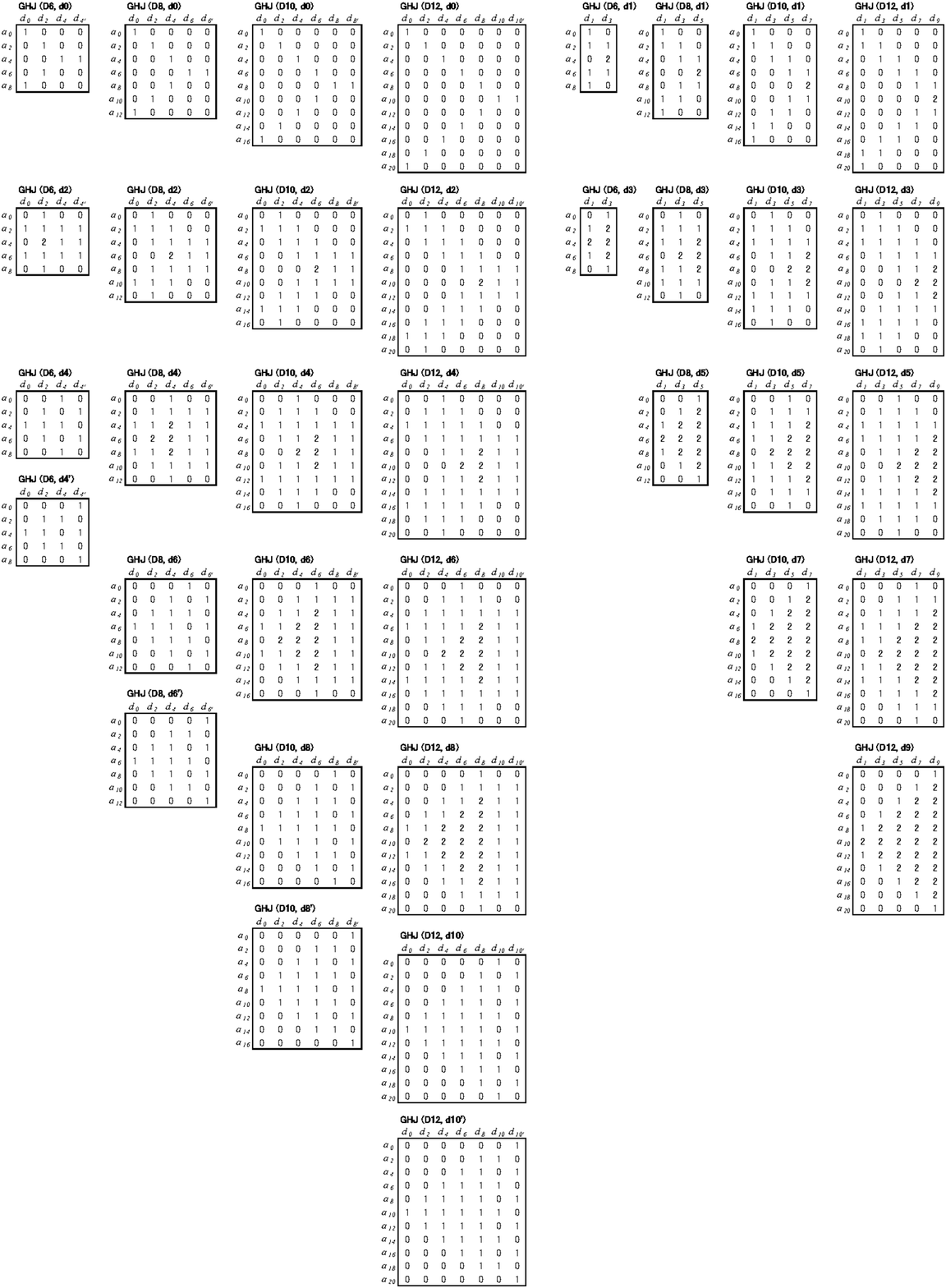}
\caption{The incidence matrices of 
the principal graphs of GHJ$(D_{even})$.}
\label{GHJ(D_even)-PG}
\end{figure}

\begin{figure}[H]
\centering
\includegraphics[width=105mm,clip]{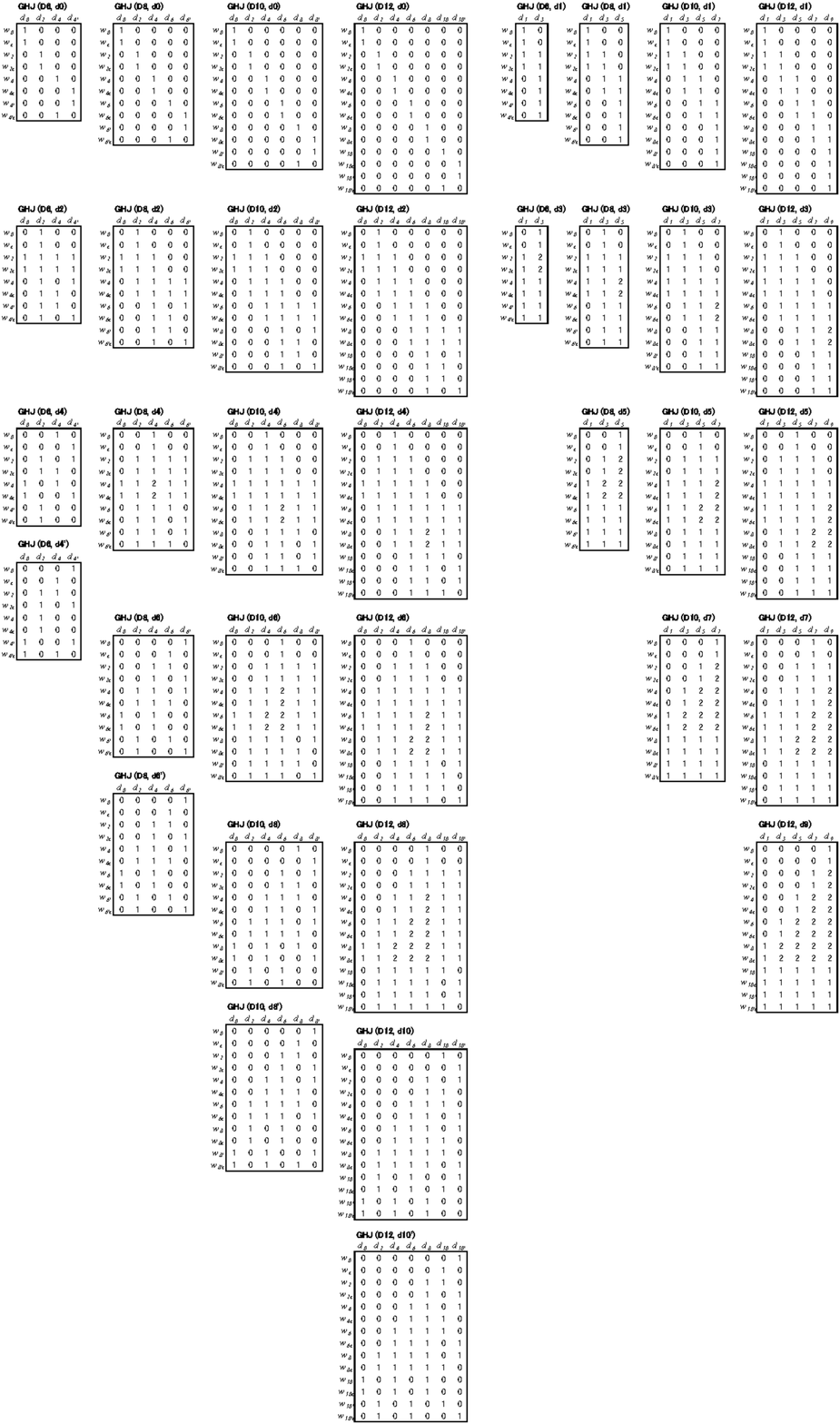}
\caption{The incidence matrices of 
the dual principal graphs of GHJ$(D_{even})$.}
\label{GHJ(D_even)-DPG}
\end{figure}

The incidence matrices of the (dual) principal graphs of 
GHJ$(D_{even},*_K)$ are also given in Figures 
\ref{GHJ(D_even)-PG} and \ref{GHJ(D_even)-DPG}. 

%%%%%%%%%%%%%%%%%%%%
%%%    D_{2n}    %%%
%%%%%%%%%%%%%%%%%%%%
\begin{figure}[H]
\unitlength 0.75mm
\thicklines
\begin{center}
\begin{picture}(100,20)(0,-5)
\multiput(10,0)(20,0){2}{\circle*{2}}
\multiput(0,10)(20,0){6}{\circle*{2}}
\multiput(70,0)(20,0){2}{\circle*{2}}
\put(90,10){\circle*{2}}
\path(0,10)(10,0)(20,10)(30,0)(40,10)
\path(60,10)(70,0)(80,10)(90,0)(100,10) 
\path(90,10)(90,0)
%\put(-9,14){\makebox(0,0){$*_A=$}}
\put(0,14){\makebox(0,0){$d_0$}}
\put(10,-4){\makebox(0,0){$d_1$}}
\put(20,14){\makebox(0,0){$d_2$}}
\put(30,-4){\makebox(0,0){$d_3$}}
\put(40,14){\makebox(0,0){$d_4$}}
\put(50,5){\makebox(0,0){$\cdots$}}
\put(60,14){\makebox(0,0){$d_{2n-6}$}}
\put(69,-5){\makebox(0,0){$d_{2n-5}$}}
\put(79,14){\makebox(0,0){$d_{2n-4}$}}
\put(92,-5){\makebox(0,0){$d_{2n-3}$}}
\put(92,14){\makebox(0,0){$d_{2n-2}$}}
\put(106,14){\makebox(0,0){$d_{2n-2}'$}}
\end{picture}
\caption{The label of vertices of the Dynkin diagram $D_{2n}$.}
\label{label-D2n}
\end{center}
%%%%%%%%%%%%%%%%%%%%
\begin{center}
\begin{picture}(80,90)(30,-13)
\path(-20,-17)(-20,80)(150,80)(150,-17)(-20,-17)
\drawline(3,70)(20,70)(20,30)(0,30)(0,67)
\drawline(0,50)(20,50)
\multiput(0,30)(0,20){2}{\circle*{1}}
\put(0,70){\makebox(0,0){$*$}}
\put(10,74){\makebox(0,0){$A_{4n-3}$}}
\put(10,54){\makebox(0,0){$D_{2n}$}}
\put(10,26){\makebox(0,0){$D_{2n}$}}
\put(-5,50){\makebox(0,0){$d_0$}}
\put(-5,30){\makebox(0,0){$d_k$}}
\put(40,70){\makebox(0,0){$\cdots \quad \longrightarrow$}}
\put(40,50){\makebox(0,0){$\cdots \quad \longrightarrow$}}
\put(40,30){\makebox(0,0){$\cdots \quad \longrightarrow$}}
\put(60,70){\makebox(0,0){$N$}}
\put(60,60){\makebox(0,0){$\cap$}}
\put(59,50){\makebox(0,0){${}^{\exists} P$}}
\put(60,40){\makebox(0,0){$\cap$}}
\put(60,30){\makebox(0,0){$M$}}
\put(110,63){\makebox(0,0){$(N\subset P) \cong (N\subset N\rtimes\Z_2) $}}
\put(110,57){\makebox(0,0){index~=~2}}
\put(110,40){\makebox(0,0){$P\neq M$}}
\put(50,15){\makebox(0,0){The principal graph}}
\put(80,15){\makebox(0,0){$\not\cong$}}
\put(115,15){\makebox(0,0){The dual principal graph}}
\put(5,7){\makebox(0,0){The number}}
\put(6,2){\makebox(0,0){of even vertices}}
\put(55,5){\makebox(0,0){$2n-1$}}
\put(80,5){\makebox(0,0){$\neq$}}
\put(105,5){\makebox(0,0){$2n+2$}}
\put(5,-7){\makebox(0,0){The fusion rule}}
\put(6,-12){\makebox(0,0){of even vertices}}
\put(55,-5){\makebox(0,0){$A_{4n-3}^{\rm even}$}}
\put(80,-5){\makebox(0,0){$\neq$}}
\put(105,-5){\makebox(0,0){${}_{D_{2n}}Z_{D_{2n}}^{\rm even}$}}
\put(57,-12){\makebox(0,0){commutative}}
\put(110,-12){\makebox(0,0){non-commutative}}
\end{picture}
\caption{ }
\label{GHJ-D2n}
\end{center}
\end{figure}
%%%%%%%%%%%%%%%%%%%%  

\begin{figure}[H]
\centering
\includegraphics[width=105mm,clip]{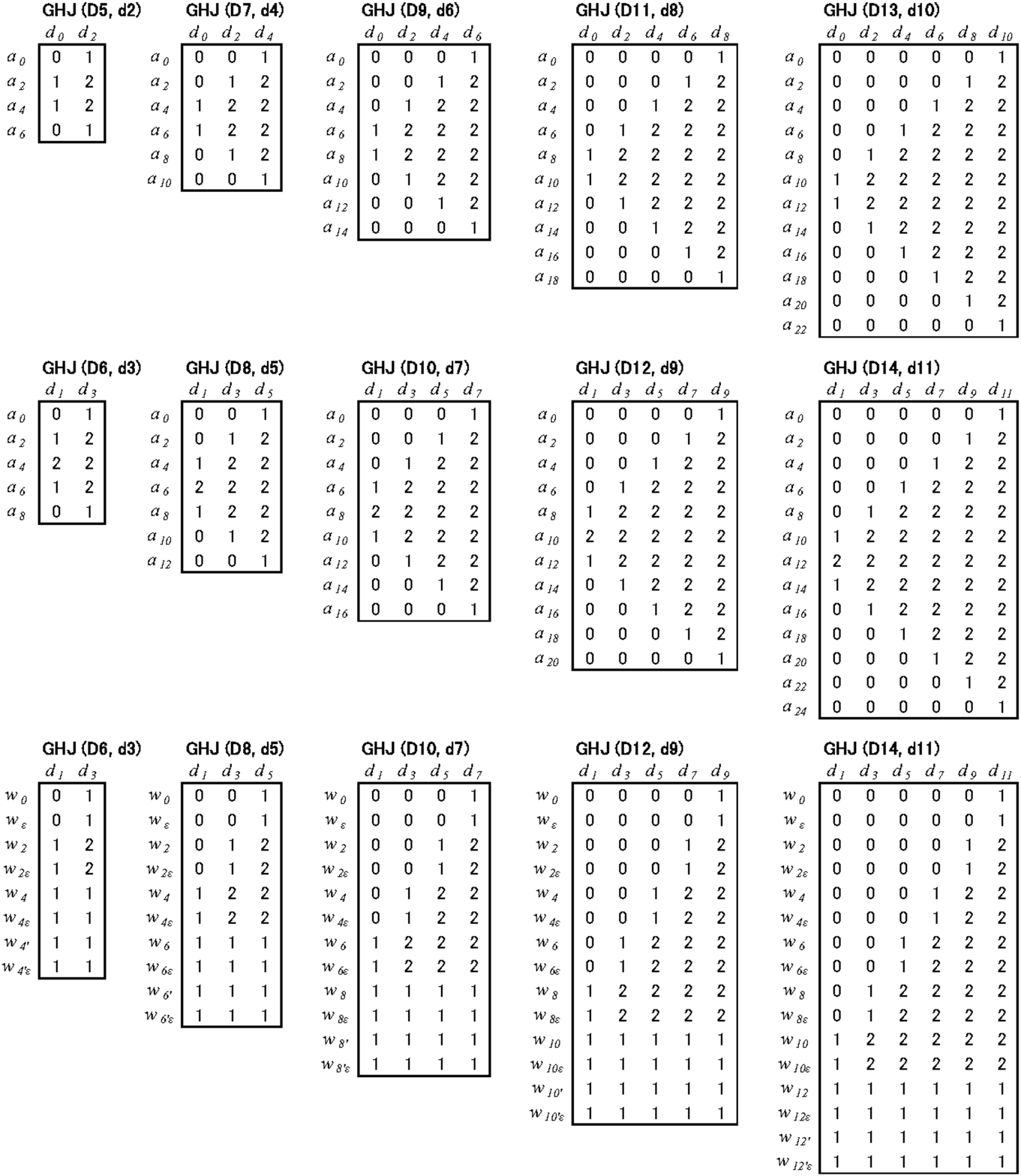}\\
{The incidence matrices of the (dual) principal graphs 
of GHJ$(D,*={\rm triple~ point})$.}
\caption{}
\label{GHJ(D-TriplePt)-PG+DPG}
\end{figure}

\begin{example}
From these computations for $\GHJ(D_{n},*_K)$ as above, 
the (dual) principal graphs of $\GHJ(D_{n},*={\rm triple~ point})$
can be obtained for general $n$ as 
in Figure \ref{GHJ(D-TriplePt)-PG+DPG}.
\end{example}

%\newpage
%\vspace{3mm}
\noindent
%%%%%%%%%%%%%%%%%%%%%%%%%%%%%%%%%%%%%%%%%%%%%%%%%%%%%%%%%%%%%
%%%%% section 4.4 %%%%%%%%%%%%%%%%%%%%%%%%%%%%%%%%%%%%%%%%%%%
%%%%%%%%%%%%%%%%%%%%%%%%%%%%%%%%%%%%%%%%%%%%%%%%%%%%%%%%%%%%%
{\sl 4.4. Goodman-de la Harpe-Jones subfactors of type $E_6$}
\medskip

We label the vertices of the Dynkin diagram $E_6$ 
by $e_0,e_1,e_2,\cdots,e_5$
as in Figure \ref{label-E6}.

\begin{figure}[H]
\centering
\includegraphics[width=140mm,clip]{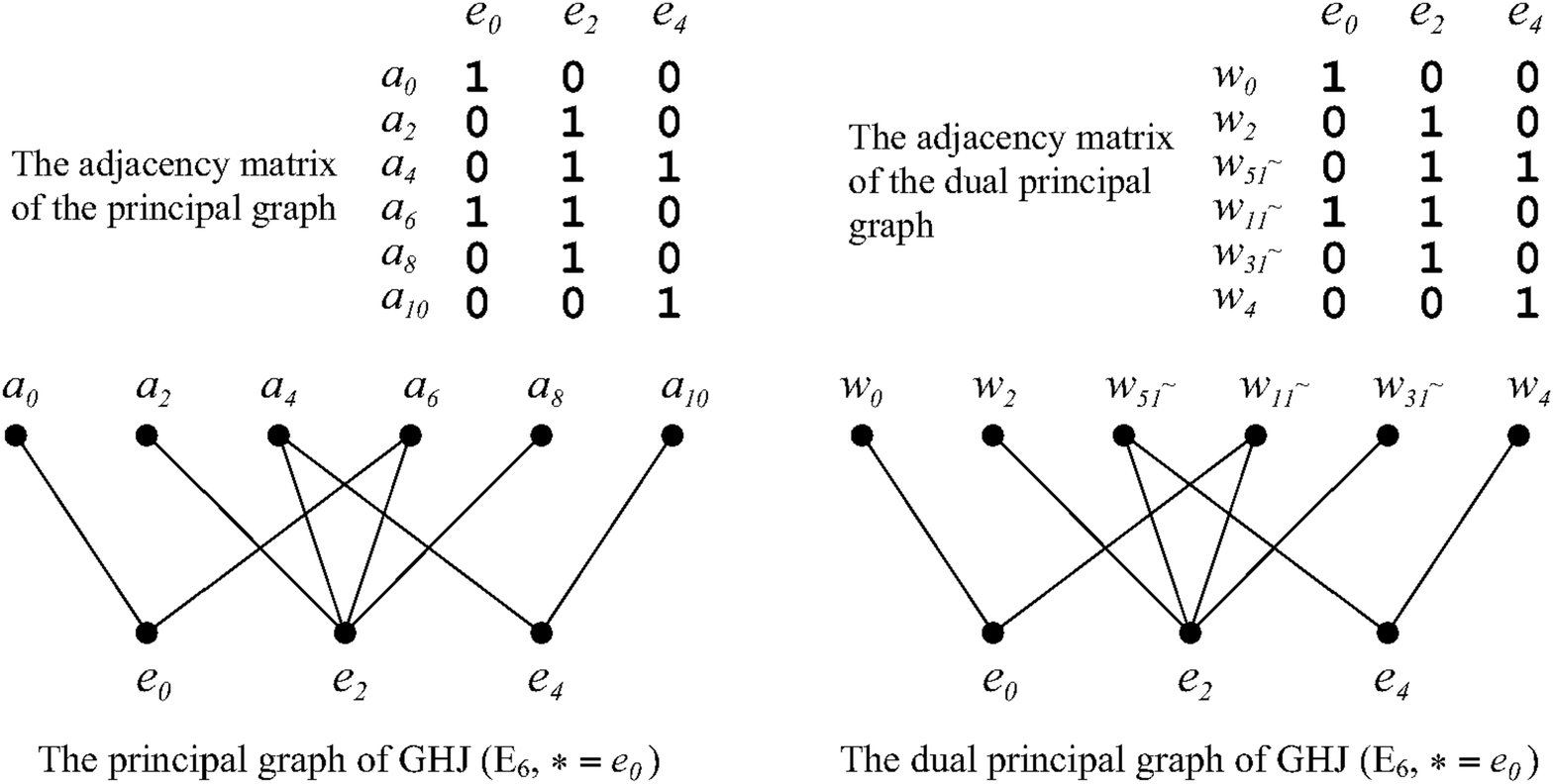}
\caption{The (dual) principal graph of GHJ$(E_6, *=e_0)$.}
\label{GHJ(E6-e0)}
\end{figure}

\begin{figure}[H]
\centering
\includegraphics[width=140mm,clip]{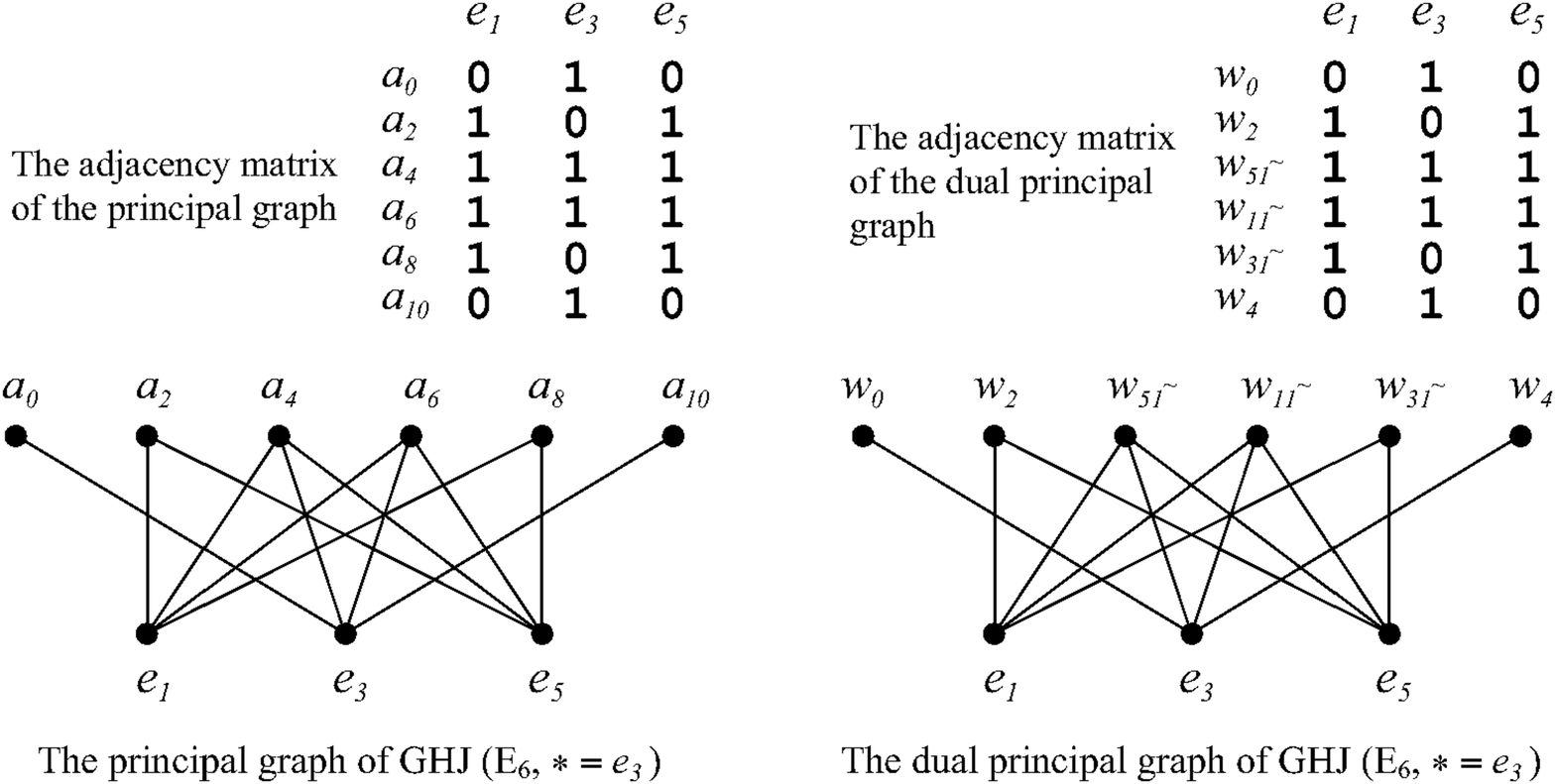}
\caption{The (dual) principal graph of GHJ$(E_6, *=e_3)$.}
\label{GHJ(E6-e3)}
\end{figure}

The Goodman-de la Harpe-Jones subfactor $\GHJ(E_6,*_K=e_0)$,  
has index $3+\sqrt{3}$ and it has the same principal and 
dual principal graph. But the fusion rules of the two graphs are different. 
This subfactor is known as the example which has the smallest index 
among such subfactors.   
The (dual) principal graphs of $\GHJ(E_6,*_K)$ are given 
in Figures \ref{GHJ(E6-e0)}`\ref{GHJ(E6-e3)}. 

\begin{figure}[H]
%\centering\includegraphics[width=140mm,clip]{qsym-e6.eps}
\centering\includegraphics[height=75mm,clip]{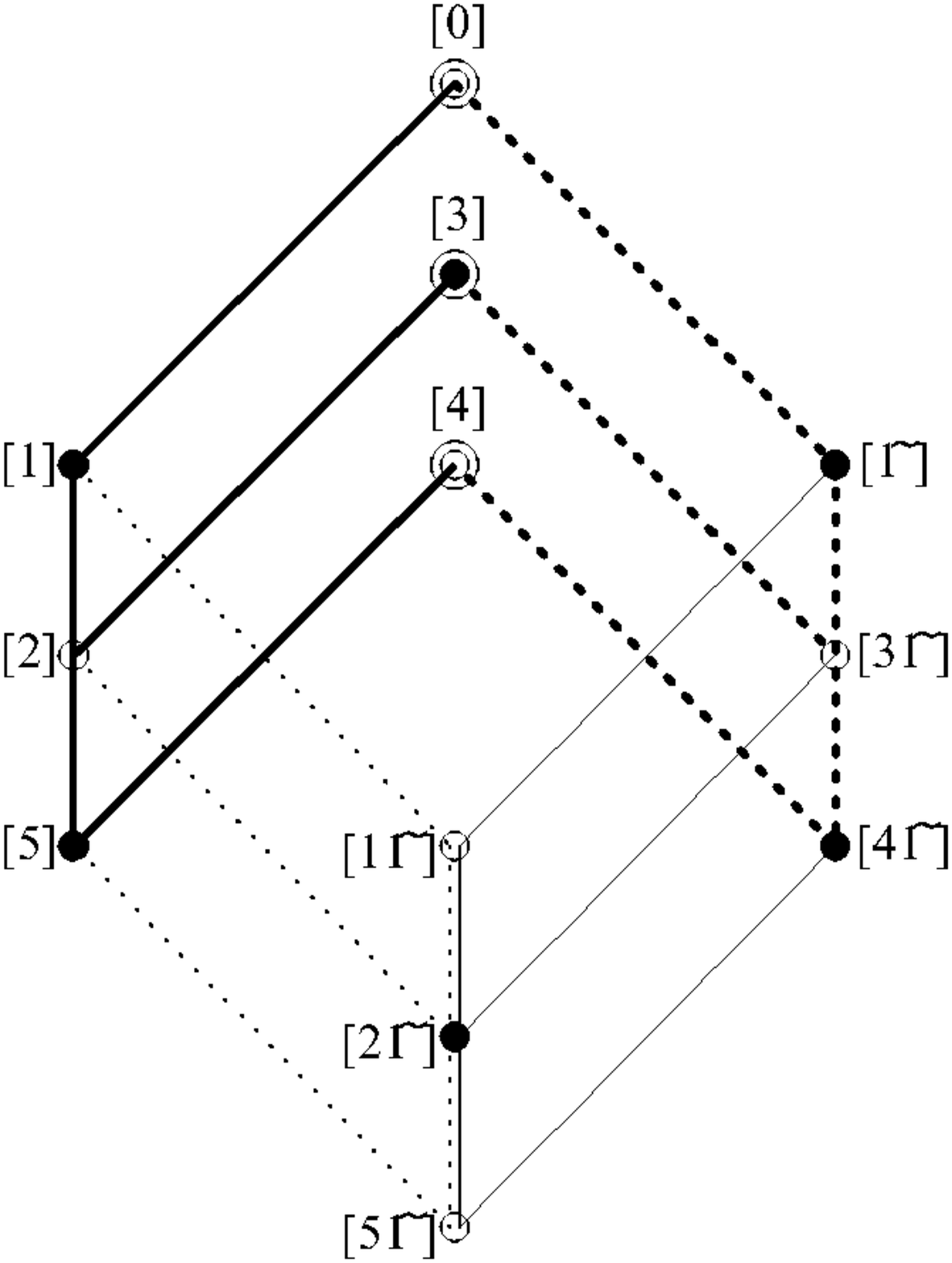}
\caption{Chiral symmetry for the Coxeter graph $E_6$}
\label{qsym-e6}
\end{figure}

\begin{figure}[H]
\centering\includegraphics[width=140mm,clip]{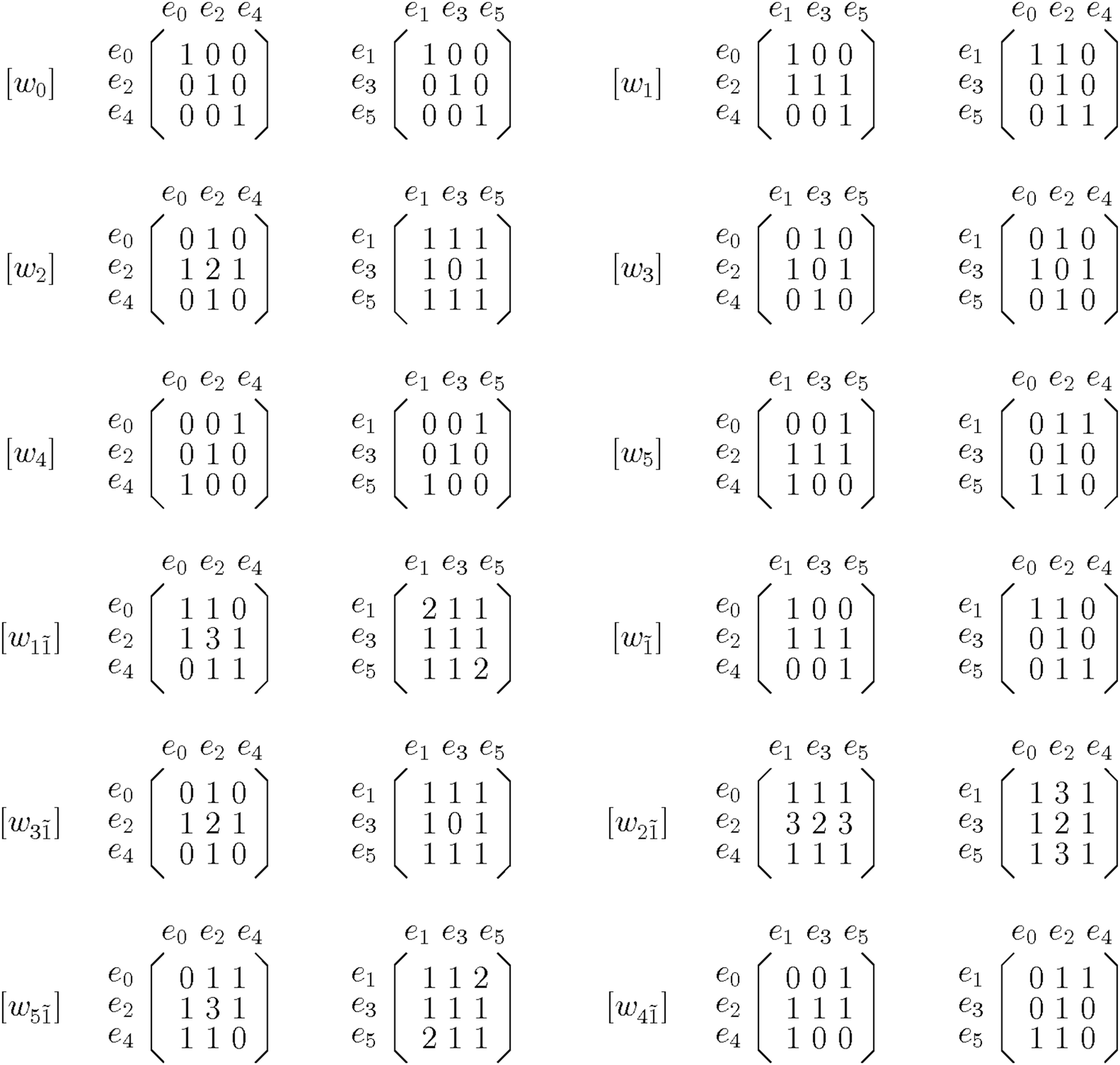}
\caption{The incidence matrices of the vertical edges 
of $E_6$-$E_6$ connections}
\label{v-edge-E6}
\end{figure}

If the vertex $*_K\neq e_0, e_4$, 
$\GHJ(E_6,*_K)$ has nontrivial intermediate subfactor 
as in Figure \ref{GHJ-E6}
because we have the decomposition of connections 
${}_Ae_k{}_E={}_Ae_0{}_E \cdot {}_E[w_k]{}_E$ for $k=1,2,3,5$. 
Here ${}_E[w_k]{}_E$ is the $E_6$-$E_6$ connection corresponding
to the vertex $[k]$ as in Figures \ref{qsym-e6} and \ref{v-edge-E6}.

%%%%%%%%%%%%%%%%%%%%
%%%%%   E_6   %%%%%%
%%%%%%%%%%%%%%%%%%%%
\begin{figure}[H]
\unitlength 0.75mm
\thicklines
\begin{center}
\begin{picture}(40,20)(0,-5)
\multiput(10,0)(10,0){3}{\circle*{2}}
\multiput(0,10)(20,0){3}{\circle*{2}}
\path(0,10)(10,0)(20,10)(30,0)(40,10)
\path(20,10)(20,0)
\put(0,14){\makebox(0,0){$e_0$}}
\put(10,-4){\makebox(0,0){$e_1$}}
\put(20,14){\makebox(0,0){$e_2$}}
\put(20,-4){\makebox(0,0){$e_3$}}
\put(30,-4){\makebox(0,0){$e_5$}}
\put(40,14){\makebox(0,0){$e_4$}}
\end{picture}
\caption{The label of vertices of the Dynkin diagram $E_6$.}
\label{label-E6}
\end{center}
%%%%%%%%%%%%%%%%%%%%
\thicklines
\begin{center}
\begin{picture}(80,82)(30,-4)
\path(-20,-8)(-20,80)(150,80)(150,-8)(-20,-8)
\drawline(3,70)(20,70)(20,30)(0,30)(0,67)
\drawline(0,50)(20,50)
\multiput(0,30)(0,20){2}{\circle*{1}}
\put(0,70){\makebox(0,0){$*$}}
\put(10,74){\makebox(0,0){$A_{11}$}}
\put(10,54){\makebox(0,0){$E_6$}}
\put(10,26){\makebox(0,0){$E_6$}}
\put(-5,50){\makebox(0,0){$e_0$}}
\put(-5,30){\makebox(0,0){$e_k$}}
\put(40,70){\makebox(0,0){$\cdots \quad \longrightarrow$}}
\put(40,50){\makebox(0,0){$\cdots \quad \longrightarrow$}}
\put(40,30){\makebox(0,0){$\cdots \quad \longrightarrow$}}
\put(60,70){\makebox(0,0){$N$}}
\put(60,60){\makebox(0,0){$\cap$}}
\put(59,50){\makebox(0,0){${}^{\exists} P$}}
\put(60,40){\makebox(0,0){$\cap$}}
\put(60,30){\makebox(0,0){$M$}}
\put(110,63){\makebox(0,0){$(N\subset P) \cong \GHJ(E_6,*_K=e_0) $}}
\put(110,57){\makebox(0,0){index~=~$3+\sqrt{3}$}}
\put(110,40){\makebox(0,0){$P\neq M$}}
\put(50,15){\makebox(0,0){The principal graph}}
\put(80,15){\makebox(0,0){$\cong$}}
\put(115,15){\makebox(0,0){The dual principal graph}}
\put(5,7){\makebox(0,0){The fusion rule}}
\put(6,2){\makebox(0,0){of even vertices}}
\put(55,5){\makebox(0,0){$A_{11}^{\rm even}$}}
\put(80,5){\makebox(0,0){$\neq$}}
\put(105,5){\makebox(0,0){${}_{E_6}Z_{E_6}^{\rm even}$}}
\put(57,-2){\makebox(0,0){commutative}}
\put(110,-2){\makebox(0,0){commutative}}
\end{picture}
\caption{ }
\label{GHJ-E6}
\end{center}
\end{figure}
%%%%%%%%%%%%%%%%%%%%  

%\newpage
%\vspace{3mm}
\noindent
%%%%%%%%%%%%%%%%%%%%%%%%%%%%%%%%%%%%%%%%%%%%%%%%%%%%%%%%%%%%%
%%%%% section 4.5 %%%%%%%%%%%%%%%%%%%%%%%%%%%%%%%%%%%%%%%%%%%
%%%%%%%%%%%%%%%%%%%%%%%%%%%%%%%%%%%%%%%%%%%%%%%%%%%%%%%%%%%%%
{\sl 4.5. Goodman-de la Harpe-Jones subfactors of type $E_7$}
\medskip

We label the vertices of the Dynkin diagram $E_7$ 
by $e_0,e_1,e_2,\cdots,e_6$
as in Figure \ref{label-E7}.

\begin{figure}[H]
\centering
\includegraphics[width=140mm,clip]{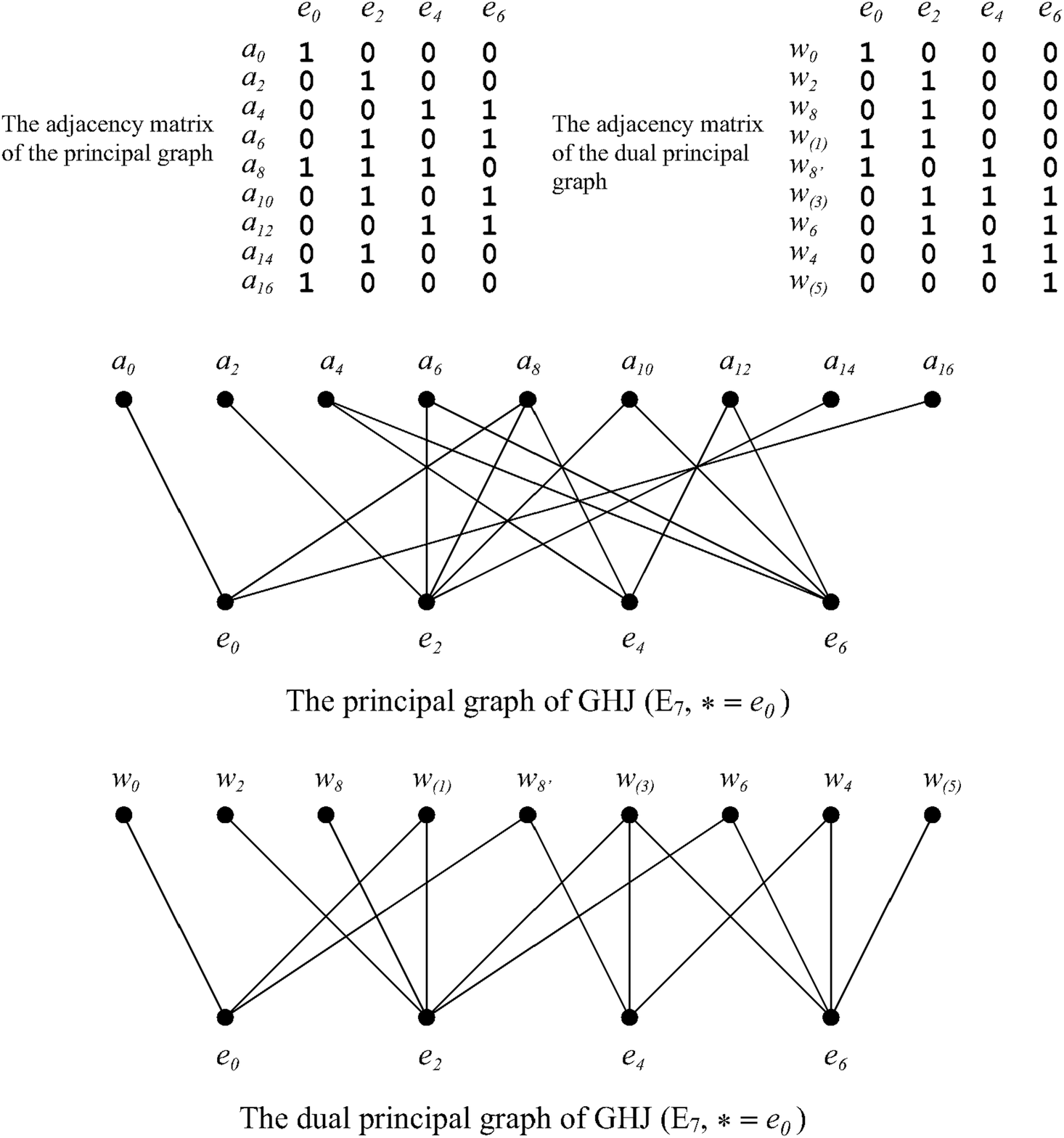}
\caption{The (dual) principal graph of GHJ$(E_7, *=e_0)$.}
\label{GHJ(E7-e0)}
\end{figure}

\begin{figure}[H]
\centering
\includegraphics[width=140mm,clip]{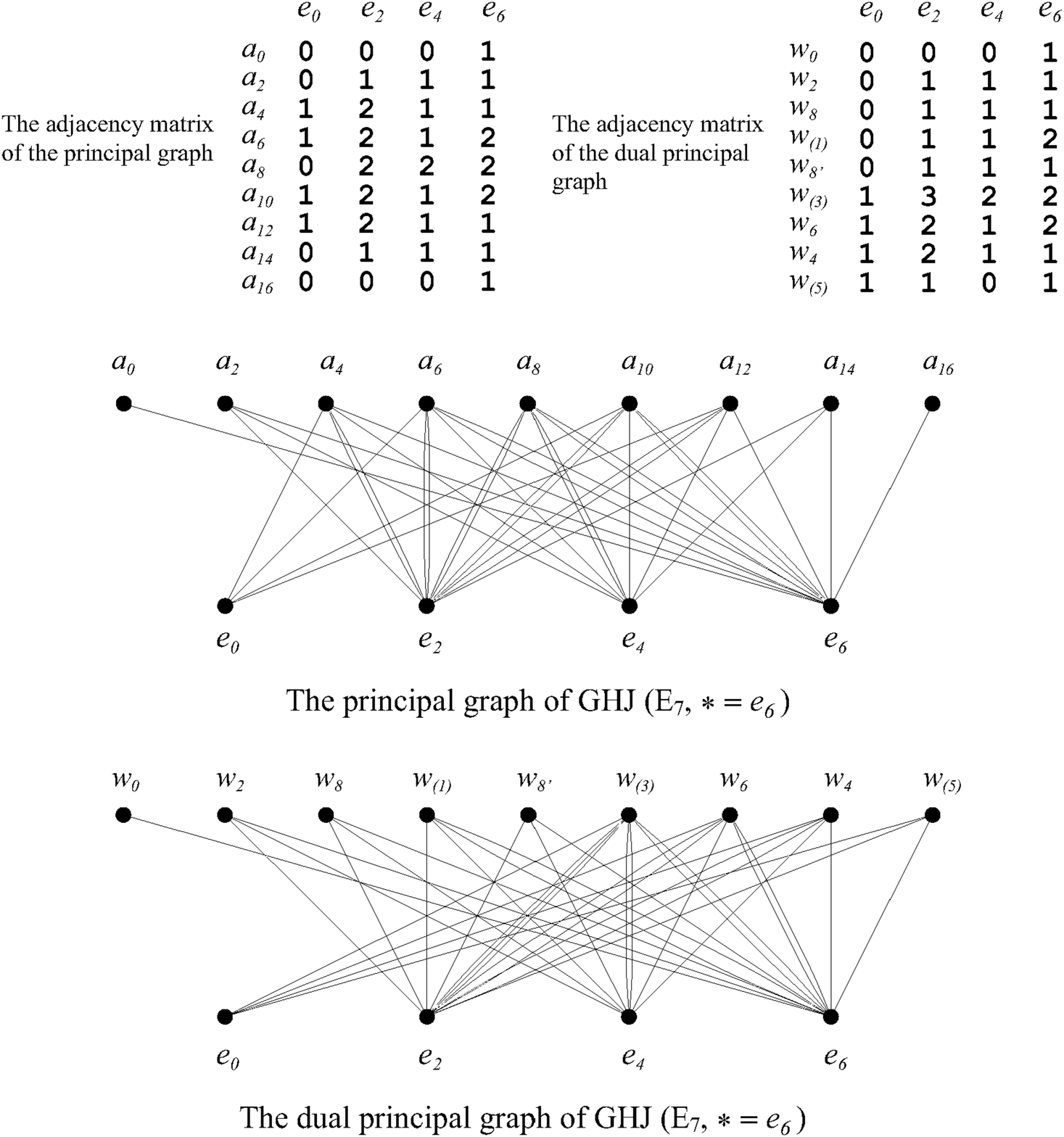}
\caption{The (dual) principal graph of the GHJ subfactor 
corresponding to $(E_7, *=e_6)$.}
\label{GHJ(E7-e6)}
\end{figure}

The Goodman-de la Harpe-Jones subfactor $\GHJ(E_7,*_K=e_0)$,  
has index $\Frac{|A_{17}|}{|E_7|}$ which is approximately $7.759$.
Here $|A_{17}|$ and $|E_7|$ represents the ``total mass'' of the graph 
$A_{17}$ and $E_7$ respectively, i.e. the sum of squares of 
normalized Perron-Frobenius eigenvalues over all the vertices of 
the graph. 
The (dual) principal graphs of $\GHJ(E_7,*_K)$ are given 
in Figures \ref{GHJ(E7-e0)}`\ref{GHJ(E7-e6)}. 

\begin{figure}[H]
%\centering\includegraphics[width=140mm,clip]{qsym-e7.eps}
%\centering\includegraphics[height=150mm,clip]{qsym-e7.eps}
\centering\includegraphics[height=80mm,clip]{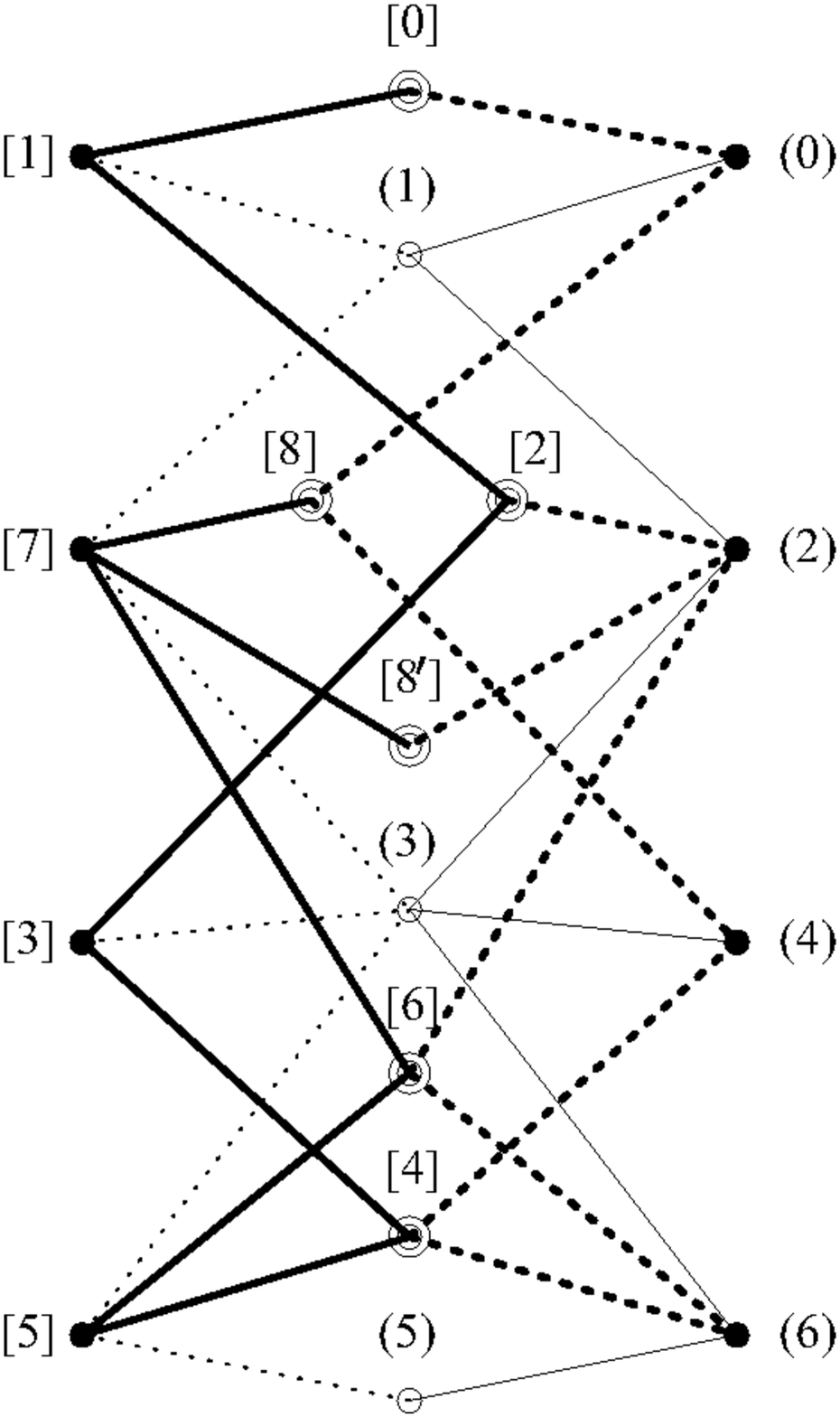}
\caption{Chiral symmetry for the Coxeter graph $E_7$}
\label{qsym-e7}
\end{figure}

\begin{figure}[H]
\centering\includegraphics[width=120mm,clip]{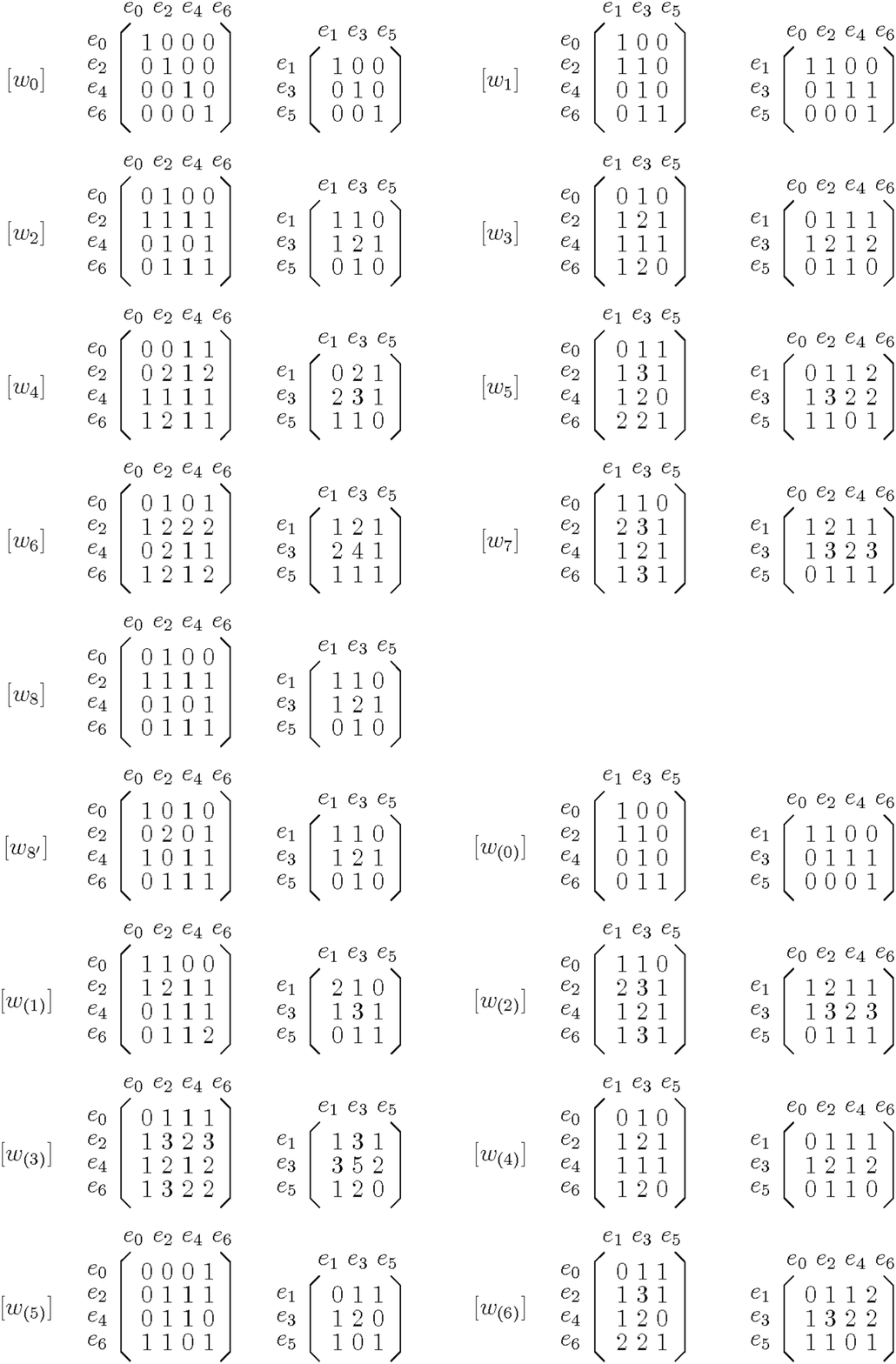}
\caption{The incidence matrices of the vertical edges 
of $E_7$-$E_7$ connections}
\label{v-edge-E7}
\end{figure}

If the vertex $*_K\neq e_0,e_4,e_5$, 
$\GHJ(E_7,*_K)$ has nontrivial intermediate subfactor 
as in Figure \ref{GHJ-E7}
because we have the decomposition of connections 
${}_Ae_k{}_E={}_Ae_0{}_E \cdot {}_E[w_k]{}_E$ for $k=1,2,3$
and ${}_Ae_6{}_E={}_Ae_0{}_E \cdot {}_E[w_{(5)}]{}_E$. 
Here ${}_E[w_k]{}_E$ is the $E_7$-$E_7$ connection corresponding
to the vertex $[k]$ $(k=1,2,3)$ and $(5)$ 
as in Figures \ref{qsym-e7} and \ref{v-edge-E7}.

%%%%%%%%%%%%%%%%%%%%
%%%%%   E_7   %%%%%%
%%%%%%%%%%%%%%%%%%%%
\begin{figure}[H]
\unitlength 0.75mm
\thicklines
\begin{center}
\begin{picture}(50,20)(0,-5)
\multiput(10,0)(20,0){3}{\circle*{2}}
\multiput(0,10)(20,0){3}{\circle*{2}}
\put(30,10){\circle*{2}}
\path(0,10)(10,0)(20,10)(30,0)(40,10)(50,0)
\path(30,10)(30,0)
\put(0,14){\makebox(0,0){$e_0$}}
\put(10,-4){\makebox(0,0){$e_1$}}
\put(20,14){\makebox(0,0){$e_2$}}
\put(30,-4){\makebox(0,0){$e_3$}}
\put(30,14){\makebox(0,0){$e_4$}}
\put(40,14){\makebox(0,0){$e_6$}}
\put(50,-4){\makebox(0,0){$e_5$}}
\end{picture}
\caption{The label of vertices of the Dynkin diagram $E_7$.}
\label{label-E7}
\end{center}
%%%%%%%%%%%%%%%%%%%%
\begin{center}
\begin{picture}(80,92)(30,-13)
\path(-20,-17)(-20,80)(150,80)(150,-17)(-20,-17)
\drawline(3,70)(20,70)(20,30)(0,30)(0,67)
\drawline(0,50)(20,50)
\multiput(0,30)(0,20){2}{\circle*{1}}
\put(0,70){\makebox(0,0){$*$}}
\put(10,74){\makebox(0,0){$A_{17}$}}
\put(10,54){\makebox(0,0){$E_7$}}
\put(10,26){\makebox(0,0){$E_7$}}
\put(-5,50){\makebox(0,0){$e_0$}}
\put(-5,30){\makebox(0,0){$e_k$}}
\put(40,70){\makebox(0,0){$\cdots \quad \longrightarrow$}}
\put(40,50){\makebox(0,0){$\cdots \quad \longrightarrow$}}
\put(40,30){\makebox(0,0){$\cdots \quad \longrightarrow$}}
\put(60,70){\makebox(0,0){$N$}}
\put(60,60){\makebox(0,0){$\cap$}}
\put(59,50){\makebox(0,0){${}^{\exists} P$}}
\put(60,40){\makebox(0,0){$\cap$}}
\put(60,30){\makebox(0,0){$M$}}
\put(110,63){\makebox(0,0){$(N\subset P) \cong \GHJ(E_7,*_K=e_0) $}}
\put(110,57){\makebox(0,0){index~$\fallingdotseq$~7.759}}
\put(110,40){\makebox(0,0){$P\neq M$}}
\put(50,15){\makebox(0,0){The principal graph}}
\put(80,15){\makebox(0,0){$\not\cong$}}
\put(115,15){\makebox(0,0){The dual principal graph}}
\put(5,7){\makebox(0,0){The number}}
\put(6,2){\makebox(0,0){of even vertices}}
\put(55,5){\makebox(0,0){$9$}}
\put(80,5){\makebox(0,0){$=$}}
\put(105,5){\makebox(0,0){$9$}}
\put(5,-7){\makebox(0,0){The fusion rule}}
\put(6,-12){\makebox(0,0){of even vertices}}
\put(55,-5){\makebox(0,0){$A_{17}^{\rm even}$}}
\put(80,-5){\makebox(0,0){$\neq$}}
\put(105,-5){\makebox(0,0){${}_{E_7}Z_{E_7}^{\rm even}$}}
\put(57,-12){\makebox(0,0){commutative}}
\put(110,-12){\makebox(0,0){commutative}}
\end{picture}
\caption{ }
\label{GHJ-E7}
\end{center}
\end{figure}
%%%%%%%%%%%%%%%%%%%%  

%\newpage
%\vspace{3mm}
\noindent
%%%%%%%%%%%%%%%%%%%%%%%%%%%%%%%%%%%%%%%%%%%%%%%%%%%%%%%%%%%%%
%%%%% section 4.6 %%%%%%%%%%%%%%%%%%%%%%%%%%%%%%%%%%%%%%%%%%%
%%%%%%%%%%%%%%%%%%%%%%%%%%%%%%%%%%%%%%%%%%%%%%%%%%%%%%%%%%%%%
{\sl 4.6. Goodman-de la Harpe-Jones subfactors of type $E_8$}
\medskip

We label the vertices of the Dynkin diagram $E_8$ 
by $e_0,e_1,e_2,\cdots,e_7$
as in Figure \ref{label-E8}.

\begin{figure}[H]
\centering
\includegraphics[width=140mm,clip]{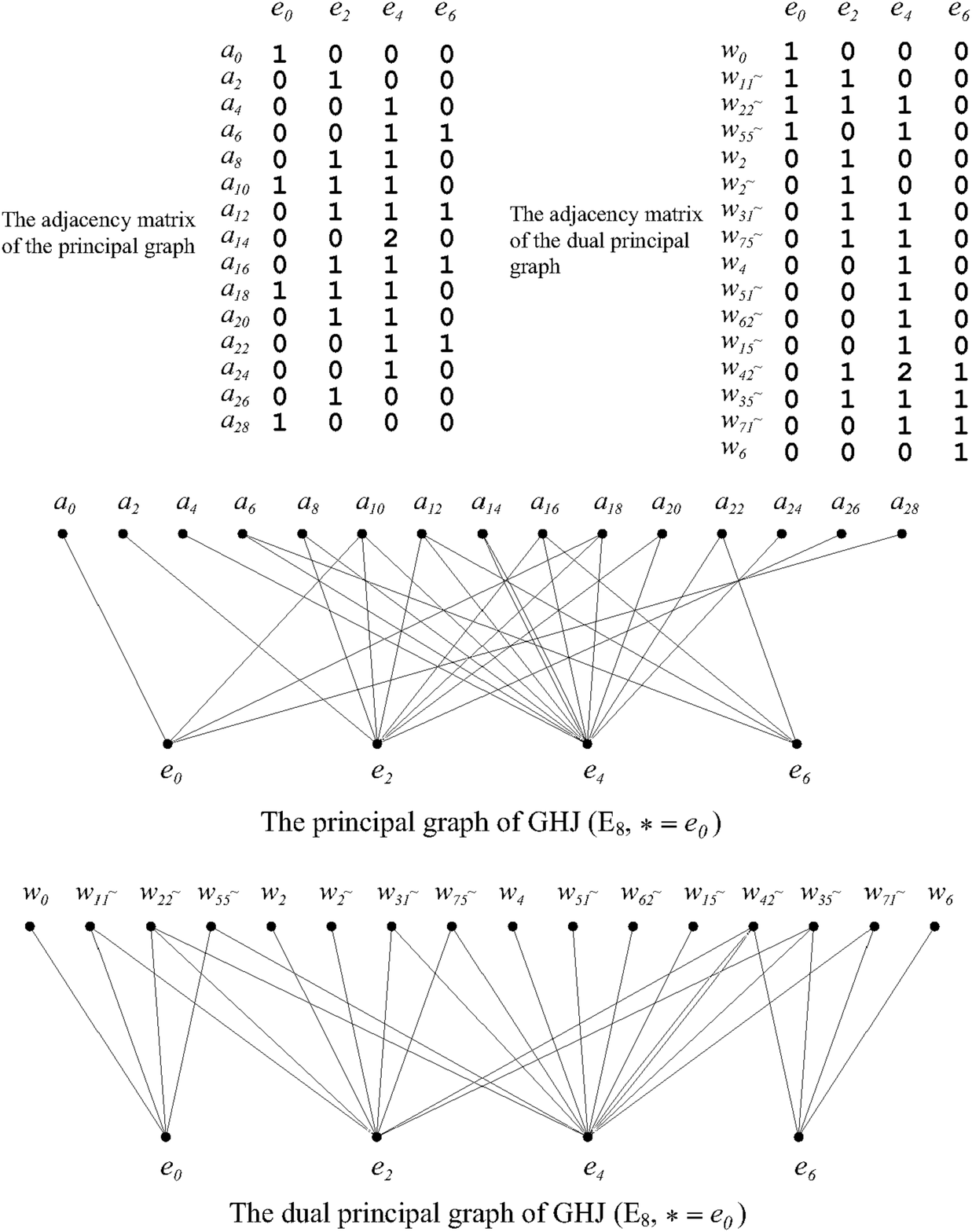}
\caption{The (dual) principal graph of the GHJ subfactor 
corresponding to $(E_8, *=e_0)$.}
\label{GHJ(E8-e0)}
\end{figure}

\begin{figure}[H]
\centering
\includegraphics[width=140mm,clip]{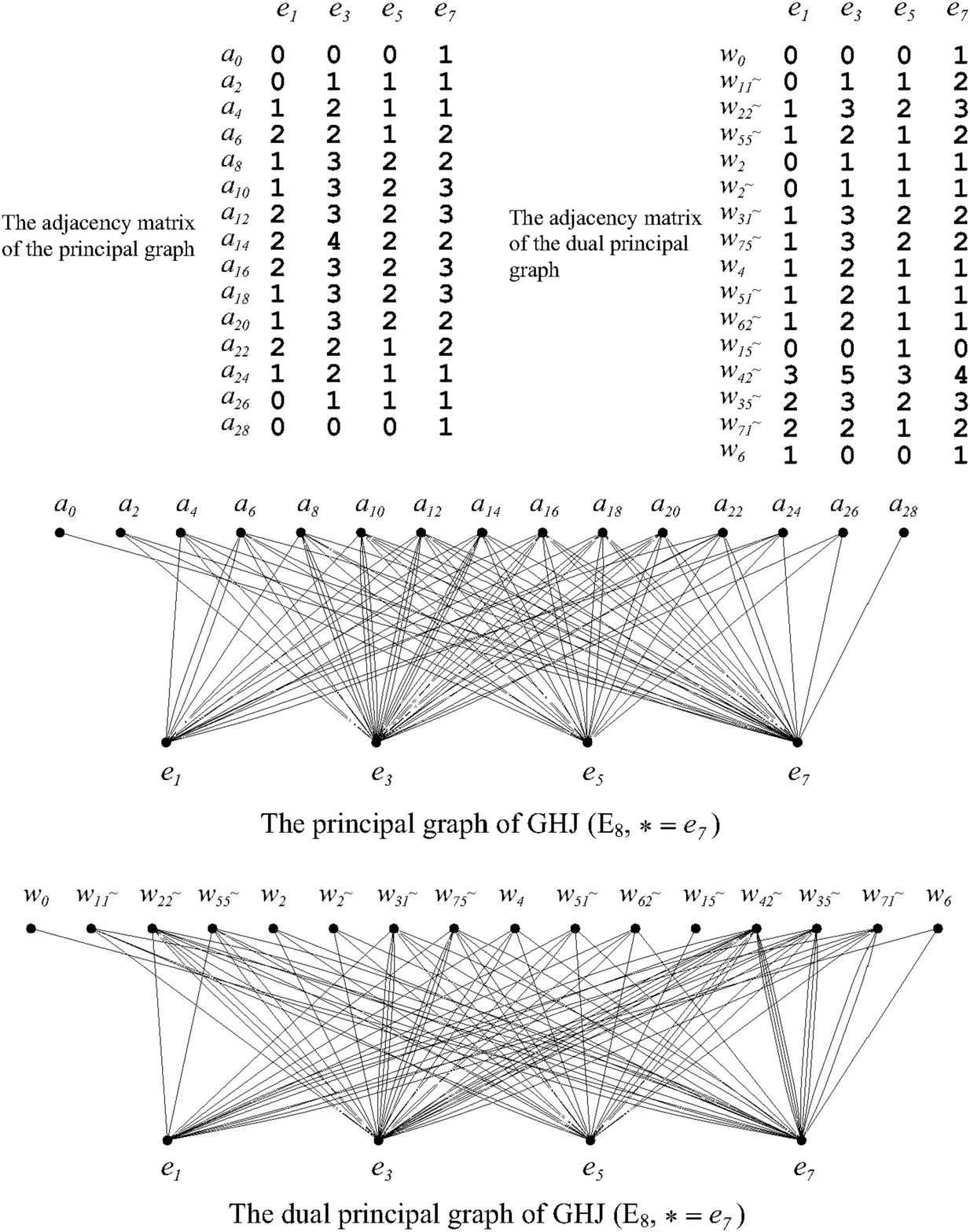}
\caption{The (dual) principal graph of the GHJ subfactor 
corresponding to $(E_8, *=e_7)$.}
\label{GHJ(E8-e7)}
\end{figure}

The Goodman-de la Harpe-Jones subfactor $\GHJ(E_8,*_K=e_0)$,  
has index $\Frac{|A_{29}|}{|E_8|}$ which is approximately $19.48$.
Here $|A_{29}|$ and $|E_8|$ represents the ``total mass'' of the graph 
$A_{29}$ and $E_8$ respectively. 
The (dual) principal graphs of $\GHJ(E_8,*_K)$ are given 
in Figures \ref{GHJ(E8-e0)}`\ref{GHJ(E8-e7)}. 

\begin{figure}[H]
\centering\includegraphics[width=120mm,clip]{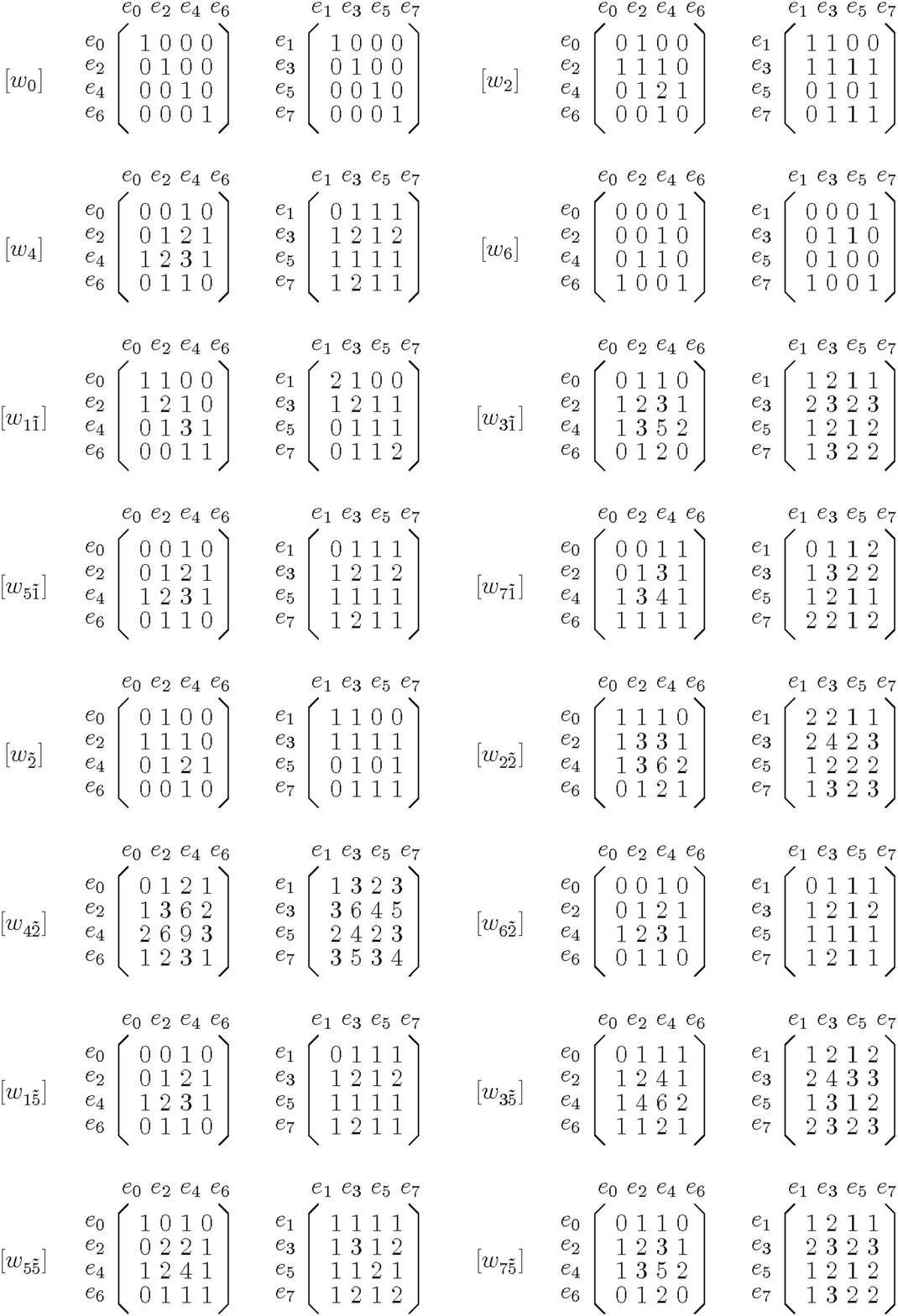}
\caption{The incidence matrices of the vertical edges 
of $E_8$-$E_8$ even connections}
\label{v-edge-E8even}
\end{figure}

\begin{figure}[H]
\centering\includegraphics[width=120mm,clip]{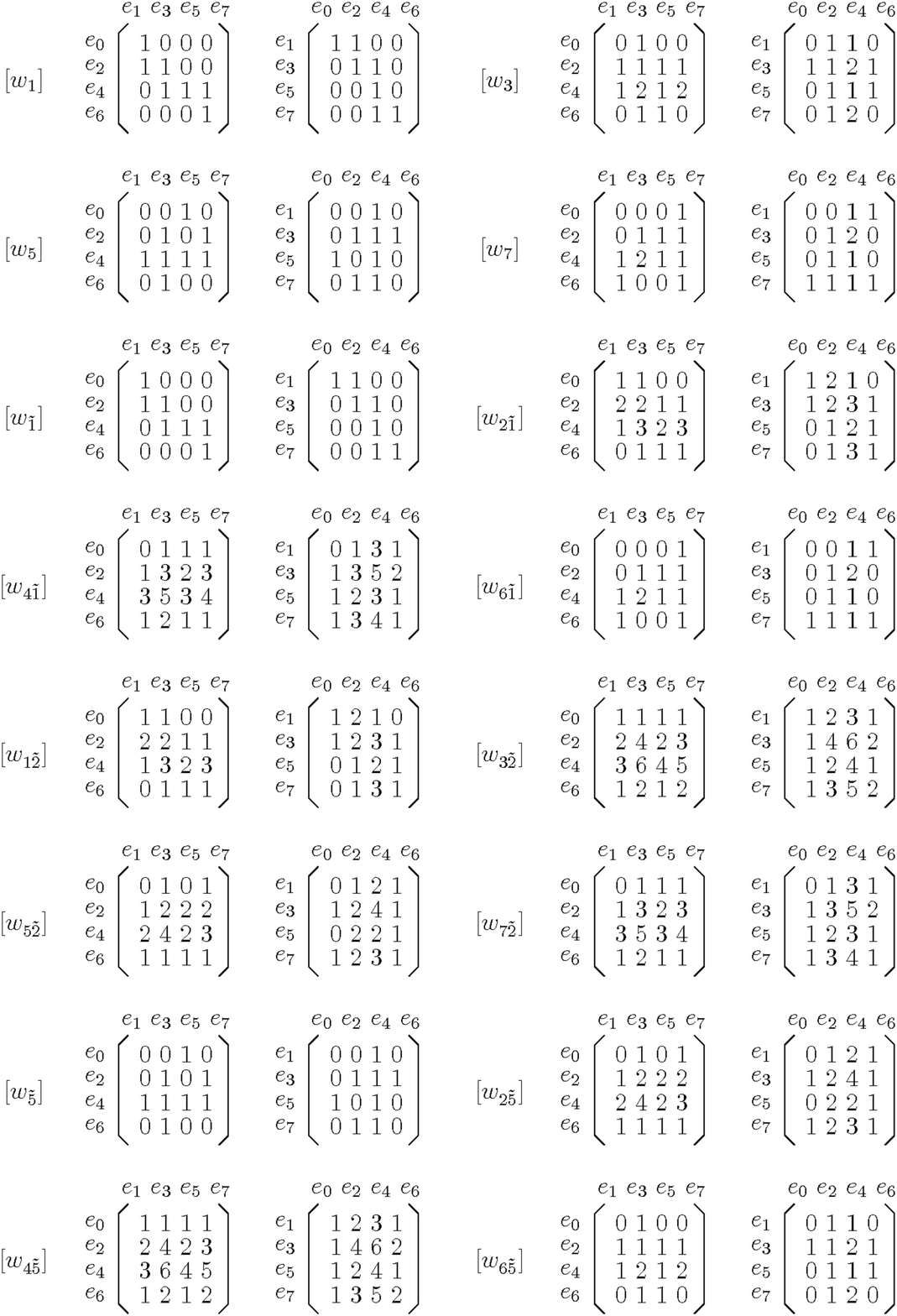}
\caption{The incidence matrices of the vertical edges 
of $E_8$-$E_8$ odd connections}
\label{v-edge-E8odd}
\end{figure}

If the vertex $*_K\neq e_0$, 
$\GHJ(E_7,*_K)$ has nontrivial intermediate subfactor 
as in Figure \ref{GHJ-E8}
because we have the decomposition of connections 
${}_Ae_k{}_E={}_Ae_0{}_E \cdot {}_E[w_k]{}_E$ for $k\neq 0$. 
Here ${}_E[w_k]{}_E$ is the $E_8$-$E_8$ connection corresponding
to the vertex $[k]$ as in Figures \ref{qsym-e8}, 
\ref{v-edge-E8even} and \ref{v-edge-E8odd}.

%%%%%%%%%%%%%%%%%%%%
%%%%%   E_8   %%%%%%
%%%%%%%%%%%%%%%%%%%%
\begin{figure}[H]
\unitlength 0.75mm
\thicklines
\begin{center}
\begin{picture}(60,20)(0,-5)
\multiput(10,0)(20,0){3}{\circle*{2}}
\multiput(0,10)(20,0){4}{\circle*{2}}
\put(40,0){\circle*{2}}
\path(0,10)(10,0)(20,10)(30,0)(40,10)(50,0)(60,10)
\path(40,10)(40,0)
\put(0,14){\makebox(0,0){$e_0$}}
\put(10,-4){\makebox(0,0){$e_1$}}
\put(20,14){\makebox(0,0){$e_2$}}
\put(30,-4){\makebox(0,0){$e_3$}}
\put(40,14){\makebox(0,0){$e_4$}}
\put(40,-4){\makebox(0,0){$e_5$}}
\put(50,-4){\makebox(0,0){$e_7$}}
\put(60,14){\makebox(0,0){$e_6$}}
\end{picture}
\caption{The label of vertices of the Dynkin diagram $E_8$.}
\label{label-E8}
\end{center}
%%%%%%%%%%%%%%%%%%%%
\begin{center}
\begin{picture}(80,92)(30,-13)
\path(-20,-17)(-20,80)(150,80)(150,-17)(-20,-17)
\drawline(3,70)(20,70)(20,30)(0,30)(0,67)
\drawline(0,50)(20,50)
\multiput(0,30)(0,20){2}{\circle*{1}}
\put(0,70){\makebox(0,0){$*$}}
\put(10,74){\makebox(0,0){$A_{29}$}}
\put(10,54){\makebox(0,0){$E_8$}}
\put(10,26){\makebox(0,0){$E_8$}}
\put(-5,50){\makebox(0,0){$e_0$}}
\put(-5,30){\makebox(0,0){$e_k$}}
\put(40,70){\makebox(0,0){$\cdots \quad \longrightarrow$}}
\put(40,50){\makebox(0,0){$\cdots \quad \longrightarrow$}}
\put(40,30){\makebox(0,0){$\cdots \quad \longrightarrow$}}
\put(60,70){\makebox(0,0){$N$}}
\put(60,60){\makebox(0,0){$\cap$}}
\put(59,50){\makebox(0,0){${}^{\exists} P$}}
\put(60,40){\makebox(0,0){$\cap$}}
\put(60,30){\makebox(0,0){$M$}}
\put(110,63){\makebox(0,0){$(N\subset P) \cong \GHJ(E_8,*_K=e_0) $}}
\put(110,57){\makebox(0,0){index~$\fallingdotseq$~19.48}}
\put(110,40){\makebox(0,0){$P\neq M$}}
\put(50,15){\makebox(0,0){The principal graph}}
\put(80,15){\makebox(0,0){$\not\cong$}}
\put(115,15){\makebox(0,0){The dual principal graph}}
\put(5,7){\makebox(0,0){The number}}
\put(6,2){\makebox(0,0){of even vertices}}
\put(55,5){\makebox(0,0){$15$}}
\put(80,5){\makebox(0,0){$\neq$}}
\put(105,5){\makebox(0,0){$16$}}
\put(5,-7){\makebox(0,0){The fusion rule}}
\put(6,-12){\makebox(0,0){of even vertices}}
\put(55,-5){\makebox(0,0){$A_{29}^{\rm even}$}}
\put(80,-5){\makebox(0,0){$\neq$}}
\put(105,-5){\makebox(0,0){${}_{E_8}Z_{E_8}^{\rm even}$}}
\put(57,-12){\makebox(0,0){commutative}}
\put(110,-12){\makebox(0,0){commutative}}
\end{picture}
\caption{ }
\label{GHJ-E8}
\end{center}
\end{figure}
%%%%%%%%%%%%%%%%%%%%  

%%%%%%%%%%%%%%%%%%%%%%%%%%%%%%%%%%%%%%%%%%%%%%%%%%%%%%%%%%%%%%%%%%%%%%
\section{An application to subequivalence on paragroups}
%%%%%%%%%%%%%%%%%%%%%%%%%%%%%%%%%%%%%%%%%%%%%%%%%%%%%%%%%%%%%%%%%%%%%%

Let $K$ be one of the Dynkin diagrams $D_{2n}(n\geq 3),E_6,E_8$
and $A_l$ the Dynkin diagram of type $A$ 
with the same Perron-Frobenius eigenvalue as $K$.
We can choose the vertex $*_K$ so that the GHJ subfactor  
$\GHJ(K,*_K)$ does not have index 2.
Let $N\subset M$ be the GHJ subfactor $\GHJ(K,*_K)$ chosen as above, then 
the fusion algebra of $N$-$N$ bimodules is isomorphic to $A_l^{\rm even}$ and 
the fusion algebra of $M$-$M$ bimodules is isomorphic to ${}_KZ_K^{\rm even}$.
Because ${}_KZ_K^{\rm even}$ contains $K^{\rm even}$ as its strict fusion
subalgebra, the paragroup of type $K$ becomes a {\sl strictly subequivalent} 
to that of type $A_l$. 
%Here, if a fusion algerba $\A$ is subequivalent but not equivalent to $\B$,
%we say that $\A$ is {\sl strictly subequivalent} to $\B$.
Here we use the terminology {\sl strictly subequivalent} in the sense
that a fusion algerba $\A$ is subequivalent but not equivalent to $\B$.
And in such a case, we denote $\A \succ \B$.

In the case of $D_4$, we can choose the direct sum of 3 connections for 
$\GHJ(D_4,*_K)$ ($*_K=d_0, d_2, d_2'$) as a connection 
for subequivalence between $A_5$ and $D_4$ paragroups. 

Hence we get the following subequivalence of paragroups.

\begin{theorem}~{\rm
The paragroups of Jones' type $A$ subfactors have the following 
strictly subequivalent paragroups. \\
\hspace{33mm}
$A_{4n-3} \succ D_{2n}$ $(n\geq 2)$,~~
$A_{11} \succ E_6$,~~ $A_{29} \succ E_8$.
}\end{theorem}

%%%%%%%%%%%%%%%%%%%%%%%%%%%%%%%%%%%%%%%%%%%%%%%%%%%%%%%%%%
%%%%%%%%%%%%%%%%%%%%%%%%%%%%%%%%%%%%%%%%%%%%%%%%%%%%%%%%%%
%%%%%%%%%%%%%%%%%%%%%%% FIGURES %%%%%%%%%%%%%%%%%%%%%%%%%%
%%%%%%%%%%%%%%%%%%%%%%%%%%%%%%%%%%%%%%%%%%%%%%%%%%%%%%%%%%
%%%%%%%%%%%%%%%%%%%%%%%%%%%%%%%%%%%%%%%%%%%%%%%%%%%%%%%%%%

%%%%%%%%%%%%%%%%%%%%%%%
%%% Essential Paths %%%
%%%%%%%%%%%%%%%%%%%%%%%

%%% EP-A4 %%%%%%%%%%%%%

%%% EP-A5 %%%%%%%%%%%%%
\begin{figure}[H]
\centering\includegraphics[width=140mm,clip]{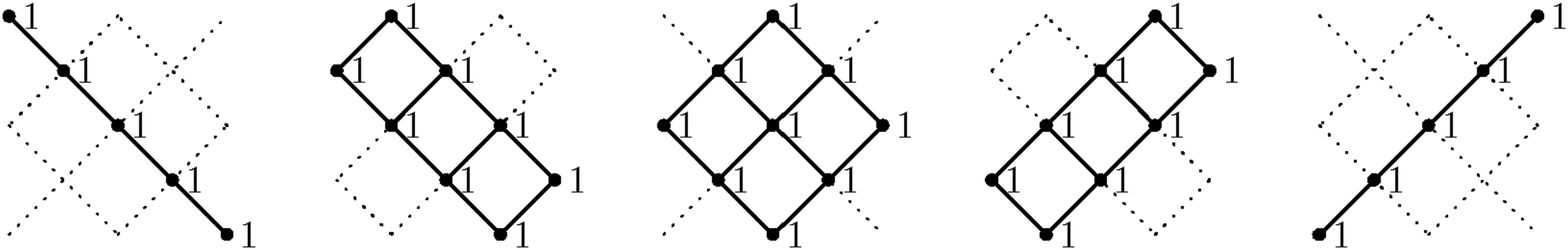}
\caption{Essential paths on the Coxeter graph $A_5$}
\label{EP-A5}
\end{figure}

%%% EP-D4 %%%%%%%%%%%%%
\begin{figure}[H]
\centering\includegraphics[width=90mm,clip]{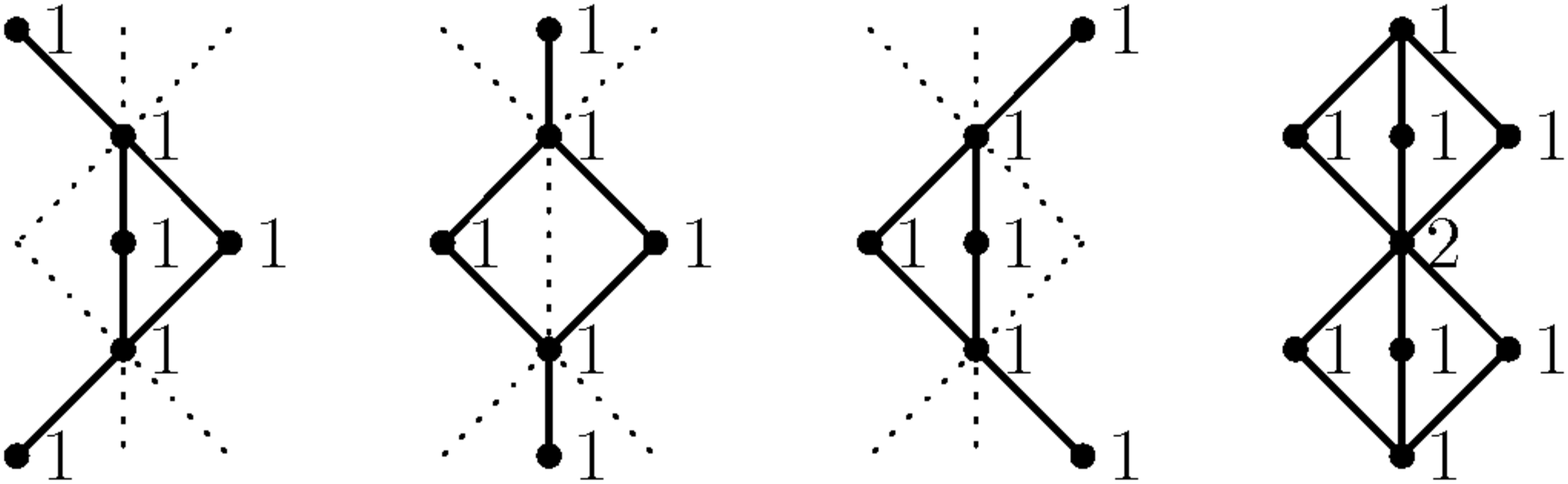}
\caption{Essential paths on the Coxeter graph $D_4$}
\label{EP-D4}
\end{figure}

%%% EP-D5 %%%%%%%%%%%%%
\begin{figure}[H]
\centering\includegraphics[width=140mm,clip]{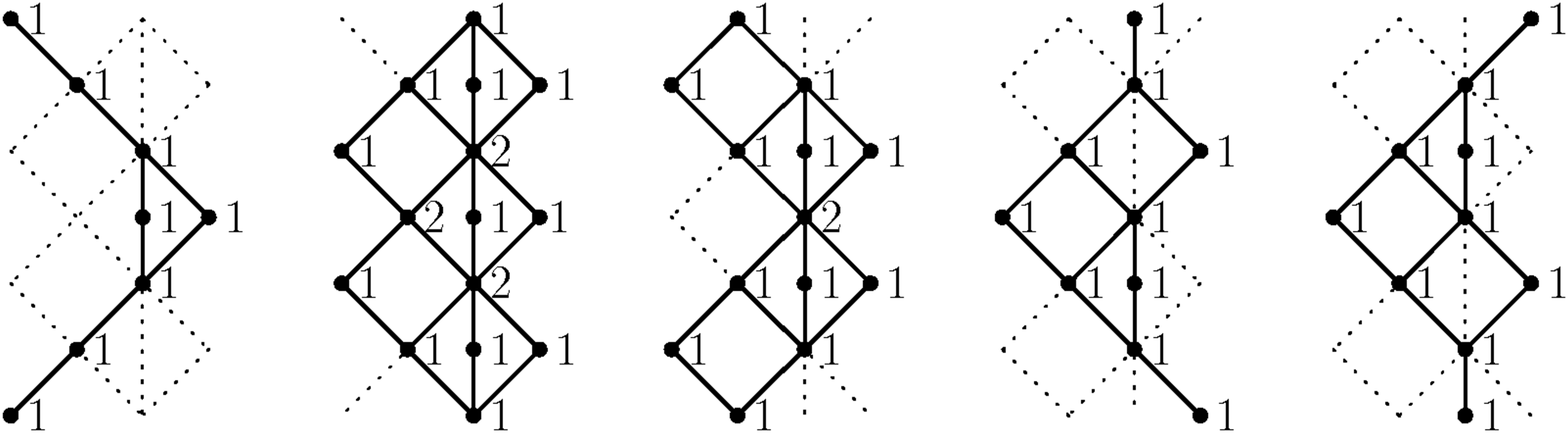}
\caption{Essential paths on the Coxeter graph $D_5$}
\label{EP-D5}
\end{figure}

%%% EP-D6 %%%%%%%%%%%%%
\begin{figure}[H]
\centering\includegraphics[width=140mm,clip]{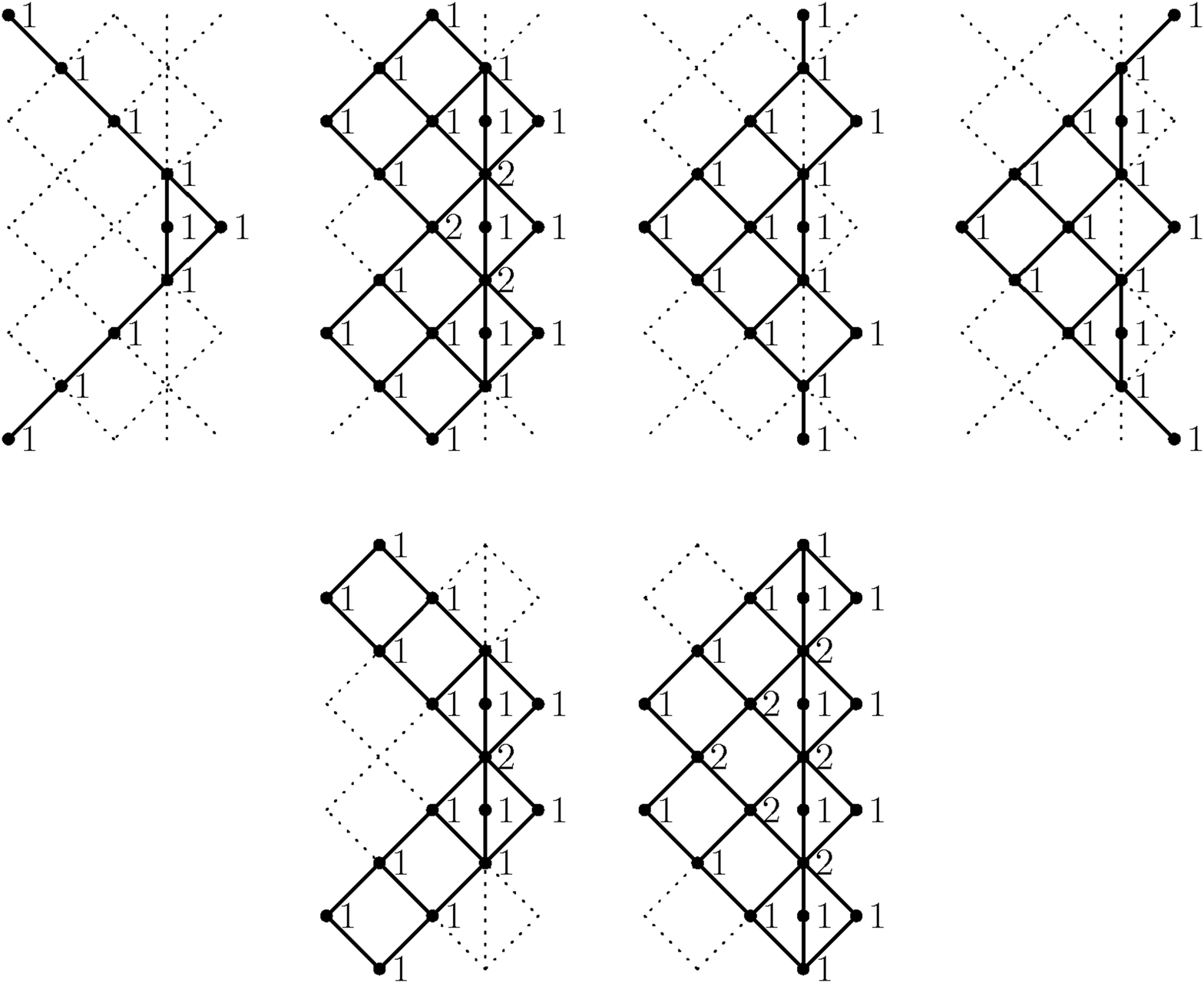}
\caption{Essential paths on the Coxeter graph $D_6$}
\label{EP-D6}
\end{figure}

%%% EP-E6 %%%%%%%%%%%%%
\begin{figure}[H]
\centering\includegraphics[width=100mm,clip]{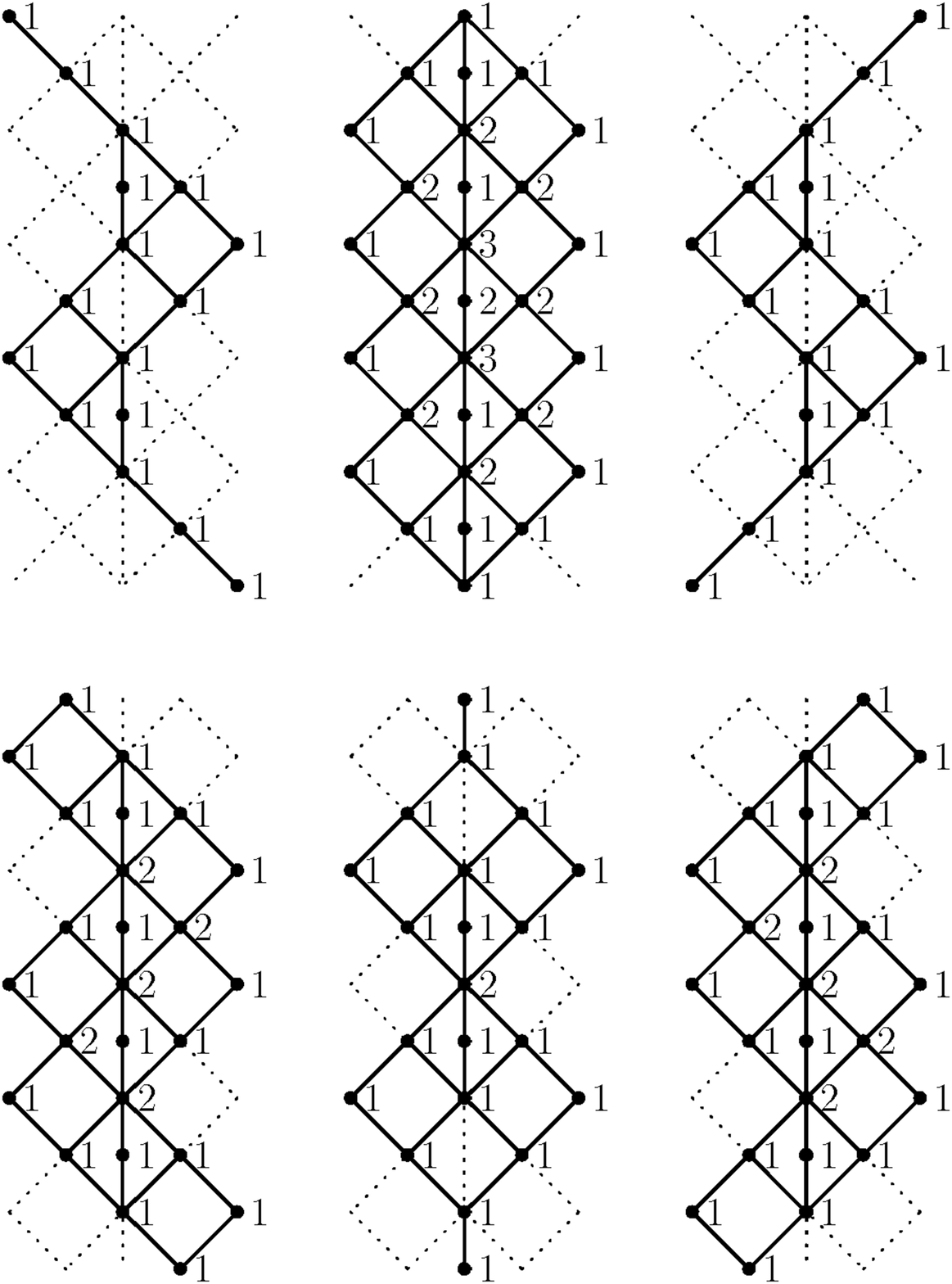}
\caption{Essential paths on the Coxeter graph $E_6$}
\label{EP-E6}
\end{figure}

%%% EP-E7-1 %%%%%%%%%%%%%
\begin{figure}[H]
%\centering\includegraphics[width=140mm,clip]{EP-E7-1.eps}
\centering\includegraphics[height=77mm,clip]{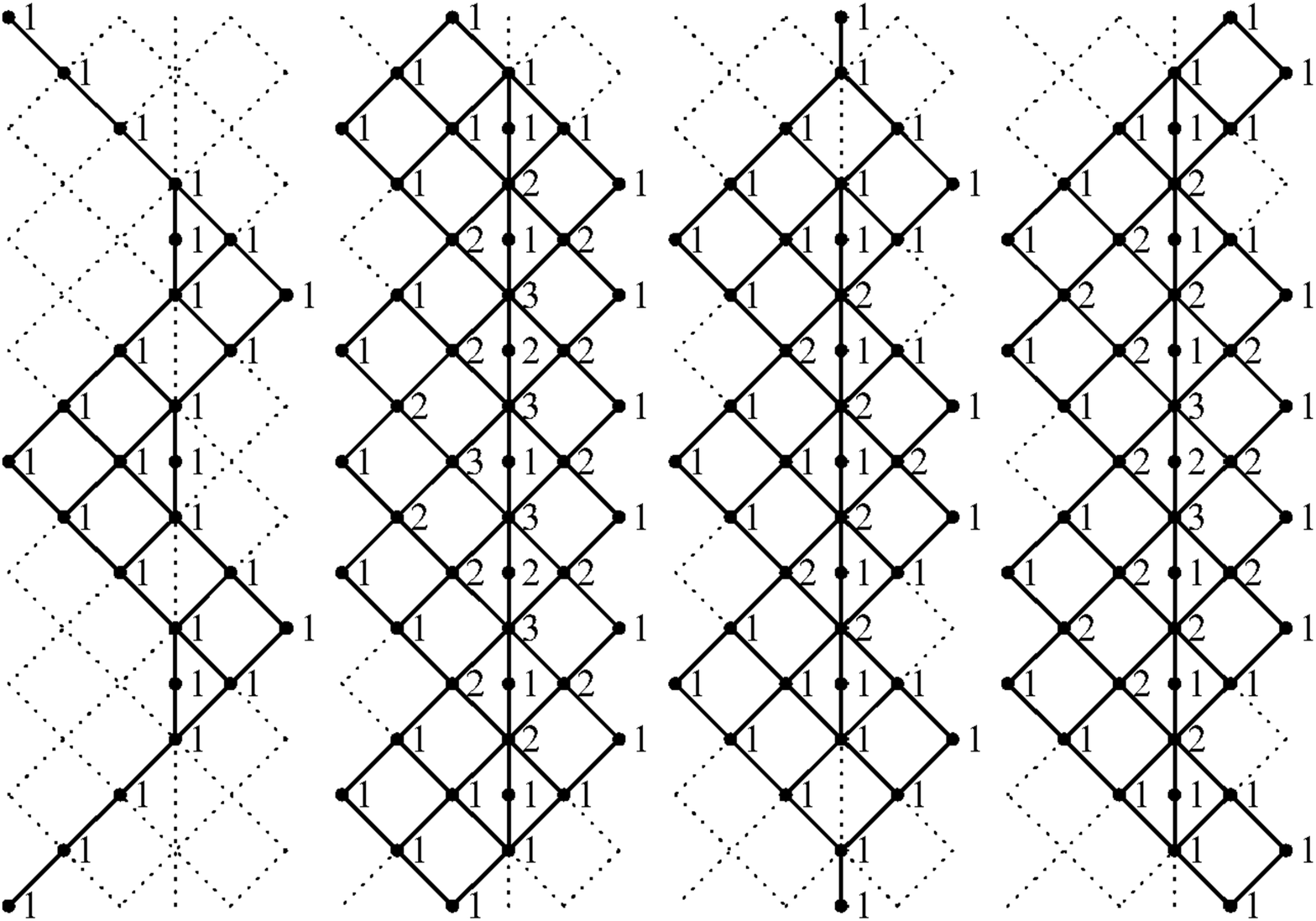}
\caption{Essential paths on the Coxeter graph $E_7$ (1)}
\label{EP-E7-1}
\end{figure}

%%% EP-E7-2 %%%%%%%%%%%%%
\begin{figure}[H]
%\centering\includegraphics[width=140mm,clip]{EP-E7-2.eps}
\centering\includegraphics[height=77mm,clip]{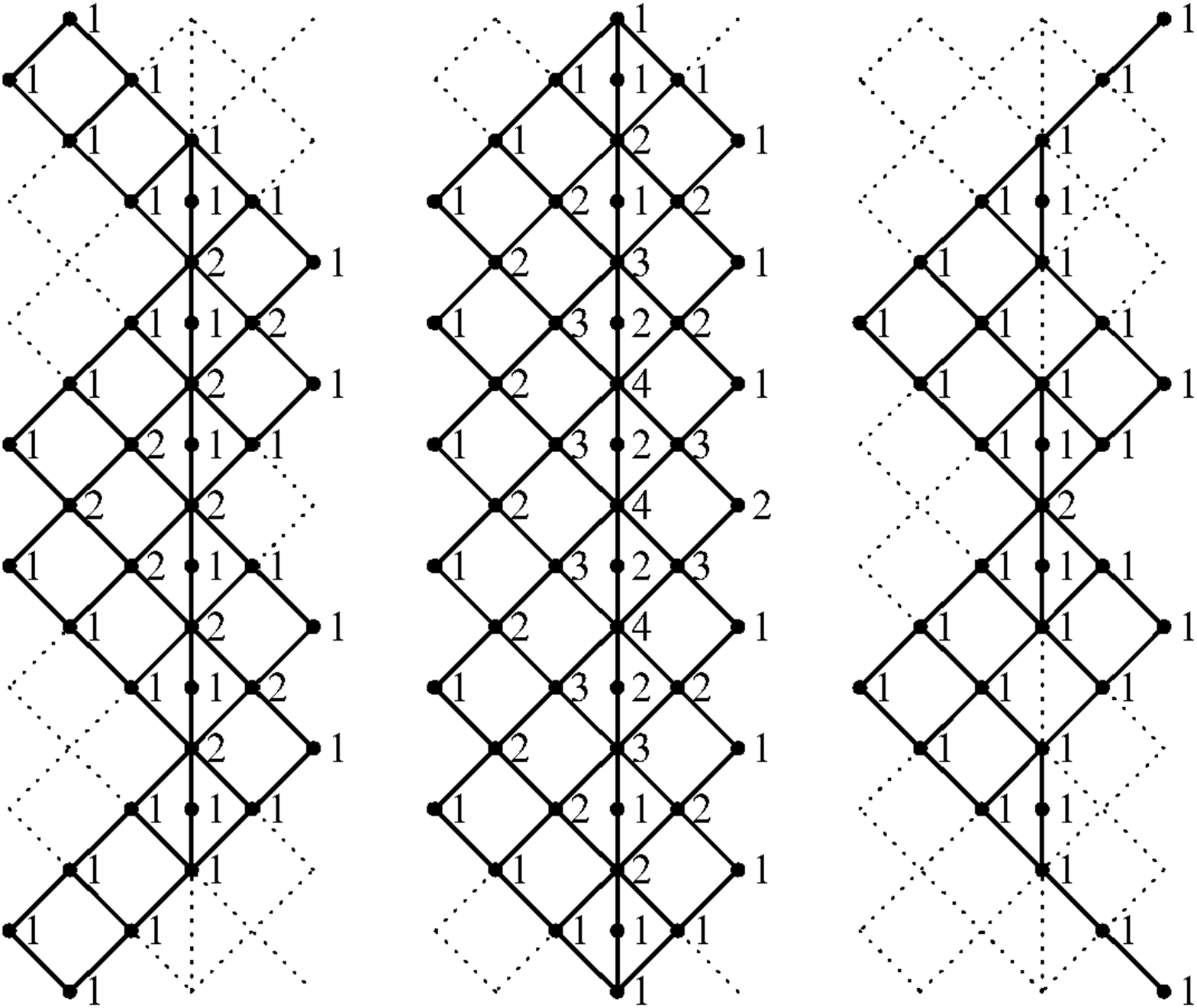}
\caption{Essential paths on the Coxeter graph $E_7$ (2)}
\label{EP-E7-2}
\end{figure}

%%% EP-E8-1 %%%%%%%%%%%%%
\begin{figure}[H]
\centering\includegraphics[width=140mm,clip]{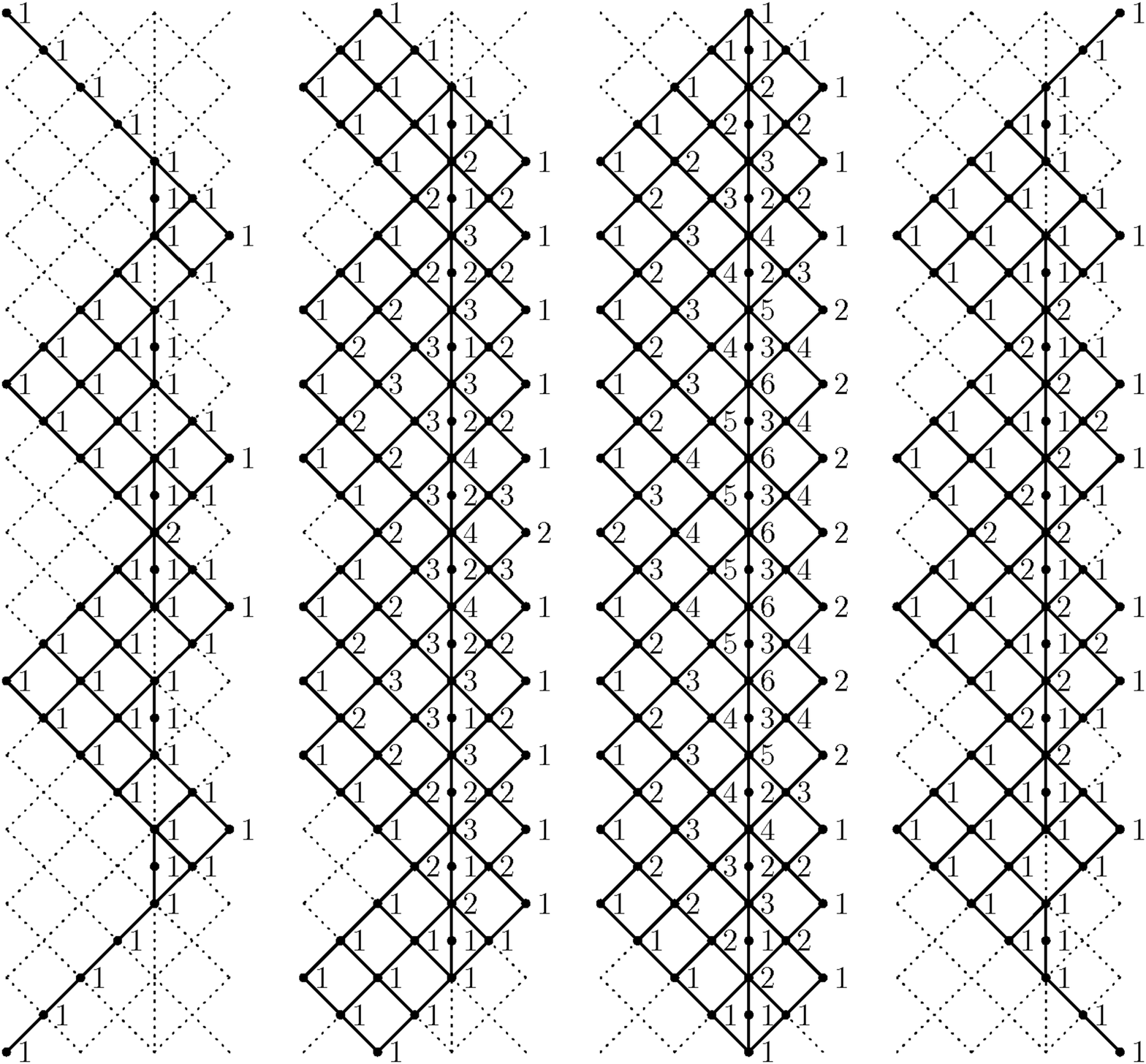}
\caption{Essential paths on the Coxeter graph $E_8$ (1)}
\label{EP-E8-1}
\end{figure}

%%% EP-E8-2 %%%%%%%%%%%%%

%%%%%%%%%%%%%%%%%%%%%%%%%%%%%%%%%%%%%%%%%%%%%%%%%%%%%%%%%%%%
%%%%%%%%%%%%%%%%%%%%%%%  Q-SYM  %%%%%%%%%%%%%%%%%%%%%%%%%%%% 
%%%%%%%%%%%%%%%%%%%%%%%%%%%%%%%%%%%%%%%%%%%%%%%%%%%%%%%%%%%%

%%% Q-SYM A %%%%%%%%%%%%%

%%% Q-SYM D_{even} %%%%%%%%%%%%%

%%% Q-SYM D_{odd} %%%%%%%%%%%%%

%%% Q-SYM E_6 %%%%%%%%%%%%%

%%% Q-SYM E_7 %%%%%%%%%%%%%

%%% Q-SYM E_8 %%%%%%%%%%%%%

%%%%%%%%%%%%%%%%%%%%%%%%%%%%%%%%%%%%%%%%%%%%%%%%%%
%%% The vertical edges of $K$-$K$ connections %%%%
%%%%%%%%%%%%%%%%%%%%%%%%%%%%%%%%%%%%%%%%%%%%%%%%%%

%%%%% Graphs for connections on $A_3$ %%%%%%%%%%%%

%%%%% Graphs for connections on $A_4$ %%%%%%%%%%%%
\begin{figure}[H]
\centering\includegraphics[width=100mm,clip]{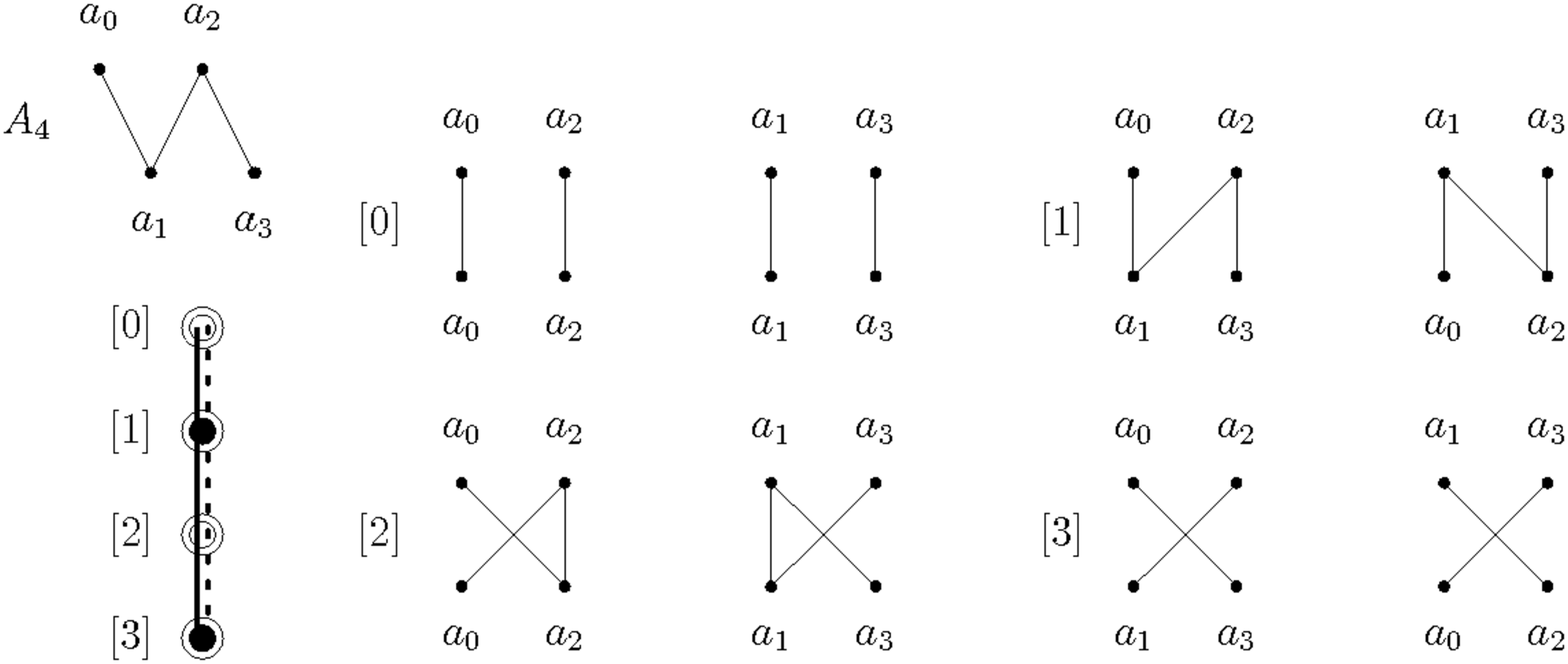}
\caption{Vertical graphs for connections on 
the Coxeter graph $A_4$}
\label{v-edge-A4}
\end{figure}

%%%%% Graphs for connections on $A_5$ %%%%%%%%%%%%
\begin{figure}[H]
\centering\includegraphics[width=100mm,clip]{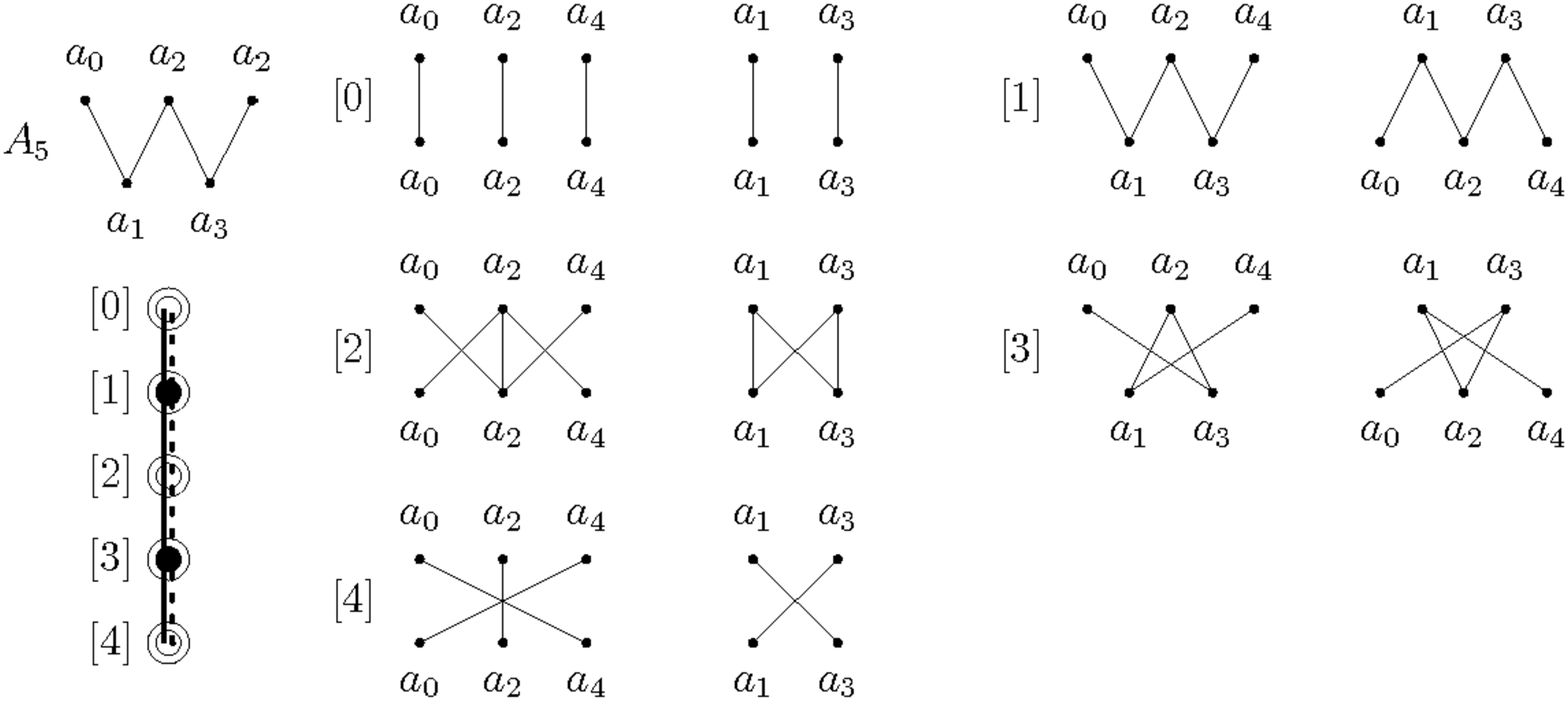}
\caption{Vertical graphs for connections on 
the Coxeter graph $A_5$}
\label{v-edge-A5}
\end{figure}

%%%%% Graphs for connections on $A_6$ %%%%%%%%%%%%
\begin{figure}[H]
\centering\includegraphics[width=130mm,clip]{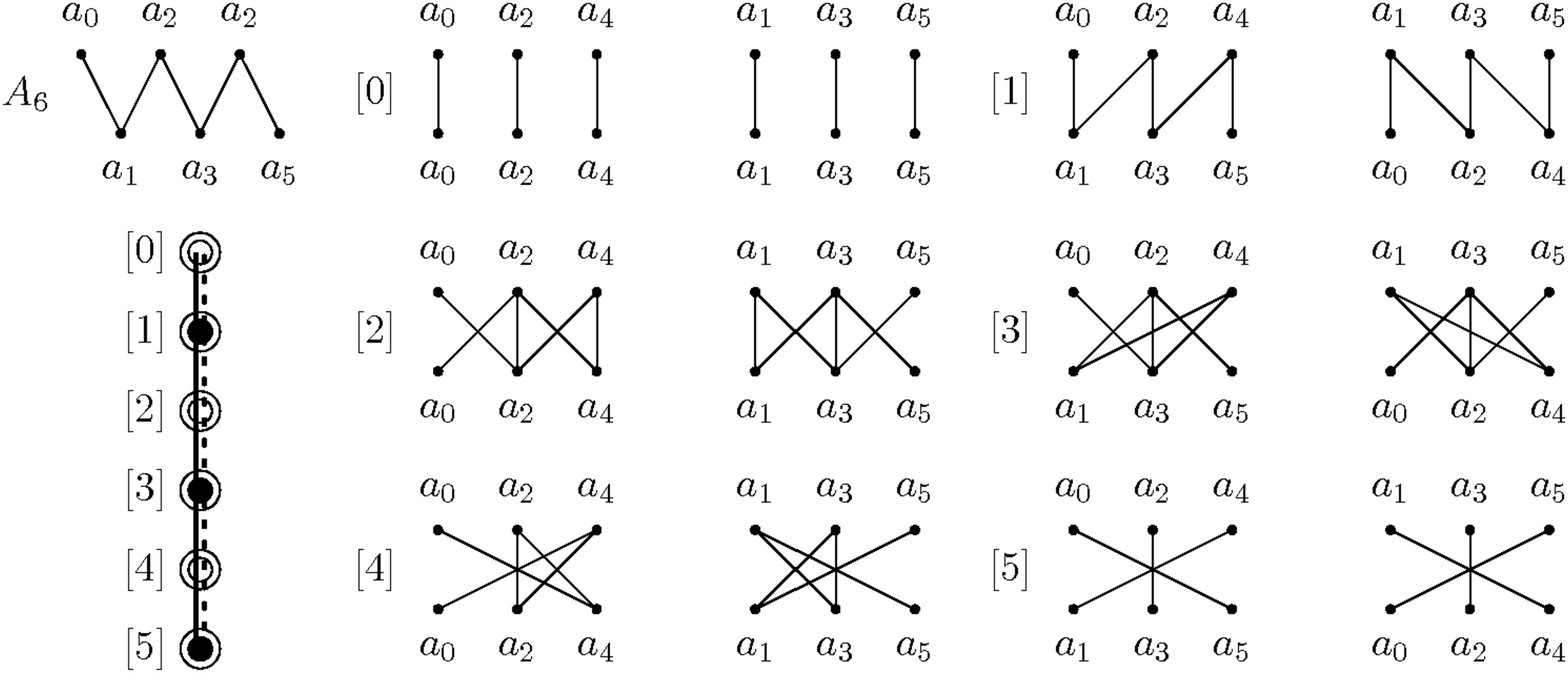}
\caption{Vertical graphs for connections on 
the Coxeter graph $A_6$}
\label{v-edge-A6}
\end{figure}

%%%%% Graphs for connections on $D_4$ %%%%%%%%%%%%

%%%%% Graphs for connections on $D_5$ %%%%%%%%%%%%

\vspace{-3mm}
%%% Graphs for connections on $D_6$ %%%%%%%%%%%%%

%%%%% vedge-e6 %%%%%%%%%%%%%%%%%%%%%%%%%%%%%%%%%%%

%%%%% vedge-e7 %%%%%%%%%%%%%%%%%%%%%%%%%%%%%%%%%%%

%%%%% vedge-e8e %%%%%%%%%%%%%%%%%%%%%%%%%%%%%%%%%%

%%%%% vedge-e8o %%%%%%%%%%%%%%%%%%%%%%%%%%%%%%%%%%

%%%%%%%%%%%%%%%%%%%%%%%%%%%%%%%%%%%%%%%%%%%%%%%%%%
%%%%%%%%%%%%%%%%%%%%%%%%%%%%%%%%%%%%%%%%%%%%%%%%%%
%%% The fusion tables of Ocneanu graphs       %%%%
%%%%%%%%%%%%%%%%%%%%%%%%%%%%%%%%%%%%%%%%%%%%%%%%%%
%%%%%%%%%%%%%%%%%%%%%%%%%%%%%%%%%%%%%%%%%%%%%%%%%%

%%%%% Fusion-Table-E6 %%%%%%%%%%%%%%%%%%%%%%%%%%%%

%%%%% Fusion-Table-E7 %%%%%%%%%%%%%%%%%%%%%%%%%%%%

%%%%% Fusion-Table-E8 %%%%%%%%%%%%%%%%%%%%%%%%%%%%

%%%%%%%%%%%%%%%%%%%%%%%%%%%%%%%%%%%%%%%%%%%%%%%%%%
%%% The (dual) pincipal graphs of             %%%%
%%%     Goodman-de la Harpe-Jones subfactors  %%%%
%%%%%%%%%%%%%%%%%%%%%%%%%%%%%%%%%%%%%%%%%%%%%%%%%%

%%%%%%%%%%%%%%%%%%%%%%%%%%%%%%%%%%%%%%%%%%%%%%%%%%
%%% The (dual) pincipal graphs of             %%%%
%%%     Goodman-de la Harpe-Jones subfactors  %%%%
%%%%%%%%%%%%%%%%%%%%%%%%%%%%%%%%%%%%%%%%%%%%%%%%%%
%%%       D5                                  %%%%
%%%%%%%%%%%%%%%%%%%%%%%%%%%%%%%%%%%%%%%%%%%%%%%%%%

%%%%% GHJ(D5-d1) %%%%%%%%%%%%%%%%%%%%%%%%%%%%%%%%%

%%%%% GHJ(D5-d2) %%%%%%%%%%%%%%%%%%%%%%%%%%%%%%%%%
\begin{figure}[H]
\centering
\includegraphics[width=140mm,clip]{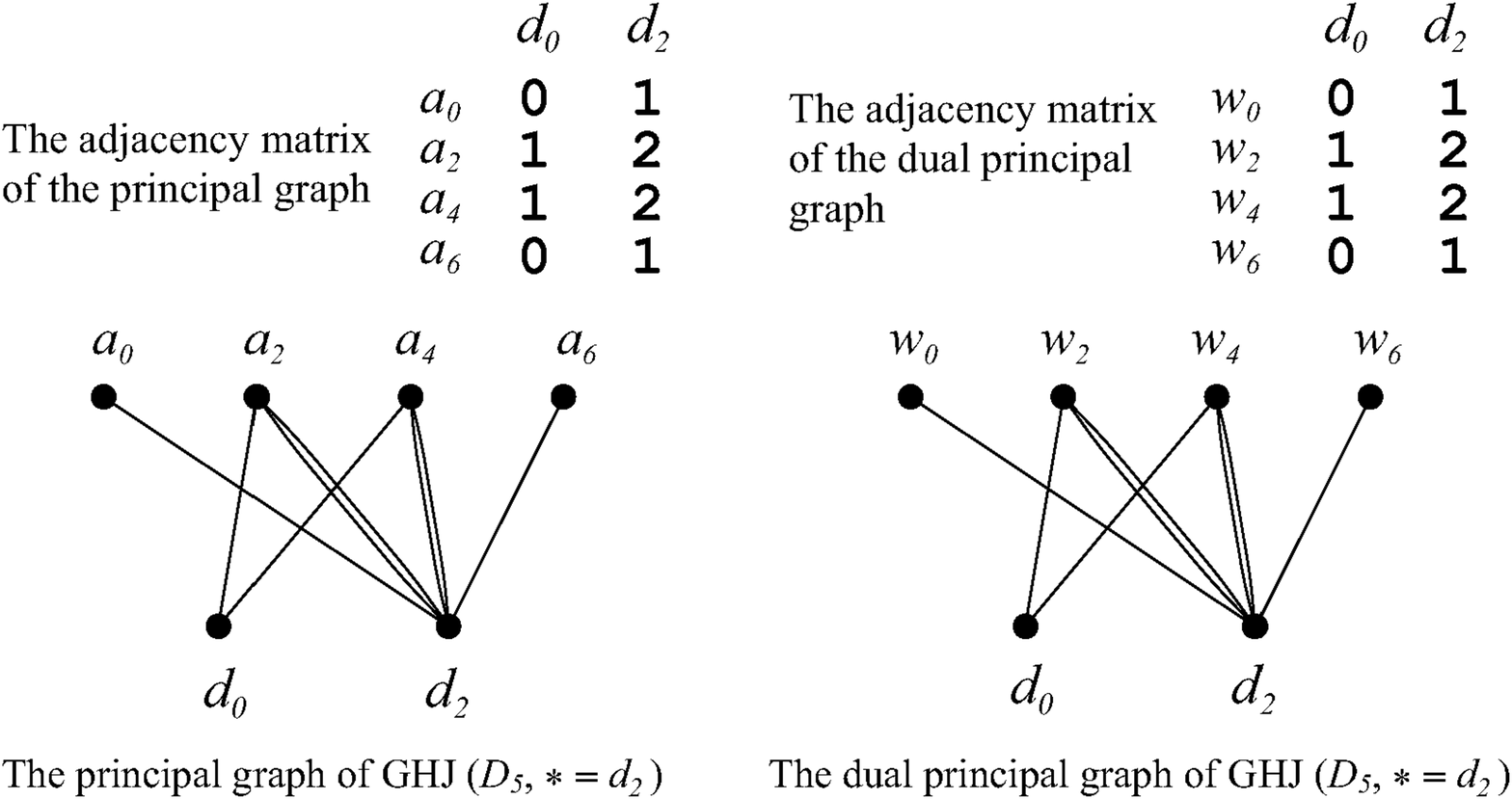}
\caption{The (dual) principal graph of GHJ$(D_5, *=d_2)$.}
\label{GHJ(D5-d2)}
\end{figure}

%%%%% GHJ(D5-d3) %%%%%%%%%%%%%%%%%%%%%%%%%%%%%%%%%
\begin{figure}[H]
\centering
\includegraphics[width=140mm,clip]{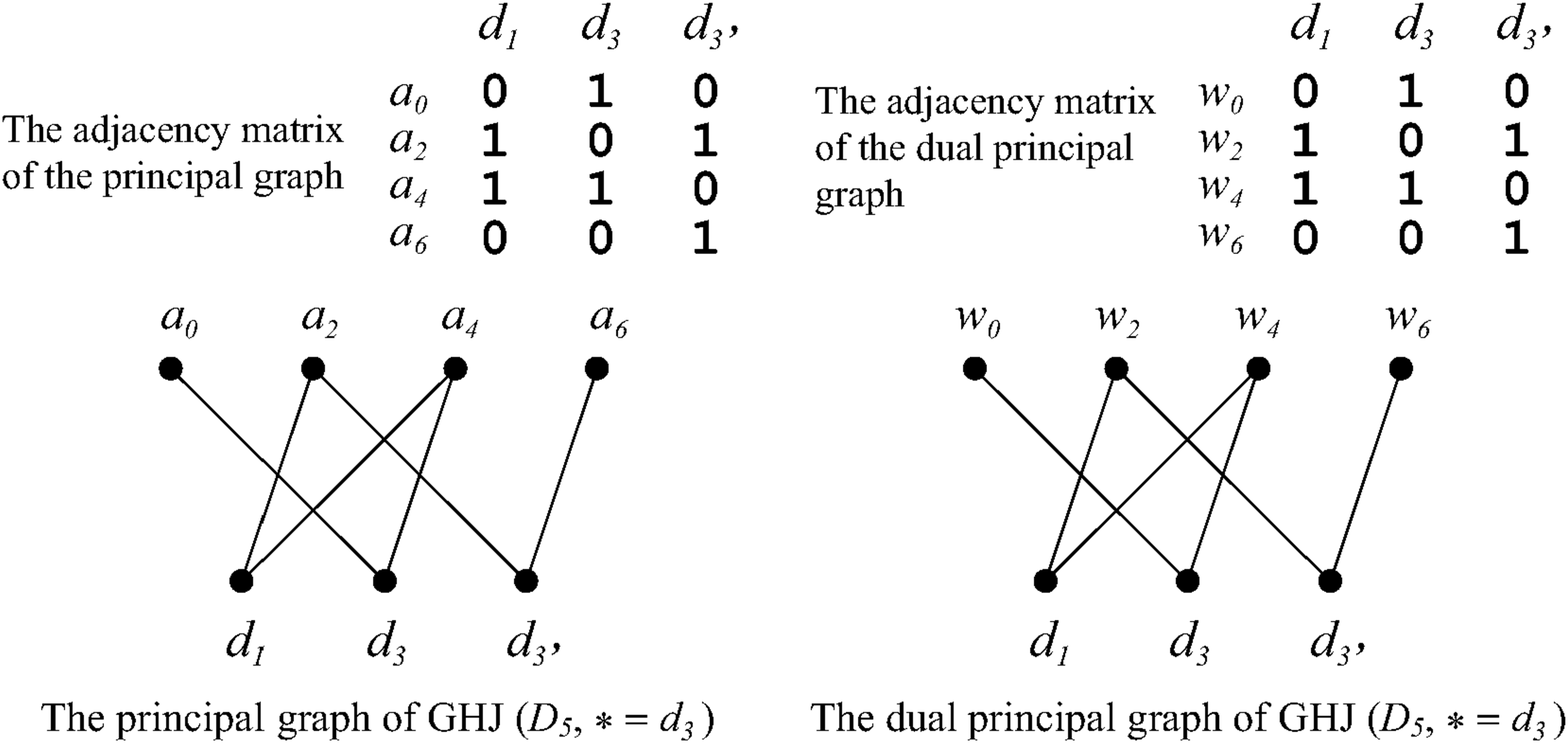}
\caption{The (dual) principal graph of GHJ$(D_5, *=d_3)$.}
\label{GHJ(D5-d3)}
\end{figure}

%%%%%%%%%%%%%%%%%%%%%%%%%%%%%%%%%%%%%%%%%%%%%%%%%%
%%% The (dual) pincipal graphs of             %%%%
%%%     Goodman-de la Harpe-Jones subfactors  %%%%
%%%%%%%%%%%%%%%%%%%%%%%%%%%%%%%%%%%%%%%%%%%%%%%%%%
%%%       D7                                  %%%%
%%%%%%%%%%%%%%%%%%%%%%%%%%%%%%%%%%%%%%%%%%%%%%%%%%

%%%%% GHJ(D7-d1) %%%%%%%%%%%%%%%%%%%%%%%%%%%%%%%%%
\begin{figure}[H]
\centering
\includegraphics[width=140mm,clip]{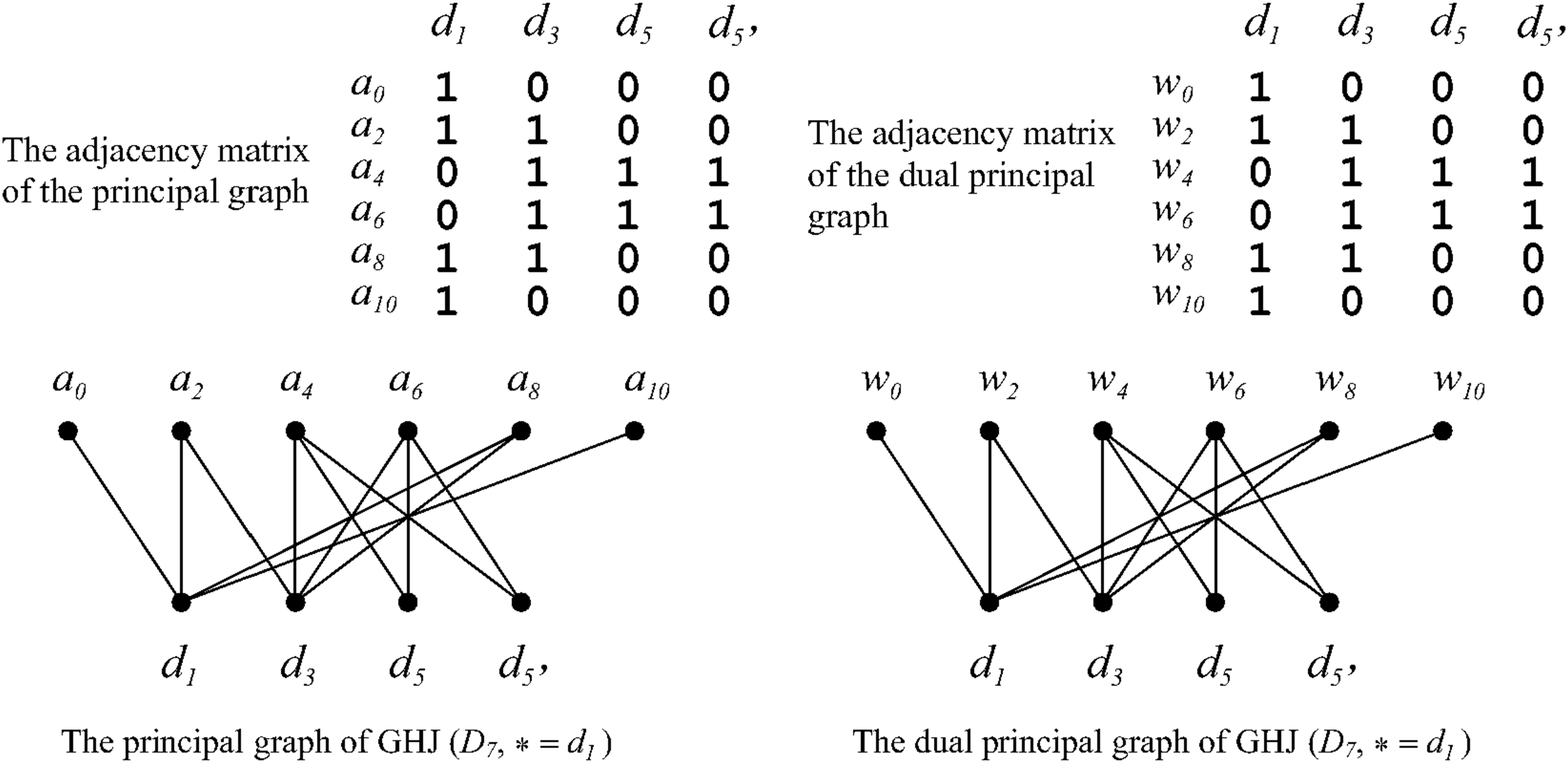}
\caption{The (dual) principal graph of GHJ$(D_7, *=d_1)$.}
\label{GHJ(D7-d1)}
\end{figure}

%%%%% GHJ(D7-d2) %%%%%%%%%%%%%%%%%%%%%%%%%%%%%%%%%
\begin{figure}[H]
\centering
\includegraphics[width=140mm,clip]{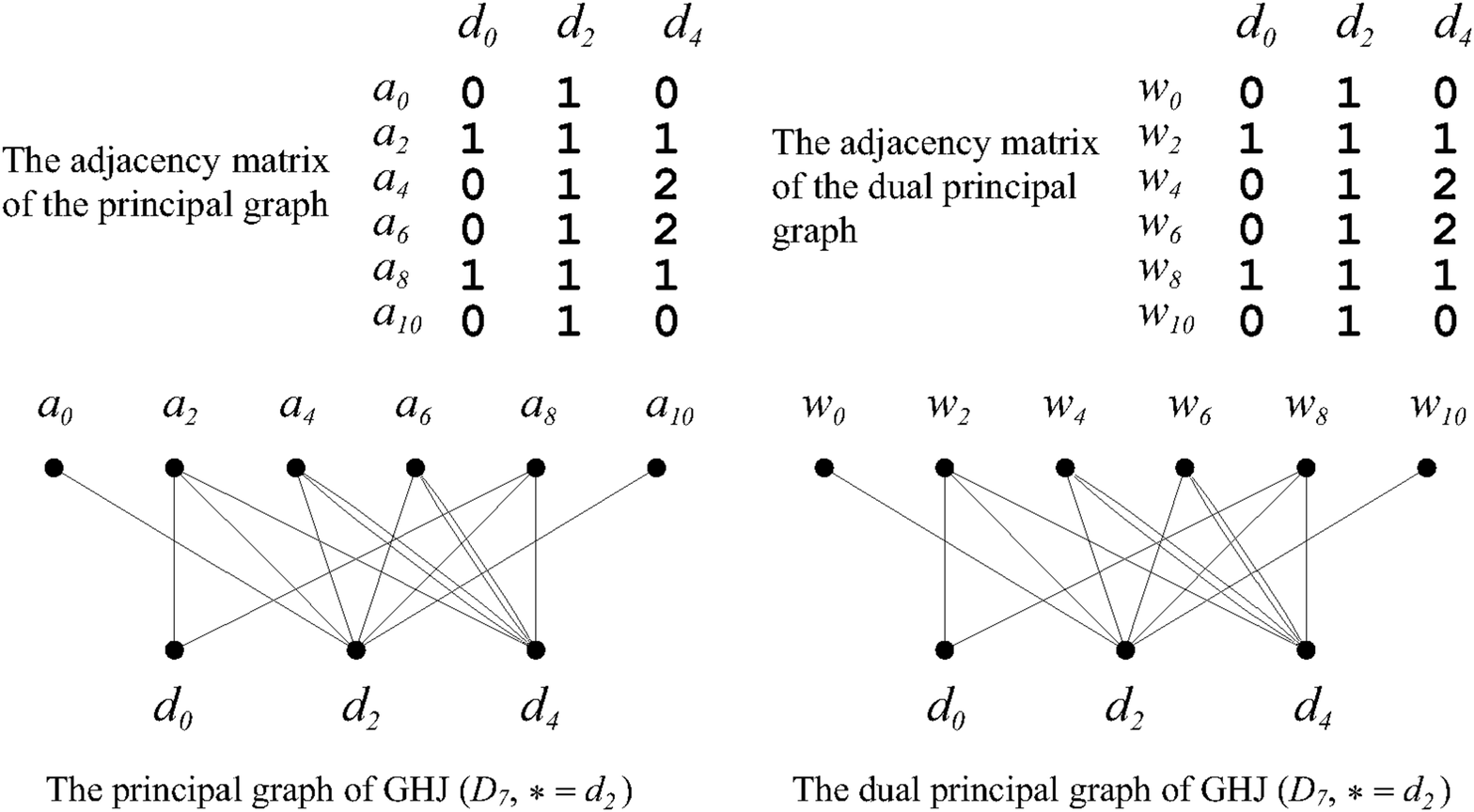}
\caption{The (dual) principal graph of GHJ$(D_7, *=d_2)$.}
\label{GHJ(D7-d2)}
\end{figure}

%%%%% GHJ(D7-d3) %%%%%%%%%%%%%%%%%%%%%%%%%%%%%%%%%
\begin{figure}[H]
\centering
\includegraphics[width=140mm,clip]{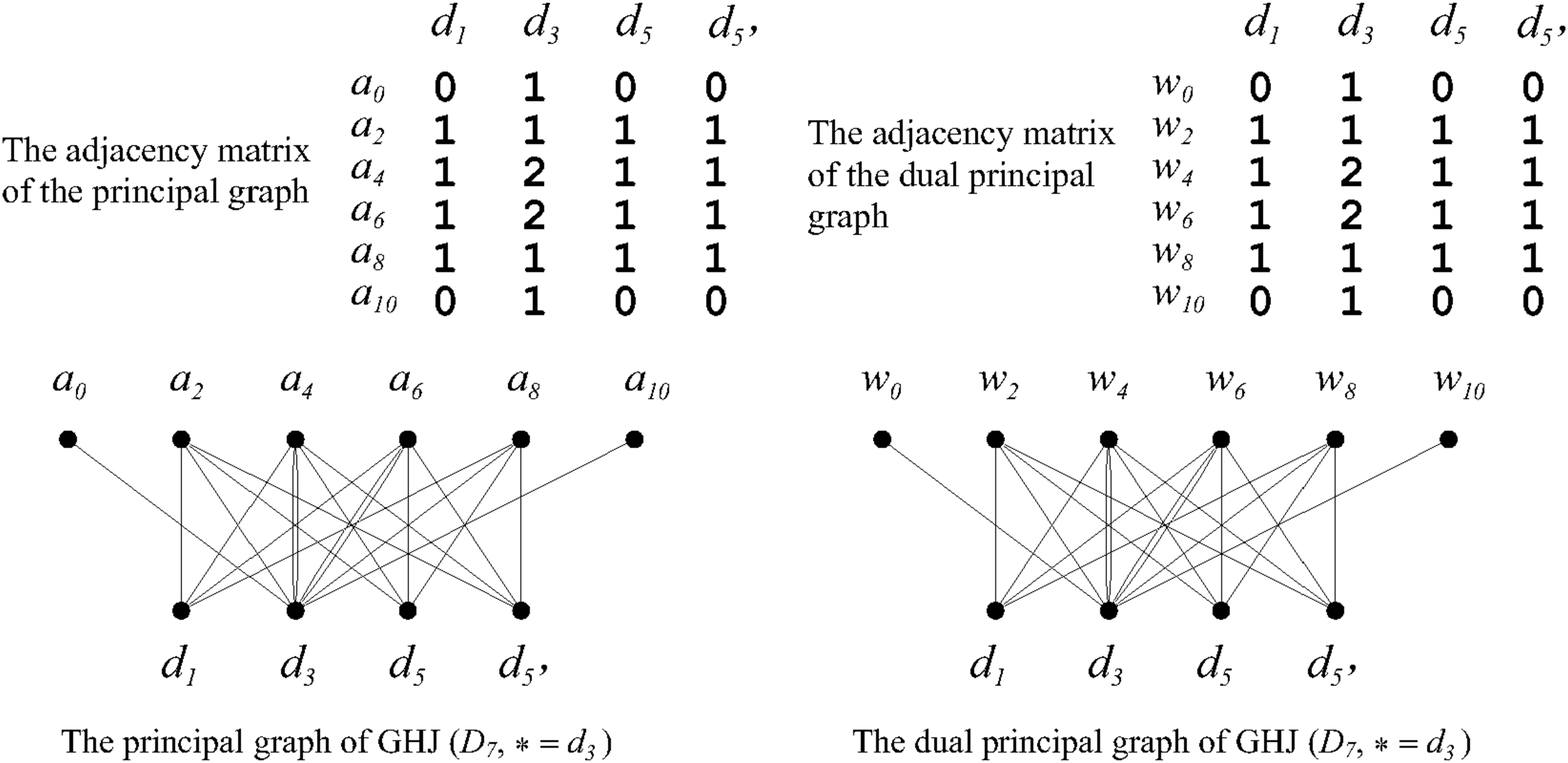}
\caption{The (dual) principal graph of GHJ$(D_7, *=d_3)$.}
\label{GHJ(D7-d3)}
\end{figure}

%%%%% GHJ(D7-d4) %%%%%%%%%%%%%%%%%%%%%%%%%%%%%%%%%
\begin{figure}[H]
\centering
\includegraphics[width=140mm,clip]{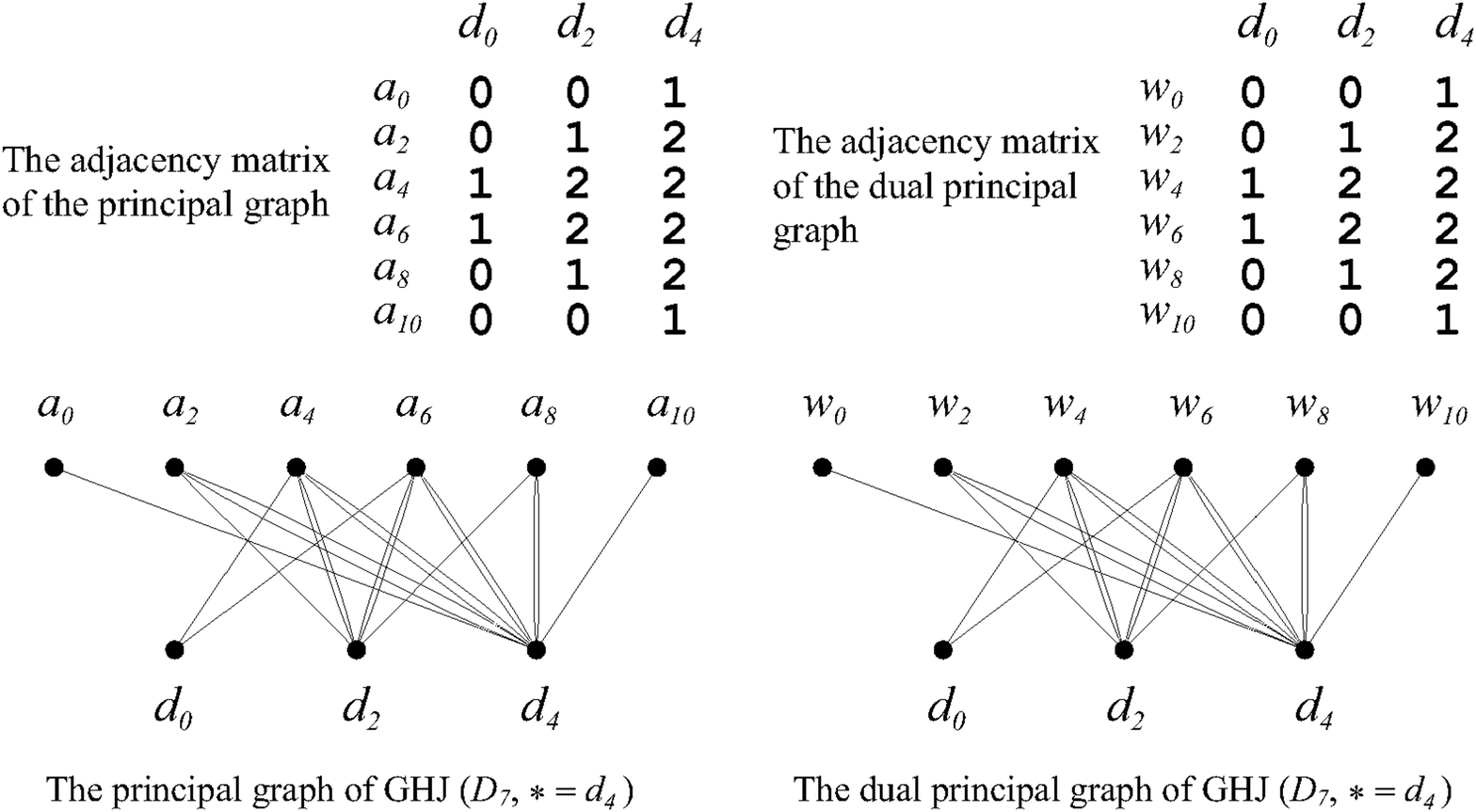}
\caption{The (dual) principal graph of GHJ$(D_7, *=d_4)$.}
\label{GHJ(D7-d4)}
\end{figure}

%%%%% GHJ(D7-d5) %%%%%%%%%%%%%%%%%%%%%%%%%%%%%%%%%
\begin{figure}[H]
\centering
\includegraphics[width=140mm,clip]{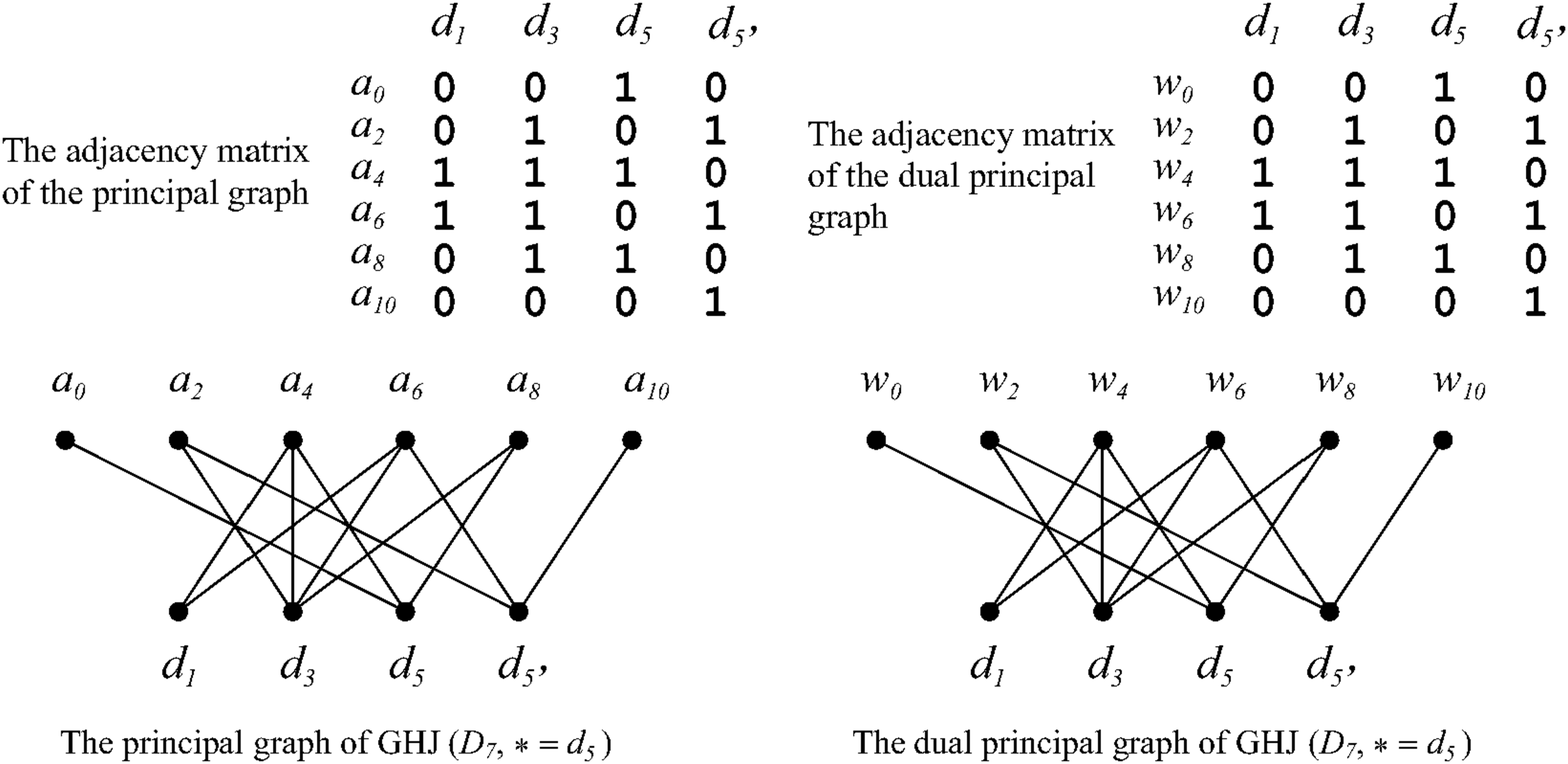}
\caption{The (dual) principal graph of GHJ$(D_7, *=d_5)$.}
\label{GHJ(D7-d5)}
\end{figure}

%%%%%%%%%%%%%%%%%%%%%%%%%%%%%%%%%%%%%%%%%%%%%%%%%%
%%% The (dual) pincipal graphs of             %%%%
%%%     Goodman-de la Harpe-Jones subfactors  %%%%
%%%%%%%%%%%%%%%%%%%%%%%%%%%%%%%%%%%%%%%%%%%%%%%%%%
%%%       D9                                  %%%%
%%%%%%%%%%%%%%%%%%%%%%%%%%%%%%%%%%%%%%%%%%%%%%%%%%

%%%%% GHJ(D9-d1) %%%%%%%%%%%%%%%%%%%%%%%%%%%%%%%%%
\begin{figure}[H]
\centering
\includegraphics[width=130mm,clip]{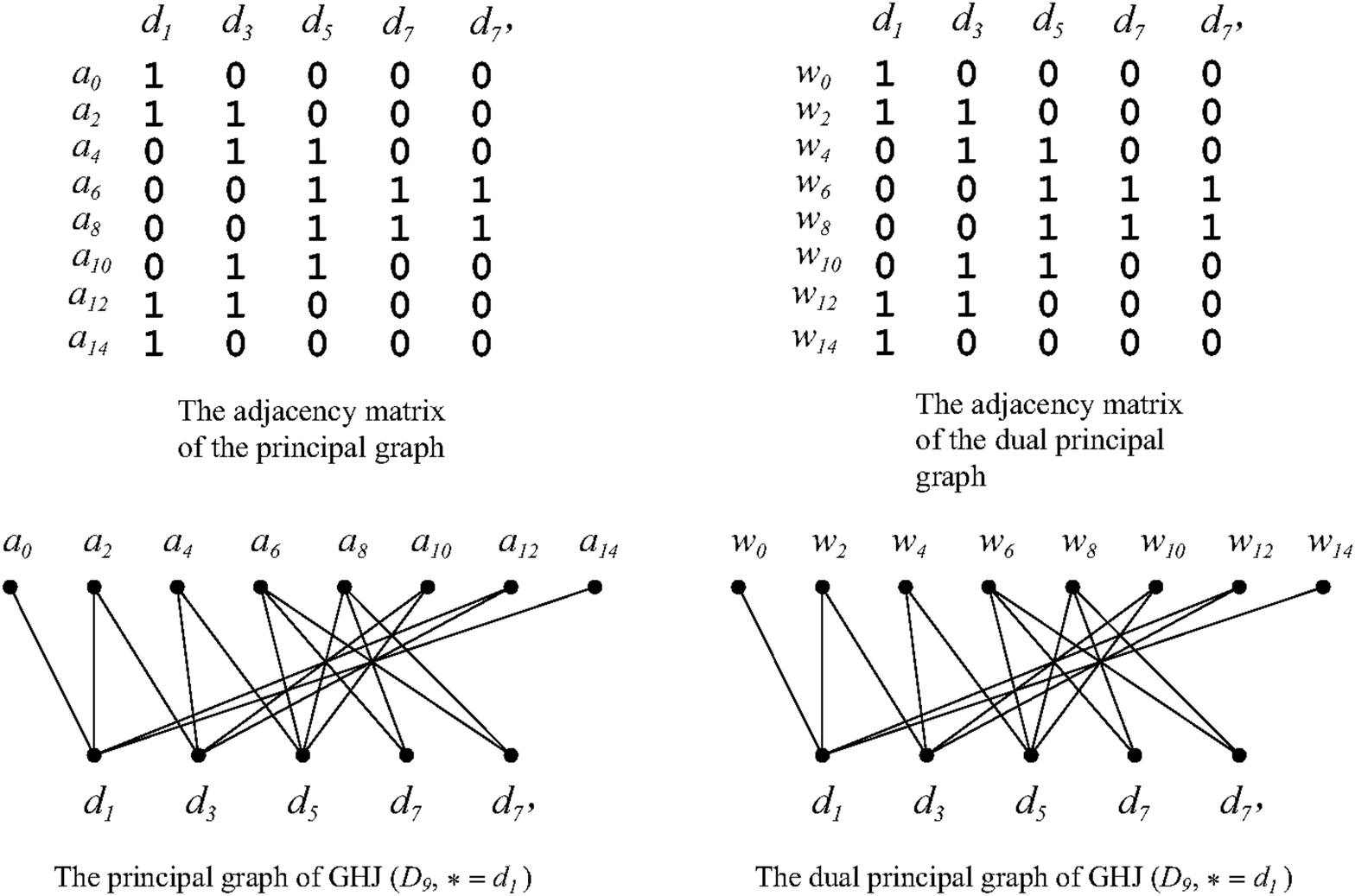}
\caption{The (dual) principal graph of GHJ$(D_9, *=d_1)$.}
\label{GHJ(D9-d1)}
\end{figure}

%%%%% GHJ(D9-d2) %%%%%%%%%%%%%%%%%%%%%%%%%%%%%%%%%
\begin{figure}[H]
\centering
\includegraphics[width=130mm,clip]{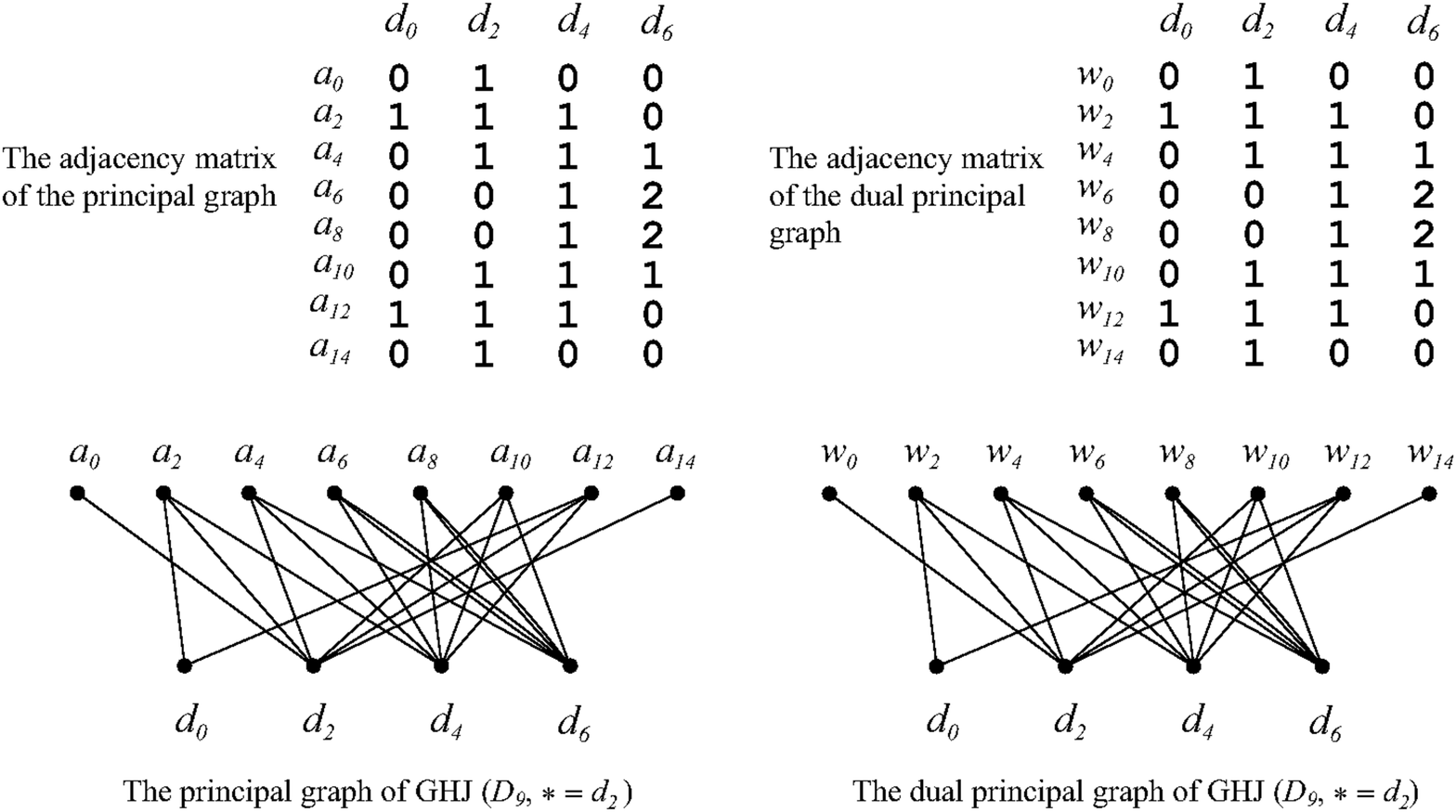}
\caption{The (dual) principal graph of GHJ$(D_9, *=d_2)$.}
\label{GHJ(D9-d2)}
\end{figure}

%%%%% GHJ(D9-d3) %%%%%%%%%%%%%%%%%%%%%%%%%%%%%%%%%
\begin{figure}[H]
\centering
\includegraphics[width=130mm,clip]{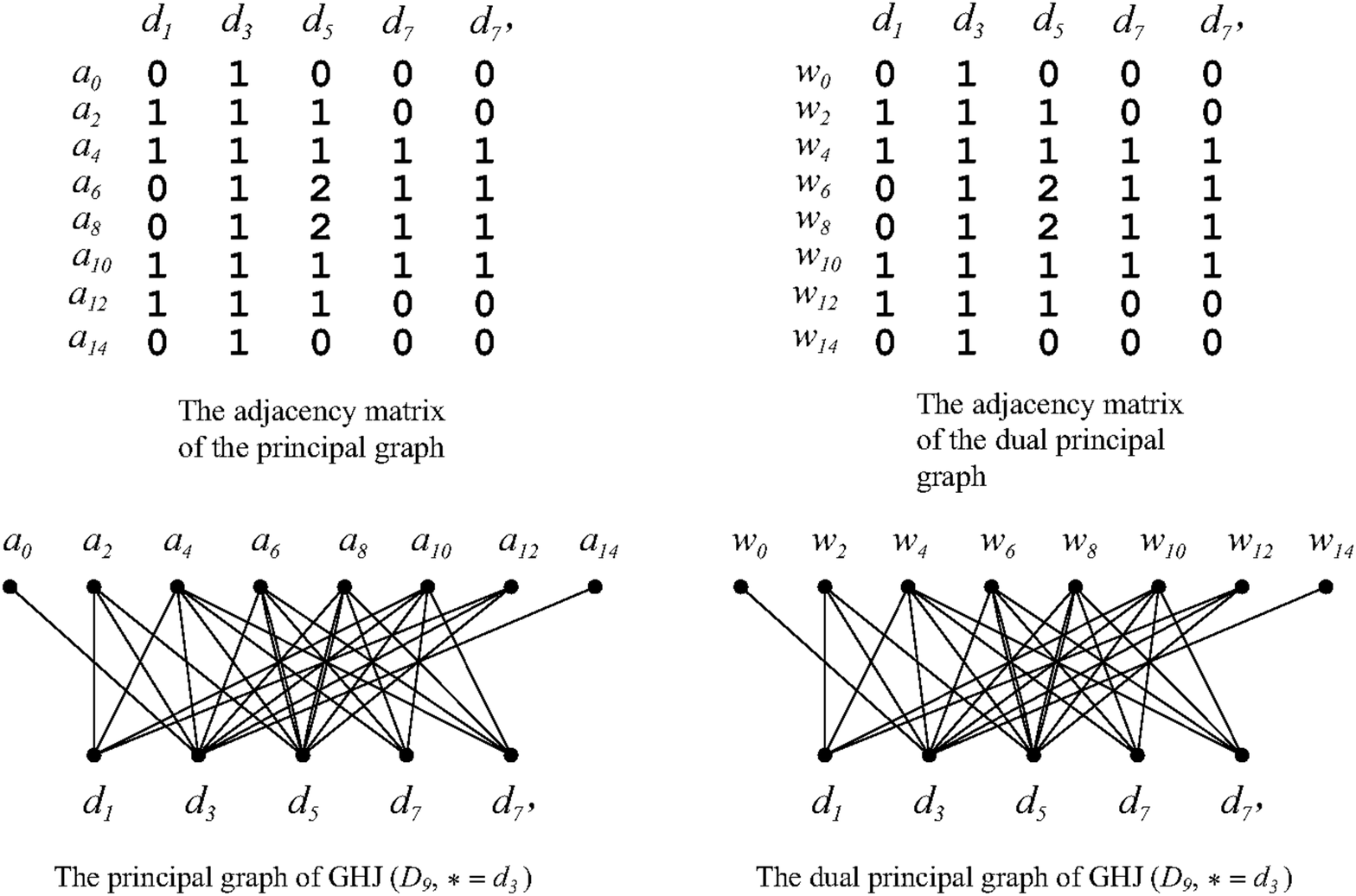}
\caption{The (dual) principal graph of GHJ$(D_9, *=d_3)$.}
\label{GHJ(D9-d3)}
\end{figure}

%%%%% GHJ(D9-d4) %%%%%%%%%%%%%%%%%%%%%%%%%%%%%%%%%
\begin{figure}[H]
\centering
\includegraphics[width=130mm,clip]{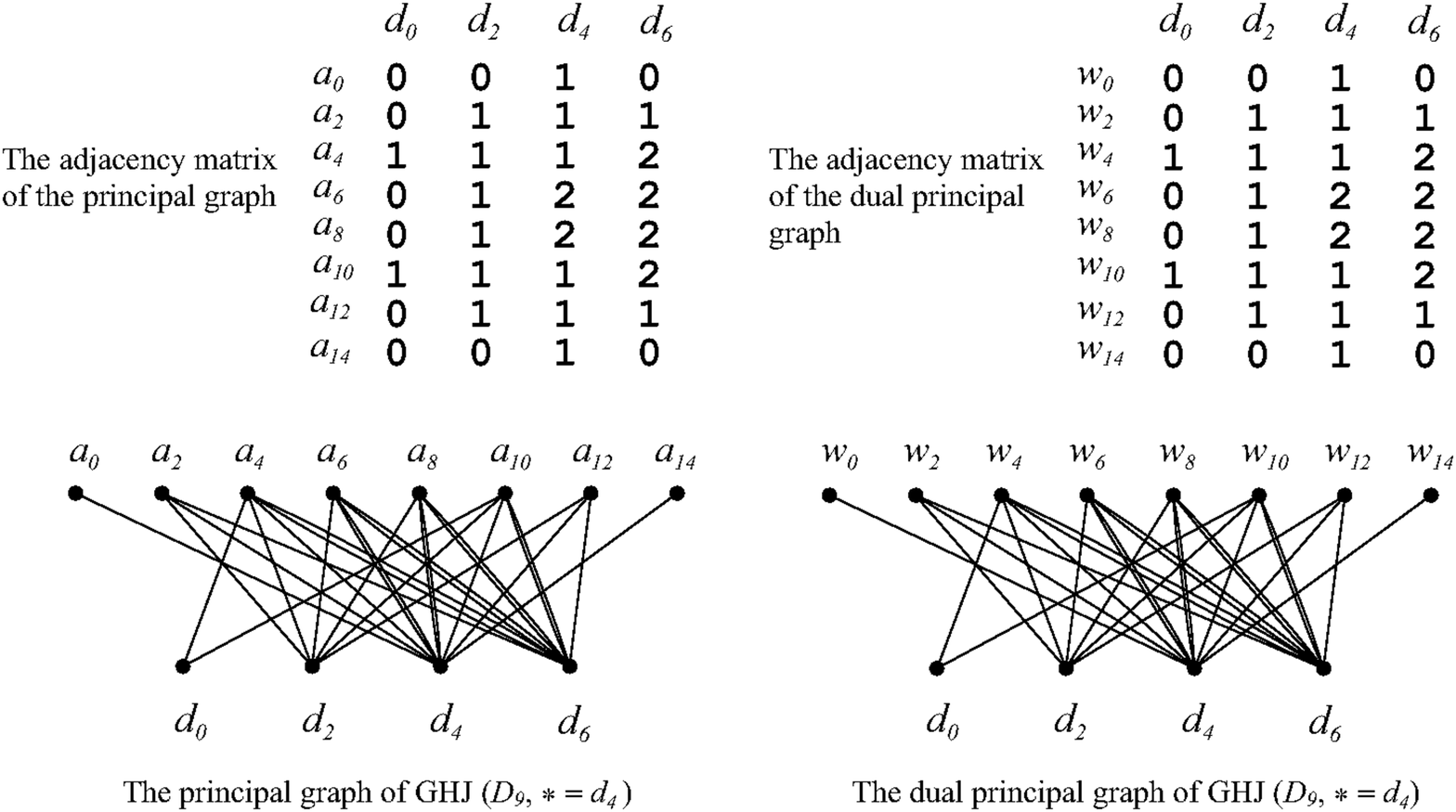}
\caption{The (dual) principal graph of GHJ$(D_9, *=d_4)$.}
\label{GHJ(D9-d4)}
\end{figure}

%%%%% GHJ(D9-d5) %%%%%%%%%%%%%%%%%%%%%%%%%%%%%%%%%
\begin{figure}[H]
\centering
\includegraphics[width=130mm,clip]{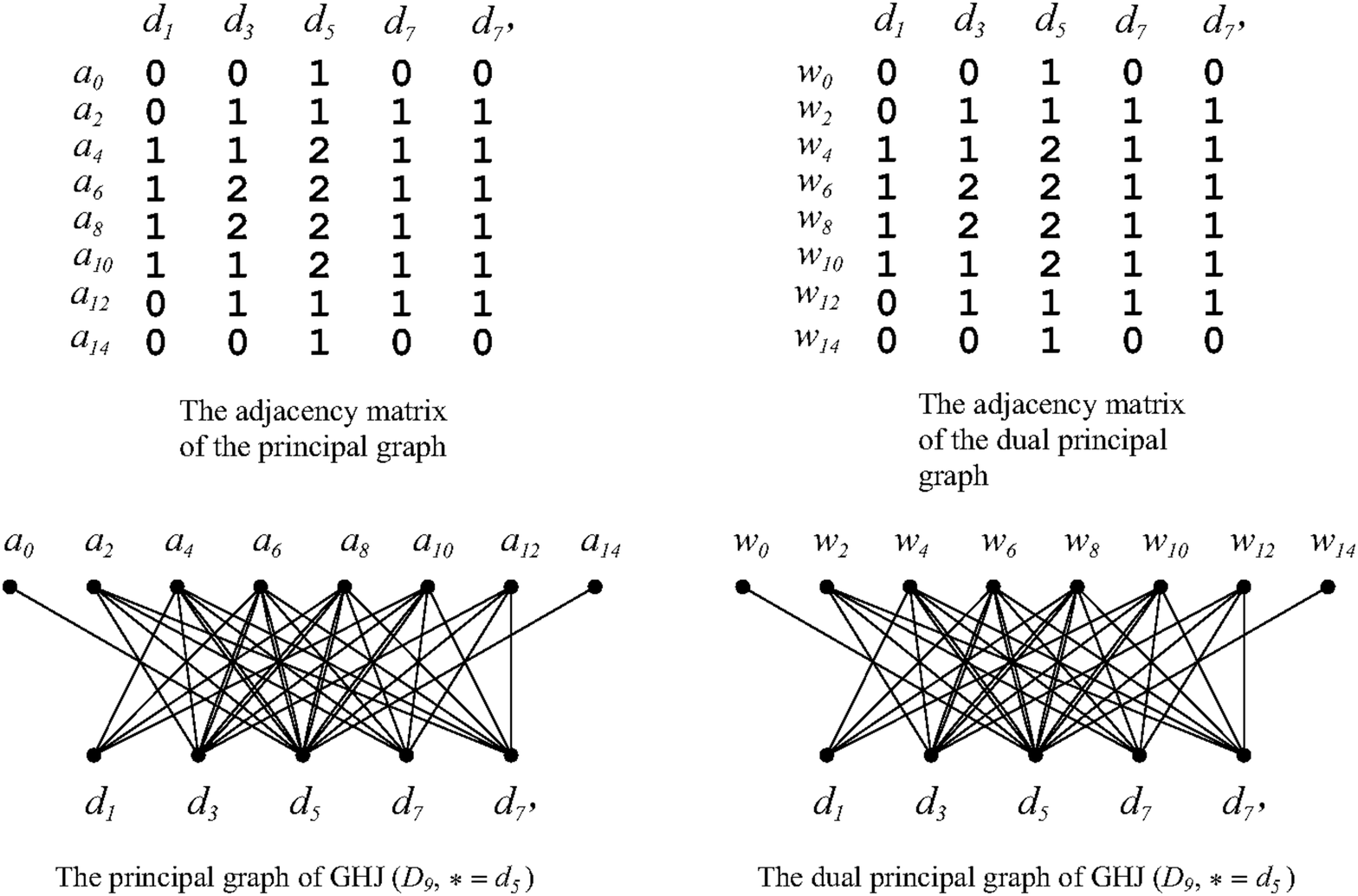}
\caption{The (dual) principal graph of GHJ$(D_9, *=d_5)$.}
\label{GHJ(D9-d5)}
\end{figure}

%%%%% GHJ(D9-d6) %%%%%%%%%%%%%%%%%%%%%%%%%%%%%%%%%
\begin{figure}[H]
\centering
\includegraphics[width=130mm,clip]{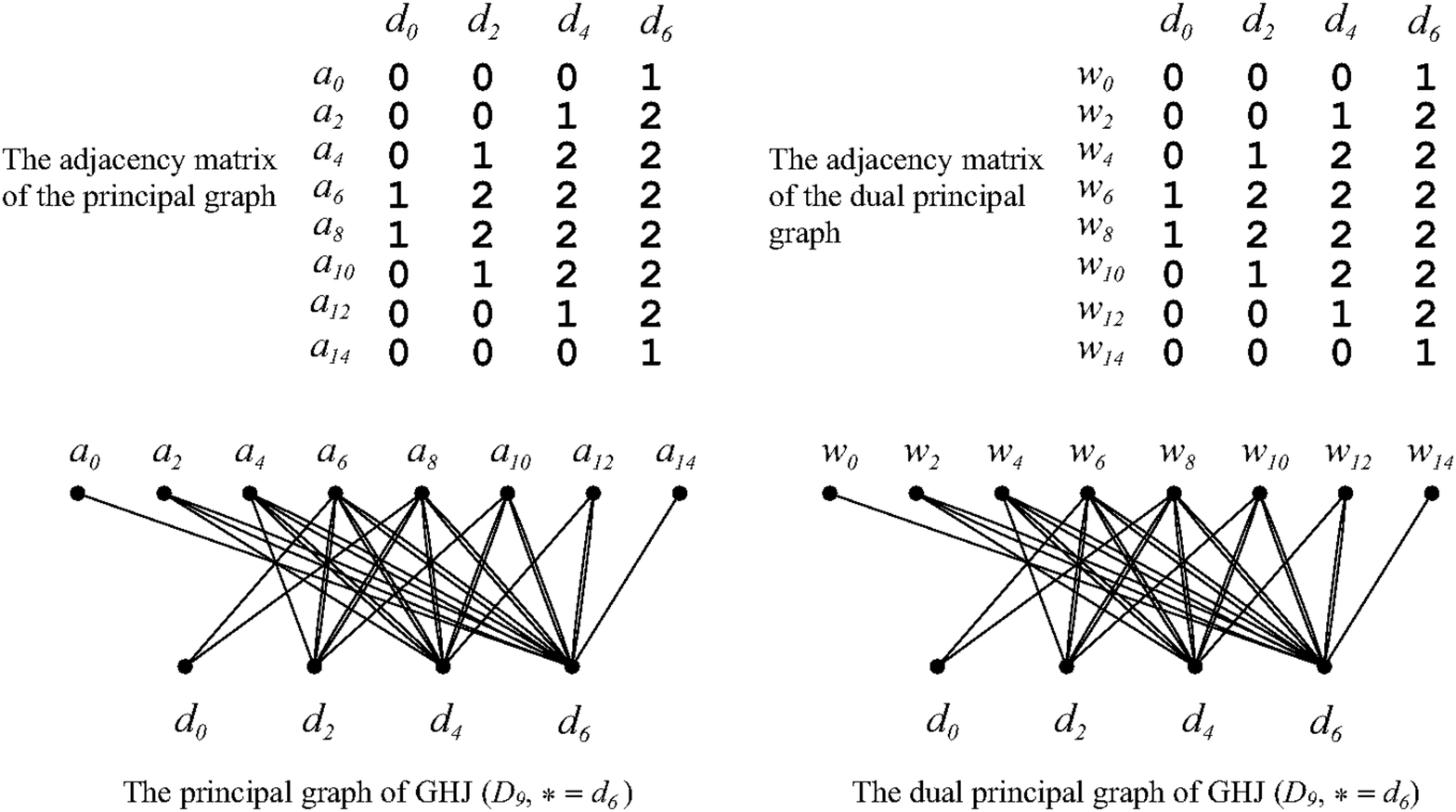}
\caption{The (dual) principal graph of GHJ$(D_9, *=d_6)$.}
\label{GHJ(D9-d6)}
\end{figure}

%%%%% GHJ(D9-d7) %%%%%%%%%%%%%%%%%%%%%%%%%%%%%%%%%
\begin{figure}[H]
\centering
\includegraphics[width=130mm,clip]{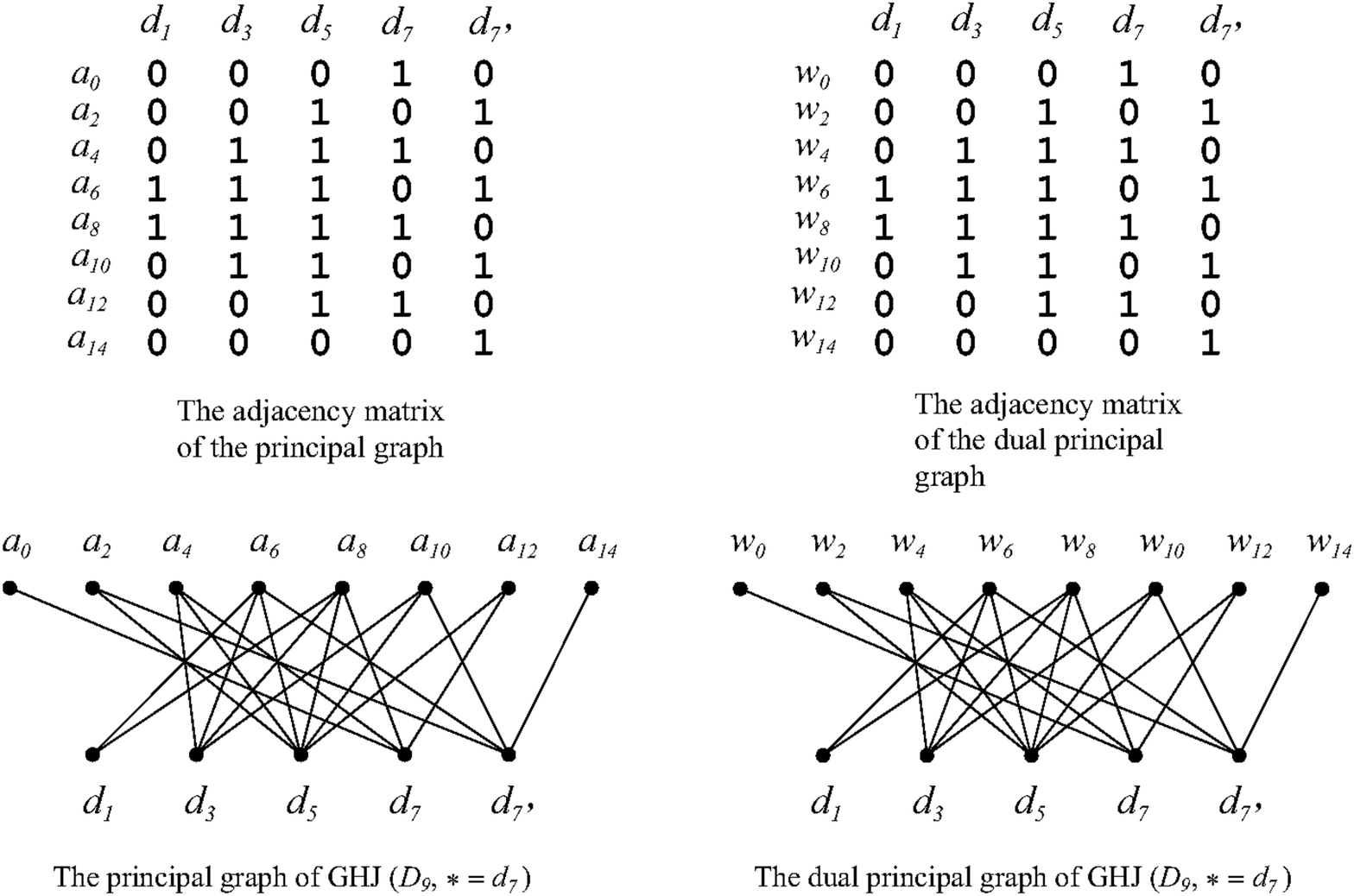}
\caption{The (dual) principal graph of GHJ$(D_9, *=d_7)$.}
\label{GHJ(D9-d7)}
\end{figure}

%%%%%%%%%%%%%%%%%%%%%%%%%%%%%%%%%%%%%%%%%%%%%%%%%%
%%% The (dual) pincipal graphs of             %%%%
%%%     Goodman-de la Harpe-Jones subfactors  %%%%
%%%%%%%%%%%%%%%%%%%%%%%%%%%%%%%%%%%%%%%%%%%%%%%%%%
%%%       D11                                  %%%%
%%%%%%%%%%%%%%%%%%%%%%%%%%%%%%%%%%%%%%%%%%%%%%%%%%

%%%%% GHJ(D11-d1) %%%%%%%%%%%%%%%%%%%%%%%%%%%%%%%%%
\begin{figure}[H]
\centering
\includegraphics[width=129mm,clip]{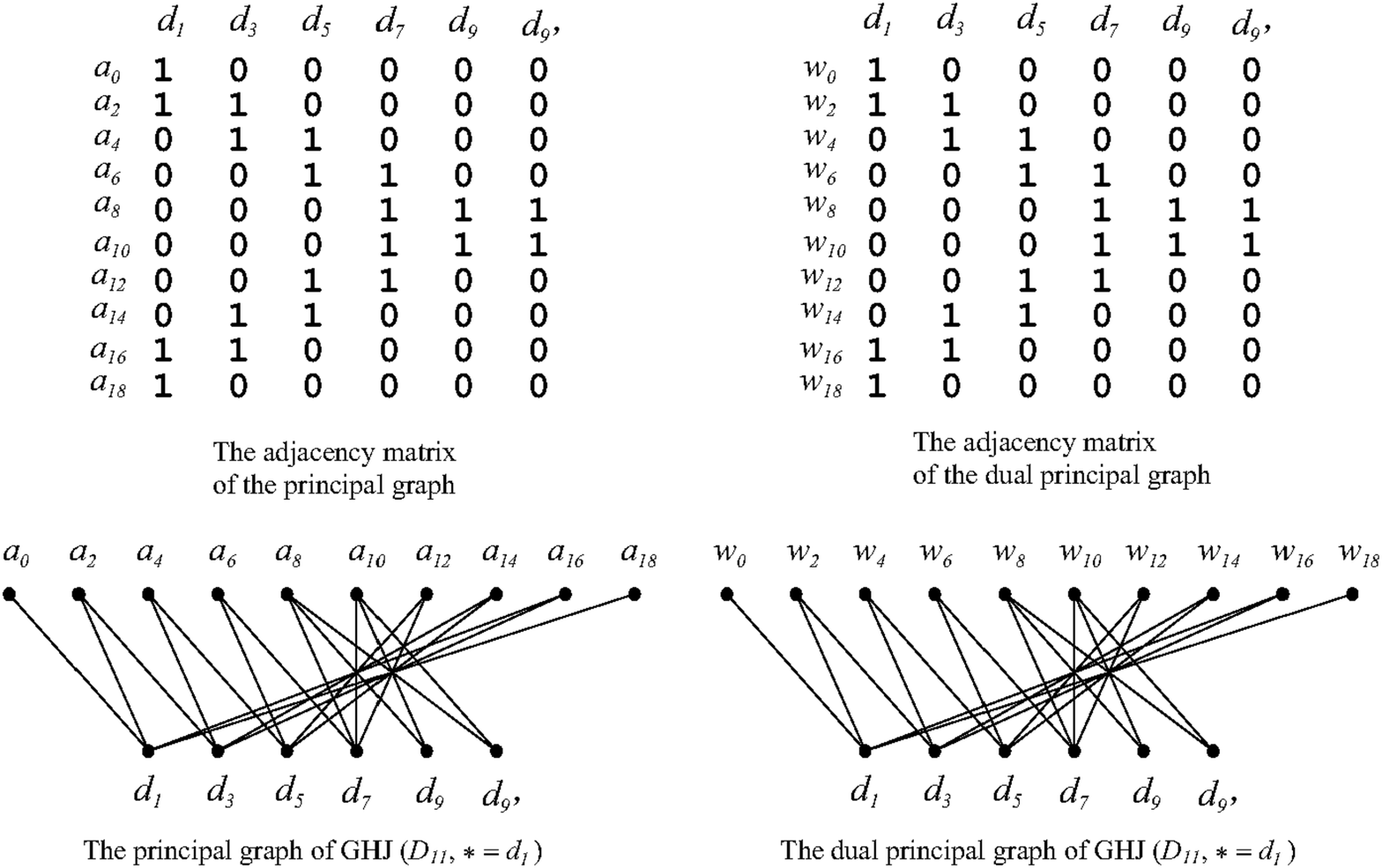}
\caption{The (dual) principal graph of GHJ$(D_{11}, *=d_1)$.}
\label{GHJ(D11-d1)}
\end{figure}

%%%%% GHJ(D11-d2) %%%%%%%%%%%%%%%%%%%%%%%%%%%%%%%%%
\begin{figure}[H]
\centering
\includegraphics[width=130mm,clip]{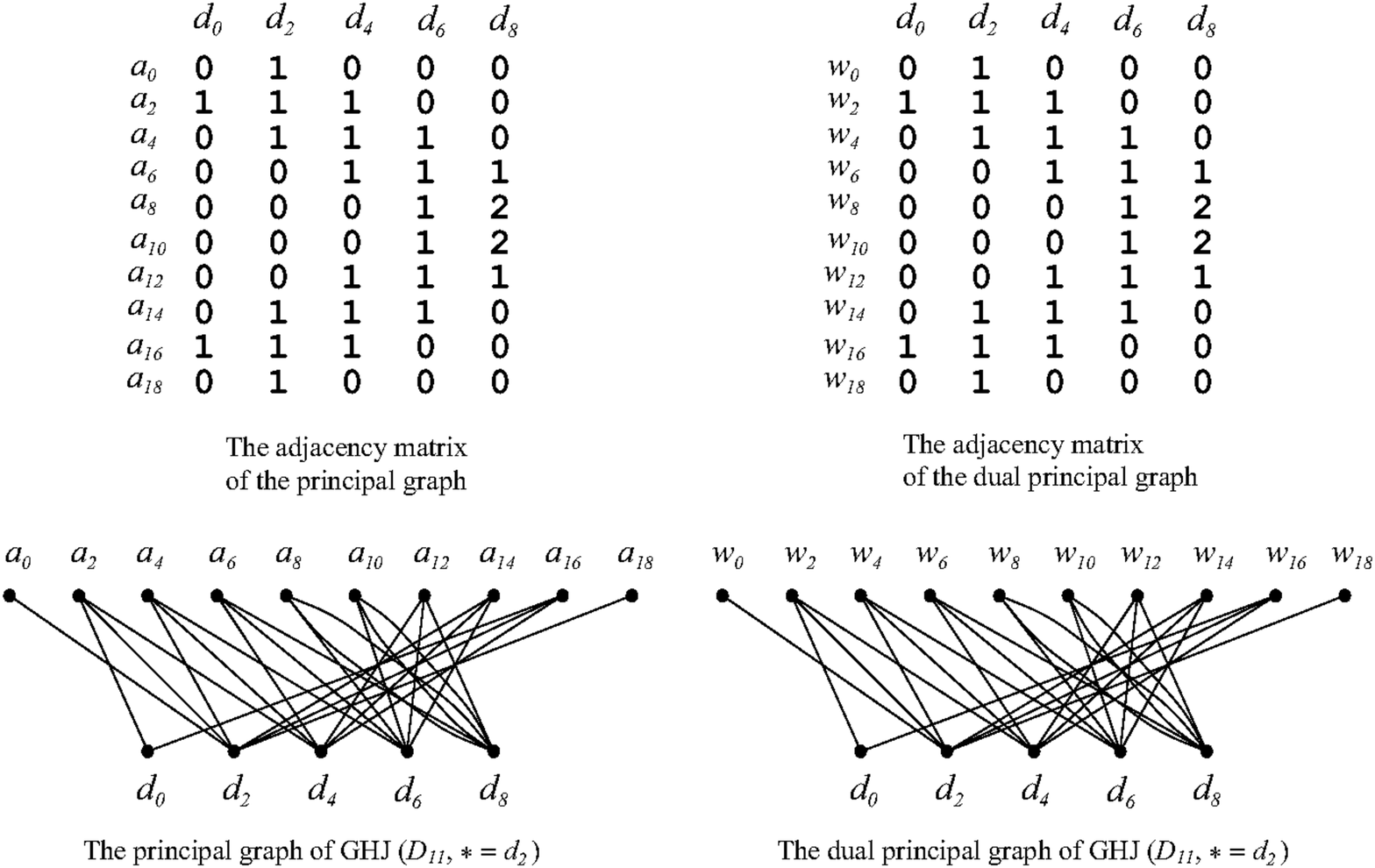}
\caption{The (dual) principal graph of GHJ$(D_{11}, *=d_2)$.}
\label{GHJ(D11-d2)}
\end{figure}

%%%%% GHJ(D11-d3) %%%%%%%%%%%%%%%%%%%%%%%%%%%%%%%%%
\begin{figure}[H]
\centering
\includegraphics[width=130mm,clip]{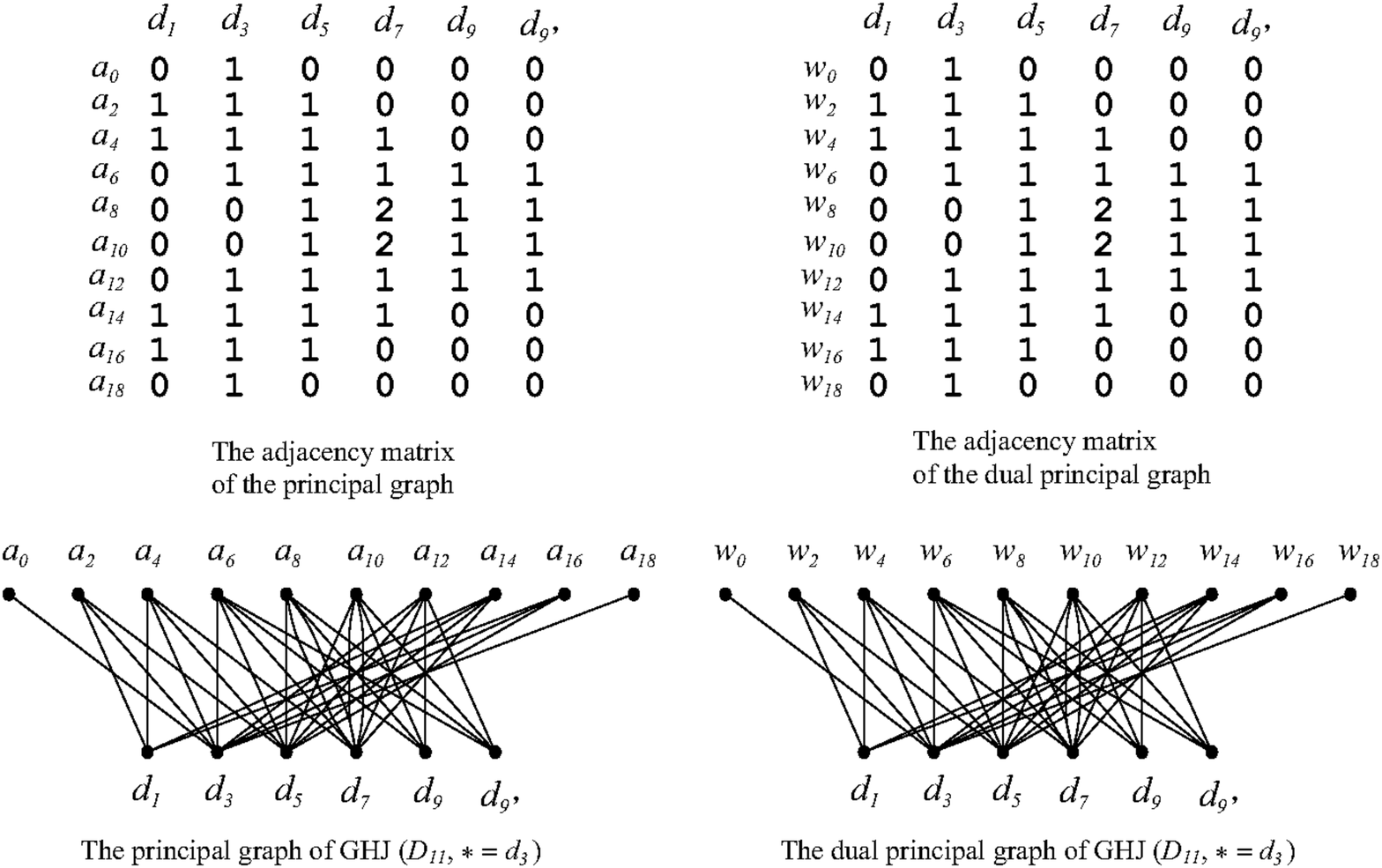}
\caption{The (dual) principal graph of GHJ$(D_{11}, *=d_3)$.}
\label{GHJ(D11-d3)}
\end{figure}

%%%%% GHJ(D11-d4) %%%%%%%%%%%%%%%%%%%%%%%%%%%%%%%%%
\begin{figure}[H]
\centering
\includegraphics[width=130mm,clip]{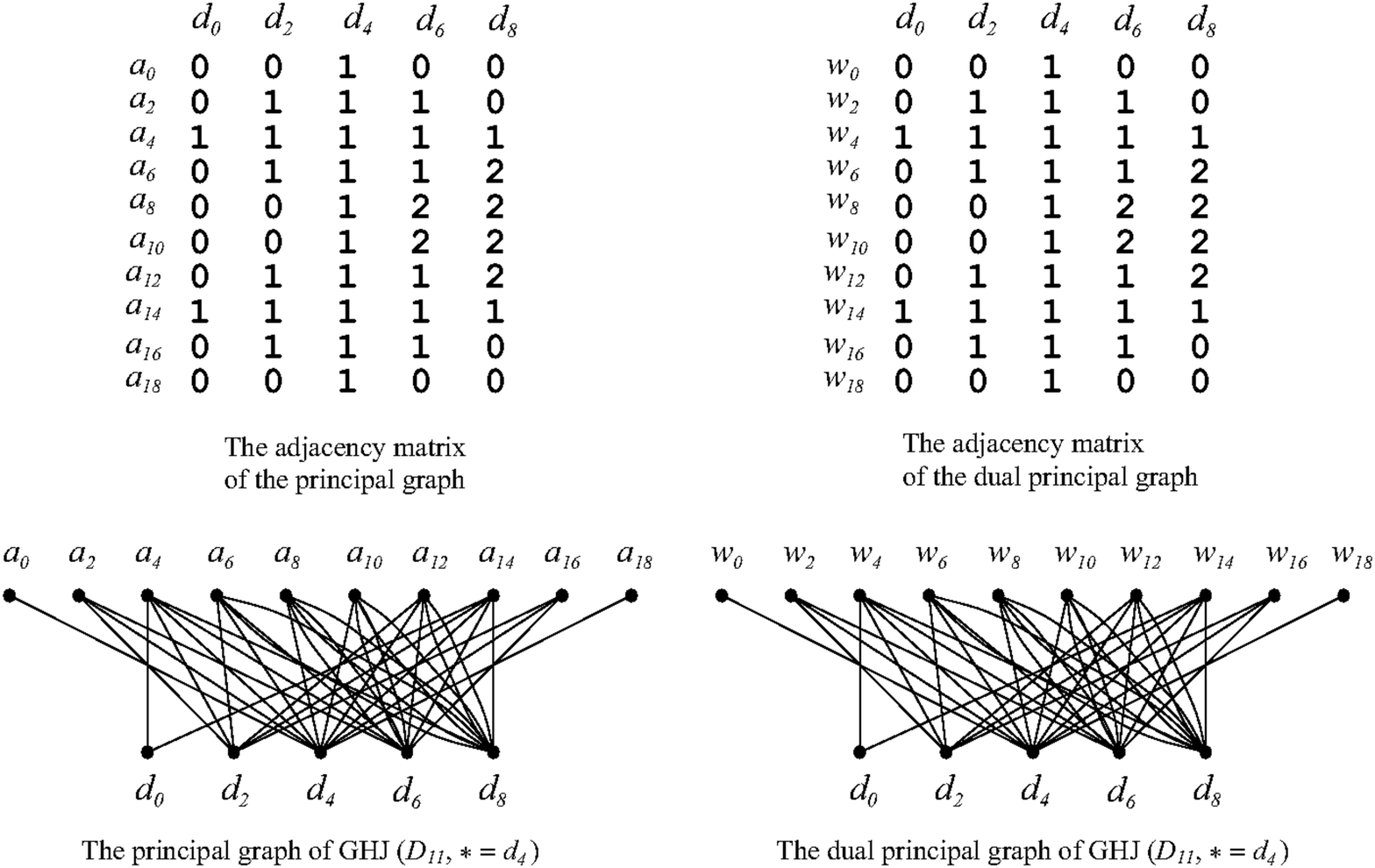}
\caption{The (dual) principal graph of GHJ$(D_{11}, *=d_4)$.}
\label{GHJ(D11-d4)}
\end{figure}

%%%%% GHJ(D11-d5) %%%%%%%%%%%%%%%%%%%%%%%%%%%%%%%%%
\begin{figure}[H]
\centering
\includegraphics[width=130mm,clip]{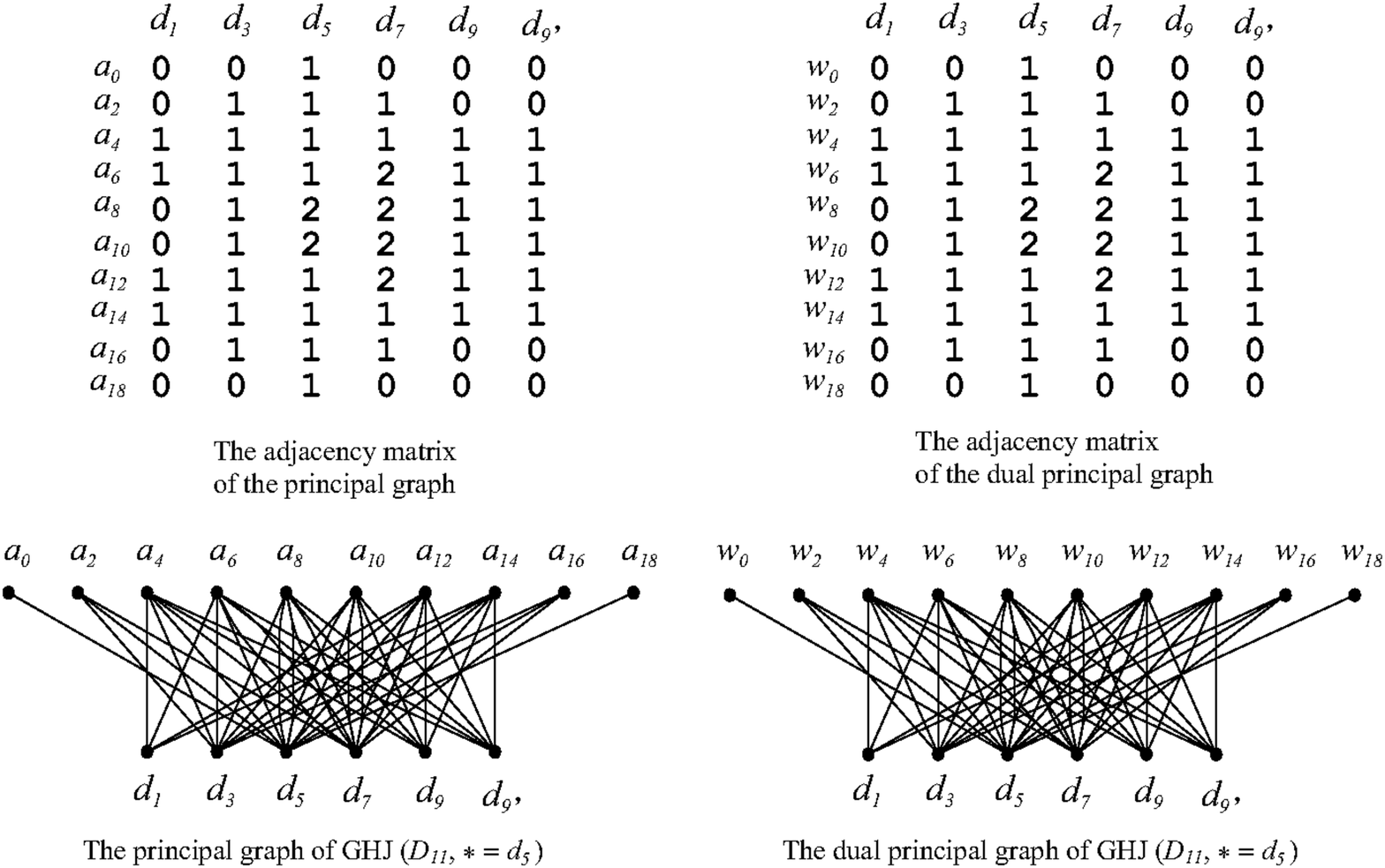}
\caption{The (dual) principal graph of GHJ$(D_{11}, *=d_5)$.}
\label{GHJ(D11-d5)}
\end{figure}

%%%%% GHJ(D11-d6) %%%%%%%%%%%%%%%%%%%%%%%%%%%%%%%%%
\begin{figure}[H]
\centering
\includegraphics[width=130mm,clip]{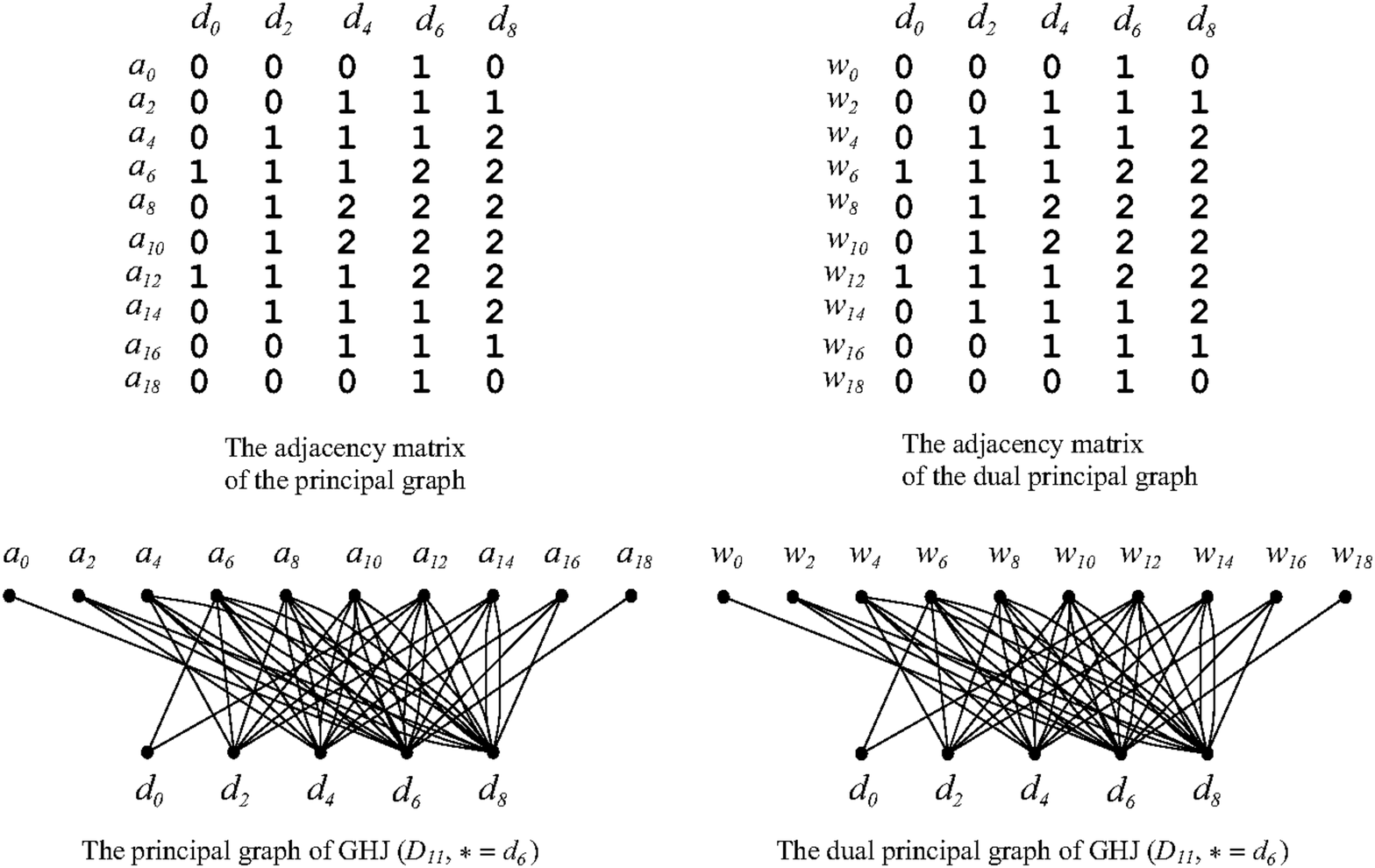}
\caption{The (dual) principal graph of GHJ$(D_{11}, *=d_6)$.}
\label{GHJ(D11-d6)}
\end{figure}

%%%%% GHJ(D9-d7) %%%%%%%%%%%%%%%%%%%%%%%%%%%%%%%%%
\begin{figure}[H]
\centering
\includegraphics[width=130mm,clip]{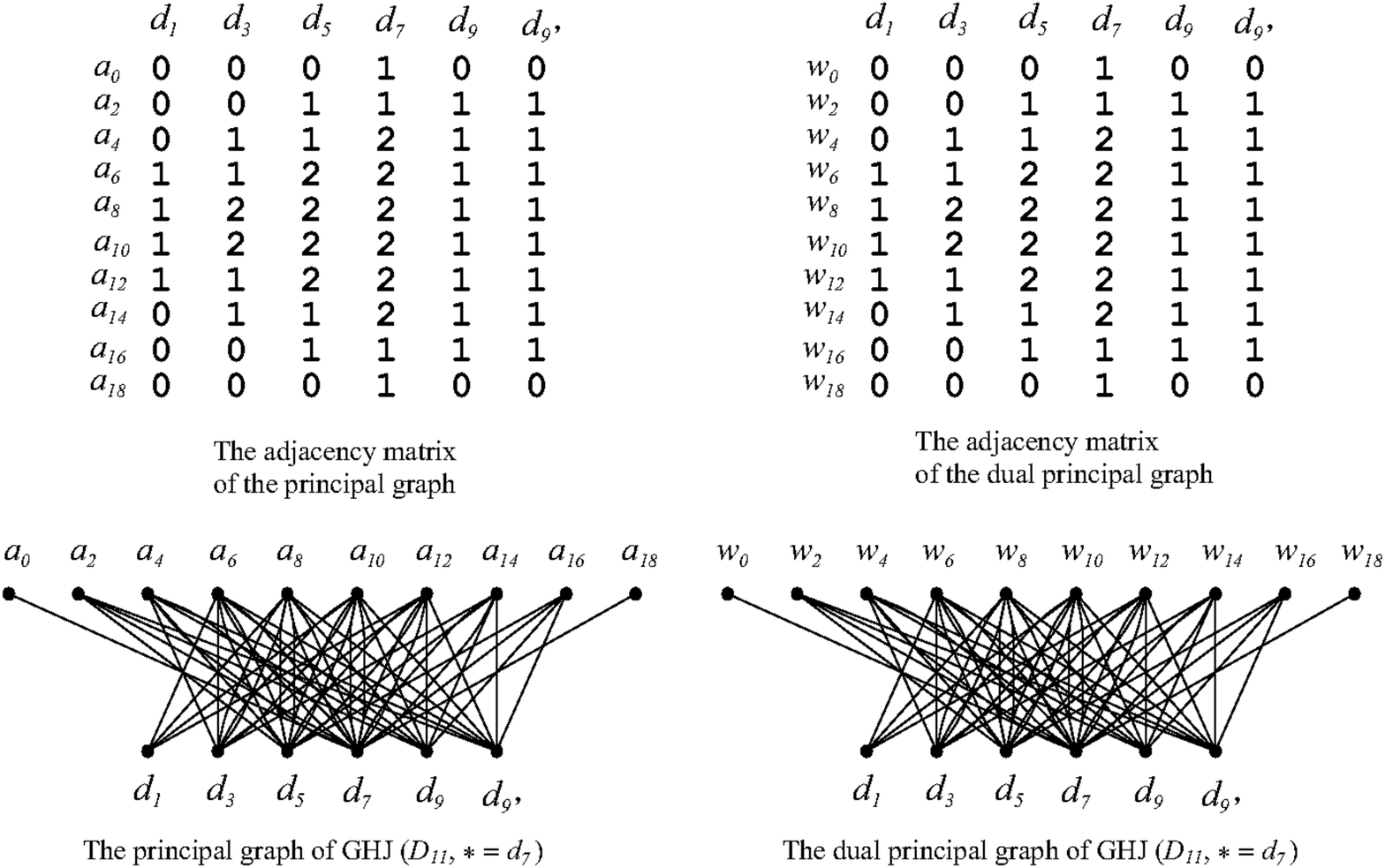}
\caption{The (dual) principal graph of GHJ$(D_{11}, *=d_7)$.}
\label{GHJ(D11-d7)}
\end{figure}

%%%%% GHJ(D11-d8) %%%%%%%%%%%%%%%%%%%%%%%%%%%%%%%%%

%%%%%%%%%%%%%%%%%%%%%%%%%%%%%%%%%%%%%%%%%%%%%%%%%%
%%% The (dual) pincipal graphs of             %%%%
%%%     Goodman-de la Harpe-Jones subfactors  %%%%
%%%%%%%%%%%%%%%%%%%%%%%%%%%%%%%%%%%%%%%%%%%%%%%%%%
%%%       D6                                  %%%%
%%%%%%%%%%%%%%%%%%%%%%%%%%%%%%%%%%%%%%%%%%%%%%%%%%

%%%%% GHJ(D6-d1) %%%%%%%%%%%%%%%%%%%%%%%%%%%%%%%%%

%%%%% GHJ(D6-d2) %%%%%%%%%%%%%%%%%%%%%%%%%%%%%%%%%
\begin{figure}[H]
\centering
\includegraphics[width=130mm,clip]{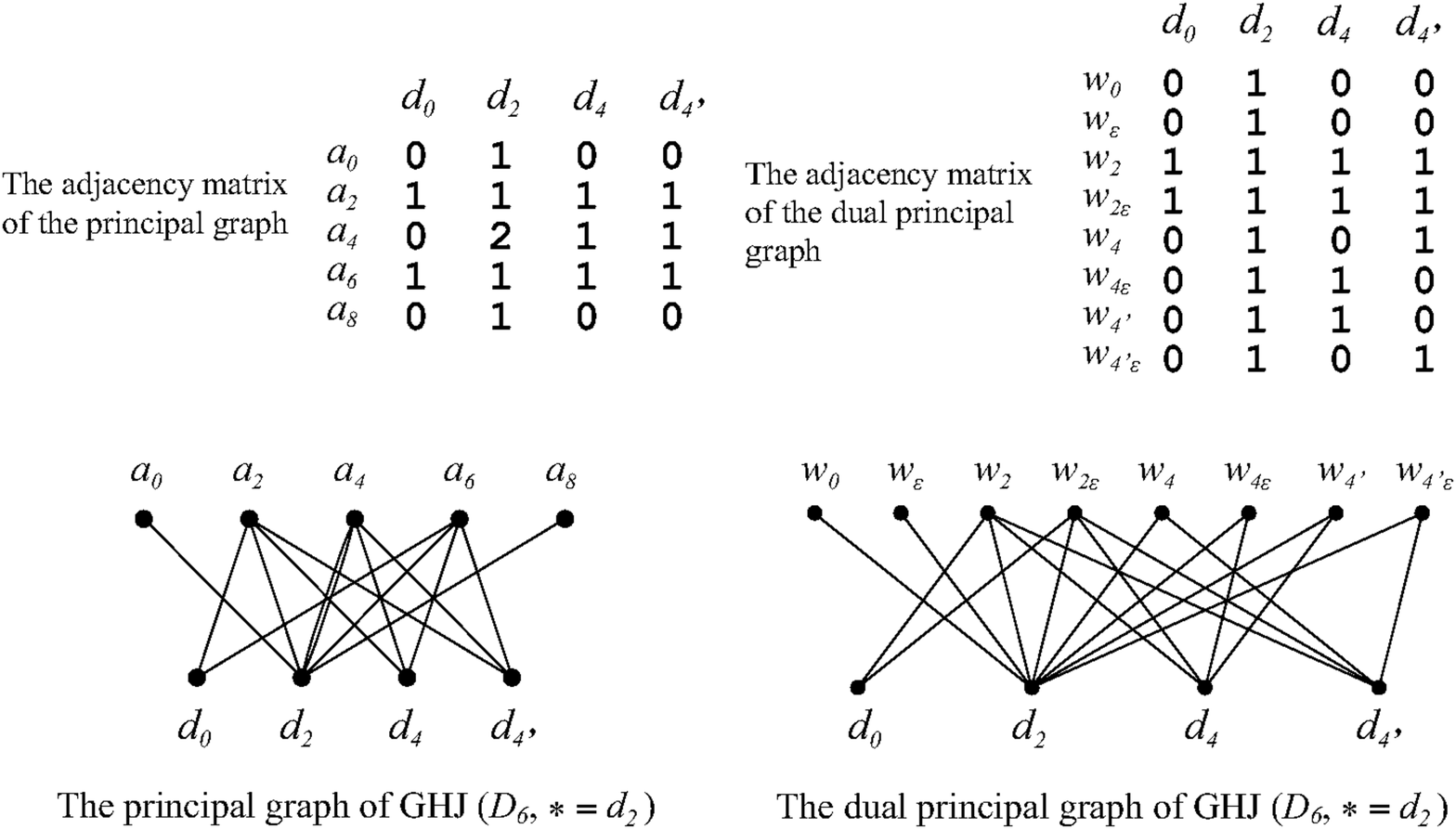}
\caption{The (dual) principal graph of GHJ$(D_6, *=d_2)$.}
\label{GHJ(D6-d2)}
\end{figure}

%%%%% GHJ(D6-d3) %%%%%%%%%%%%%%%%%%%%%%%%%%%%%%%%%
\begin{figure}[H]
\centering
\includegraphics[width=130mm,clip]{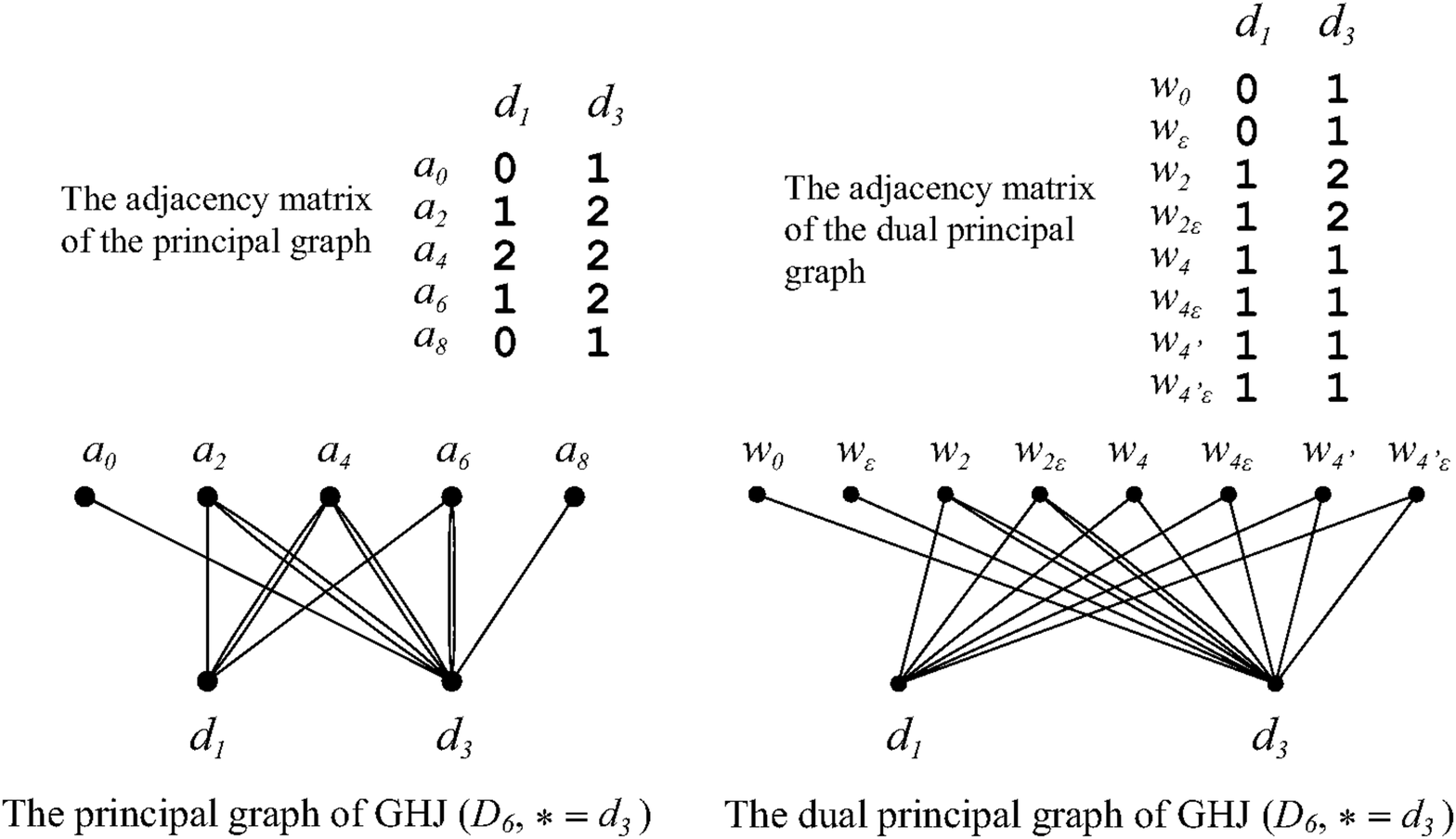}
\caption{The (dual) principal graph of GHJ$(D_6, *=d_3)$.}
\label{GHJ(D6-d3)}
\end{figure}

%%%%% GHJ(D6-d4) %%%%%%%%%%%%%%%%%%%%%%%%%%%%%%%%%
\begin{figure}[H]
\centering
\includegraphics[width=130mm,clip]{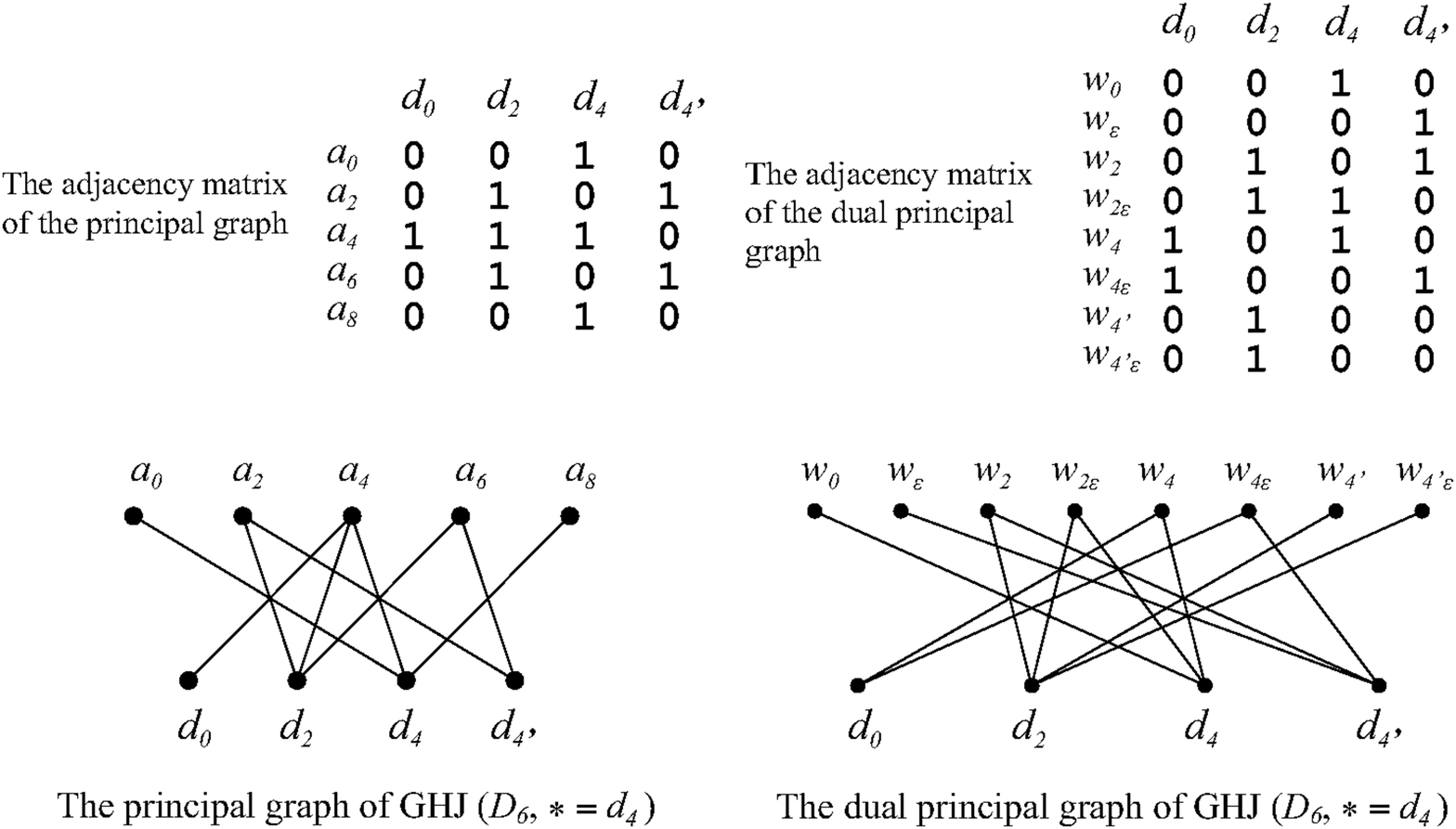}
\caption{The (dual) principal graph of GHJ$(D_6, *=d_4)$.}
\label{GHJ(D6-d4)}
\end{figure}

%%%%%%%%%%%%%%%%%%%%%%%%%%%%%%%%%%%%%%%%%%%%%%%%%%
%%% The (dual) pincipal graphs of             %%%%
%%%     Goodman-de la Harpe-Jones subfactors  %%%%
%%%%%%%%%%%%%%%%%%%%%%%%%%%%%%%%%%%%%%%%%%%%%%%%%%
%%%       D8                                  %%%%
%%%%%%%%%%%%%%%%%%%%%%%%%%%%%%%%%%%%%%%%%%%%%%%%%%

%%%%% GHJ(D8-d1) %%%%%%%%%%%%%%%%%%%%%%%%%%%%%%%%%
\begin{figure}[H]
\centering
\includegraphics[width=120mm,clip]{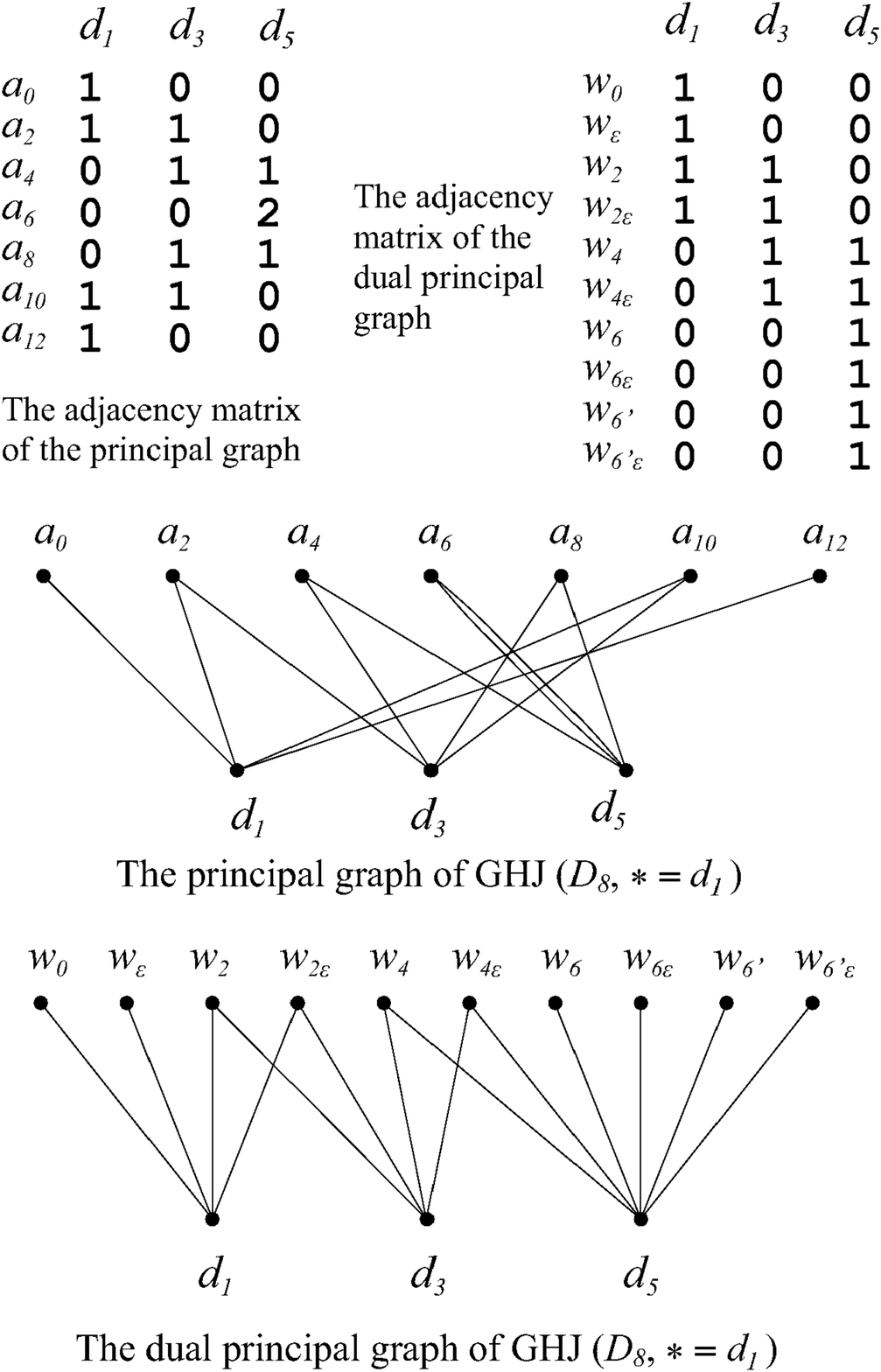}
\caption{The (dual) principal graph of GHJ$(D_8, *=d_1)$.}
\label{GHJ(D8-d1)}
\end{figure}

%%%%% GHJ(D8-d2) %%%%%%%%%%%%%%%%%%%%%%%%%%%%%%%%%
\begin{figure}[H]
\centering
\includegraphics[width=130mm,clip]{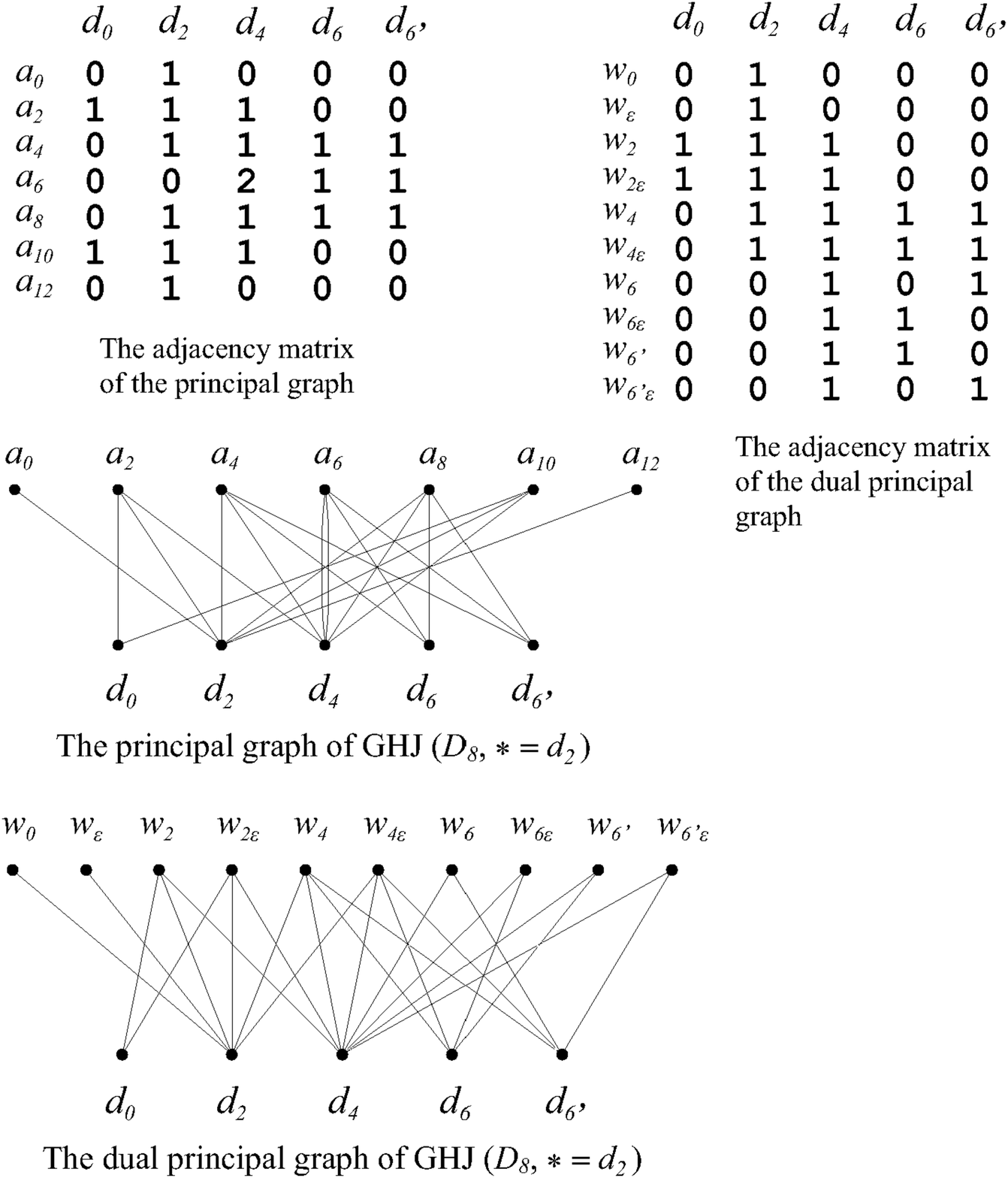}
\caption{The (dual) principal graph of GHJ$(D_8, *=d_2)$.}
\label{GHJ(D8-d2)}
\end{figure}

%%%%% GHJ(D8-d3) %%%%%%%%%%%%%%%%%%%%%%%%%%%%%%%%%
\begin{figure}[H]
\centering
\includegraphics[width=120mm,clip]{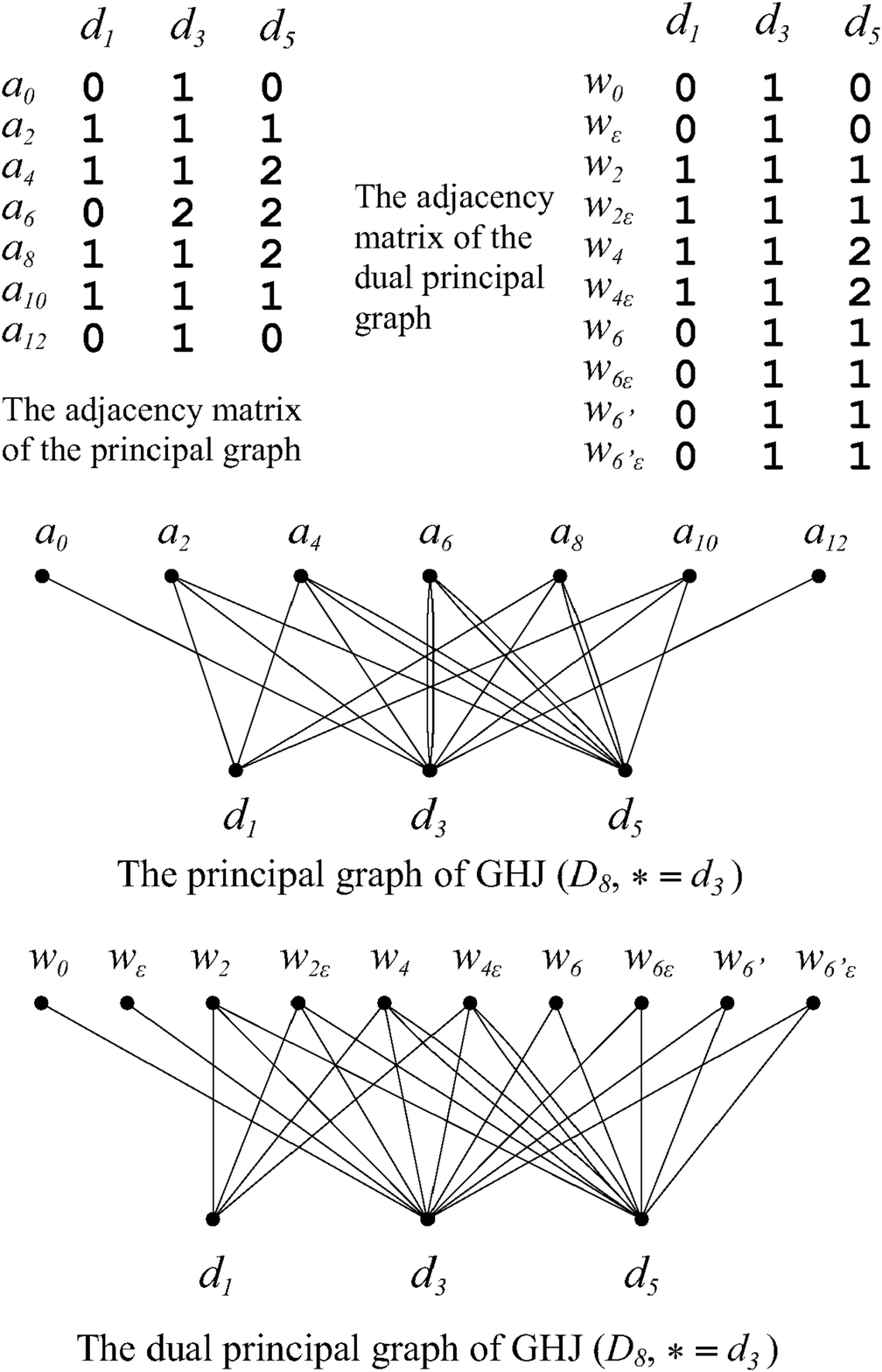}
\caption{The (dual) principal graph of GHJ$(D_8, *=d_3)$.}
\label{GHJ(D8-d3)}
\end{figure}

%%%%% GHJ(D8-d4) %%%%%%%%%%%%%%%%%%%%%%%%%%%%%%%%%
\begin{figure}[H]
\centering
\includegraphics[width=130mm,clip]{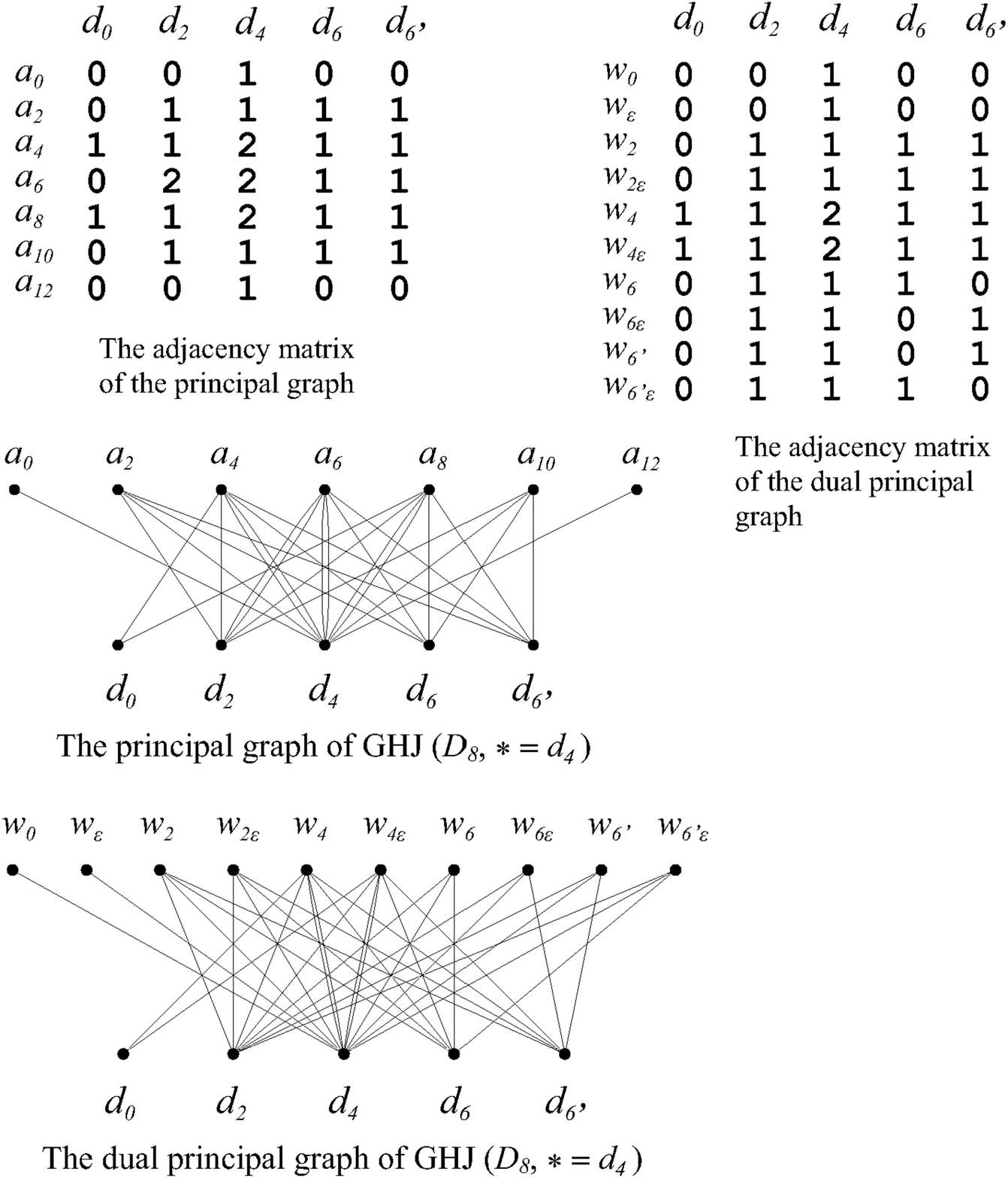}
\caption{The (dual) principal graph of GHJ$(D_8, *=d_4)$.}
\label{GHJ(D8-d4)}
\end{figure}

%%%%% GHJ(D8-d5) %%%%%%%%%%%%%%%%%%%%%%%%%%%%%%%%%
\begin{figure}[H]
\centering
\includegraphics[width=120mm,clip]{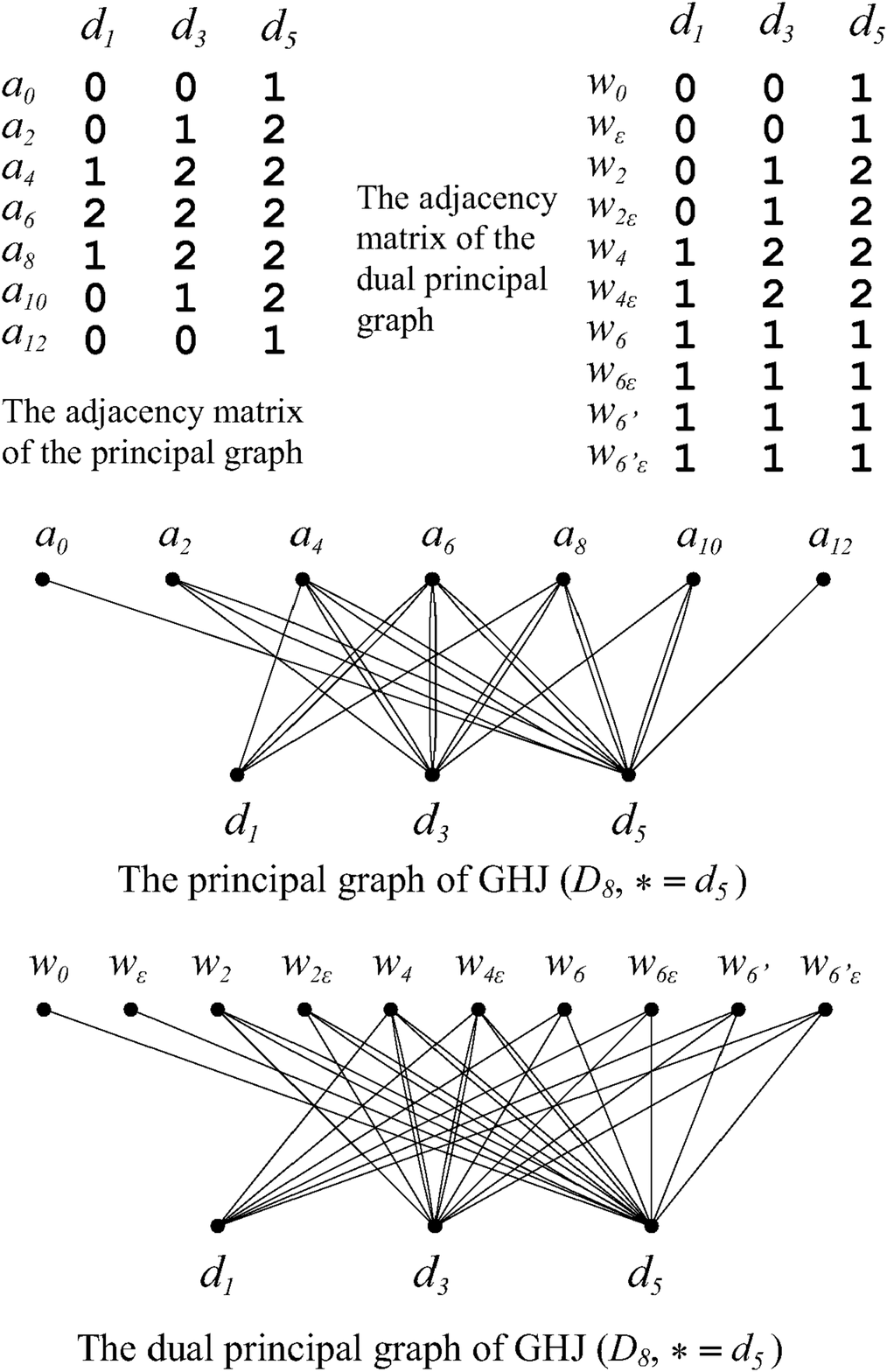}
\caption{The (dual) principal graph of GHJ$(D_8, *=d_5)$.}
\label{GHJ(D8-d5)}
\end{figure}

%%%%% GHJ(D8-d6) %%%%%%%%%%%%%%%%%%%%%%%%%%%%%%%%%
\begin{figure}[H]
\centering
\includegraphics[width=130mm,clip]{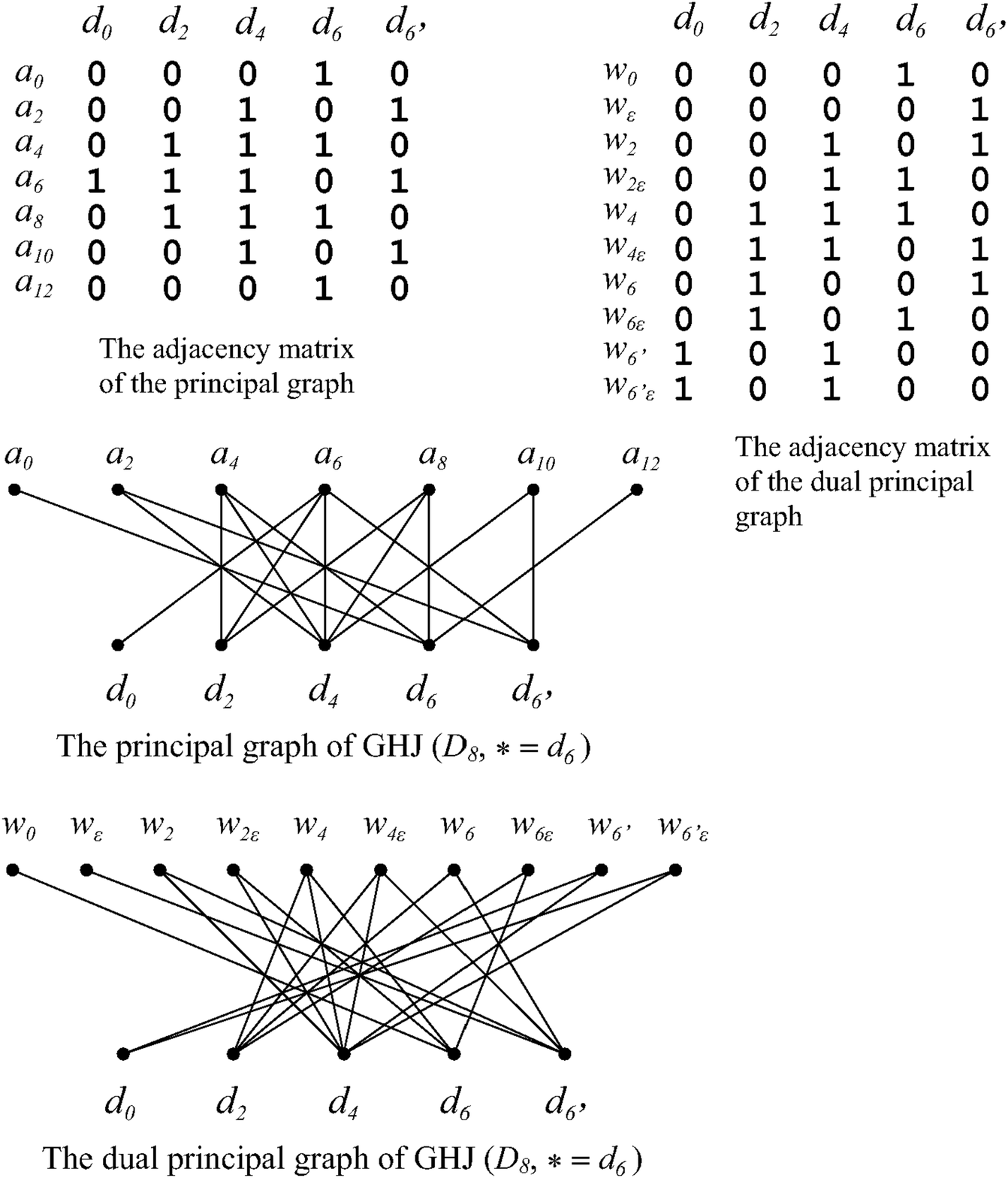}
\caption{The (dual) principal graph of GHJ$(D_8, *=d_6)$.}
\label{GHJ(D8-d6)}
\end{figure}

%%%%%%%%%%%%%%%%%%%%%%%%%%%%%%%%%%%%%%%%%%%%%%%%%%
%%% The (dual) pincipal graphs of             %%%%
%%%     Goodman-de la Harpe-Jones subfactors  %%%%
%%%%%%%%%%%%%%%%%%%%%%%%%%%%%%%%%%%%%%%%%%%%%%%%%%
%%%       D10                                  %%%%
%%%%%%%%%%%%%%%%%%%%%%%%%%%%%%%%%%%%%%%%%%%%%%%%%%

%%%%% GHJ(D10-d1) %%%%%%%%%%%%%%%%%%%%%%%%%%%%%%%%%
\begin{figure}[H]
\centering
\includegraphics[width=130mm,clip]{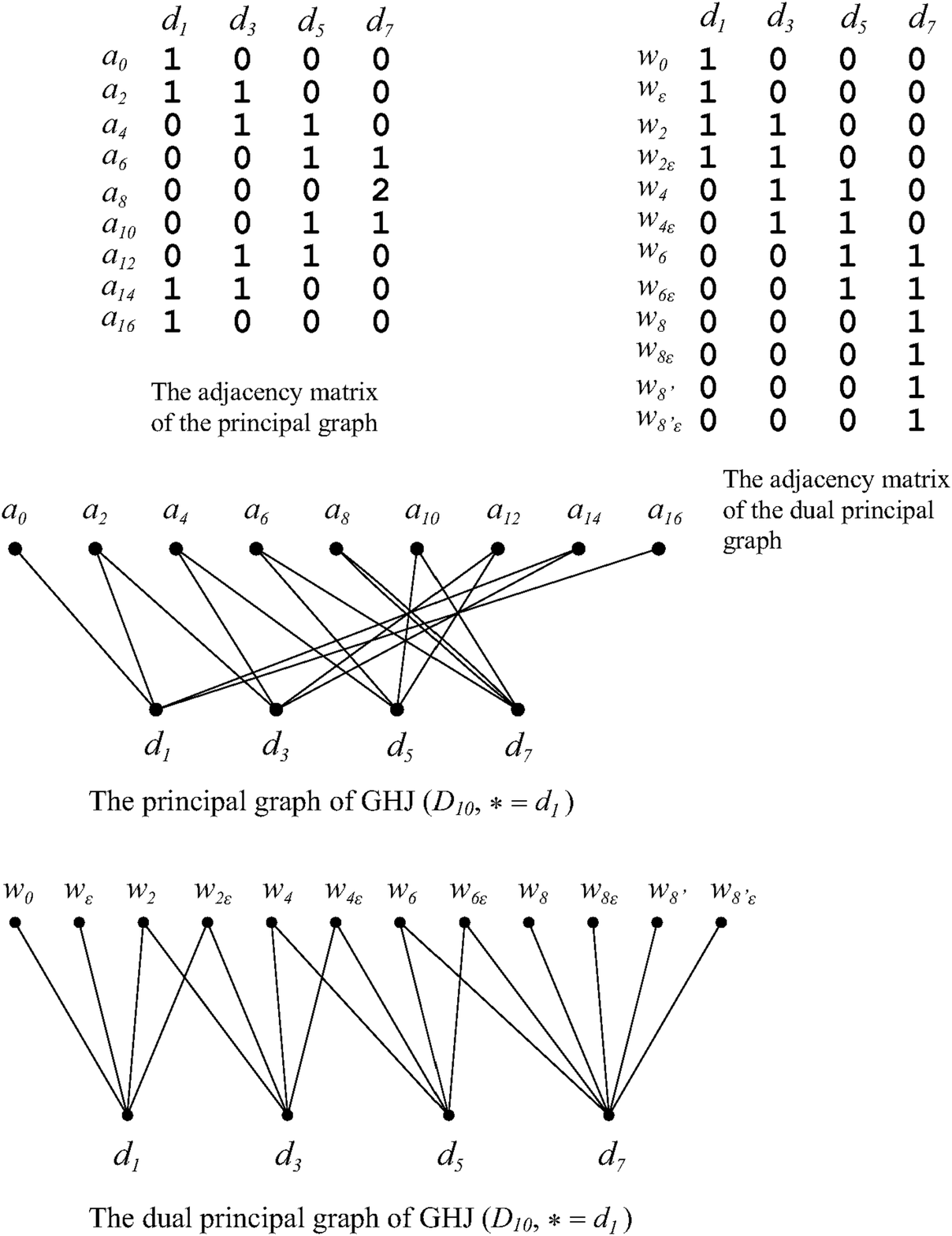}
\caption{The (dual) principal graph of GHJ$(D_{10}, *=d_1)$.}
\label{GHJ(D10-d1)}
\end{figure}

%%%%% GHJ(D10-d2) %%%%%%%%%%%%%%%%%%%%%%%%%%%%%%%%%
\begin{figure}[H]
\centering
\includegraphics[width=130mm,clip]{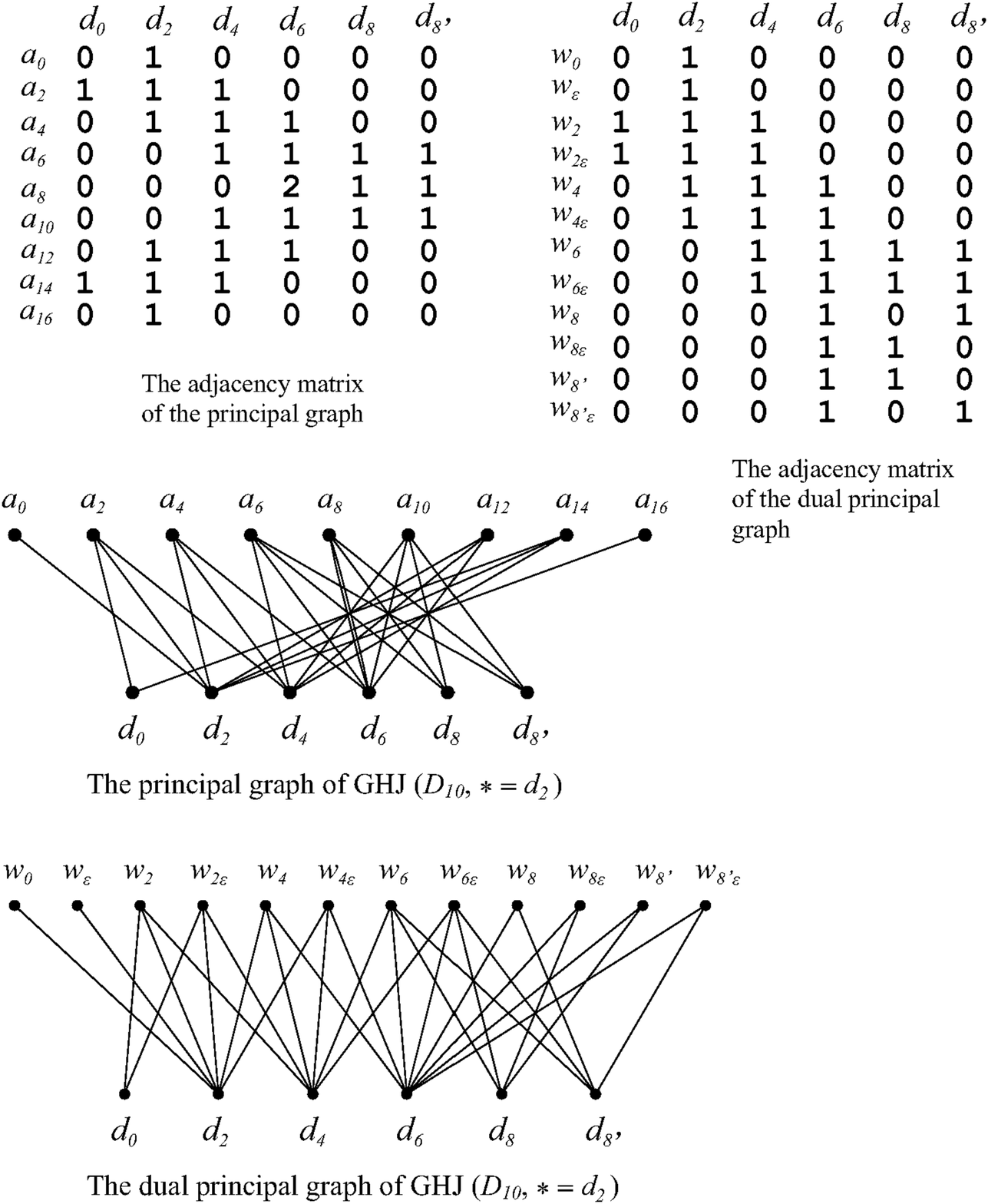}
\caption{The (dual) principal graph of GHJ$(D_{10}, *=d_2)$.}
\label{GHJ(D10-d2)}
\end{figure}

%%%%% GHJ(D10-d3) %%%%%%%%%%%%%%%%%%%%%%%%%%%%%%%%%
\begin{figure}[H]
\centering
\includegraphics[width=130mm,clip]{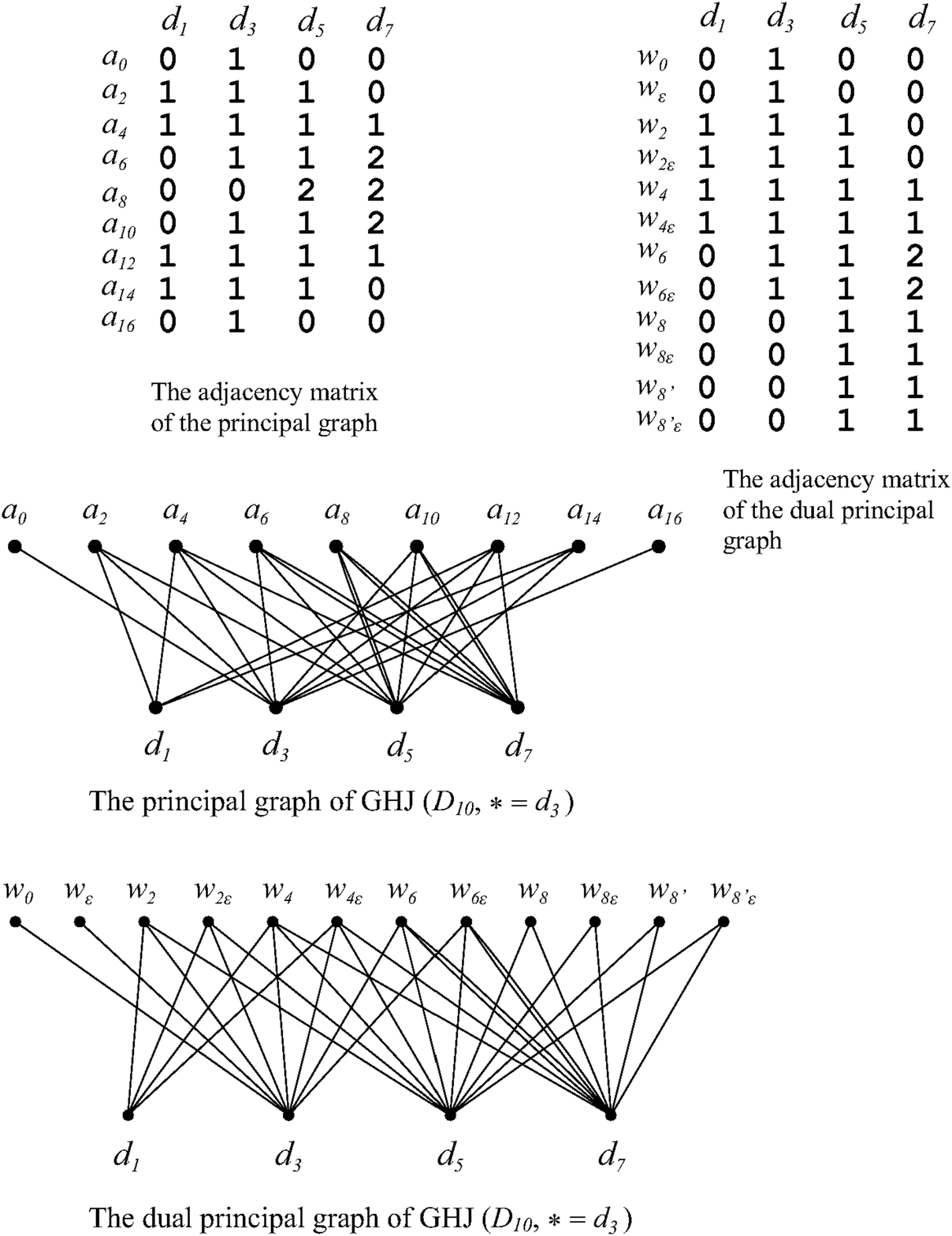}
\caption{The (dual) principal graph of GHJ$(D_{10}, *=d_3)$.}
\label{GHJ(D10-d3)}
\end{figure}

%%%%% GHJ(D10-d4) %%%%%%%%%%%%%%%%%%%%%%%%%%%%%%%%%
\begin{figure}[H]
\centering
\includegraphics[width=130mm,clip]{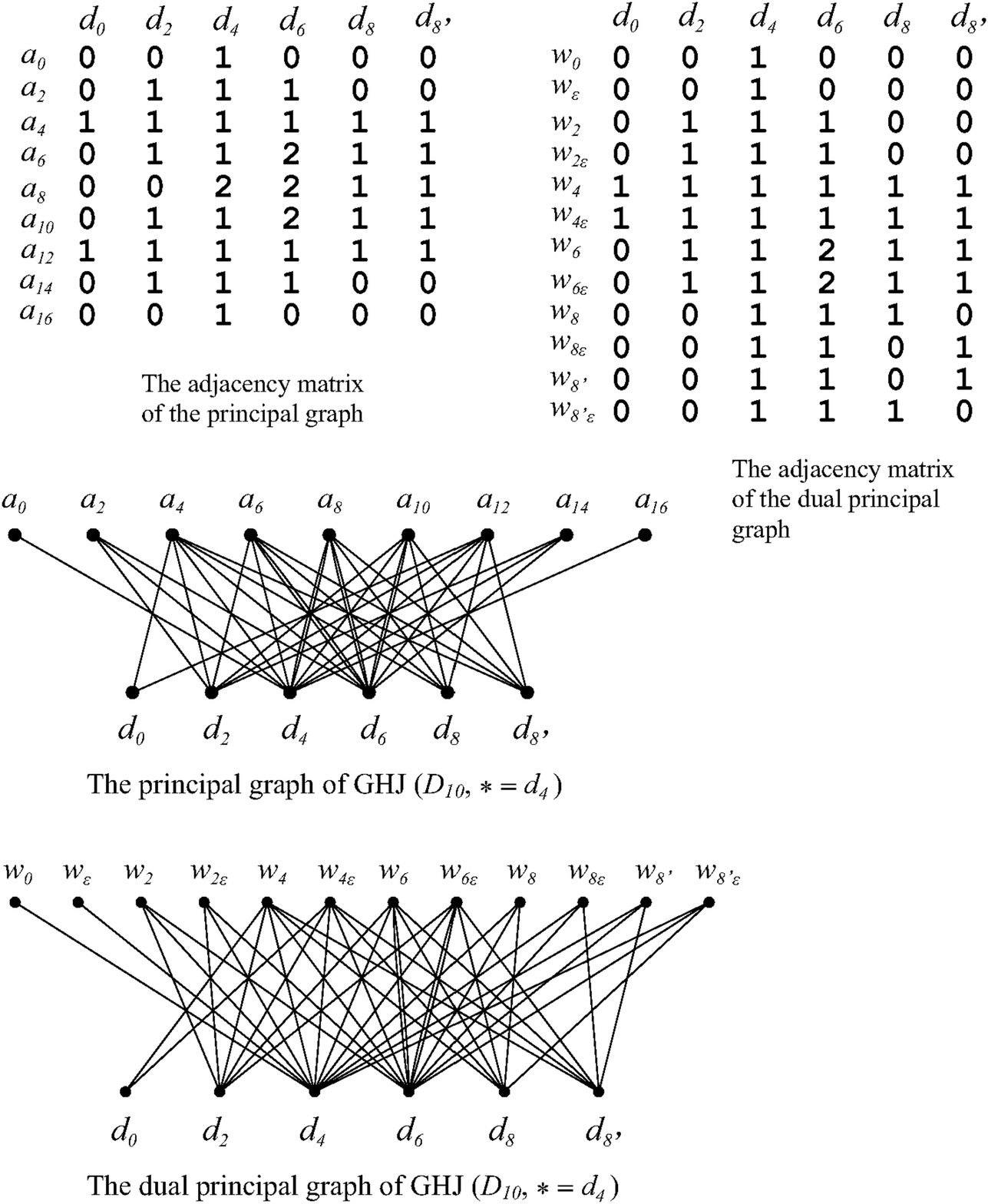}
\caption{The (dual) principal graph of GHJ$(D_{10}, *=d_4)$.}
\label{GHJ(D10-d4)}
\end{figure}

%%%%% GHJ(D10-d5) %%%%%%%%%%%%%%%%%%%%%%%%%%%%%%%%%
\begin{figure}[H]
\centering
\includegraphics[width=130mm,clip]{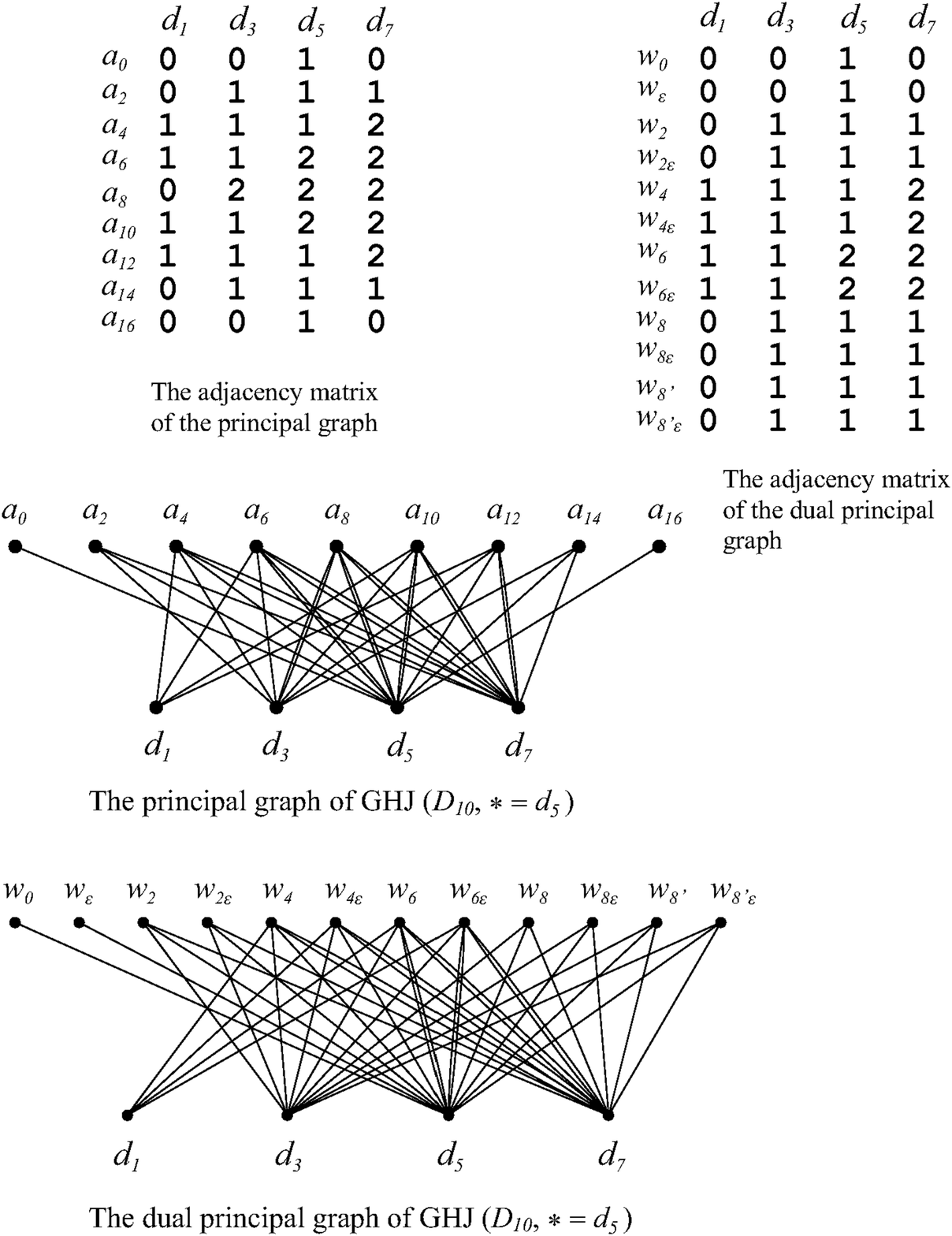}
\caption{The (dual) principal graph of GHJ$(D_{10}, *=d_5)$.}
\label{GHJ(D10-d5)}
\end{figure}

%%%%% GHJ(D10-d6) %%%%%%%%%%%%%%%%%%%%%%%%%%%%%%%%%
\begin{figure}[H]
\centering
\includegraphics[width=130mm,clip]{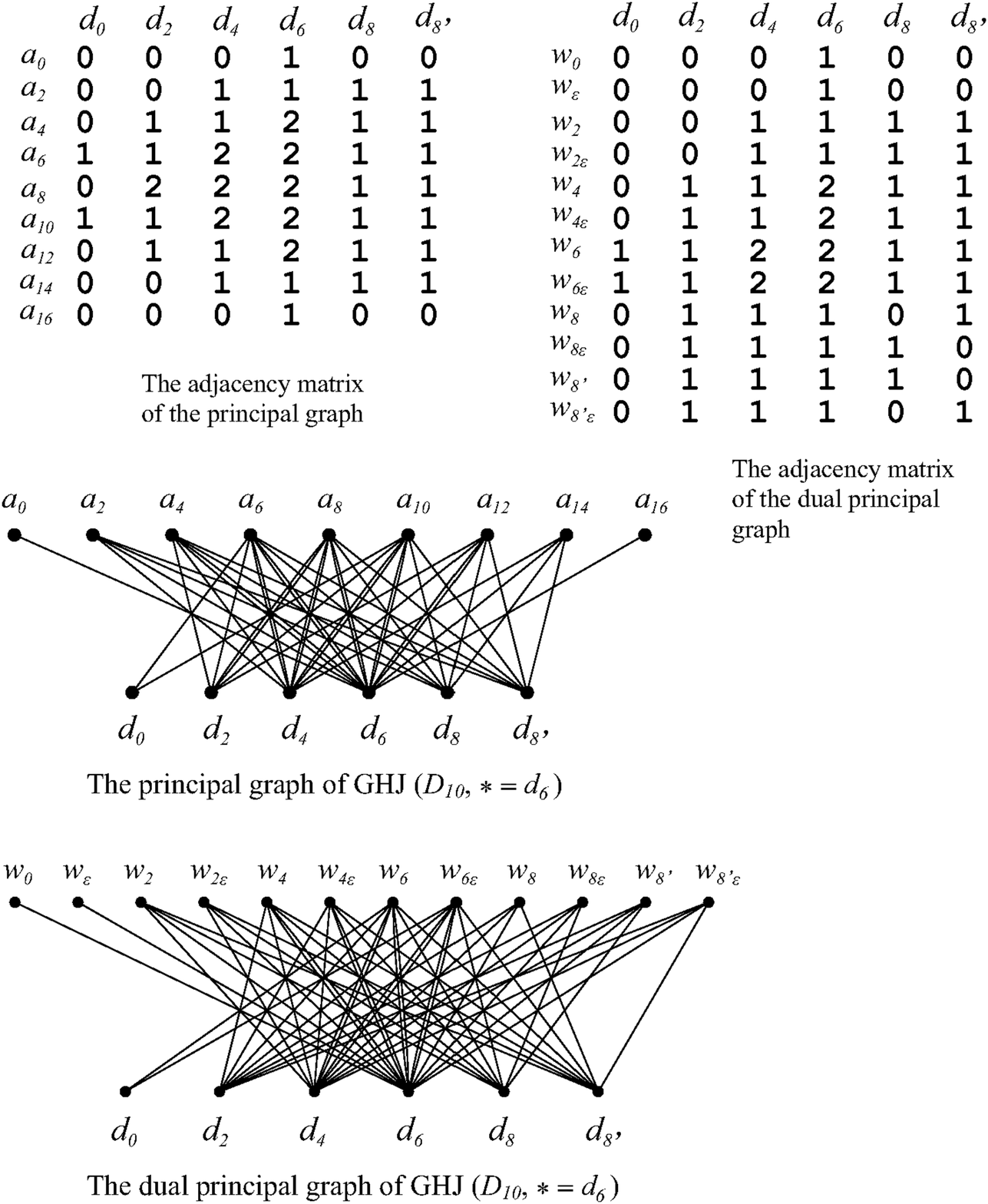}
\caption{The (dual) principal graph of GHJ$(D_{10}, *=d_6)$.}
\label{GHJ(D10-d6)}
\end{figure}

%%%%% GHJ(D10-d7) %%%%%%%%%%%%%%%%%%%%%%%%%%%%%%%%%
\begin{figure}[H]
\centering
\includegraphics[width=130mm,clip]{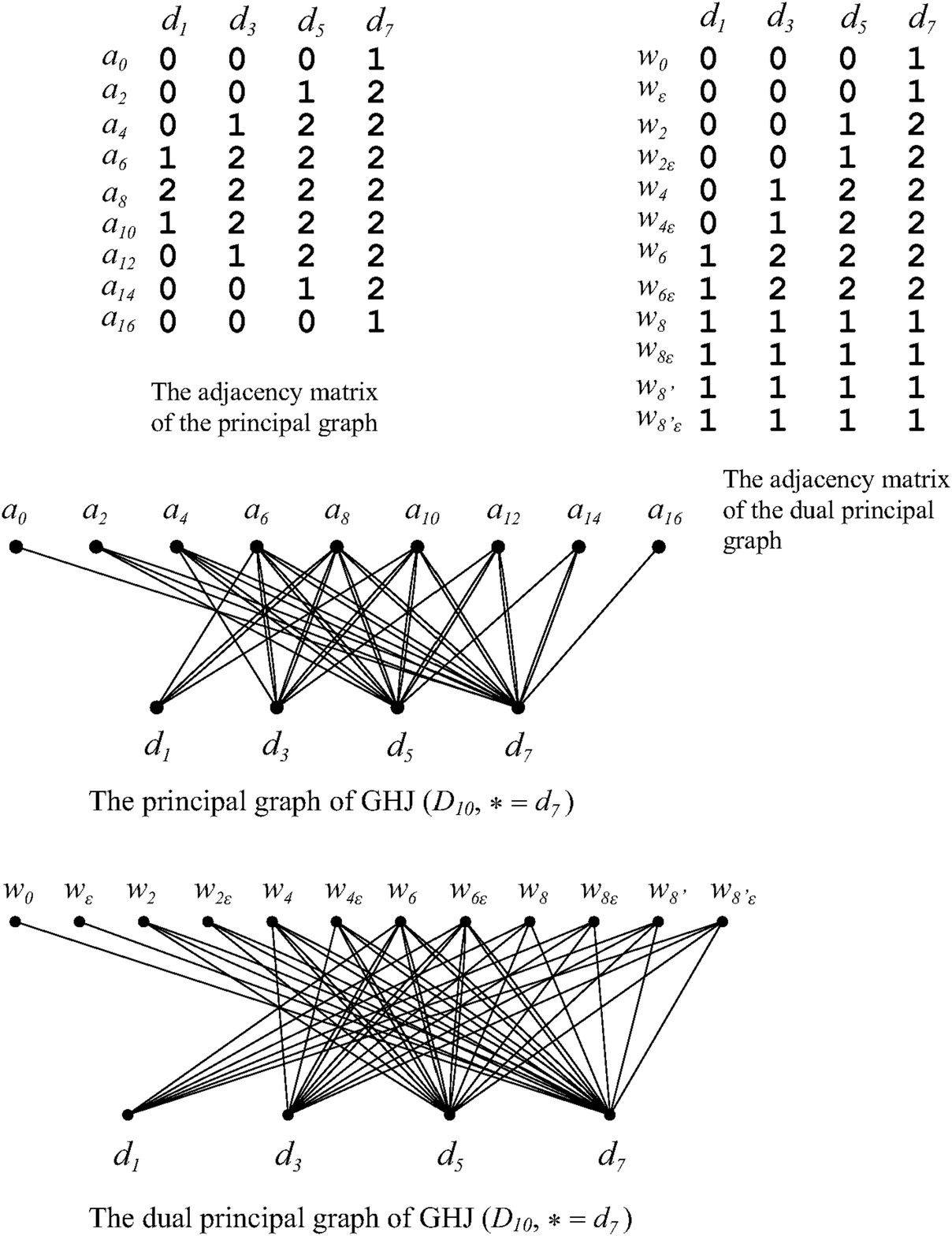}
\caption{The (dual) principal graph of GHJ$(D_{10}, *=d_7)$.}
\label{GHJ(D10-d7)}
\end{figure}

%%%%% GHJ(D10-d8) %%%%%%%%%%%%%%%%%%%%%%%%%%%%%%%%%
\begin{figure}[H]
\centering
\includegraphics[width=130mm,clip]{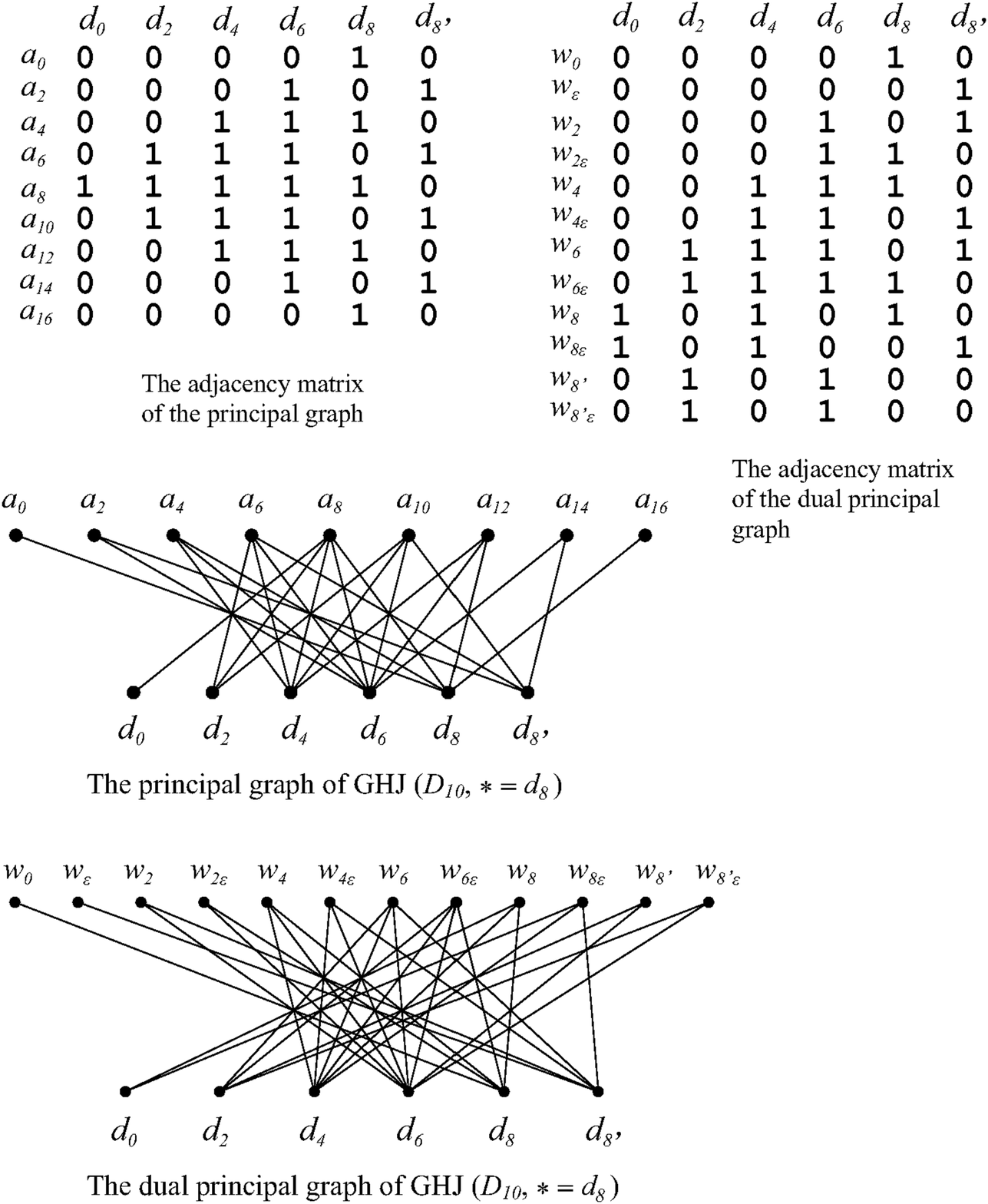}
\caption{The (dual) principal graph of GHJ$(D_{10}, *=d_8)$.}
\label{GHJ(D10-d8)}
\end{figure}

%%%%%%%%%%%%%%%%%%%%%%%%%%%%%%%%%%%%%%%%%%%%%%%%%%
%%% The (dual) pincipal graphs of             %%%%
%%%     Goodman-de la Harpe-Jones subfactors  %%%%
%%%%%%%%%%%%%%%%%%%%%%%%%%%%%%%%%%%%%%%%%%%%%%%%%%
%%%       D12                                  %%%%
%%%%%%%%%%%%%%%%%%%%%%%%%%%%%%%%%%%%%%%%%%%%%%%%%%

%%%%% GHJ(D12-d1) %%%%%%%%%%%%%%%%%%%%%%%%%%%%%%%%%
\begin{figure}[H]
\centering
\includegraphics[width=130mm,clip]{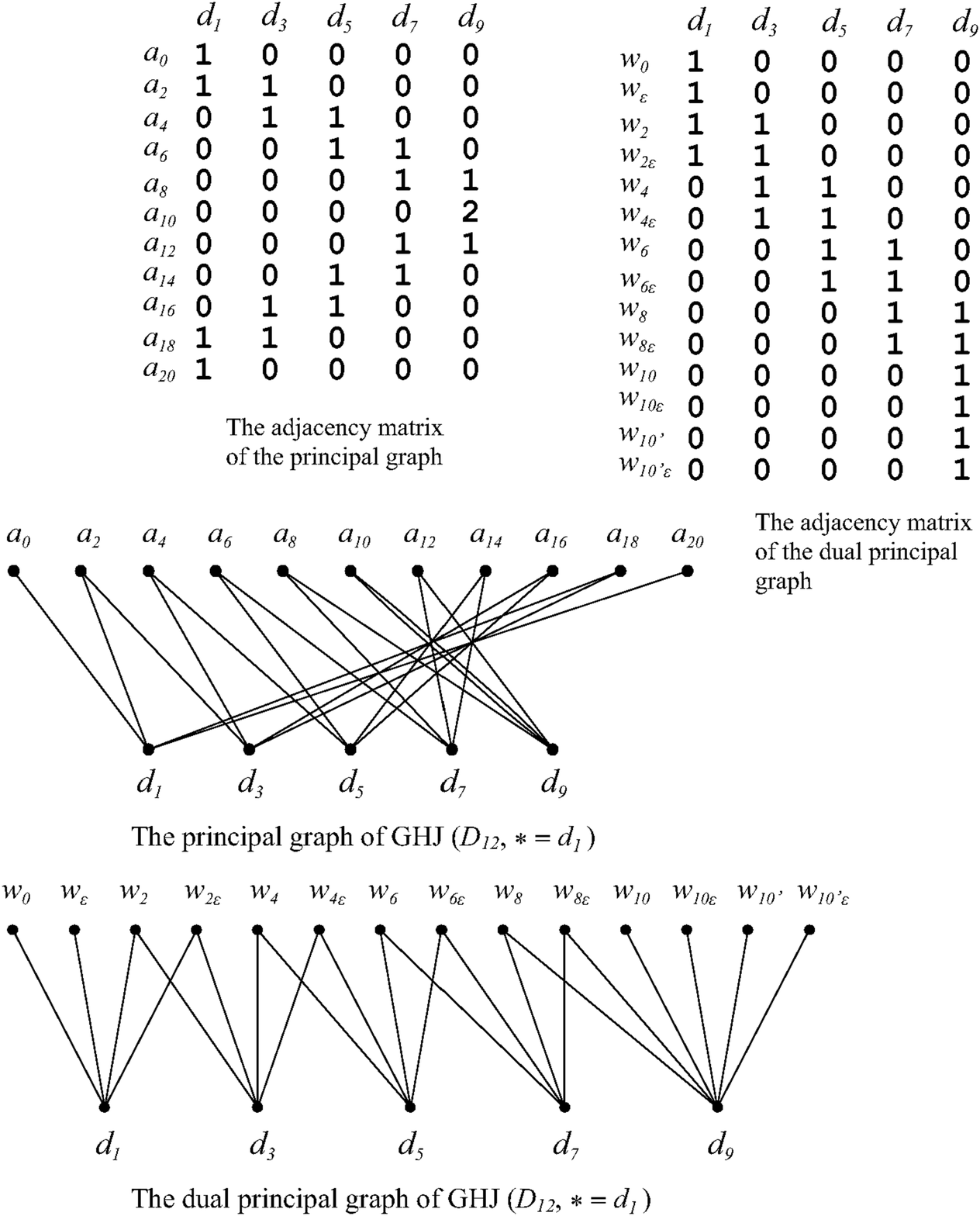}
\caption{The (dual) principal graph of GHJ$(D_{12}, *=d_1)$.}
\label{GHJ(D12-d1)}
\end{figure}

%%%%% GHJ(D12-d2) %%%%%%%%%%%%%%%%%%%%%%%%%%%%%%%%%
\begin{figure}[H]
\centering
\includegraphics[width=130mm,clip]{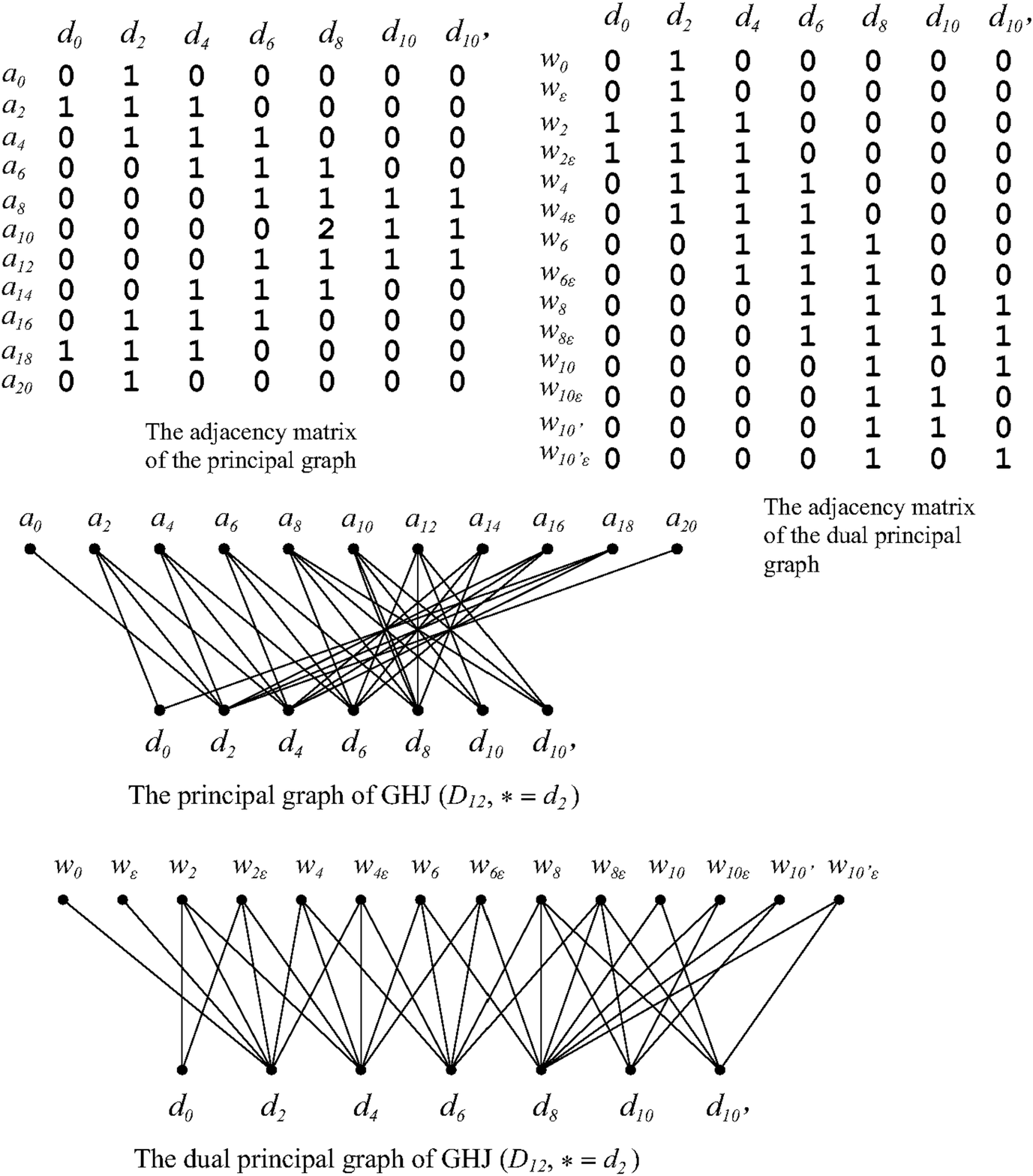}
\caption{The (dual) principal graph of GHJ$(D_{12}, *=d_2)$.}
\label{GHJ(D12-d2)}
\end{figure}

%%%%% GHJ(D12-d3) %%%%%%%%%%%%%%%%%%%%%%%%%%%%%%%%%
\begin{figure}[H]
\centering
\includegraphics[width=130mm,clip]{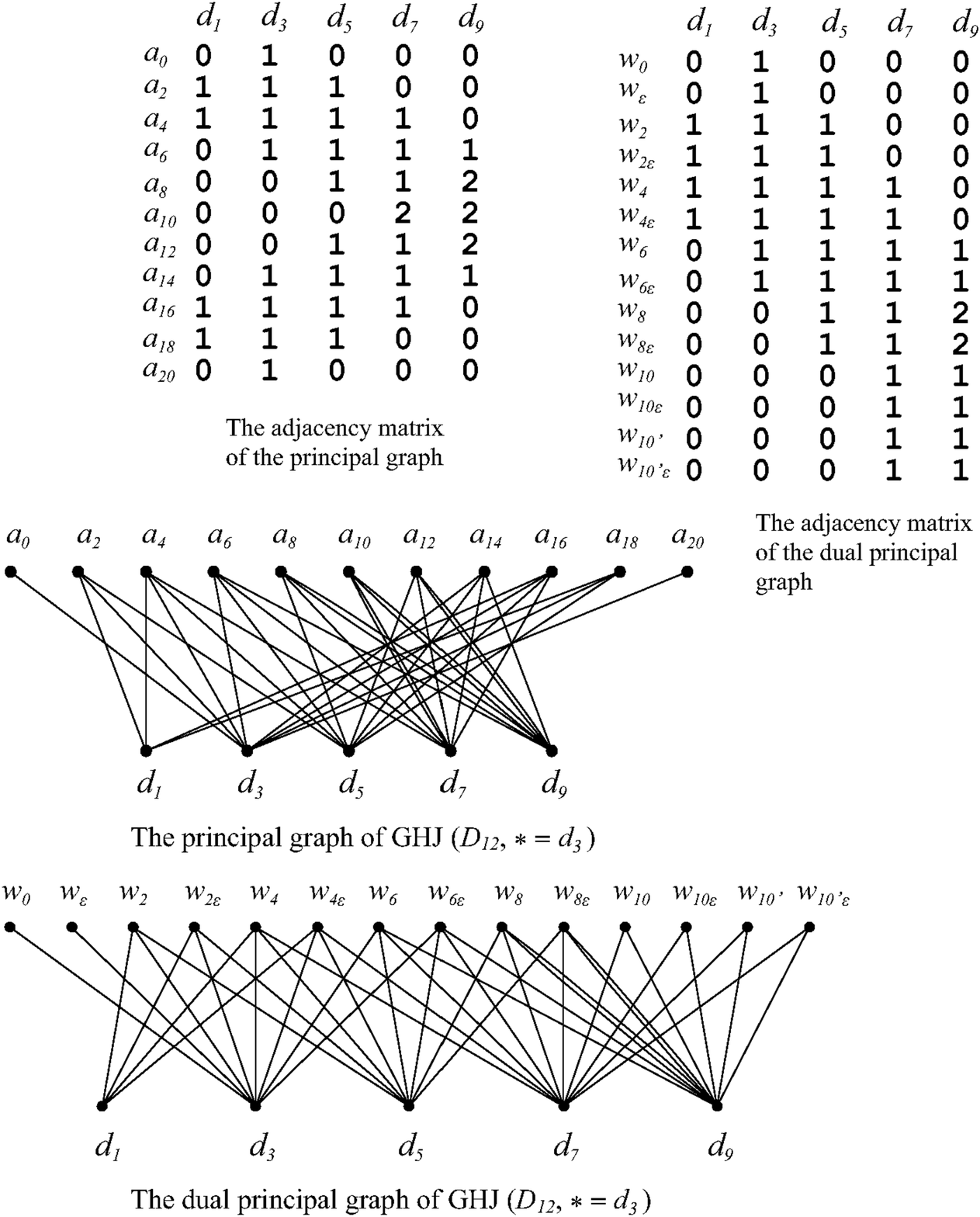}
\caption{The (dual) principal graph of GHJ$(D_{12}, *=d_3)$.}
\label{GHJ(D12-d3)}
\end{figure}

%%%%% GHJ(D12-d4) %%%%%%%%%%%%%%%%%%%%%%%%%%%%%%%%%
\begin{figure}[H]
\centering
\includegraphics[width=130mm,clip]{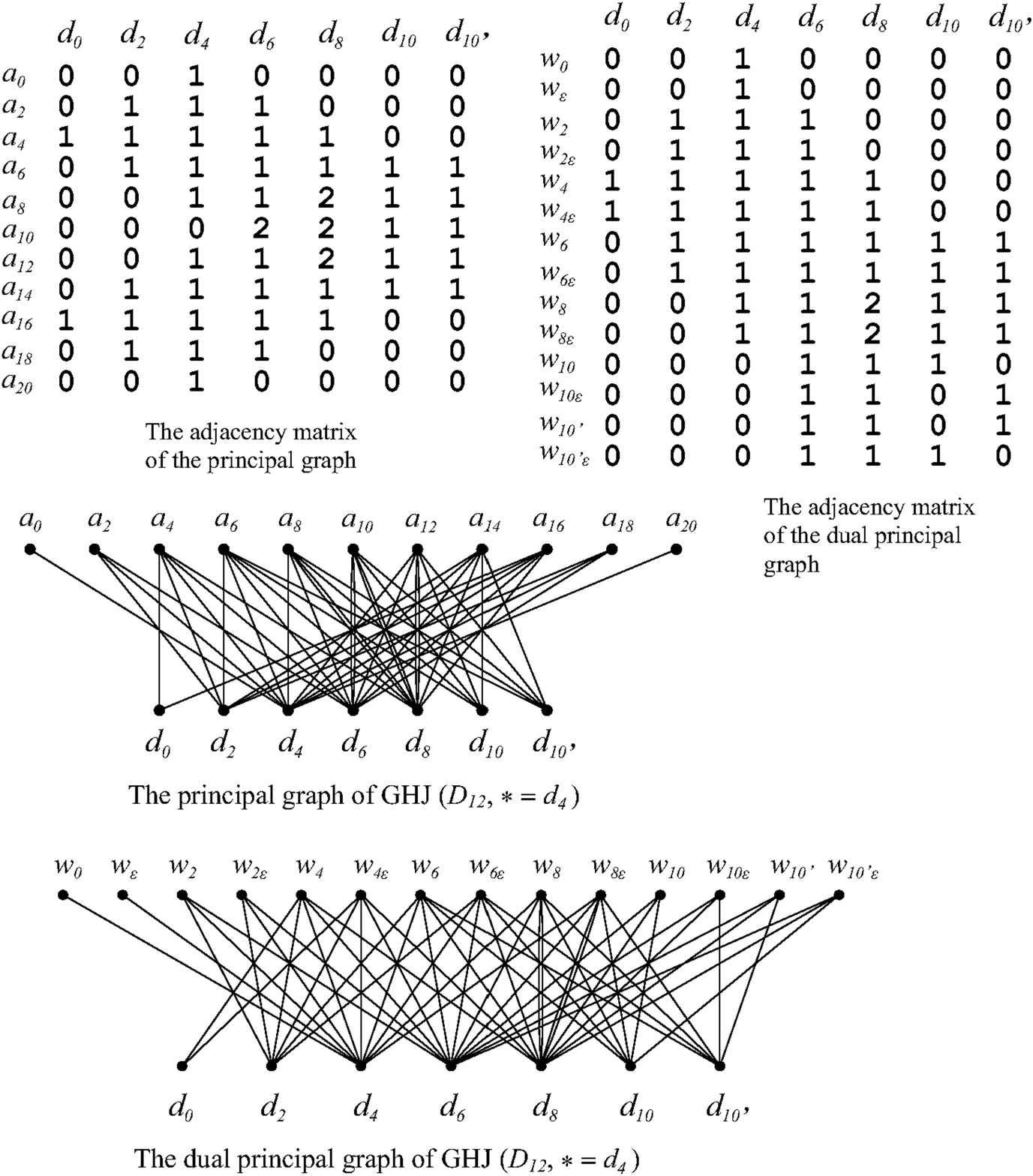}
\caption{The (dual) principal graph of GHJ$(D_{12}, *=d_4)$.}
\label{GHJ(D12-d4)}
\end{figure}

%%%%% GHJ(D12-d5) %%%%%%%%%%%%%%%%%%%%%%%%%%%%%%%%%
\begin{figure}[H]
\centering
\includegraphics[width=130mm,clip]{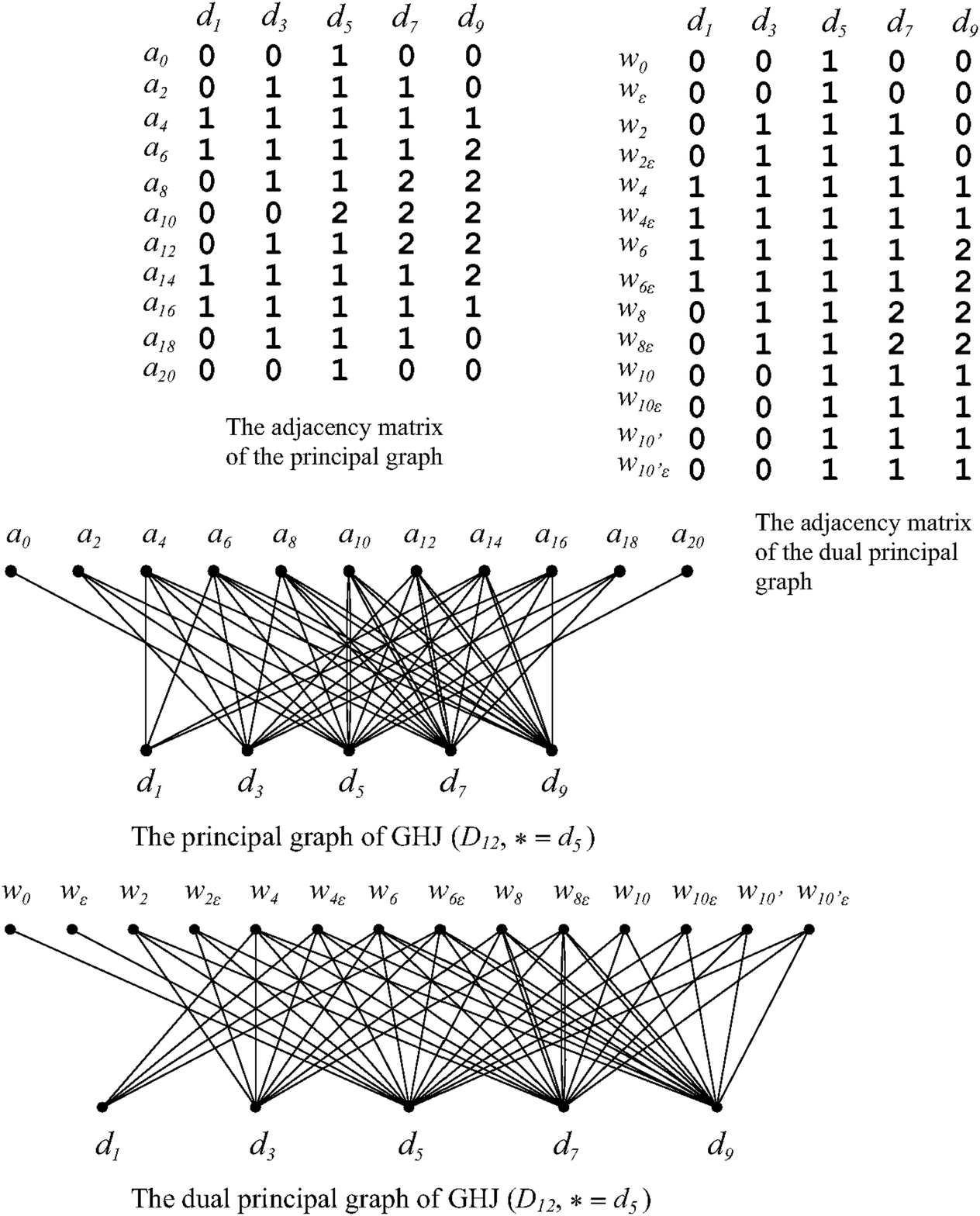}
\caption{The (dual) principal graph of GHJ$(D_{12}, *=d_5)$.}
\label{GHJ(D12-d5)}
\end{figure}

%%%%% GHJ(D12-d6) %%%%%%%%%%%%%%%%%%%%%%%%%%%%%%%%%
\begin{figure}[H]
\centering
\includegraphics[width=130mm,clip]{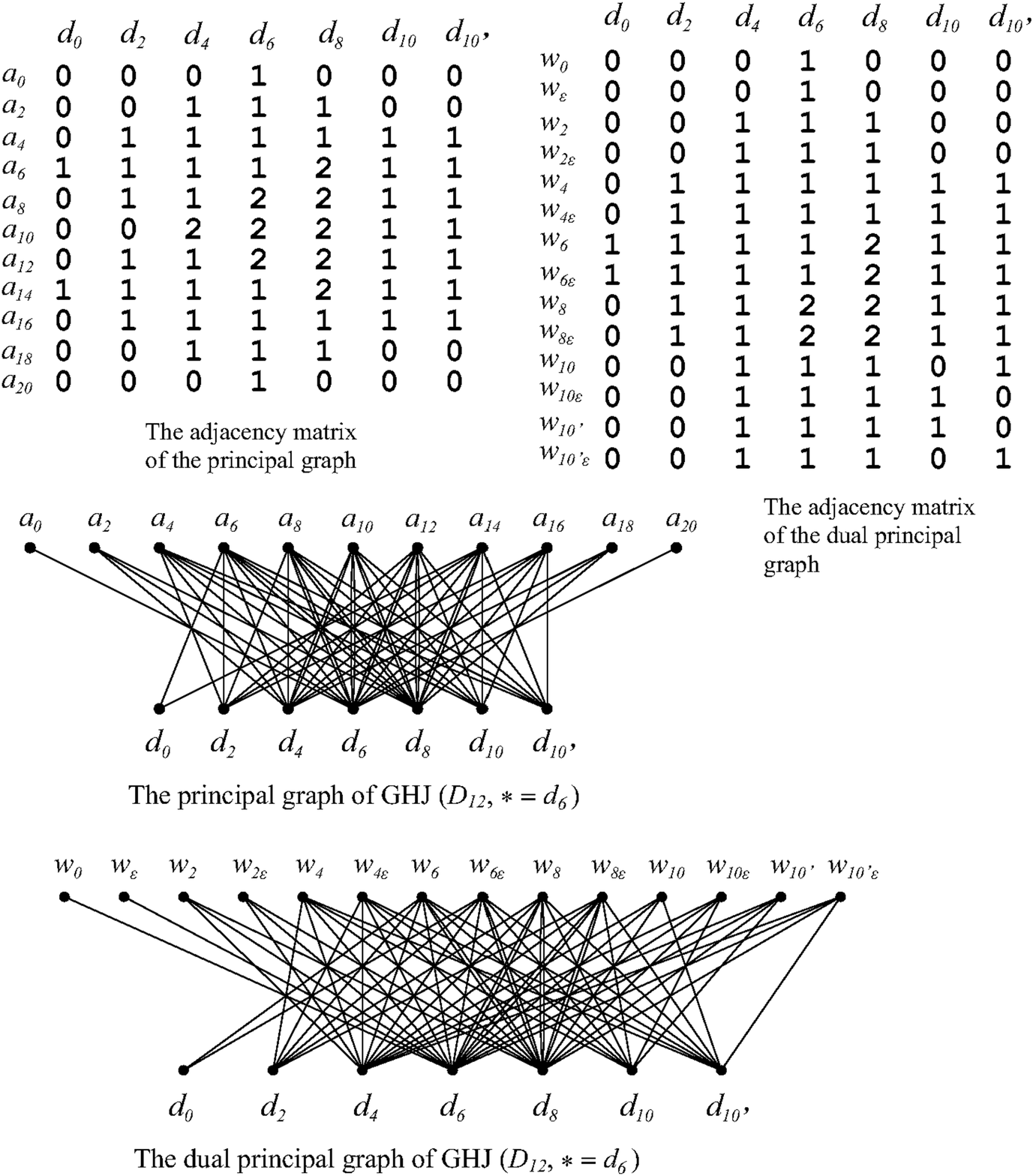}
\caption{The (dual) principal graph of GHJ$(D_{12}, *=d_6)$.}
\label{GHJ(D12-d6)}
\end{figure}

%%%%% GHJ(D12-d7) %%%%%%%%%%%%%%%%%%%%%%%%%%%%%%%%%
\begin{figure}[H]
\centering
\includegraphics[width=130mm,clip]{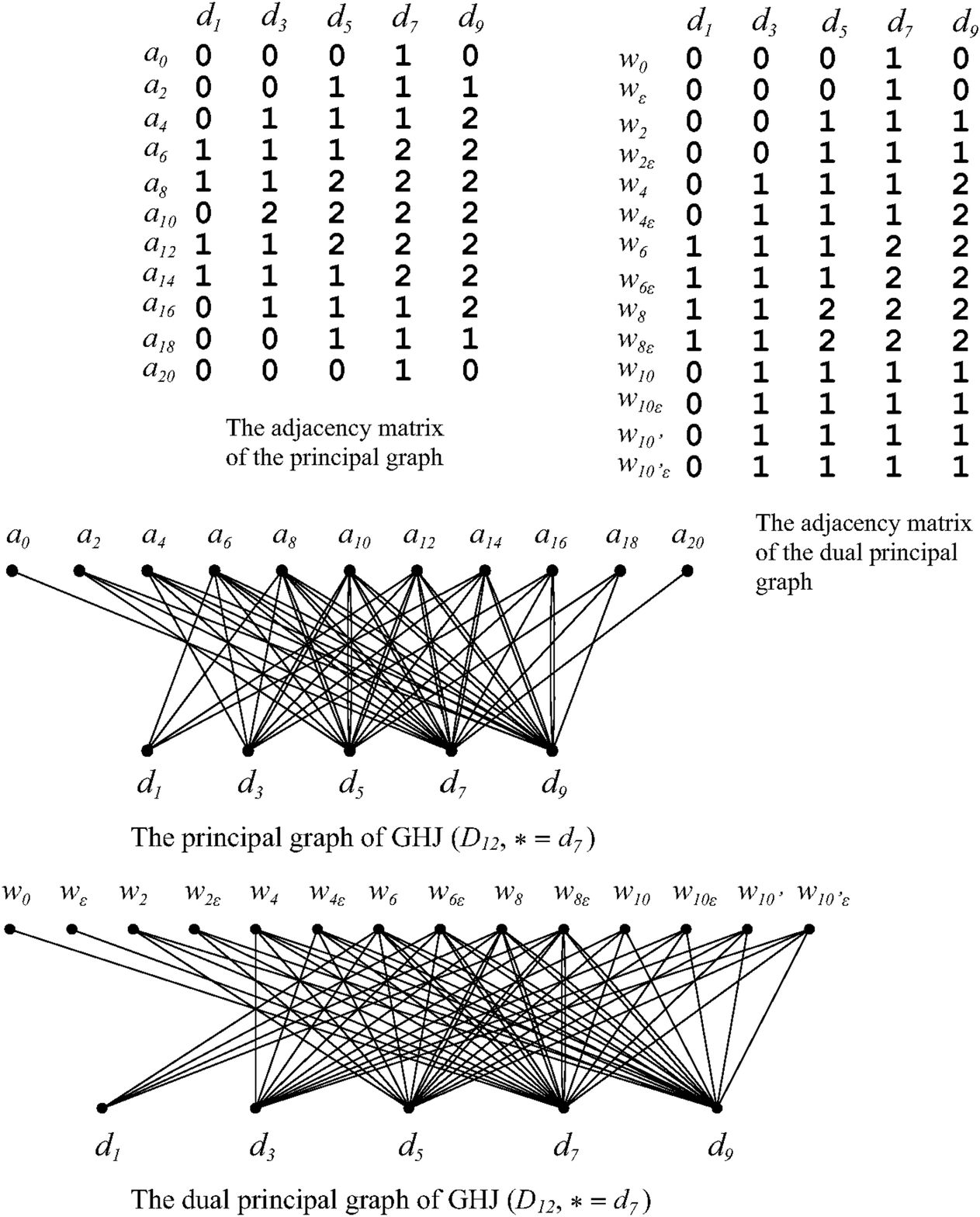}
\caption{The (dual) principal graph of GHJ$(D_{12}, *=d_7)$.}
\label{GHJ(D12-d7)}
\end{figure}

%%%%% GHJ(D12-d8) %%%%%%%%%%%%%%%%%%%%%%%%%%%%%%%%%
\begin{figure}[H]
\centering
\includegraphics[width=130mm,clip]{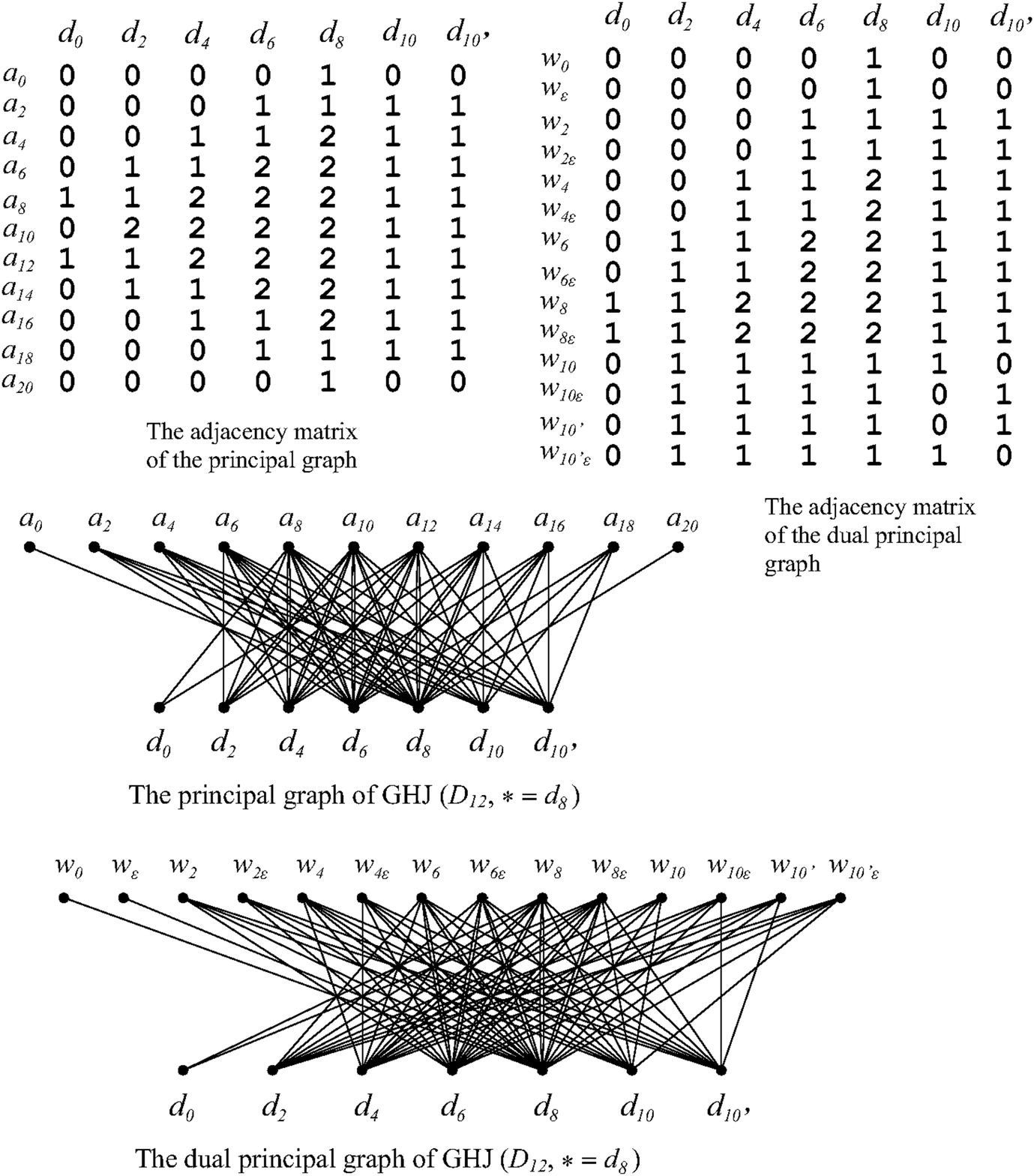}
\caption{The (dual) principal graph of GHJ$(D_{12}, *=d_8)$.}
\label{GHJ(D12-d8)}
\end{figure}

%%%%% GHJ(D12-d9) %%%%%%%%%%%%%%%%%%%%%%%%%%%%%%%%%
\begin{figure}[H]
\centering
\includegraphics[width=130mm,clip]{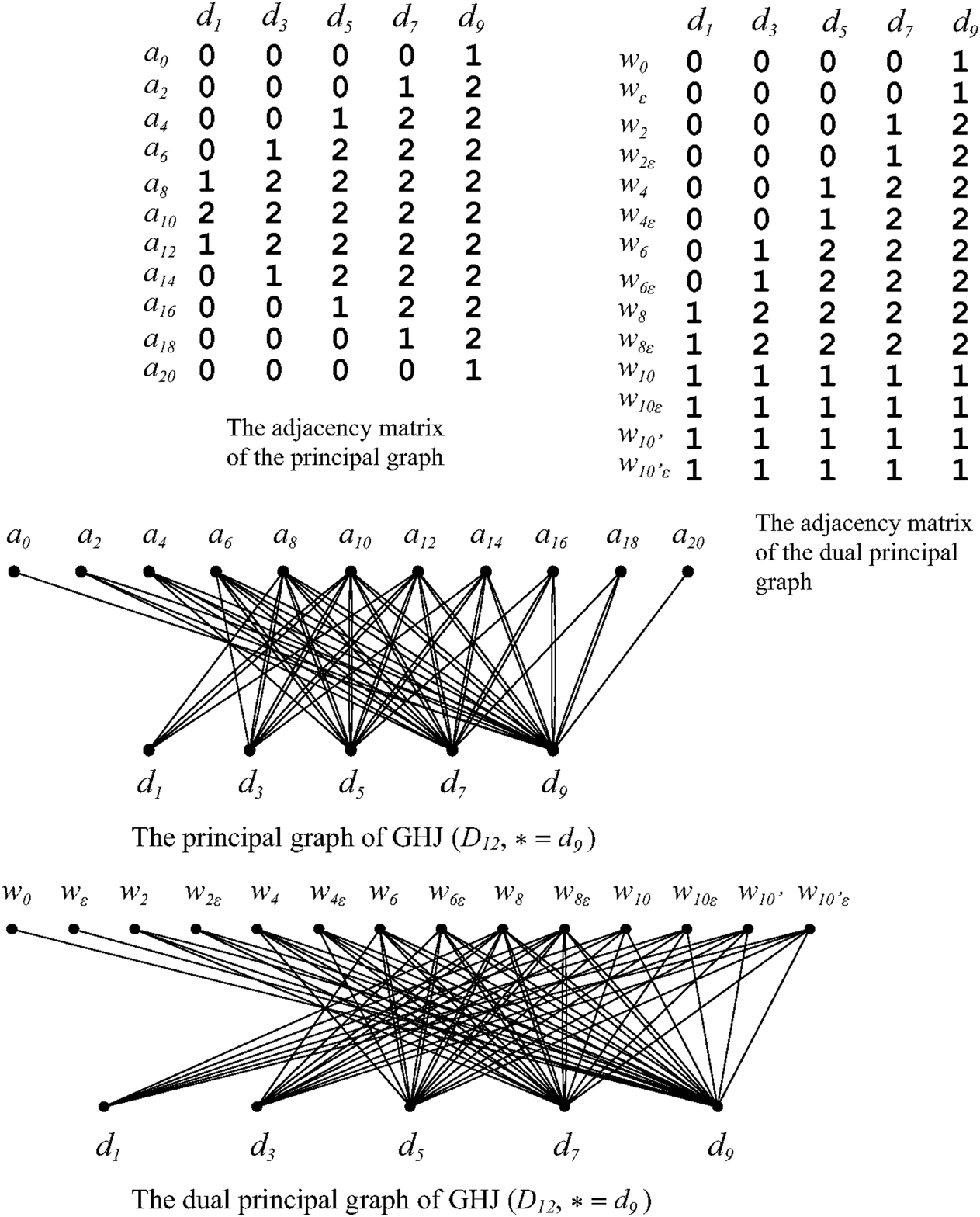}
\caption{The (dual) principal graph of GHJ$(D_{12}, *=d_9)$.}
\label{GHJ(D12-d9)}
\end{figure}

%%%%% GHJ(D12-d10) %%%%%%%%%%%%%%%%%%%%%%%%%%%%%%%%%

%%%%%%%%%%%%%%%%%%%%%%%%%%%%%%%%%%%%%%%%%%%%%%%%%%
%%% The (dual) pincipal graphs of             %%%%
%%%     Goodman-de la Harpe-Jones subfactors  %%%%
%%%%%%%%%%%%%%%%%%%%%%%%%%%%%%%%%%%%%%%%%%%%%%%%%%
%%%  Adj-Matrices PG+DPG for type D           %%%%
%%%%%%%%%%%%%%%%%%%%%%%%%%%%%%%%%%%%%%%%%%%%%%%%%%

%%%%%%%%%%%%%%%%%%%%%%%%%%%%%%%%%%%%%%%%%%%%%%%%%%

%%%%%%%%%%%%%%%%%%%%%%%%%%%%%%%%%%%%%%%%%%%%%%%%%%

%%%%%%%%%%%%%%%%%%%%%%%%%%%%%%%%%%%%%%%%%%%%%%%%%%

%%%%%%%%%%%%%%%%%%%%%%%%%%%%%%%%%%%%%%%%%%%%%%%%%%

%%%%%%%%%%%%%%%%%%%%%%%%%%%%%%%%%%%%%%%%%%%%%%%%%%

%%%%%%%%%%%%%%%%%%%%%%%%%%%%%%%%%%%%%%%%%%%%%%%%%%
%%% The (dual) pincipal graphs of             %%%%
%%%     Goodman-de la Harpe-Jones subfactors  %%%%
%%%%%%%%%%%%%%%%%%%%%%%%%%%%%%%%%%%%%%%%%%%%%%%%%%
%%%       E6                                  %%%%
%%%%%%%%%%%%%%%%%%%%%%%%%%%%%%%%%%%%%%%%%%%%%%%%%%

%%%%% GHJ(E6-e0) %%%%%%%%%%%%%%%%%%%%%%%%%%%%%%%%%

%\vspace{-10mm}
%%%%% GHJ(E6-e1) %%%%%%%%%%%%%%%%%%%%%%%%%%%%%%%%%
\begin{figure}[H]
\centering
\includegraphics[width=140mm,clip]{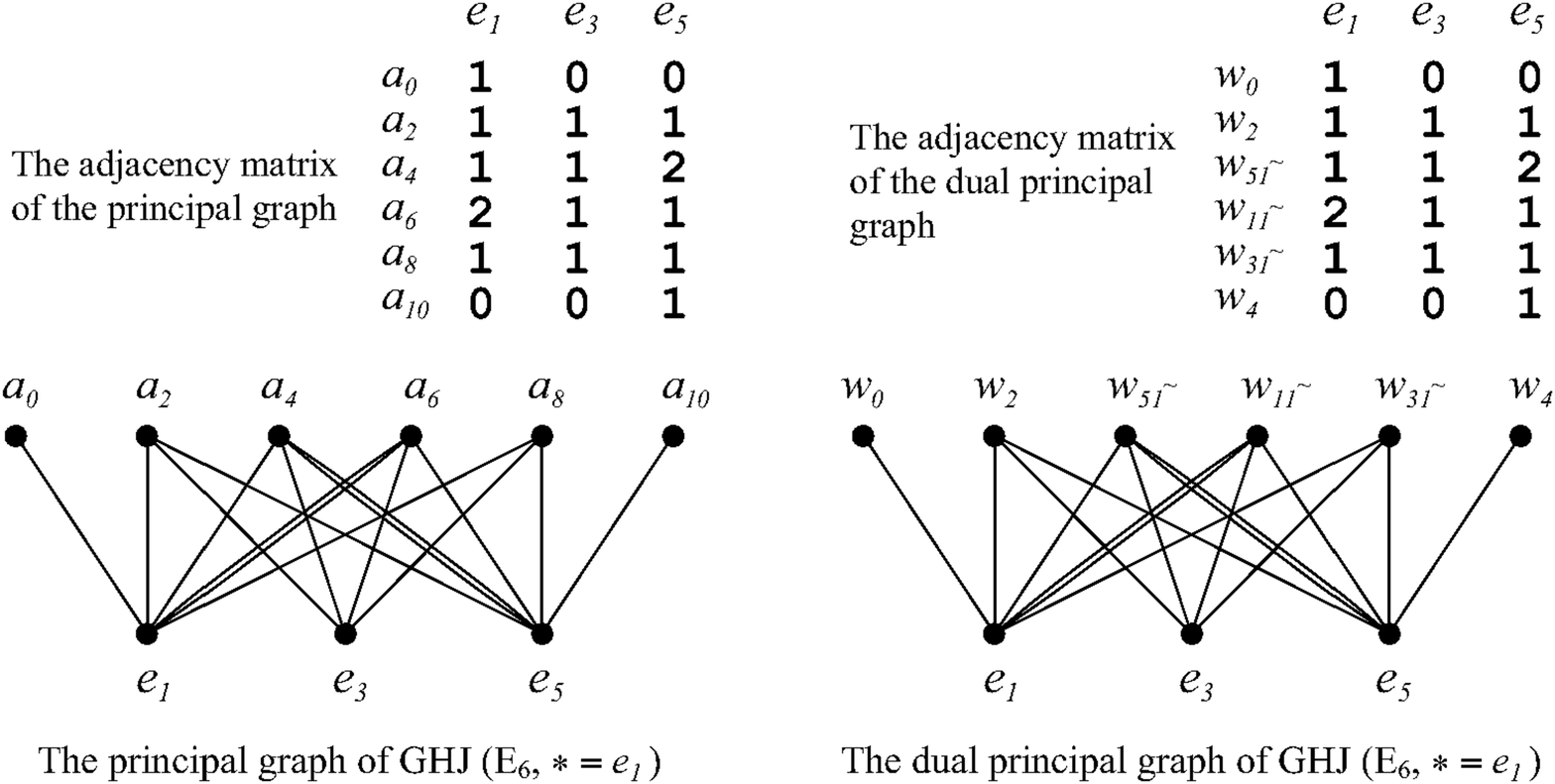}
\caption{The (dual) principal graph of GHJ$(E_6, *=e_1)$.}
\label{GHJ(E6-e1)}
\end{figure}

%%%%% GHJ(E6-e2) %%%%%%%%%%%%%%%%%%%%%%%%%%%%%%%%%
\begin{figure}[H]
\centering
\includegraphics[width=140mm,clip]{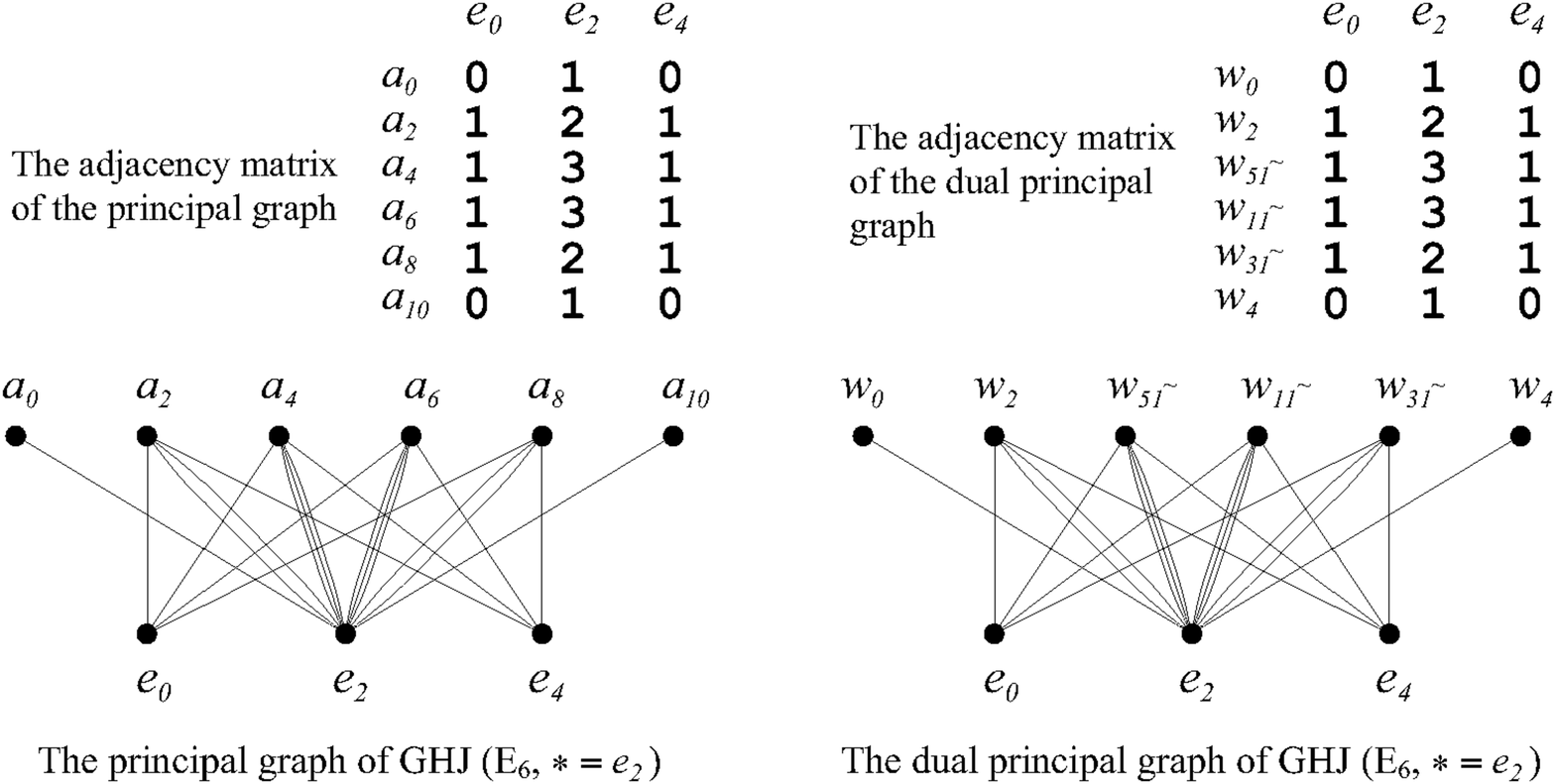}
\caption{The (dual) principal graph of GHJ$(E_6, *=e_2)$.}
\label{GHJ(E6-e2)}
\end{figure}

%%%%% GHJ(E6-e3) %%%%%%%%%%%%%%%%%%%%%%%%%%%%%%%%%

%%%%%%%%%%%%%%%%%%%%%%%%%%%%%%%%%%%%%%%%%%%%%%%%%%
%%% The (dual) pincipal graphs of             %%%%
%%%     Goodman-de la Harpe-Jones subfactors  %%%%
%%%%%%%%%%%%%%%%%%%%%%%%%%%%%%%%%%%%%%%%%%%%%%%%%%
%%%       E7                                  %%%%
%%%%%%%%%%%%%%%%%%%%%%%%%%%%%%%%%%%%%%%%%%%%%%%%%%

%%%%% GHJ(E7-e0) %%%%%%%%%%%%%%%%%%%%%%%%%%%%%%%%%

%%%%% GHJ(E7-e1) %%%%%%%%%%%%%%%%%%%%%%%%%%%%%%%%%
\begin{figure}[H]
\centering
\includegraphics[width=140mm,clip]{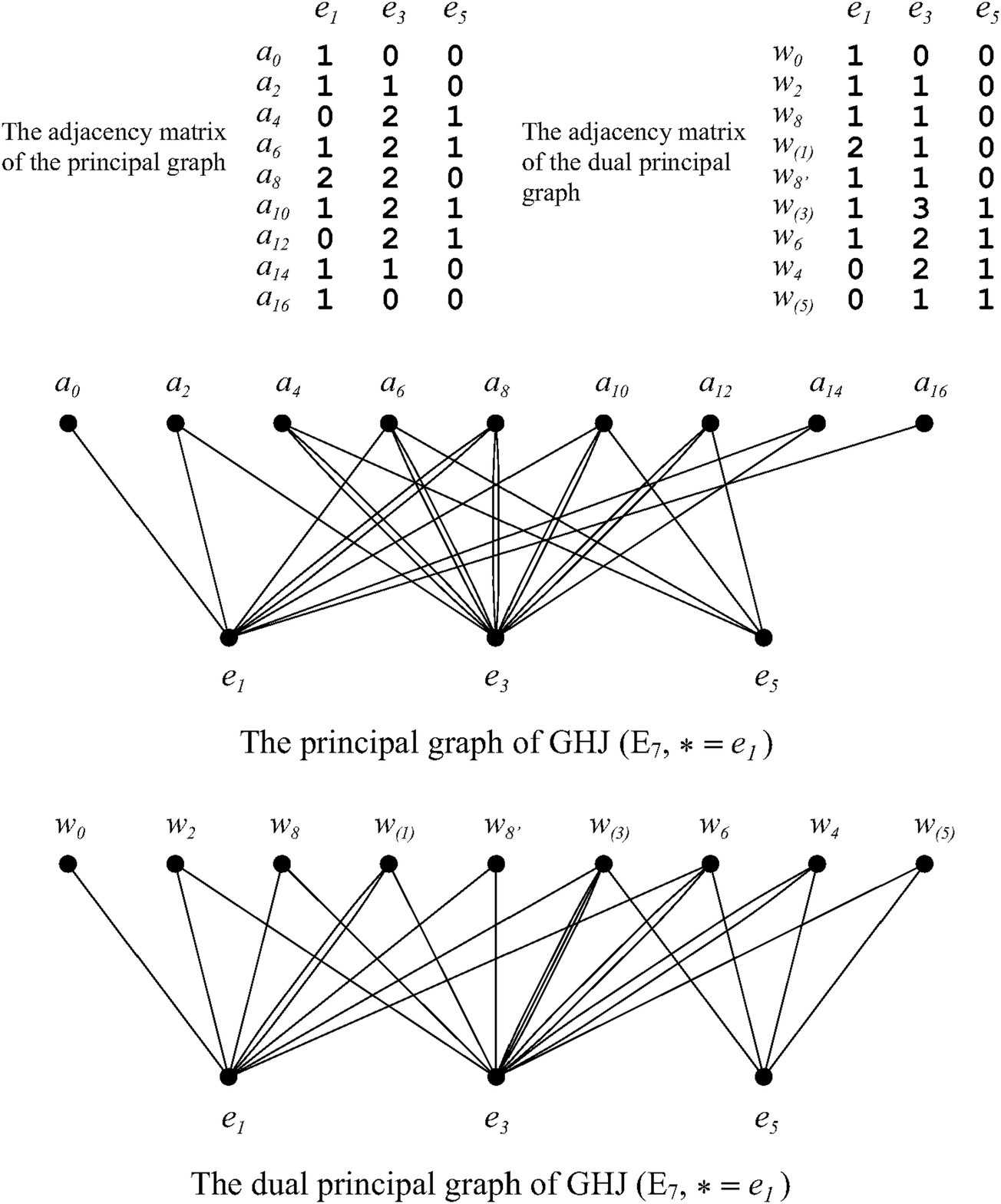}
\caption{The (dual) principal graph of the GHJ subfactor 
corresponding to $(E_7, *=e_1)$.}
\label{GHJ(E7-e1)}
\end{figure}

%%%%% GHJ(E7-e2) %%%%%%%%%%%%%%%%%%%%%%%%%%%%%%%%%
\begin{figure}[H]
\centering
\includegraphics[width=140mm,clip]{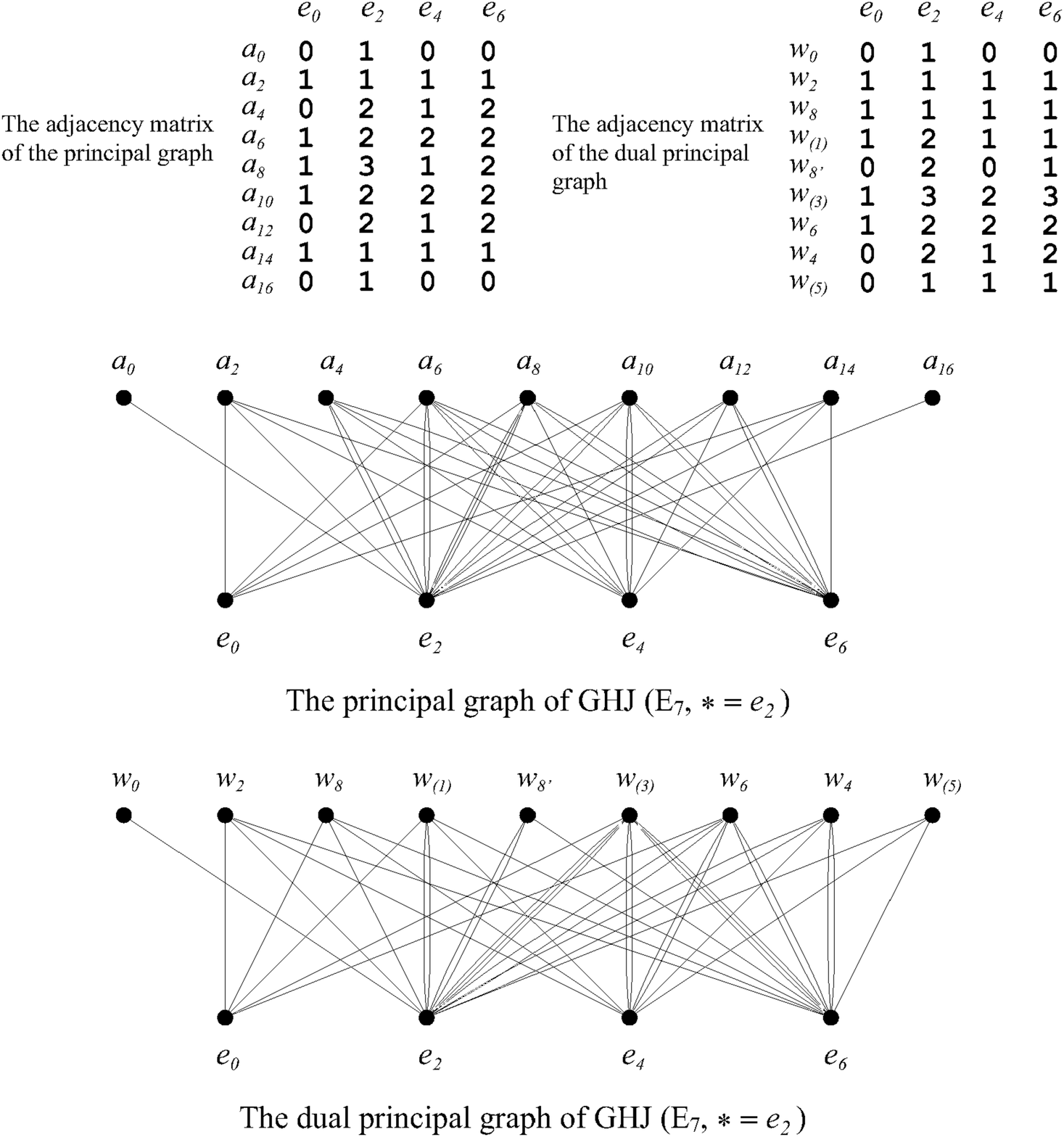}
\caption{The (dual) principal graph of the GHJ subfactor 
corresponding to $(E_7, *=e_2)$.}
\label{GHJ(E7-e2)}
\end{figure}

%%%%% GHJ(E7-e3) %%%%%%%%%%%%%%%%%%%%%%%%%%%%%%%%%
\begin{figure}[H]
\centering
\includegraphics[width=140mm,clip]{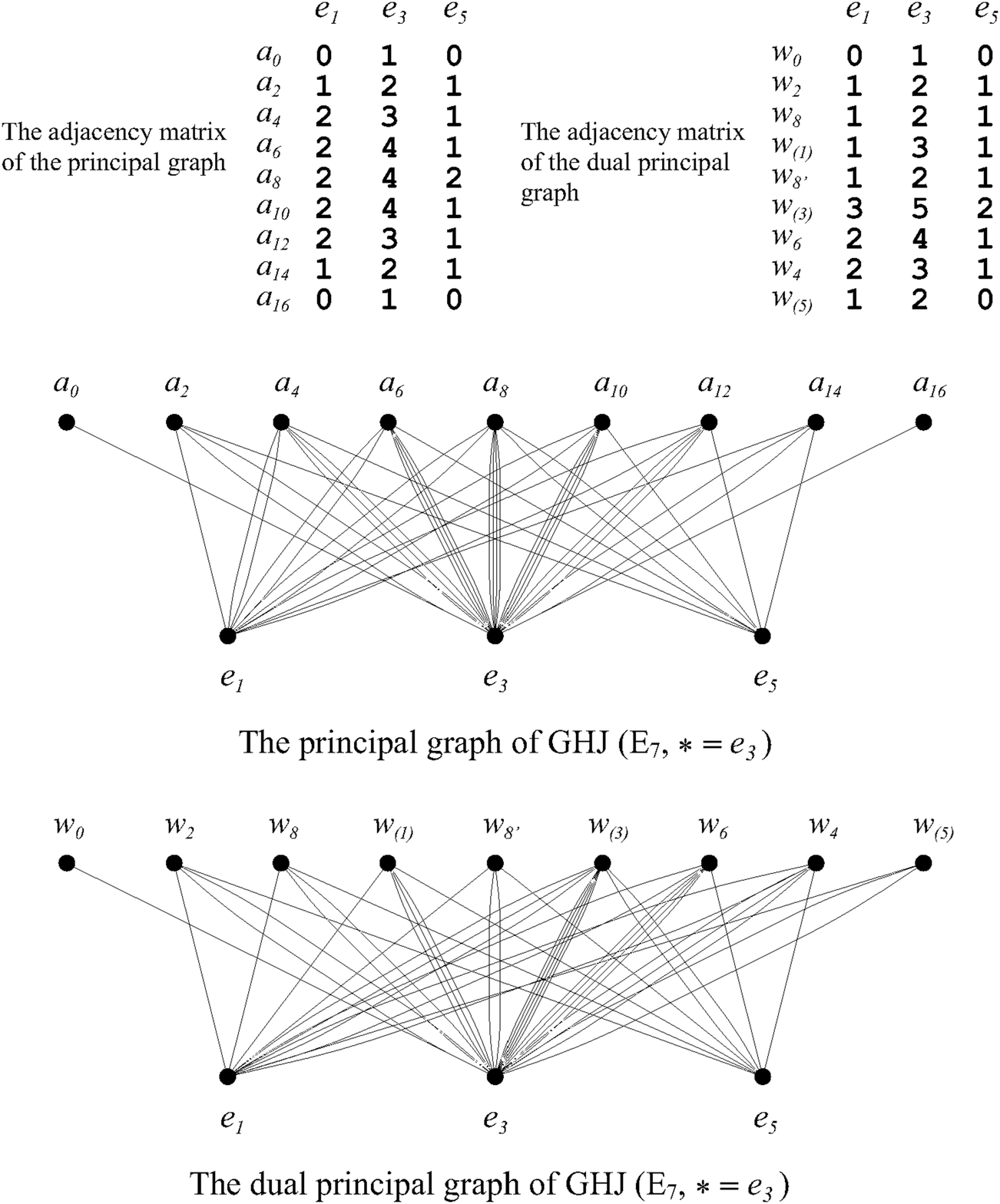}
\caption{The (dual) principal graph of the GHJ subfactor 
corresponding to $(E_7, *=e_3)$.}
\label{GHJ(E7-e3)}
\end{figure}

%%%%% GHJ(E7-e4) %%%%%%%%%%%%%%%%%%%%%%%%%%%%%%%%%
\begin{figure}[H]
\centering
\includegraphics[width=140mm,clip]{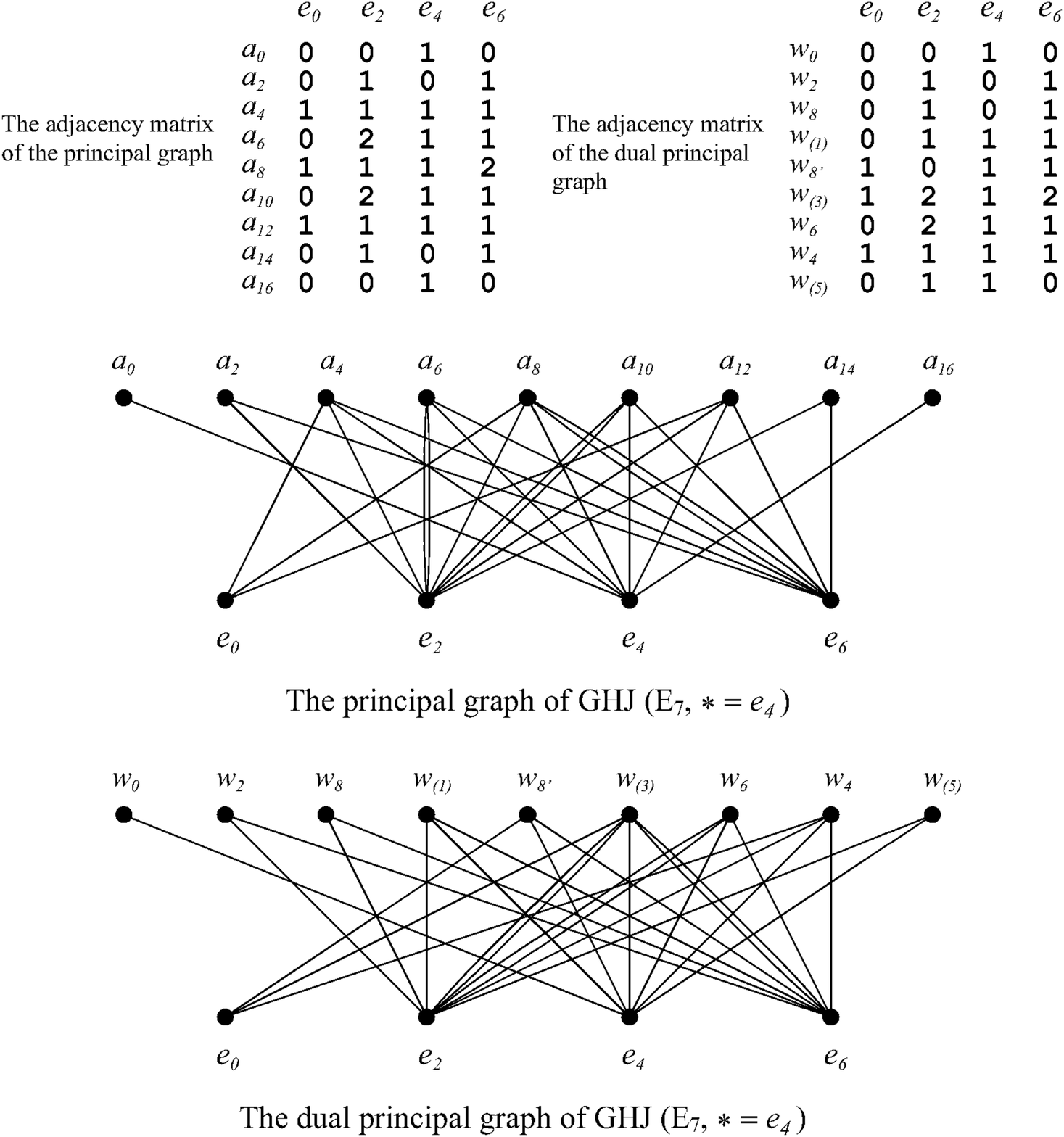}
\caption{The (dual) principal graph of the GHJ subfactor 
corresponding to $(E_7, *=e_4)$.}
\label{GHJ(E7-e4)}
\end{figure}

%%%%% GHJ(E7-e5) %%%%%%%%%%%%%%%%%%%%%%%%%%%%%%%%%
\begin{figure}[H]
\centering
\includegraphics[width=140mm,clip]{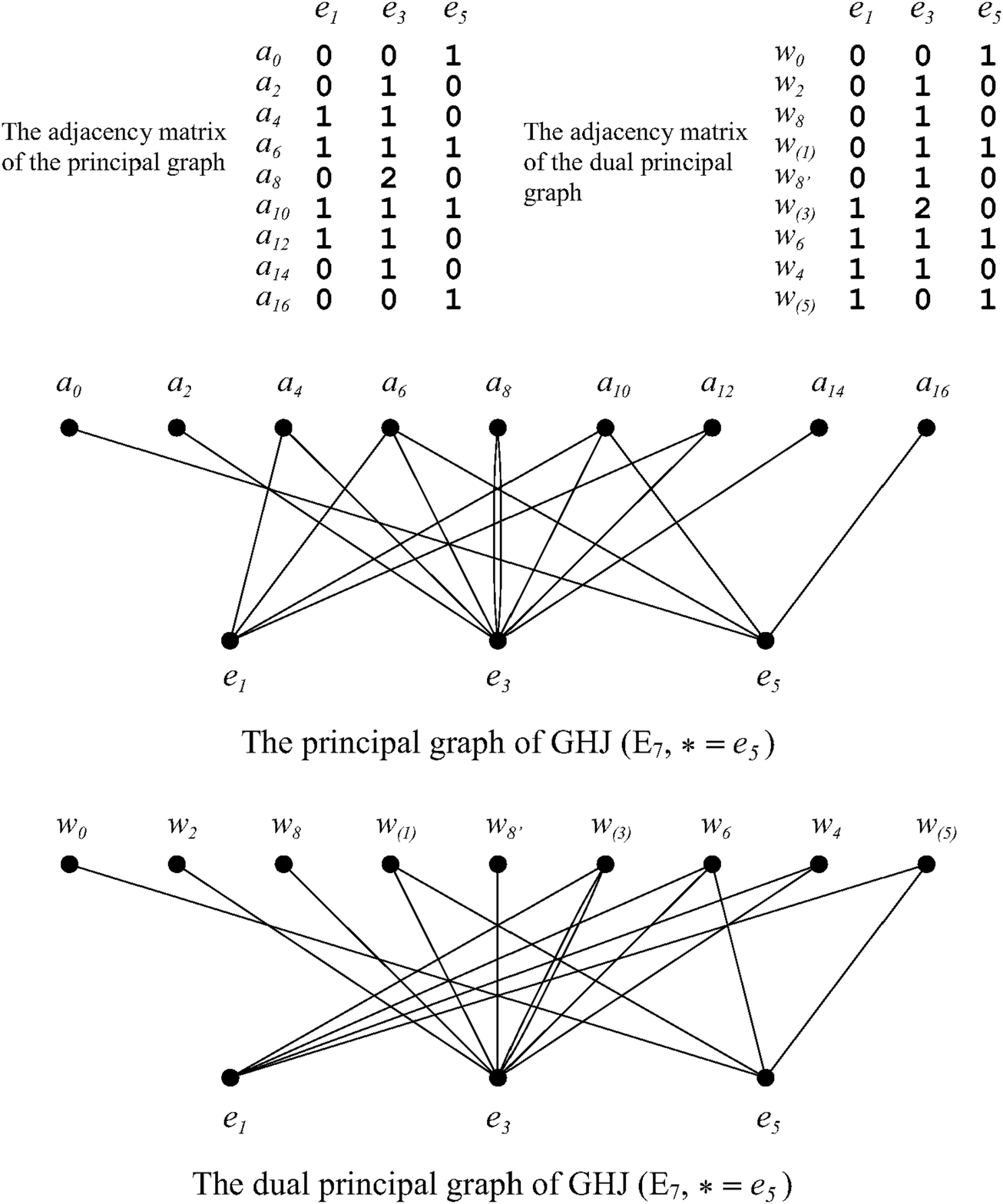}
\caption{The (dual) principal graph of the GHJ subfactor 
corresponding to $(E_7, *=e_5)$.}
\label{GHJ(E7-e5)}
\end{figure}

%%%%% GHJ(E7-e6) %%%%%%%%%%%%%%%%%%%%%%%%%%%%%%%%%

%%%%%%%%%%%%%%%%%%%%%%%%%%%%%%%%%%%%%%%%%%%%%%%%%%
%%% The (dual) pincipal graphs of             %%%%
%%%     Goodman-de la Harpe-Jones subfactors  %%%%
%%%%%%%%%%%%%%%%%%%%%%%%%%%%%%%%%%%%%%%%%%%%%%%%%%
%%%       E8                                  %%%%
%%%%%%%%%%%%%%%%%%%%%%%%%%%%%%%%%%%%%%%%%%%%%%%%%%

%%%%% GHJ(E8-e0) %%%%%%%%%%%%%%%%%%%%%%%%%%%%%%%%%

%%%%% GHJ(E8-e1) %%%%%%%%%%%%%%%%%%%%%%%%%%%%%%%%%
\begin{figure}[H]
\centering
\includegraphics[width=140mm,clip]{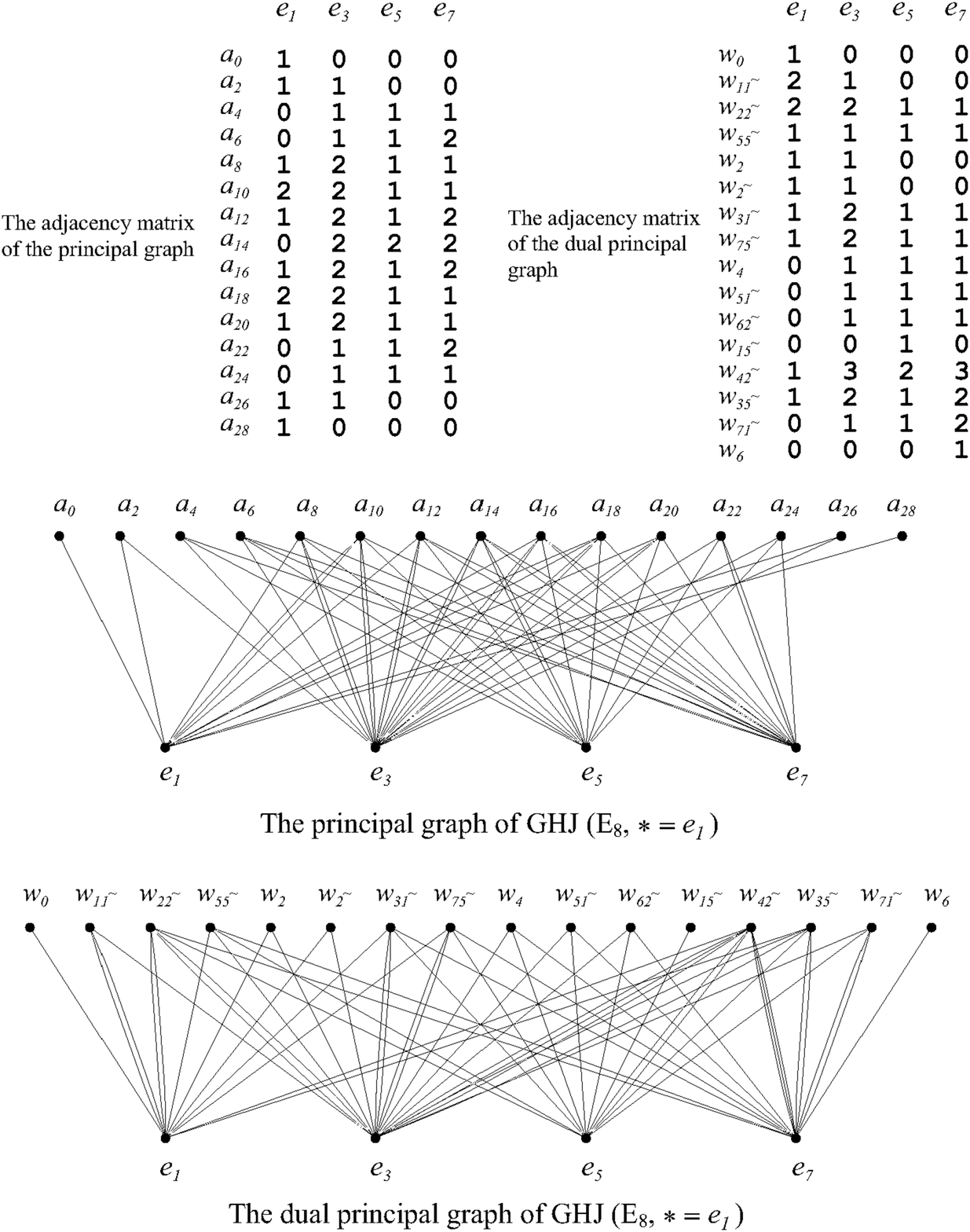}
\caption{The (dual) principal graph of the GHJ subfactor 
corresponding to $(E_8, *=e_1)$.}
\label{GHJ(E8-e1)}
\end{figure}

%%%%% GHJ(E8-e2) %%%%%%%%%%%%%%%%%%%%%%%%%%%%%%%%%
\begin{figure}[H]
\centering
\includegraphics[width=140mm,clip]{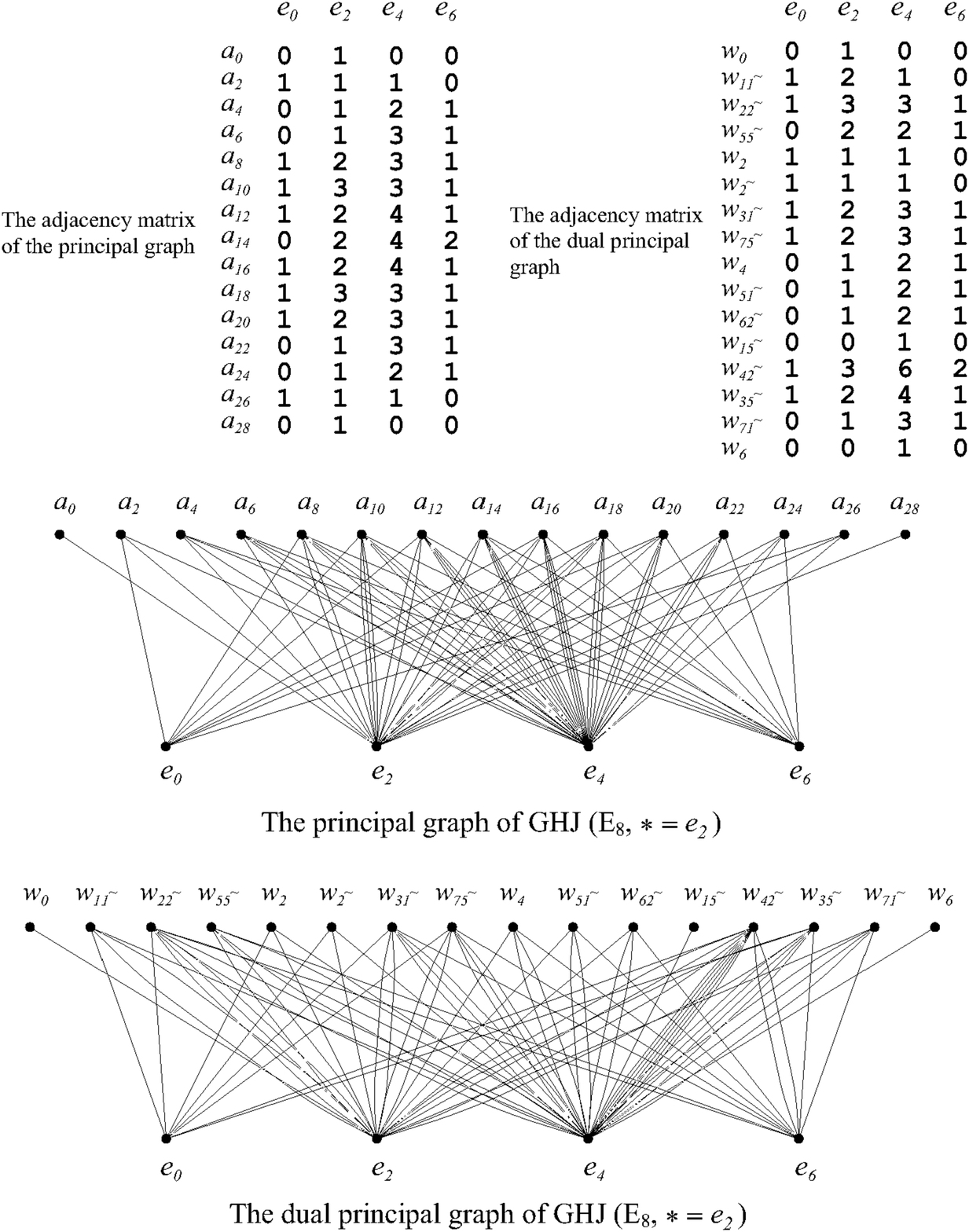}
\caption{The (dual) principal graph of the GHJ subfactor 
corresponding to $(E_8, *=e_2)$.}
\label{GHJ(E8-e2)}
\end{figure}

%%%%% GHJ(E8-e3) %%%%%%%%%%%%%%%%%%%%%%%%%%%%%%%%%
\begin{figure}[H]
\centering
\includegraphics[width=140mm,clip]{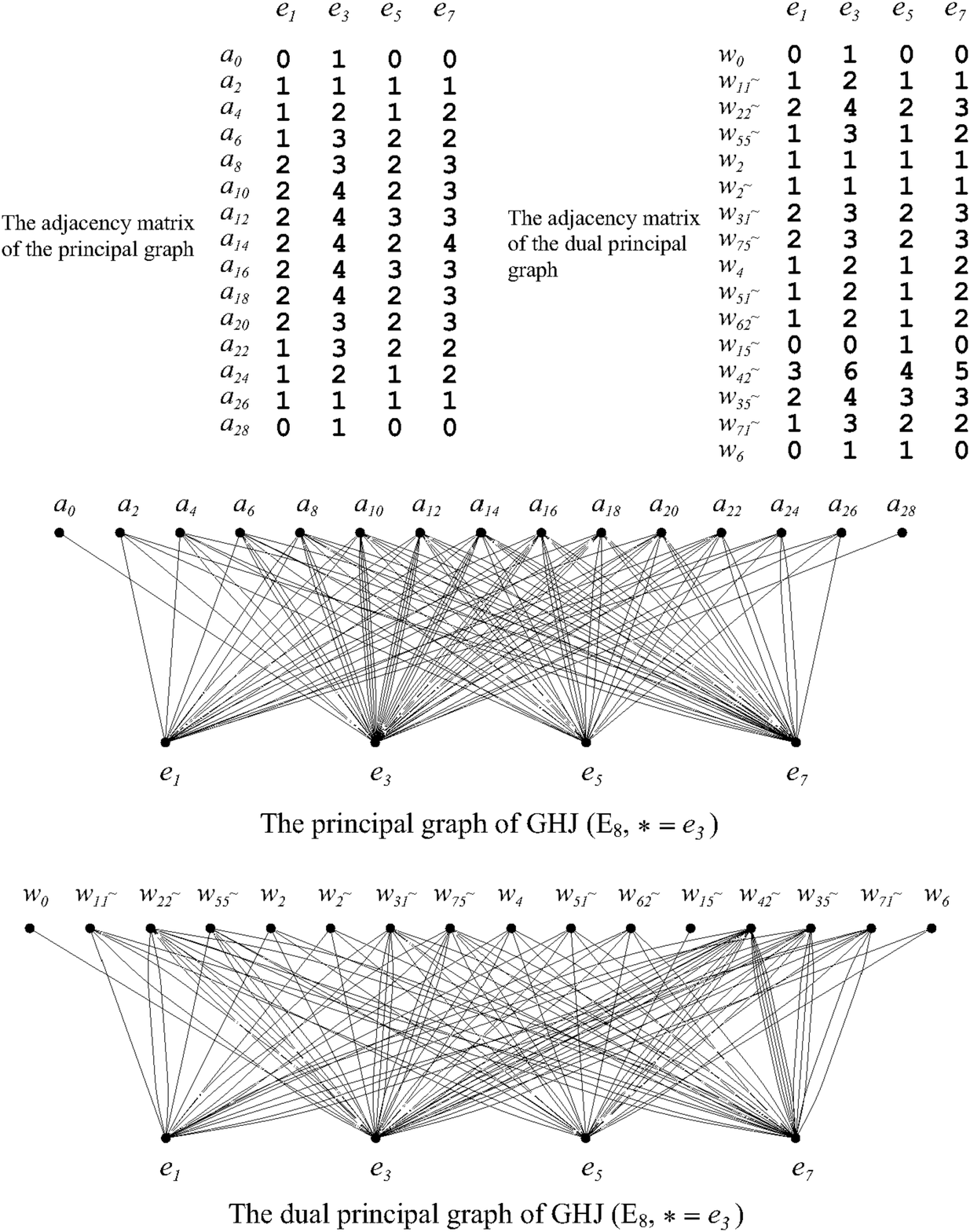}
\caption{The (dual) principal graph of the GHJ subfactor 
corresponding to $(E_8, *=e_3)$.}
\label{GHJ(E8-e3)}
\end{figure}

%%%%% GHJ(E8-e4) %%%%%%%%%%%%%%%%%%%%%%%%%%%%%%%%%
\begin{figure}[H]
\centering
\includegraphics[width=140mm,clip]{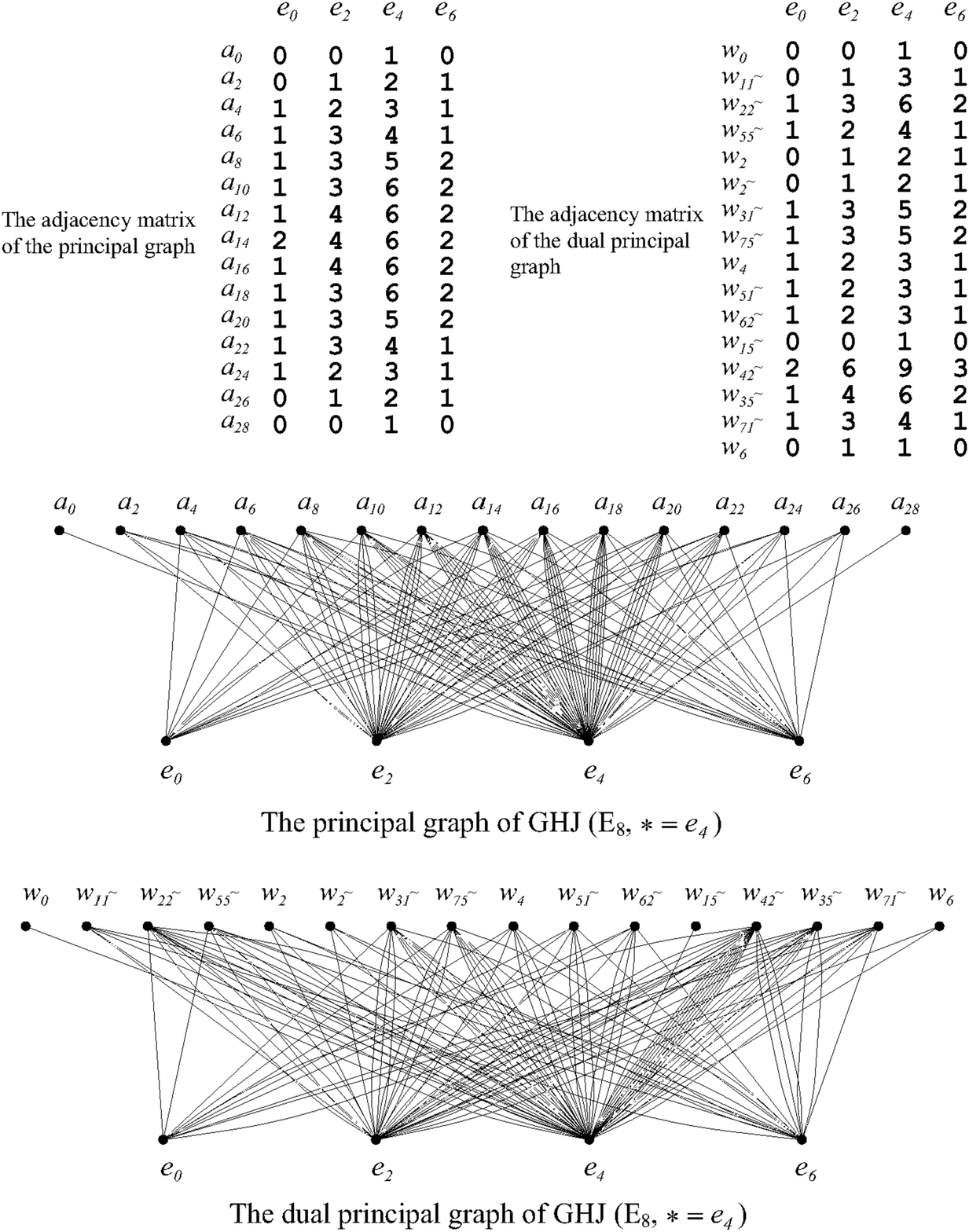}
\caption{The (dual) principal graph of the GHJ subfactor 
corresponding to $(E_8, *=e_4)$.}
\label{GHJ(E8-e4)}
\end{figure}

%%%%% GHJ(E8-e5) %%%%%%%%%%%%%%%%%%%%%%%%%%%%%%%%%
\begin{figure}[H]
\centering
\includegraphics[width=140mm,clip]{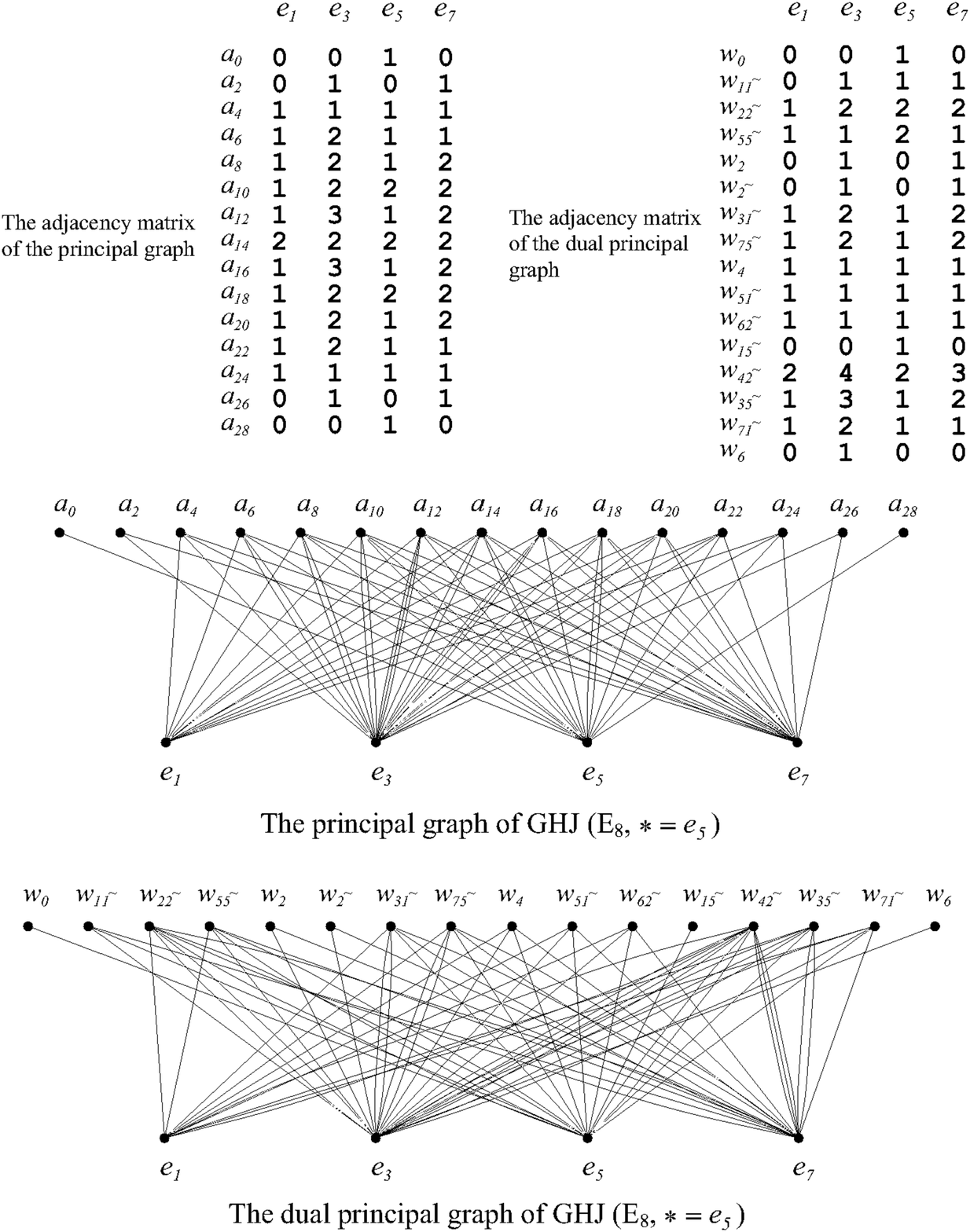}
\caption{The (dual) principal graph of the GHJ subfactor 
corresponding to $(E_8, *=e_5)$.}
\label{GHJ(E8-e5)}
\end{figure}

%%%%% GHJ(E8-e6) %%%%%%%%%%%%%%%%%%%%%%%%%%%%%%%%%
\begin{figure}[H]
\centering
\includegraphics[width=140mm,clip]{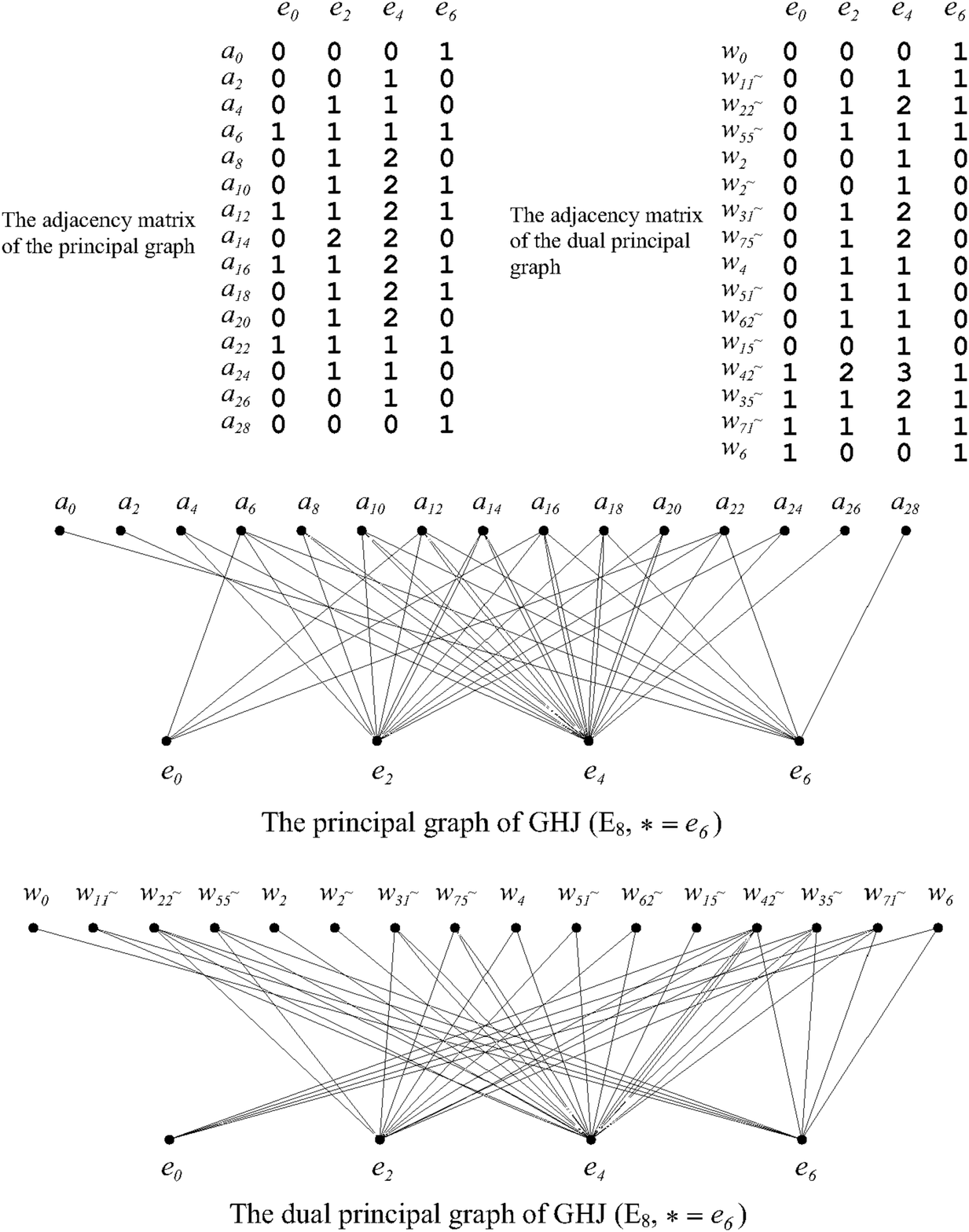}
\caption{The (dual) principal graph of the GHJ subfactor 
corresponding to $(E_8, *=e_6)$.}
\label{GHJ(E8-e6)}
\end{figure}

%%%%% GHJ(E8-e7) %%%%%%%%%%%%%%%%%%%%%%%%%%%%%%%%%

%%%%%%%%%%%%%%%%%%%%%%%%%%%%%%%%%%%%%%%%%%%%%%%%%%%%%%%%%%
%%%%%%%%%%%%%%%%%%%%%%%%%%%%%%%%%%%%%%%%%%%%%%%%%%%%%%%%%%
%%%%%%%%%%%%%%%%%%%% Bibliography %%%%%%%%%%%%%%%%%%%%%%%%
%%%%%%%%%%%%%%%%%%%%%%%%%%%%%%%%%%%%%%%%%%%%%%%%%%%%%%%%%%
%%%%%%%%%%%%%%%%%%%%%%%%%%%%%%%%%%%%%%%%%%%%%%%%%%%%%%%%%%

%%%%%%%%%%%%%%%%%%%%%%%%%%%%%%%%%%%%%%%%%%%%%%%%%%%%%%%%%%%%%%%%%%%%
%%%%%%%%%%%%%%%%%%%%%%%%%%%%%%%%%%%%%%%%%%%%%%%%%%%%%%%%%%%%%%%%%%%%
%%%%%%%%%%%%%%%%%%%%%%%%%%%%%%%%%%%%%%%%%%%%%%%%%%%%%%%%%%%%%%%%%%%%

\begin{thebibliography}{99}

\bibitem{AH}
Asaeda, M. and U. Haagerup, {Exotic subfactors of finite depth
with Jones indices $(5+\sqrt{13})/2$ and $(5+\sqrt{17})/2$}, 
Comm. Math. Phys. {\bf 202} (1999), 1--63.

\bibitem{BiN}
Bion-Nadal, J., 
{Subfactor of the hyperfinite II$_1$ factor with 
Coxeter graph $E_6$ as invariant},
J. Operator Theory {\bf 28} (1992), 27--50.

\bibitem{Bs}
Bisch, D., 
{Principal graphs of subfactors with small Jones index},
Math. Ann. {\bf 311} (1998), 223--231.

\bibitem{Bs2}
Bisch, D., 
{On the existence of central sequences in subfactors}, 
Trans. Amer. Math. Soc. {\bf 321} (1990), 117--128.

\bibitem{BE1}
B\"ockenhauer, J. and D. E. Evans,
{Modular Invariants, Graphs and $\a$-Induction for Nets
of Subfactors I}, Comm. Math. Phys. {\bf 197} (1998), 361--386.

\bibitem{BE2}
B\"ockenhauer, J. and D. E. Evans,
{Modular Invariants, Graphs and $\a$-Induction for Nets
of Subfactors II}, Comm. Math. Phys. {\bf 200} (1999), 57--103.

\bibitem{BE3}
B\"ockenhauer, J. and D. E. Evans,
{Modular Invariants, Graphs and $\a$-Induction for Nets
of Subfactors III}, Comm. Math. Phys. {\bf 205} (1999), 183--228. 

%\bibitem{BEK1}
%B\"ockenhauer, J., D. E. Evans and Y. Kawahigashi, 
%{\it On $\alpha$-induction, chiral generators
%and modular invariants for subfactors},
%Comm. Math. Phys. {\bf 208} (1999), 429--487.

\bibitem{BEK2}
B\"ockenhauer, J., D. E. Evans and Y. Kawahigashi, 
{\it Chiral structure of modular invariants for subfactors},
Comm. Math. Phys. {\bf 210} (2000), 733--784.

%\bibitem{BEK3}
%B\"ockenhauer, J., Evans, D. E. and Kawahigashi, Y. (preprint 2000).
%Longo-Rehren subfactors arising from $\alpha$-induction.

\bibitem{EK-book} 
Evans, D. E. and Y. Kawahigashi,
{Quantum symmetries on operator algebras},
Oxford University Press, Oxford, 1998.

\bibitem{CIZ}
Cappelli, A., C. Itzykson, and J.-B. Zuber,
{The $A$-$D$-$E$ classification of minimal and
$A^{(1)}_1$ conformal invariant theories},
Comm. Math. Phys. {\bf 113} (1987), 1--26.


\bibitem{GHJ} 
Goodman, F., P. de la Harpe and V. F. R. Jones,
{Coxeter graphs and towers of algebras}, MSRI Publications,
14, Springer, Berlin, 1989.

\bibitem{G}
Goto, S., 
{On Ocneanu's theory of double triangle algebras for subfactors 
and classification of irreducible connections on the Dynkin diagrams},
Expos. Math. {\bf 28} (2010), 218--253.

\bibitem{I1}
Izumi, M.,
{Application of fusion rules to classification of subfactors}, 
Publ. Res. Inst. Math. Sci. {\bf 27} (1991), 953--994.

\bibitem{I2}
Izumi, M.,
{On flatness of the Coxeter graph $E_8$},
Pacific J. Math. {\bf 166} (1994), 305--327.

\bibitem{J1}
Jones, V. F. R.,
{Index for subfactors},
Invent. Math. {\bf 72} (1983), 1--15.

\bibitem{Ka1}
Kawahigashi, Y., 
{On flatness of Ocneanu's connections on the Dynkin diagrams
and classification of subfactors}, 
J. Funct. Anal. {\bf 127} (1995), 63--107.

\bibitem{Ka2}
Kawahigashi, Y.,
{Classification of paragroup actions on subfactors},
Publ. Res. Inst. Math. Sci. {\bf 31} (1995), 481--517.

\bibitem{O1} 
Ocneanu, A.,
{Quantized group string algebras and Galois theory for algebras},
in ``Operator algebras and applications, Vol. 2 (Warwick, 1987),''
London Math. Soc. Lect. Note Series Vol. 136,
Cambridge University Press, (1988), 119--172.

\bibitem{Oc}
Ocneanu, A., 
{Paths on Coxeter diagrams:
from Platonic solids and singularities 
to minimal models and subfactors},
(Notes recorded by S. Goto), in
{\em Lectures on operator theory},
(ed. B. V. Rajarama Bhat et al.), 
The Fields Institute Monographs, Providence, Rhode Island:
AMS Publications. (2000), 243--323.

\bibitem{Ok} Okamoto, S.,
{Invariants for subfactors arising from Coxeter graphs},
Current Topics in Operator Algebras, World Scientific Publishing, 
(1991), 84--103.

\bibitem{SV}
Sunder, V. S. and A. K. Vijayarajan,
{On the non-occurrence of the Coxeter graphs 
$\beta_{2n+1}$, $E_7$, $D_{2n+1}$ as
principal graphs of an inclusion of II$_1$ factors},
Pacific J. Math. {\bf 161} (1993), 185--200.

\bibitem{Xu} Xu, F.,
{New braided endomorphisms from conformal inclusions},
Comm. Math. Phys. {\bf 192} (1998), 349--403.



\end{thebibliography}
\end{document}